\documentclass[11pt,a4paper]{amsart}
\usepackage[all]{xy}
\usepackage{amssymb,longtable}
\usepackage{epsfig}

\makeatletter\@addtoreset{equation}{section} \makeatother

\newtheorem{theorem}[equation]{Theorem}

\newtheorem{lemma}[equation]{Lemma}

\newtheorem{conjecture}[equation]{Conjecture}

\theoremstyle{definition}

\newtheorem{definition}[equation]{Definition}

\theoremstyle{remark}
\newtheorem{remark}[equation]{Remark}

\author{Dimitra Kosta} \thanks{I would like to thank Ivan Cheltsov,
Constantin Shramov and Jihun Park for valuable comments and corrections
and the Postech Mathematics Institute for the hospitality during
the workshop on Groups of Birational Automorphisms in April 2009.
This paper was completed under financial support provided by IKY
(Greek State Scholarship Foundation).}

\title{Del Pezzo surfaces with Du Val singularities}

\begin{document}

\begin{abstract}
In this paper we show that del Pezzo surfaces of degree 1 with Du
Val singular points of type
\begin{eqnarray*}
\mathbb{A}_4 \text{, } \mathbb{A}_4 + \mathbb{A}_4 \text{,
}\mathbb{A}_4 + \mathbb{A}_3 \text{, } \mathbb{A}_4 + 2
\mathbb{A}_1
\text{, } \mathbb{A}_4 + \mathbb{A}_1 \text{, }\mathbb{A}_3+ 4 \mathbb{A}_1  \text{, }\\
 \mathbb{A}_3+ 3
\mathbb{A}_1 \text{, } 2\mathbb{A}_3+ 2 \mathbb{A}_1  \text{, }
\mathbb{A}_3+ 2 \mathbb{A}_1  \text{, } \mathbb{A}_3+ \mathbb{A}_1
\text{, }2 \mathbb{A}_3 \text{, } \mathbb{A}_3 \text{, }
\end{eqnarray*}
 admit a K\"ahler-Einstein metric.
Moreover we are going to compute global log canonical thresholds
of del Pezzo surfaces of degree 1 with Du Val singularities, and of
del Pezzo surfaces of Picard group $\mathbb{Z}$ with Du Val
 singularities.
\end{abstract}

\maketitle\thispagestyle{empty}%

\section{Introduction}

A result of Demailly-Kollar \cite{DK01} has recently drawn a lot
of attention to global log canonical thresholds of Fano varieties,
which are algebraic counterparts of the $\alpha$-invariant of Tian
for smooth Fano varieties (see \cite[Appendix~A]{Fano3folds}). At
first, we are going to introduce some definitions from singularity
theory while more details could be found in  the classical
reference \cite{Kollar-sing}.

Suppose that $X$ is a normal variety and $D= \sum d_i D_i$ is a
$\mathbb{Q}$-divisor on $X$ such that $K_X + D$ is
$\mathbb{Q}$-Cartier and let $f : Y \to X$ be a birational
morphism, where $Y$ is normal. We can write
$$
K_Y  \sim_{\mathbb{Q}} f^*(K_X+D) + \sum a(X, D, E) E \text{ .}
$$
\begin{definition}
 The discrepancy of the log pair $(X, D)$ is the number
 $$
 \text{discrep}(X, D) = \text{inf}_E \left\lbrace a(F, D, E)  |
 E \text{ is exceptional divisor over }X \right\rbrace \text{ .}
 $$
 The total discrepancy of the log pair $(X, D)$ is the number
$$
 \text{totaldiscrep}(X, D) = \text{inf}_E \left\lbrace a(F, D, E)  |
 E \text{ is divisor over }X \right\rbrace \text{ .}
 $$

\end{definition}

We say that the log pair $K_X + D$ is

\begin{itemize}
 \item Kawamata log terminal (or log terminal) iff $\text{totaldiscrep}(X, D)> -1$
 \item log canonical  iff  $\text{discrep}(X, D) \geq -1$.
\end{itemize}

Assume, now, that $X$ is a variety with log terminal singularities,
let $Z \subset X$ be a closed subvariety and let $D$ be an
effective $\mathbb{Q}$-Cartier $\mathbb{Q}$-divisor on $X$. Then
the number
$$
 \text{lct}_{Z}(X, D) = \text{sup} \left\lbrace \lambda \in \mathbb{Q} | \text{ the log pair }
 (X, \lambda D) \text{ is log canonical along } Z \right\rbrace
 $$
is called the log canonical threshold of $D$ along $Z$ and is a
positive rational number. For $Z=X$ we use the notation
$\text{lct}(X, D)$, instead of $\text{lct}_{Z}(X, D)$.
 $$
 \text{lct}(X, D) = \text{sup} \left\lbrace \lambda \in \mathbb{Q} | \text{ the log pair }
 (X, \lambda D) \text{ is log canonical  } \right\rbrace \text{ .}
 $$

 Suppose, moreover, that $X$ is a Fano variety.
\begin{definition}
 The global log canonical threshold of  $X$ is the number
 $$
 \text{lct}(X) = \text{inf} \left\lbrace \text{lct}(X,D) | D \text{ effective  divisor on
}
 X \text{ such that } D \sim_{\mathbb{Q}} -K_X \right\rbrace \text{ .}
 $$
\end{definition}

Another way to see the global log canonical threshold is to take
the $\text{inf}_n \left\lbrace \text{lct}_n (X) \right\rbrace $,
where
$$
\text{lct}_n (X) = \text{inf} \left\lbrace
\text{lct}(X,\frac{1}{n} D) | D \text{ effective }
\mathbb{Q}\text{-divisor on }
 X \text{ such that } D \in |-nK_X| \right\rbrace
$$
In the case when $X$ is a del Pezzo surface with Du Val
singularities, such that $K_X^2=1$ the number $\text{lct}_1 (X)$
was computed in \cite{Park2002}.

In particular, global log canonical thresholds are related to the
existence of K\"ahler-Einstein metrics on Fano varieties, as we
can see in the following result due to \cite{DK01}, \cite{Nadel},
\cite{Tian90} .

\begin{theorem}
\label{Kollar-Demailly}
 Let $X$ be an $n$-dimensional Fano variety with at most quotient singularities.
 The variety $X$ has a K\"ahler-Einstein metric if the inequality holds
 $$
\mathrm{lct}(X) > \frac{n}{n+1} \text{ .}
 $$
\end{theorem}

For the rest of this paper we are going to assume that $X$ is a
del Pezzo surface with at most Du Val singular
points\footnote[1]{All varieties are assumed to be projective,
normal and defined over $\mathbb{C}$.}. The problem of existence
of K\"ahler-Einstein metrics on smooth del Pezzo surfaces was
completely settled by Tian in \cite{Tian90}

Moreover the following is due to \cite{Cheltsov-lct}.

\begin{theorem}
\label{theorem:dP1}  Let $X$ be a smooth del Pezzo surface. Then
$$
\mathrm{lct}\big(X\big)=\left\{%
\aligned
&1/3\ \mathrm{when}\ X\cong\mathbb{F}_{1}\ \mathrm{or}\ K_{X}^{2}\in\big\{7,9\big\},\\%
&1/2\ \mathrm{when}\ X\cong\mathbb{P}^{1}\times\mathbb{P}^{1}\ \mathrm{or}\ K_{X}^{2}\in\big\{5,6\big\},\\%
&2/3\ \mathrm{when}\ K_{X}^{2}=4\ \mathrm{or}\ X\ \mathrm{is\ a\ cubic\ surface\ in}\ \mathbb{P}^{3}\ \mathrm{with\ an\ Eckardt\ point},\\%
&3/4\ \mathrm{when}\ X\ \mathrm{is\ a\ cubic\ surface\ in}\ \mathbb{P}^{3}\ \mathrm{without\ Eckardt\ points},\\%
&3/4\ \mathrm{when}\ K_{X}^{2}=2\ \mathrm{and}\ |-K_{X}|\ \mathrm{has\ a\ tacnodal\ curve},\\%
&5/6\ \mathrm{when}\ K_{X}^{2}=2\ \mathrm{and}\ |-K_{X}|\ \mathrm{has\ no\ tacnodal\ curves},\\%
&5/6\ \mathrm{when}\ K_{X}^{2}=1\ \mathrm{and}\ |-K_{X}|\ \mathrm{has\ a\ cuspidal\ curve},\\%
&1\ \mathrm{when}\ K_{X}^{2}=1\ \mathrm{and}\ |-K_{X}|\ \mathrm{has\ no\ cuspidal\ curves}.\\%
\endaligned\right.%
$$
\end{theorem}

If now $S_3 \subset \mathbb{P}^3$ is a singular cubic surface with
Du Val singularities and $S_3$ admits a K\"ahler-Einstein metric,
then according to \cite{Ding-Tian92} it can only have points of
type $\mathbb{A}_1$ or $\mathbb{A}_2$.

Moreover on a del Pezzo surface $S_2$ of degree 2 with only
$\mathbb{A}_1$ or $\mathbb{A}_2$ singularities a K\"ahler-Einstein
metric exists due to \cite{Kollar-Ghigi}. In their method they
consider $S_2$ as a double cover of $ \mathbb{P}^2$ and use a
K\"ahler-Einstein metric on $\mathbb{P}^2$ to construct a
K\"ahler-Einstein metric on $S_2$. A del Pezzo surface $S_1$ of
degree 1 can be realised as a double cover of the cone
$\mathbb{P}(1,1,2)$, however $\mathbb{P}(1,1,2)$ does not admit a
K\"ahler-Einstein metric. Thus one cannot apply the same idea to
prove existence of a K\"ahler-Einstein metric on $S_1$.

Due to \cite{Cheltsov-lct} we have the following.

\begin{theorem}
\label{dpezzoA1A2} Let $X$ be a del Pezzo surface with only Du Val
singularities only of type $\mathbb{A}_1$ or $\mathbb{A}_2$ such
that $K_{X}^{2}=1$.~Then
$$
\mathrm{lct}\big(X\big)=\left\{%
\aligned
&1\hspace{0.5cm} \mathrm{when}\ |-K_{X}|\ \mathrm{does\ not\ have\ cuspidal\ curves},\\%
&2/3\ \mathrm{when}\ |-K_{X}|\ \mathrm{has\ a\ cuspidal\ curve}\
C\ \mathrm{such\ that}\ \mathrm{Sing}(C)= \mathbb{A}_2,\\%
&3/4\ \mathrm{when}\ |-K_{X}|\ \mathrm{has\ a\ cuspidal\ curve}\
C\ \mathrm{such\ that}\ \mathrm{Sing}(C)=\mathbb{A}_1\\%
& \hspace{0.75cm} \mathrm{ and\ no\ cuspidal\ curve}\ C\ \mathrm{such\ that}\ \mathrm{Sing}(C)=\mathbb{A}_2,\\%
&5/6\ \mathrm{in\ the\ remaining\ cases}.\\%
\endaligned\right.%
$$
\end{theorem}

By Theorem~\ref{dpezzoA1A2} and Theorem~\ref{Kollar-Demailly} we
get that on every del Pezzo surface of degree 1 that has at most
ordinary double points a K\"ahler-Einstein metric exists.

The main purpose of this paper is to prove the following result.
\begin{theorem}
\label{Dimitra1} Let $X$ be a degree 1 del Pezzo surface having
the following type of Du Val singular points:
\begin{eqnarray*}
\mathbb{A}_4 \text{, } \mathbb{A}_4 + \mathbb{A}_4 \text{,
}\mathbb{A}_4 + \mathbb{A}_3 \text{, } \mathbb{A}_4 + 2 \mathbb{A}_1
\text{, } \mathbb{A}_4 + \mathbb{A}_1 \text{, }\mathbb{A}_3+ 4 \mathbb{A}_1  \text{, }\\
 \mathbb{A}_3+ 3
\mathbb{A}_1 \text{, } 2\mathbb{A}_3+ 2 \mathbb{A}_1  \text{, }
\mathbb{A}_3+ 2 \mathbb{A}_1  \text{, } \mathbb{A}_3+ \mathbb{A}_1
\text{, }2 \mathbb{A}_3 \text{, } \mathbb{A}_3 \text{ .}
\end{eqnarray*}
 Then $X$
admits a K\"ahler-Einstein metric.
\end{theorem}

Furthermore from Table 1 to Table~\ref{degree1A4A3} we give a list
of all global log canonical thresholds for
\begin{itemize}
 \item del Pezzo surfaces of degree 1 with Du Val
 singularities,
 \item del Pezzo surfaces of Picard group $\mathbb{Z}$ with Du Val
 singularities.\footnote[2]{Global log canonical thresholds of
  cubic surfaces with Du Val singularities were computed in
  \cite{Cheltsov-cubic} and  $\text{lct}(\mathbb{P}(1,1,2)) = \frac{1}{4}$ (see \cite{Fano3folds}). }

\end{itemize}

We see that Table~\ref{degree1A4A3} together with
Theorem~\ref{Kollar-Demailly} imply the existence of a
K\"ahler-Einstein metric on every del Pezzo surface of degree 1
that has the singularities mentioned in Theorem~\ref{Dimitra1}.

\begin{remark}
In  \cite{Tian90} and \cite{Shi} the invariant $\alpha_{m,2}(X)$
was introduced. One can see that
$\alpha_{m,2}(X)\geq\text{lct}(X)$ and $\alpha_{m,2}(X)$ goes to
$\text{lct}(X)$ as $m$ goes to $+\infty$. However it never reaches
$\text{lct} (X)$ if there are only finitely many
$\mathbb{Q}$-divisors $D\sim_{\mathbb{Q}} -K_{X}$, such that
$\text{lct} (X)= \text{lct} (X,D)$. From the proofs of
Lemma~\ref{A5}, Lemma~\ref{A6} and Theorem~\ref{dpezzoA1A2} it
follows that this is exactly the case when $\text{lct}(X) =
\frac{2}{3}$ and $X$ is a del Pezzo of degree 1 with Du Val
singularities. It follows from \cite{Tian90} and \cite{Shi} that a
K\"ahler-Einstein metric exists on a smooth del Pezzo surface $X$
if $\alpha_{m,2}(X)>\frac{2}{3}$. It is expected that the same is
true in case $X$ is an orbifold del Pezzo surface.
\end{remark}

Therefore we expect to have the following result.

\begin{conjecture}
Let $X$ be a degree 1 del Pezzo surface having only Du Val
singularities of type $\mathbb{A}_n$, for $n \leq 6$, then $X$
admits a K\"ahler-Einstein metric.
\end{conjecture}

Apart from their connection to the existence of  K\"ahler-Einstein
metrics global log canonical thresholds have a birational
application. For the following result we refer the reader to
\cite{Pa01} and \cite{Cheltsov-cubic}.

\begin{theorem}
\label{theorem:Park-Cheltsov} Let $V, \bar{V}$ be two varieties
and $Z$ be a smooth curve. Suppose that there is a commutative
diagram
\begin{equation}
\label{equation:commutative-diagram} \xymatrix{
&V\ar@{->}[d]_{\pi}\ar@{-->}[rr]^{\rho}&&\ \bar{V}\ar@{->}[d]^{\bar{\pi}}&\\%
&Z\ar@{=}[rr]&&Z&}
\end{equation}
such that $\pi$ and $\bar{\pi}$ are flat morphisms, and $\rho$ is
a birational map that induces an isomorphism
\begin{equation}
\label{equation:isomorphism}
\rho\big\vert_{V\setminus F}\colon V\setminus F\longrightarrow\bar{V}\setminus\bar{F},%
\end{equation}
where $F$ and $\bar{F}$ are scheme fibers of $\pi$ and $\bar{\pi}$
over a point $O\in Z$, respectively. Suppose that
\begin{itemize}
\item the~varieties $V$ and $\bar{V}$ have terminal $\mathbb{Q}$-factorial singularities,%
\item the~divisors $-K_{V}$ and $-K_{\bar{V}}$ are $\pi$-ample and $\bar{\pi}$-ample, respectively,%
\item the~fibers $F$ and $\bar{F}$ are irreducible.
\end{itemize}
Then $\rho$ is an isomorphism if one of the~following conditions
hold:
\begin{itemize}
\item the~varieties $F$ and $\bar{F}$ have log terminal singularities, and $\mathrm{lct}(F)+\mathrm{lct}(\bar{F})>1$;%
\item the~variety $F$ has log terminal singularities, and $\mathrm{lct}(F)\geqslant 1$.%
\end{itemize}
\end{theorem}

\section{Preliminaries}

Let $X$ be a del Pezzo surface with Du Val singularities. Assume
that the global log canonical threshold is
$$
\text{lct} (X)< \omega \leq 1 \text{ , }
$$
where $\omega $ is a positive rational number such that $\omega
\leq \frac{1}{K_X^2} $. This means that there is an effective
$\mathbb{Q}$-divisor $D$, with $D \sim_{\mathbb{Q}} -K_X$, such
that the log pair $(X, \lambda D)$ is not log canonical for some
rational number $\lambda < \omega \leq 1$.

  \begin{remark}
  Suppose that $Z$ is an effective $\mathbb{Q}$-divisor  on $X$ such that $(X, \lambda Z)$ is log
canonical and $Z \sim_{\mathbb{Q}} -K_X$. Then if $(X, \lambda D)$
not log canonical
  also the pair
  $$
  ( X, \frac{1}{1-\alpha } (\lambda D - \alpha \lambda Z ) )
  $$
  is not log canonical, where $\alpha \in \mathbb{Q}$ such that $ 0\leq \alpha < 1$.
 \end{remark}

Denote now by $\text{LCS} (X, \lambda D)$ the locus of log
canonical singularities, that is the set of all points of $X$
where the pair $(X, \lambda D)$ is not Kawamata log terminal. The
following Connectedness Theorem can be found in
\cite{Kollar-flips} (Ch.17).

\begin{theorem}
\label{connectedness}
 If $-(K_X+\lambda D)$ is ample, then the log canonical locus $\text{LCS} (X, \lambda D)$
is connected.
\end{theorem}

From the way log canonicity is defined for the log pair $(X,
\lambda D)$, one should understand all resolutions of
singularities of  the log pair $(X, \lambda D)$. However instead
we will use the following condition on multiplicity that follows
from \cite{Kollar-Mori98} (4.5).

\begin{remark}
\label{mult>1}
 For a smooth point $P$ of the surface $X$ the condition that the pair $(X, \lambda D)$
 is not log canonical implies that $ \text{mult}_P D > 1 \text{ .}$
\end{remark}

The following Lemma is going to be a key ingredient of the proof
of Theorem~\ref{Dimitra1}.

\begin{lemma}
\label{SingPoint} The pair $(X, \lambda D)$ is log canonical
everywhere except for a Du Val singular point $P$, where $(X,
\lambda D)$ is not log canonical.
\end{lemma}

\begin{proof}
By Theorem~\ref{connectedness} the log canonical locus $\text{LCS}
(X, \lambda D)$ is connected, since
$$-(K_X + \lambda D) \sim_{\mathbb{Q}} -(1 - \lambda) K_X$$ is ample.
Suppose now that there is an irreducible curve $C$ on the surface
$X$, such that $C \subset \text{LCS} (X, \lambda D)$. Then $C \in
\text{Supp} (D)$ and we can write $D = m C + \Omega$, where $m$ is
a rational number $m \lambda \geq 1$ and $\Omega$ is an effective
$\mathbb{Q}$-divisor such that $C \not \in \text{Supp} (\Omega)$.

But then Remark~\ref{mult>1} implies that
$$
K_X^2 = D \cdot (-K_X) = m  C \cdot (-K_X) +  \Omega \cdot (-K_X)
> m  C \cdot (-K_X) \geq \frac{1}{\lambda} \text{deg}C > \frac{1}{\omega} >K_X^2 \text{ ,}
$$
which is a contradiction. Therefore the log canonical locus is
zero-dimensional and there is a point $P \in D$ where the log pair
$(X, \lambda D)$ is not log canonical. Moreover we can assume that
$P$ is not a smooth point of $X$. This follows from
\cite{Cheltsov-lct}, where the case of smooth del Pezzo surfaces
is treated.

\end{proof}

The following theorem is known as adjunction or inversion of
adjunction ( \cite{Kollar-flips} ).
\begin{theorem}
\label{adjunction} Let $X$ be normal and $S \subset X$ be an
irreducible Cartier divisor. Let $B$ be an effective
$\mathbb{Q}$-Cartier $\mathbb{Q}$-divisor and assume that $K_X + S
+ B$ is $\mathbb{Q}$-Cartier and $S$ is Kawamata log terminal such
that $S \not \subseteq \text{Supp}B$. Then
$$(X, S + B)\ \mathrm{ is\ log\ canonical\ near}\ S \iff (S,B|S)\ \mathrm{ is\ log\ canonical } \text{.}$$
\end{theorem}

Throughout this paper we are going to refer to
Theorem~\ref{adjunction} simply as adjunction.

\section{Del Pezzo surfaces of degree 1 with exactly one Du Val singularity}

All possible combinations of Du Val singularities on a del Pezzo
surface of degree 1 are given in \cite{Urabe}.

\subsection{Del Pezzo surfaces of degree 1 with exactly one
$\mathbb{A}_3$ type singularity}

In this section we will prove the following.

\begin{lemma}
\label{A3} Let $X$ be a del Pezzo surface with exactly one Du Val
singularity of type $\mathbb{A}_3$ and $K_X^2=1$. Then the global
log canonical threshold of $X$ is
$$
\mathrm{lct} (X) = 1 \text{ .}
$$
\end{lemma}

\begin{proof}
Let $Z$ be the unique curve in the linear system $|-K_X|$ that contains $P$, where $P$ 
is the  Du Val singular point of type $\mathbb{A}_3$.
We take the minimal resolution $\pi: \tilde{X} \to X$ of $X$,
that contracts the exceptional curves $E_1, E_2, E_3$ to the singular point $P$.
The following diagram shows how the exceptional curves intersect each other.

\bigskip

$\mathbb{A}_3 \text{   }$ \xymatrix{ {\bullet}^{E_1} \ar@{-}[r] &
{\bullet}^{E_2} \ar@{-}[r] & {\bullet}^{E_3} }
\bigskip

\noindent Then
$$
\tilde{Z} \sim_{\mathbb{Q}} \pi^*(Z)-E_1-E_2-E_3\text{ .}
$$
Suppose that $\mathrm{lct}(X) < \mathrm{lct}_1(X) = 1$, then there
exists an effective $\mathbb{Q}$-divisor $D\in X$ such that $D
\sim_{\mathbb{Q}} -K_X$ and the log pair $(X, \lambda D)$ is not
log canonical, where $\lambda < 1$. According to Lemma~\ref{SingPoint}
the pair $(X, \lambda D)$ is log canonical everywhere except for 
a singular point $P\in X$, at which point $(X, \lambda D)$ is not log canonical.
Since the curve $Z$ is irreducible we may assume that the support of $D$
does not contain $Z$. Then
$$
\tilde{D} \sim_{\mathbb{Q}} \pi^*(D)-a_1E_1-a_2E_2-a_3E_3 \text{, }
$$
and from the inequalities
\begin{eqnarray*}
0 \leq \tilde{D} \cdot \tilde{Z} & = & 1 - a_1 - a_3\\
0 \leq E_1 \cdot \tilde{D} & = & 2a_1 - a_2\\
0 \leq E_2 \cdot \tilde{D} & = & 2a_2 - a_1 - a_3\\
0 \leq E_3 \cdot \tilde{D} & = & 2a_3 - a_2
\end{eqnarray*}
we see that
$$
2a_1 \geq a_2 \text{ , } \frac{3}{2}a_2 \geq a_3 \text{  and  }
2a_3 \geq a_2 \text{ , } \frac{3}{2} a_2 \geq a_1 \text{ .}
$$
Moreover we get the following upper bounds $a_1 \leq \frac{3}{4}
\text{, } a_2 \leq 1 \text{, } a_3 \leq \frac{3}{4} \text{ .}$ The
equivalence
$$
K_{\tilde{X}} + \lambda \tilde{D} + \lambda a_1E_1 + \lambda
a_2E_2 + \lambda a_3E_3 \sim_{\mathbb{Q}} \pi^*(K_X + \lambda D)
$$
implies that there is a point $Q\in E_1\cup E_2\cup E_3$ such that
the pair $K_{\tilde{X}} + \lambda \tilde{D} + \lambda a_1E_1 +
\lambda a_2E_2 + \lambda a_3E_3$ is not log canonical at $Q$.

\begin{itemize}
\item If the point $Q \in E_1$ and $Q\not \in E_2$ then
$$
K_{\tilde{X}} + \lambda \tilde{D} + \lambda a_1E_1
$$
is not log canonical at the point $Q$ and so is the pair
$$
K_{\tilde{X}} + \lambda \tilde{D} + E_1 \text{ .}
$$
By adjunction $(E_1, \lambda \tilde{D}|_{E_1})$ is not log
canonical at $Q$ and
$$\frac{4}{3} a_1 \geq 2a_1 - \frac{2}{3} a_1 \geq 2a_1 - a_2 = \tilde{D} \cdot E_1 \geq
 \text{mult}_Q\Bigl(\tilde{D} \cdot E_1 \Bigr)  > \frac{1}{\lambda} > 1 \text{ ,}$$
 which is false.

\item If $Q\in E_2$ but $Q\not \in E_1 \cup E_3$ then
$$
K_{\tilde{X}} + \lambda \tilde{D} + \lambda a_2E_2
$$
is not log canonical at the point $Q$ and so is the pair
$$
K_{\tilde{X}} + \lambda \tilde{D} + E_2 \text{, since } \lambda
a_2 \leq 1 \text{ .}
$$
By adjunction $(E_2, \lambda \tilde{D}|_{E_2})$ is not log
canonical at $Q$ and
$$2a_2 - \frac{a_2}{2}- \frac{a_2}{2} \geq 2 a_2 - a_1 - a_3 = \tilde{D} \cdot E_2 \geq
 \text{mult}_Q\Bigl(\tilde{D} \cdot E_2 \Bigr) > \frac{1}{\lambda}>1 \text{ ,}$$
which is false.

\item If $Q\in E_1 \cap E_2$ then the log pair
$$
K_{\tilde{X}} + \lambda \tilde{D} + \lambda a_1E_1 + \lambda
a_2E_2
$$
is not log canonical at the point $Q$ and so are the log pairs
$$
K_{\tilde{X}} + \lambda \tilde{D} + E_1 + \lambda a_2E_2 \text{
and } K_{\tilde{X}} + \lambda \tilde{D} + \lambda a_1E_1 + E_2
\text{ .}
$$
By adjunction it follows that
$$
2a_1 - a_2 = \tilde{D} \cdot E_1 \geq
\text{mult}_Q\Bigl(\tilde{D}|_{E_1} \Bigr) =
 \text{mult}_Q\Bigl(\tilde{D} \cdot E_1 \Bigr) > 1 - a_2
$$
and
$$
2a_2 - a_1 - a_3 = \tilde{D} \cdot E_2 \geq
\text{mult}_Q\Bigl(\tilde{D}|_{E_2} \Bigr) =
 \text{mult}_Q\Bigl(\tilde{D} \cdot E_2 \Bigr) > 1 - a_1 \text{ , }
$$
which implies that
$$
2a_2 - \frac{a_2}{2} \geq 2a_2 - a_3 > 1 \Rightarrow  a_2 >
\frac{2}{3} \text{ .}
$$

\end{itemize}

\bigskip

Consider now the blow-up $\sigma_1 : \tilde{X}_1 \to \tilde{X}$ of
the surface $\tilde{X}$ at the point $Q$ that contracts the
$(-1)$-curve $F_1$ to the point $Q$. Then for the strict
transforms of the exceptional divisors $E_1, E_2, E_3$ we have
\begin{eqnarray*}
\tilde{E_1} & \sim_{\mathbb{Q}} & \sigma_1 ^*(E_1) - F_1\\
\tilde{E_2} & \sim_{\mathbb{Q}} & \sigma_1 ^*(E_2) - F_1\\
\tilde{E_3} & \sim_{\mathbb{Q}} & \sigma_1 ^*(E_3) \text{ .}
\end{eqnarray*}
Let now
$$\pi_1 : \tilde{X}_1 \stackrel{\sigma_1}{\to} \tilde{X} \stackrel{\pi}{\to} X$$ be the
composition $\pi_1 = \pi \circ \sigma_1$. We have
\begin{eqnarray*}
K_{\tilde{X}_1} & \sim_{\mathbb{Q}} & \sigma_1^*(K_{\tilde{X}}) +
F_1
  \sim_{\mathbb{Q}} \sigma_1^*(\pi^*(K_X)) + F_1
  \sim_{\mathbb{Q}} \pi_1^*(K_X) + F_1
\end{eqnarray*}
and
\begin{eqnarray*}
\tilde{D}_1 & \sim_{\mathbb{Q}} & \sigma_1^*(\tilde{D}) - m_1 F1\\
 & \sim_{\mathbb{Q}} & \sigma_1^*(\pi^*(D)-a_1E_1-a_2E_2-a_3E_3) - m_1 F_1\\
 & \sim_{\mathbb{Q}} & \pi_1^*(D) - a_1 \tilde{E_1} - a_2 \tilde{E_2} - a_3 \tilde{E_3} - (a_1 + a_2 +
m_1) F_1 \text{ ,}
\end{eqnarray*}
where $m_1 = \text{mult}_{Q} \tilde{D}$. Also the strict transform
of the anticanonical curve $Z$ is
\begin{eqnarray*}
\tilde{Z}_1 & \sim_{\mathbb{Q}} & \sigma_1^*(\tilde{Z})\\
 & \sim_{\mathbb{Q}} & \sigma_1^*(\pi^*(Z) - E_1 - E_2 - E_3)\\
 & \sim_{\mathbb{Q}} & \pi_1^*(Z) - \tilde{E_1} - \tilde{E_2} - \tilde{E_3} - 2F_1 \text{ .}
\end{eqnarray*}
From the inequalities
\begin{eqnarray*}
0  \leq  \tilde{D}_1  \cdot \tilde{Z}_1  & = & 1 - a_1 - a_3\\
0  \leq  \tilde{E_1} \cdot \tilde{D}_1 & = & 2a_1 - a_2 - m_1\\
0  \leq  \tilde{E_2} \cdot \tilde{D}_1  & = & 2a_2 - a_1 - a_3 - m_1\\
0  \leq  \tilde{E_3} \cdot \tilde{D}_1  & = & 2a_3 - a_2\\
0  \leq  F_1 \cdot \tilde{D}_1  & = & m_1
\end{eqnarray*}
we get that $m_1 \leq \frac{1}{2}$. The equivalence
$$
K_{\tilde{X}_1 } + \lambda \tilde{D}_1  + \lambda a_1\tilde{E_1} +
\lambda a_2\tilde{E_2} + \lambda a_3\tilde{E_3} + (\lambda(a_1 +
a_2 + m_1) -1)F_1 \sim_{\mathbb{Q}} \pi_1^*(K_X + \lambda D)
$$
implies that there is a point $Q_1 \in F_1$ such that the pair
$$
K_{\tilde{X}_1 } + \lambda \tilde{D}_1  + \lambda a_1\tilde{E_1} +
\lambda a_2\tilde{E_2} + (\lambda(a_1 + a_2 + m_1) -1)F_1
$$
is not log canonical at $Q_1$.
\begin{itemize}
\item[(i)] Suppose $Q_1 \in \tilde{E_1} \cap F_1$, then the log
pair
$$
K_{\tilde{X}_1 } + \lambda \tilde{D}_1  + \lambda a_1\tilde{E_1} +
(\lambda(a_1 + a_2 + m_1) -1)F_1
$$
is not log canonical at $Q_1$ and so are the pairs
$$
K_{\tilde{X}_1 } + \lambda \tilde{D}_1  + \tilde{E_1} +
(\lambda(a_1 + a_2 + m_1) -1)F_1
$$
and
$$
K_{\tilde{X}_1 } + \lambda \tilde{D}_1  + \lambda a_1\tilde{E_1} +
F_1 \text{ .}
$$
By adjunction it follows that
$$
2a_1 - a_2 - m_1 = \tilde{D}_1  \cdot \tilde{E_1} \geq
 \text{mult}_{Q_1}\Bigl(\tilde{D}_1  \cdot \tilde{E_1} \Bigr) > 1 - a_1 - a_2 - m_1 + 1
\text{ ,}$$
 and
 $$
m_1 = \tilde{D}_1  \cdot F_1 \geq
 \text{mult}_{Q_1}\Bigl(\tilde{D}_1  \cdot F_1 \Bigr) > 1 - a_1\text{ ,}$$

 which is false.

\item[(ii)] Suppose $Q_1 \in F_1 \backslash (\tilde{E_1} \cup
\tilde{E_2})$ then the log pair
$$K_{\tilde{X}_1 } + \lambda \tilde{D}_1  + (\lambda(a_1 + a_2 + m_1) -1)F_1$$ is not log canonical at $Q_1$
and so is the pair
$$K_{\tilde{X}_1 } + \lambda \tilde{D}_1  + F_1 \text{  , since  }  0 \leq \lambda(a_1 + a_2 + m_1) -1
\leq 1$$ By adjunction it follows that
$$
m_1 = \tilde{D}_1  \cdot F_1 \geq
 \text{mult}_{Q_1}\Bigl(\tilde{D}_1  \cdot F_1 \Bigr) > 1 \text{ ,}$$ which is impossible
since $m_1 \leq \frac{1}{2}$.

\item[(iii)] Suppose $Q_1 \in \tilde{E_2} \cap F_1$ then the log
pair
$$K_{\tilde{X}_1 } + \lambda \tilde{D}_1  + \lambda a_2\tilde{E_2} + (\lambda(a_1 + a_2 + m_1) -1)F_1 $$ is not
log terminal at $Q_1$ and so are the pairs
$$
K_{\tilde{X}_1 } + \lambda \tilde{D}_1  + \tilde{E_2} +
(\lambda(a_1 + a_2 + m_1) - 1)F_1
$$
and
$$
K_{\tilde{X}_1 } + \lambda \tilde{D}_1  + \lambda a_2\tilde{E_2} +
F_1 \text{ .}
$$
By adjunction it follows that
$$
2a_2 - a_1 - a_3 - m_1 = \tilde{D}_1  \cdot \tilde{E_2} \geq
 \text{mult}_{Q_1}\Bigl(\tilde{D}_1  \cdot \tilde{E_2} \Bigr) > 1 - a_1 - a_2 - m_1 + 1$$
 and
   $$
m_1 = \tilde{D}_1  \cdot F_1 \geq
 \text{mult}_{Q_1}\Bigl(\tilde{D}_1  \cdot F_1 \Bigr) > 1 - a_2\text{ , }$$
 which imply that
 $$
 3a_2 - \frac{a_2}{2} \geq 3a_2 - a_3 > 2 \Rightarrow a_2 > \frac{4}{5} \text{ .}
 $$

  \end{itemize}

  Consider now the blow-up $\sigma_2: \tilde{X}_2  \to \tilde{X}_1 $
of the surface $\tilde{X}_1$ at the point $Q_1$ that contracts the
$(-1)$-curve $F_2$ to the point $Q_1$. We then have

  \begin{eqnarray*}
   K_{\tilde{X}_2 } & \sim_{\mathbb{Q}} & \pi_2^*(K_X) + \tilde{F}_1 + 2F_2\\
  \tilde{D}_2  & \sim_{\mathbb{Q}} & \pi_2^*(D) -a_1 \tilde{E_1} -a_2 \tilde{E_2} - a_3 \tilde{E_3} -
(a_1 + a_2 + m_1) \tilde{F}_1-(a_1+2a_2+m_1+m_2) F_2\\
  \tilde{Z}_2  & \sim_{\mathbb{Q}} & \pi_2^*(Z) - \tilde{E_1} - \tilde{E_2} - \tilde{E_3} - 2\tilde{F}_1
- 3F_2 \text{ ,}
  \end{eqnarray*}
and
\begin{eqnarray*}
 0 \leq \tilde{D}_2  \cdot \tilde{Z}_2  & = & 1 - a_1 - a_3\\
0 \leq \tilde{E}_1 \cdot \tilde{D}_2  & = & 2a_1 - a_2 - m_1\\
0 \leq \tilde{E}_2 \cdot \tilde{D}_2  & = & 2a_2 - a_1 - a_3 - m_1 - m_2\\
0 \leq \tilde{E}_1 \cdot \tilde{D}_2  & = & 2a_3 - a_2\\
0 \leq \tilde{F}_1 \cdot \tilde{D}_2  & = & m_1 - m_2\\
0 \leq F_2 \cdot \tilde{D}_2  & = & m_2 \text{ ,}
\end{eqnarray*}
where $m_2 = \text{mult}_{Q_1} \tilde{D}_2 $.

Because of the equivalence
  \begin{eqnarray*}
  \lefteqn{\pi_2^*(K_X + \lambda D)  \sim_{\mathbb{Q}} }\\
 & & K_{\tilde{X}_2 } + \lambda \tilde{D}_2  + \lambda a_1 \tilde{E}_1
+ \lambda a_2 \tilde{E}_2 + \lambda a_3 \tilde{E}_3 + (\lambda(a_1
+ a_2 + m_1) -1)
 \tilde{F}_1 + (\lambda(a_1 + 2a_2 + m_1 + m_2) - 2) F_2
  \end{eqnarray*}

there is a point $Q_2 \in F_2$ such that the pair
$$
K_{ \tilde{X}_2} + \lambda \tilde{D}_2  + \lambda a_1 \tilde{E}_1
+ \lambda a_2 \tilde{E}_2 + \lambda a_3 \tilde{E}_3 + (\lambda(a_1
+ a_2 + m_1) -1) \tilde{F}_1 + (\lambda(a_1 + 2a_2 + m_1 + m_2) -
2) F_2
$$
is not log canonical at $Q_2$.
\begin{itemize}
\item Suppose $Q_2 \in F_2 \cap \tilde{F}_1$, then the log pair
$$
 \left( \tilde{X}_2 , \lambda \tilde{D}_2  + (\lambda(a_1 + a_2 + m_1) -1)
 \tilde{F}_1 + (\lambda(a_1 + 2a_2 + m_1 + m_2) - 2) F_2 \right)
$$
 is not log canonical at $Q_2$ and so are the log pairs
 $$
K_{\tilde{X}_2 } + \lambda \tilde{D}_2  + \tilde{F}_1 +
(\lambda(a_1 + 2a_2 + m_1 + m_2) - 2) F_2
 $$
and
$$
K_{\tilde{X}_2 } + \lambda \tilde{D}_2  + (\lambda(a_1 + a_2 +
m_1) -1) \tilde{F}_1 +  F_2 \text{ .}
$$
By adjunction it follows that
 $$
 m_1 - m_2 = \tilde{D}_2 \cdot \tilde{F}_1 \geq \text{mult}_{Q_2} \Bigl(\tilde{D}_2 \cdot
\tilde{F}_1 \Bigr) > 1 - ( a_1 + 2a_2 + m_1 + m_2 - 2)
 $$
 and
 $$
 m_2 = \tilde{D}_2  \cdot F_2 \geq \text{mult}_{Q_2} \Bigl( \tilde{D}_2 \cdot F_2 \Bigr)
1 - ( a_1 + a_2 + m_1 - 1) \text{ ,}
$$ which is false.

\item Suppose $O \in F_2 \backslash (\tilde{F}_1 \cup
\tilde{E_2})$, then the log pair
$$
 K_{\tilde{X}_2 } + \lambda \tilde{D}_2  + (\lambda(a_1 + 2a_2 + m_1 + m_2) - 2) F_2
$$
is not log canonical at $Q_2$ and so is the log pair
$$
 K_{\tilde{X}_2 } + \lambda  \tilde{D}_2  +  F_2   \text{   , since   } \lambda(a_1 + 2a_2 + m_1 + m_2) -
2 \leq 1 \text{ .}
$$
By adjunction it follows that
$$
 m_2 = \tilde{D}_2 \cdot F_2 \geq \text{mult}_{Q_2} \Bigl(\tilde{D}_2 \cdot F_2 \Bigr) >
1  \text{ ,}
 $$
 which is false.

\item Suppose $Q_2 \in F_2 \cap \tilde{E}_2$, then the log pair
$$
 K_{\tilde{X}_2 } + \lambda \tilde{D}_2  + (\lambda(a_1 + 2a_2 + m_1 + m_2) - 2) F_2 + \lambda a_2 \tilde{E}_2
$$
 is not log canonical at $Q_2$ and so are the log pairs
 $$
 K_{\tilde{X}_2 } + \lambda \tilde{D}_2  + (\lambda(a_1 + 2a_2 + m_1 + m_2) - 2) F_2 +  \tilde{E}_2
 $$
and
$$
K_{\tilde{X}_2 } + \lambda \tilde{D}_2  + F_2 + \lambda a_2
\tilde{E}_2 \text{ .}
$$
By adjunction it follows that
 $$
 2a_2 - a_1 - a_3 - m_1 - m_2 = \tilde{D}_2 \cdot \tilde{E}_2 \geq \text{mult}_{Q_2}
\Bigl(\tilde{D}_2 \cdot \tilde{E}_2\Bigr) > 1 - ( a_1 + 2a_2 + m_1
+ m_2 - 2)
 $$
 and
 $$
 m_2 = \tilde{D}_2 \cdot F_2 \geq \text{mult}_{Q_2} \Bigl(\tilde{D}_2 \cdot F_2 \Bigr) >
1 - a_2 \text{, }
 $$
which implies that
$$
4a_2 - \frac{a_2}{2} \geq 4a_2 - a_3 > 3 \Rightarrow a_2 >
\frac{6}{7} \text{ .}
$$
\end{itemize}

Consider now the blow-up $\sigma_k: \tilde{X}_k  \to
\tilde{X}_{k-1} $ of the surface $\tilde{X}_{k-1}$ at the point
$Q_{k-1}$ that contracts the $(-1)$-curve $F_k$ to the point
$Q_{k-1}$. We then have

\begin{eqnarray*}
K_{\tilde{X}_k } & \sim_{\mathbb{Q}} & \pi_k^*(K_X) + \tilde{F}_1
+ 2\tilde{F}_2 + 3 \tilde{F}_3
+...+ (k-1) \tilde{F}_{k-1}+ k F_k\\
\tilde{D}_k  & \sim_{\mathbb{Q}} & \pi_k^*(D) -a_1 \tilde{E_1}
-a_2 \tilde{E_2} - a_3 \tilde{E_3} - (a_1
+ a_2 + m_1) \tilde{F}_1-(a_1+2a_2+m_1+m_2) \tilde{F}_2 -...\\
& & - \left( a_1 + (k-1) a_2 + m_1 + m_2+ ..+ m_{k-1} \right)
\tilde{F}_{k-1}
- \left( a_1 + k a_2 + m_1 + m_2+ ..+ m_k \right) F_k  \\
\tilde{Z}_k  & \sim_{\mathbb{Q}} & \pi_k^*(Z) - \tilde{E_1} -
\tilde{E_2} - \tilde{E_3} - 2\tilde{F}_1 - 3\tilde{F}_2 -...- k
\tilde{F}_{k-1} - (k+1) F_k \text{ ,}
\end{eqnarray*}
and
\begin{eqnarray*}
 0 \leq \tilde{D}_k  \cdot \tilde{Z}_k  & = & 1 - a_1 - a_3\\
0 \leq \tilde{E}_1 \cdot \tilde{D}_k  & = & 2a_1 - a_2 - m_1\\
0 \leq \tilde{E}_2 \cdot \tilde{D}_k  & = & 2a_2 - a_1 - a_3 - m_1 - m_2 -...- m_k\\
0 \leq \tilde{E}_1 \cdot \tilde{D}_k  & = & 2a_3 - a_2\\
0 \leq \tilde{F}_1 \cdot \tilde{D}_k  & = & m_1 - m_2\\
0 \leq \tilde{F}_2 \cdot \tilde{D}_k  & = & m_2 - m_3\\
 & . & \\
  & . & \\
   & . & \\
0 \leq \tilde{F}_{k-1} \cdot \tilde{D}_k  & = & m_{k-1} - m_k\\
0 \leq F_k \cdot \tilde{D}_k  & = & m_k \text{ ,}
\end{eqnarray*}
where $m_i = \text{mult}_{Q_{i-1}} \tilde{D}_i \text{ , for  }
i=1,...,k$.

Because of the equivalence
\begin{eqnarray*}
\lefteqn{\pi_k^*(K_X + \lambda D) \sim_{\mathbb{Q}} }\\
 & & K_{\tilde{X}_k } + \lambda \tilde{D}_k  + \lambda a_1 \tilde{E}_1
+ \lambda a_2 \tilde{E}_2 + \lambda a_3 \tilde{E}_3 + (\lambda(a_1
+ a_2 + m_1) -1)
 \tilde{F}_1  \\
& &  + (\lambda(a_1 + 2a_2 + m_1 + m_2) - 2) \tilde{F}_2 +.....\\
& & + \left( \lambda ( a_1 + (k-1) a_2 + m_1 + m_2+...+ m_{k-1}) -
(k-1) \right) \tilde{F}_{k-1} +\\
& &  + \left( \lambda (a_1 + k a_2 + m_1 + m_2+...+ m_k ) -k
\right) F_k
\end{eqnarray*}

there is a point $Q_k \in F_k$ such that the pair
\begin{eqnarray*}
 K_{\tilde{X}_k } + \lambda \tilde{D}_k + \lambda a_2 \tilde{E}_2+ \left( \lambda (a_1 + (k-1) a_2 + m_1 +
m_2+...+ m_{k-1}) - (k-1) \right) \tilde{F}_{k-1}\\
 + \left( \lambda ( a_1 + k a_2 + m_1 + m_2+...+ m_k) -k \right) F_k
\end{eqnarray*}
is not log canonical at $Q_k$.
\begin{itemize}
\item Suppose $Q_k \in F_k \cap \tilde{F}_{k-1}$, then the log
pair
\begin{eqnarray*}
 K_{\tilde{X}_k } + \lambda \tilde{D}_k + \left( \lambda (a_1 + (k-1) a_2 + m_1 + m_2+...+ m_{k-1}) - (k-1)
\right) \tilde{F}_{k-1}\\
  + \left( \lambda ( a_1 + k a_2 + m_1 + m_2+...+ m_k) -k \right) F_k
\end{eqnarray*}
 is not log canonical at $Q_2$ and so are the log pairs
 $$
K_{\tilde{X}_k } + \lambda \tilde{D}_k  + \tilde{F}_{k-1} + \left(
\lambda ( a_1 + k a_2 + m_1 + m_2+...+ m_k) -k \right) F_k
 $$
and
$$
K_{\tilde{X}_k } + \lambda \tilde{D}_k  + \left( \lambda (a_1 +
(k-1) a_2 + m_1 + m_2+...+ m_{k-1}) - (k-1) \right)
\tilde{F}_{k-1} +  F_k \text{ .}
$$
By adjunction it follows that
 $$
 m_{k-1} - m_k = \tilde{D}_k \cdot \tilde{F}_{k-1} \geq \text{mult}_{Q_k}
\Bigl(\tilde{D}_k \cdot \tilde{F}_{k-1} \Bigr) > 1 - ( a_1 + ka_2
+ m_1 + m_2 +...+m_k - k) \text{, }
 $$
  which is a contradiction. Indeed from the inequality above we have that
$$
 a_1 + ka_2 + m_1 + m_2 +...+m_{k-2} + 2 m_{k-1} > k+1
 $$
 but since $m_1 \geq m_2 \geq ...\geq m_k $ we get that
 $$
 a_1 + ka_2 + k m_1  > k+1 \text{ .}
 $$
 However the inequality $0 \leq \tilde{E}_1 \cdot \tilde{D}_k   =  2a_1 - a_2 - m_1$
 finally gives us
 $$(2k+1) a_1 > k+1 \Rightarrow a_3 < \frac{k}{2k+1} \Rightarrow a_2 \leq 2a_3 <
\frac{2k}{2k+1} \text{ .}$$

\medskip

\item Suppose $Q_k \in F_k \backslash (\tilde{F}_{k-1} \cup
\tilde{E_2})$, then the log pair
$$
 K_{\tilde{X}_k } + \lambda \tilde{D}_k  + \left( \lambda (a_1 + k a_2 + m_1 + m_2+...+ m_k) -k \right) F_k
$$
is not log canonical at $Q_k$ and so is the log pair
$$
 K_{\tilde{X}_k } + \lambda \tilde{D}_k  +  F_k   \text{   , since   } \left( \lambda (a_1 + k a_2 + m_1
+ m_2+...+ m_k) -k \right) \leq 1 \text{ .}
$$
By adjunction it follows that
$$
 m_k = \tilde{D}_k \cdot F_k \geq \text{mult}_{Q_k} \Bigl(\tilde{D}_k \cdot F_k \Bigr) >
1  \text{ ,}
 $$
 which is false, since $ \frac{1}{2} \geq m_1 \geq m_2 \geq ... \geq m_k$.

\medskip

\item Suppose $Q_k \in F_k \cap \tilde{E}_2$, then the log pair
$$
 K_{\tilde{X}_k } + \lambda \tilde{D}_k  + \left( \lambda (a_1 + k a_2 + m_1 + m_2+...+ m_k) -k \right)
F_k + \lambda a_2 \tilde{E}_2
$$
 is not log canonical at $Q_k$ and so are the log pairs
 $$
 K_{\tilde{X}_k } + \lambda \tilde{D}_k  + \left( \lambda (a_1 + k a_2 + m_1 + m_2+...+ m_k) -k \right)
F_k +  \tilde{E}_2
 $$
and
$$
K_{\tilde{X}_k } + \lambda \tilde{D}_k  + F_k + \lambda a_2
\tilde{E}_2 \text{ .}
$$
By adjunction it follows that
 $$
 2a_2 - a_1 - a_3 - m_1 - m_2-...-m_k = \tilde{D}_k \cdot \tilde{E}_2 \geq
\text{mult}_{Q_k} \Bigl(\tilde{D}_k \cdot \tilde{E}_2\Bigr) > 1 -
\left( a_1 + k a_2 + m_1 + m_2+...+ m_k -k \right)
 $$
 and
 $$
 m_k = \tilde{D}_k \cdot F_k \geq \text{mult}_{Q_k} \Bigl(\tilde{D}_k \cdot F_k \Bigr) >
1 - a_2 \text{, }
 $$
which implies that
$$
(k+2) a_2 - \frac{a_2}{2} \geq (k+2) a_2 - a_3 > k+1 \Rightarrow
a_2 > \frac{2k+2}{2k+3} \text{ .}
$$
\end{itemize}

\begin{remark}
 It remains to be shown that after the $k$-th blow up we have
 $$
 \left( \lambda (a_1 + k a_2 + m_1 + m_2+...+ m_k) -k \right) \leq 1 \text{, }
 $$
 and for this it is enough to show that
$$
 \left( a_1 + k a_2 + m_1 + m_2+...+ m_k -k \right) \leq 1 \text{ . }
 $$
 Suppose that we have blown up $k-1$ times, then $a_2 > \frac{2k}{2k+1}$.
 Let us assume on the contrary that
 \begin{eqnarray*}
  a_1 + k a_2 + m_1 + m_2+...+ m_k -k  > 1 & \Rightarrow & a_1 + k a_2 + m_1 + m_2+...+
m_k > k + 1\\
  a_1 + 2k a_1  \geq a_1 + k a_2 + k m_1  > k + 1 & \Rightarrow &  a_1 >
\frac{k+1}{2k+1}\\
  a_3 \leq 1 - a_1 < \frac{k}{2k+1} & \Rightarrow & a_2 \leq 2 a_3 < \frac{2k}{2k+1}
\text{, }
 \end{eqnarray*}
 which is a contradiction.

\end{remark}

\end{proof}

\subsection{Del Pezzo surfaces of degree 1 with  one
$\mathbb{A}_4$ type singularity}\footnote{I am grateful to Jihun Park for letting me know about a mistake on the upper bound of
$\mathrm{lct} (X)$.} In this section we will prove the following.

\begin{lemma}
\label{A4} Let $X$ be a del Pezzo surface with one Du Val
singularity of type $\mathbb{A}_4$ and $K_X^2=1$. Then the global
log canonical threshold of $X$ is
$$
\mathrm{lct} (X) = \frac{4}{5} \text{ .}
$$
\end{lemma}

\begin{proof}

Let $Z$ be the unique curve in the linear system $|-K_X|$ that contains $P$, where $P$ 
is the  Du Val singular point of type $\mathbb{A}_4$.
We take the minimal resolution $\pi: \tilde{X} \to X$ of $X$,
that contracts the exceptional curves $E_1, E_2, E_3, E_4$ to the singular point $P$.
The following diagram shows how the exceptional curves intersect each other.

\bigskip

$\mathbb{A}_4$. \xymatrix{ {\bullet}^{E_1} \ar@{-}[r] &
{\bullet}^{E_2} \ar@{-}[r] & {\bullet}^{E_3} \ar@{-}[r] &
{\bullet}^{E_4}}
\bigskip

\noindent Then
$$
\tilde{Z} \sim_{\mathbb{Q}} \pi^*(Z)-E_1-E_2-E_3-E_4\text{ .}
$$

\noindent Furthermore there exists a unique smooth irreducible element $C$
of the linear system $|- 2 K_X|$, such that $E_2 \cap E_3 \in C$. For the
strict transform of the irreducible curve $C$ we have
$$
\tilde{C} + E_1 + 2E_2 + 2E_3 + E_4 \in |-2 K_{\tilde{X}}| \text{
.}
$$
We can assume that $C \not \in \text{Supp} D$ and then
\begin{eqnarray}
\label{mult>2}
\text{mult}_{Q} \tilde{D} \leq \tilde{C} \cdot \tilde{D} = 2 - a_2
- a_3  \Rightarrow \text{mult}_{Q} \tilde{D} + a_2 + a_3  \leq 2
\text{ .}
\end{eqnarray}
Since $\frac{1}{2} C \sim_{\mathbb{Q}} -K_X$, we have that
$$
\mathrm{lct}(X) \leq \mathrm{lct}(X,\frac{1}{2}C) = \frac{4}{5}
$$

Suppose that $\mathrm{lct}(X) < \frac{4}{5}$, then there
exists an effective $\mathbb{Q}$-divisor $D\in X$ such that $D
\sim_{\mathbb{Q}} -K_X$ and the log pair $(X, \lambda D)$ is not
log canonical, where $\lambda < \frac{4}{5}$. According to Lemma~\ref{SingPoint}
the pair $(X, \lambda D)$ is log canonical everywhere except for 
a singular point $P\in X$, at which point $(X, \lambda D)$ is not log canonical.
Since the curve $Z$ is irreducible we may assume that the support of $D$
does not contain $Z$, thus $0 \leq \tilde{D} \cdot \tilde{Z}$. 
Then
$$
\tilde{D} \sim_{\mathbb{Q}} \pi_1^*(D) - a_1E_1 - a_2E_2 - a_3E_3
- a_4E_4 \text{, }
$$
and from the inequalities
\begin{eqnarray*}
0 \leq \tilde{D} \cdot \tilde{Z} & = & 1 - a_1 - a_4\\
0 \leq E_1 \cdot \tilde{D} & = & 2a_1 - a_2\\
0 \leq E_2 \cdot \tilde{D} & = & 2a_2 - a_1 - a_3\\
0 \leq E_3 \cdot \tilde{D} & = & 2a_3 - a_2 - a_4\\
0 \leq E_4 \cdot \tilde{D} & = & 2a_4 - a_3
\end{eqnarray*}
we get the following upper bounds $a_1 \leq \frac{4}{5} \text{, } a_2 \leq \frac{6}{5}
\text{, } a_3 \leq \frac{6}{5} \text{, } a_4 \leq \frac{4}{5}$.

The equivalence
$$
K_{\tilde{X}} + \lambda \tilde{D} + \lambda a_1E_1 + \lambda
a_2E_2 + \lambda a_3E_3 + \lambda a_4E_4 \sim_{\mathbb{Q}}
\pi_1^*(K_X + \lambda D)
$$
implies that there is a point $Q\in E_1\cup E_2\cup E_3\cup E_4$
such that the pair $K_{\tilde{X}} + \lambda \tilde{D} + \lambda
a_1E_1 + a_2E_2 + \lambda a_3E_3 + \lambda a_4E_4$ is not log
canonical at $Q$.

\begin{itemize}
\item If the point $Q \in E_1$ and $Q\not \in E_2$ then
$$
K_{\tilde{X}} + \lambda \tilde{D} + \lambda a_1E_1
$$
is not log canonical at the point $Q$ and so is the pair
$$
K_{\tilde{X}} + \lambda \tilde{D} + E_1 \text{ , since  } \lambda
a_1 \leq 1 \text{ .}
$$
By adjunction $(E_1, \lambda \tilde{D}|_{E_1})$ is not log
canonical at $Q$ and
$$1 \geq \frac{5}{4} a_1 \geq 2a_1 - \frac{3}{4} a_1 \geq 2a_1 - a_2 = \tilde{D} \cdot
E_1 \geq
 \text{mult}_Q\Bigl(\tilde{D} \cdot E_1 \Bigr)  > \frac{1}{\lambda} > \frac{6}{5} \text{ ,}$$
which is false.

\item If $Q\in E_2$ but $Q\not \in E_1 \cup E_3$ then
$$
K_{\tilde{X}} + \lambda \tilde{D} + \lambda a_2E_2
$$
is not log canonical at the point $Q$ and so is the pair
$$
K_{\tilde{X}} + \lambda \tilde{D} + E_2 \text{ , since  } \lambda
a_2 \leq 1 \text{ .}
$$
By adjunction $(E_2, \lambda \tilde{D}|_{E_2})$ is not log
canonical at $Q$ and
$$1 \geq \frac{5}{6} a_2 \geq 2a_2 - \frac{a_2}{2} - \frac{2}{3} a_2 \geq 2a_2 - a_1 -
a_3 = \tilde{D} \cdot E_2 \geq
 \text{mult}_Q\Bigl(\tilde{D} \cdot E_2 \Bigr) > \frac{1}{\lambda} > \frac{6}{5} \text{ ,}$$
which is false.

\item If $Q\in E_3$ but $Q\not \in E_2 \cup E_4$ then
$$
K_{\tilde{X}} + \lambda \tilde{D} + \lambda a_3E_3
$$
is not log canonical at the point $Q$ and so is the pair
$$
K_{\tilde{X}} + \lambda \tilde{D} + E_3 \text{ , since  } \lambda
a_3 \leq 1 \text{ .}
$$
By adjunction $(E_3, \lambda \tilde{D}|_{E_3})$ is not log
canonical at $Q$ and
$$1 \geq \frac{5}{6} a_3 \geq 2a_3 - \frac{2}{3} a_3 - \frac{a_3}{2} \geq 2a_3 - a_2 -
a_4 = \tilde{D} \cdot E_3 \geq
 \text{mult}_Q\Bigl(\tilde{D} \cdot E_3 \Bigr) > \frac{1}{\lambda} > \frac{6}{5} \text{ ,}$$
 which is false.

\item If $Q\in E_4$ but $Q\not \in E_3$ then
$$
K_{\tilde{X}} + \lambda \tilde{D} + \lambda a_4E_4
$$
is not log canonical at the point $Q$ and so is the pair
$$
K_{\tilde{X}} + \lambda \tilde{D} + E_4 \text{ , since  } \lambda
a_4 \leq 1 \text{ .}
$$
By adjunction $(E_4, \lambda \tilde{D}|_{E_4})$ is not log
canonical at $Q$ and
$$1 \geq \frac{5}{4} a_4 \geq 2a_4 - \frac{3}{4} a_4 \geq 2a_4 - a_3 = \tilde{D} \cdot
E_4 \geq
 \text{mult}_Q\Bigl(\tilde{D} \cdot E_4 \Bigr) > \frac{1}{\lambda} > \frac{6}{5} \text{ ,}$$
 which is false.

\item If $Q\in E_1 \cap E_2$ then the log pair
$$
K_{\tilde{X}} + \lambda \tilde{D} + \lambda a_1E_1 + \lambda
a_2E_2
$$
is not log canonical at the point $Q$ and so are the log pairs
$$
K_{\tilde{X}} + \lambda \tilde{D} + E_1 + \lambda a_2E_2 \text{
and } K_{\tilde{X}} + \lambda \tilde{D} + \lambda a_1E_1 + E_2
\text{ .}
$$
By adjunction it follows that
$$
2a_2 - a_1 - \frac{2}{3} a_2 \geq 2a_2 - a_1 - a_3 = \tilde{D}
\cdot E_2 \geq
 \text{mult}_Q\Bigl(\tilde{D} \cdot E_2 \Bigr) > \frac{1}{\lambda}- a_1 > \frac{6}{5} - a_1 \text{ ,}
$$
and
$$
2a_1 - a_2 = \tilde{D} \cdot E_1 \geq \text{mult}_Q\Bigl(\tilde{D}
\cdot E_1 \Bigr) > \frac{1}{\lambda} - a_2 > \frac{6}{5} - a_2 \text{ .}
$$
These imply that $a_1 > \frac{6}{10}$ and $a_2 > \frac{9}{10}$.

\end{itemize}

Consider now the blow-up $\sigma_1 : \tilde{X}_1 \to \tilde{X}$ of
the surface $\tilde{X}$ at the point $Q$ that contracts the
$(-1)$-curve $F_1$ to the point $Q$. Then for the strict
transforms of the exceptional divisors $E_1, E_2, E_3, E_4$ we
have
\begin{eqnarray*}
\tilde{E_1} & \sim_{\mathbb{Q}} & \sigma_1 ^*(E_1) - F_1\\
\tilde{E_2} & \sim_{\mathbb{Q}} & \sigma_1 ^*(E_2) - F_1\\
\tilde{E_3} & \sim_{\mathbb{Q}} & \sigma_1 ^*(E_3)\\
\tilde{E_4} & \sim_{\mathbb{Q}} & \sigma_1 ^*(E_4) \text{ .}
\end{eqnarray*}
Let now
$$\pi_1 : \tilde{X}_1 \stackrel{\sigma_1}{\to} \tilde{X} \stackrel{\pi}{\to} X$$ be the
composition $\pi_1 = \pi \circ \sigma_1$. We have
\begin{eqnarray*}
K_{\tilde{X}_1} & \sim_{\mathbb{Q}} & \sigma_1^*(K_{\tilde{X}}) +
F_1
  \sim_{\mathbb{Q}} \sigma_1^*(\pi^*(K_X)) + F_1
  \sim_{\mathbb{Q}} \pi_1^*(K_X) + F_1
\end{eqnarray*}
and
\begin{eqnarray*}
\tilde{D}_1 & \sim_{\mathbb{Q}} & \sigma_1^*(\tilde{D}) - m_1 F1\\
 & \sim_{\mathbb{Q}} & \sigma_1^*(\pi^*(D)-a_1E_1-a_2E_2-a_3E_3 -a_4E_4) - m_1 F_1\\
 & \sim_{\mathbb{Q}} & \pi_1^*(D) - a_1 \tilde{E_1} - a_2 \tilde{E_2} - a_3 \tilde{E_3} - a_4 \tilde{E_4} - (a_1 + a_2 +
m_1) F_1 \text{ ,}
\end{eqnarray*}
where $m_1 = \text{mult}_{Q} \tilde{D}$. Also the strict transform
of the anticanonical curve $Z$ is
\begin{eqnarray*}
\tilde{Z}_1 & \sim_{\mathbb{Q}} & \sigma_1^*(\tilde{Z})\\
 & \sim_{\mathbb{Q}} & \sigma_1^*(\pi^*(Z) - E_1 - E_2 - E_3 - E_4)\\
 & \sim_{\mathbb{Q}} & \pi_1^*(Z) - \tilde{E_1} - \tilde{E_2} - \tilde{E_3} - \tilde{E_4} - 2F_1 \text{ .}
\end{eqnarray*}
From the inequalities
\begin{eqnarray*}
0  \leq  \tilde{D}_1  \cdot \tilde{Z}_1  & = & 1 - a_1 - a_4\\
0  \leq  \tilde{E_1} \cdot \tilde{D}_1 & = & 2a_1 - a_2 - m_1\\
0  \leq  \tilde{E_2} \cdot \tilde{D}_1  & = & 2a_2 - a_1 - a_3 - m_1\\
0  \leq  \tilde{E_3} \cdot \tilde{D}_1  & = & 2a_3 - a_2 - a_4\\
0  \leq  \tilde{E_4} \cdot \tilde{D}_1  & = & 2a_4 - a_3\\
0  \leq  F_1 \cdot \tilde{D}_1  & = & m_1
\end{eqnarray*}
we get that $m_1 \leq \frac{1}{2}$. The equivalence
$$
K_{\tilde{X}_1 } + \lambda \tilde{D}_1  + \lambda a_1\tilde{E_1} +
\lambda a_2\tilde{E_2} + \lambda a_3\tilde{E_3} + \lambda a_4
\tilde{E_4} + (\lambda(a_1 + a_2 + m_1) -1)F_1 \sim_{\mathbb{Q}}
\pi_1^*(K_X + \lambda D)
$$
implies that there is a point $Q_1 \in F_1$ such that the pair
$$
K_{\tilde{X}_1 } + \lambda \tilde{D}_1  + \lambda a_1\tilde{E_1} +
\lambda a_2\tilde{E_2} + (\lambda(a_1 + a_2 + m_1) -1)F_1
$$
is not log canonical at $Q_1$.
\begin{itemize}
\item Suppose $Q_1 \in \tilde{E_1} \cap F_1$, then the log pair
$$
K_{\tilde{X}_1 } + \lambda \tilde{D}_1  + \lambda a_1\tilde{E_1} +
(\lambda(a_1 + a_2 + m_1) -1)F_1
$$
is not log canonical at $Q_1$ and so are the pairs
$$
K_{\tilde{X}_1 } + \lambda \tilde{D}_1  + \tilde{E_1} +
(\lambda(a_1 + a_2 + m_1) -1)F_1
$$
and
$$
K_{\tilde{X}_1 } + \lambda \tilde{D}_1  + \lambda a_1\tilde{E_1} +
F_1 \text{ .}
$$
By adjunction it follows that
$$
2a_1 - a_2 - m_1 = \tilde{D}_1  \cdot \tilde{E_1} \geq
 \text{mult}_{Q_1}\Bigl(\tilde{D}_1  \cdot \tilde{E_1} \Bigr) > \frac{12}{5} - a_1 - a_2 - m_1
\text{ ,}$$
 and
 $$
m_1 = \tilde{D}_1  \cdot F_1 \geq
 \text{mult}_{Q_1}\Bigl(\tilde{D}_1  \cdot F_1 \Bigr) > \frac{1}{\lambda}-a_1 >\frac{6}{5} - a_1 \text{ ,}$$

 which is false.

\item Suppose $Q_1 \in F_1 \backslash (\tilde{E_1} \cup
\tilde{E_2})$ then the log pair
$$K_{\tilde{X}_1 } + \lambda \tilde{D}_1  + (\lambda(a_1 + a_2 + m_1) -1)F_1$$ is not log canonical at $Q_1$
and so is the pair
$$K_{\tilde{X}_1 } + \lambda \tilde{D}_1  + F_1 \text{  , since  }  0 \leq \frac{5}{6}(a_1 + a_2 + m_1) -1
\leq 1$$ By adjunction it follows that
$$
m_1 = \tilde{D}_1  \cdot F_1 \geq
 \text{mult}_{Q_1}\Bigl(\tilde{D}_1  \cdot F_1 \Bigr) > \frac{1}{\lambda} > \frac{6}{5} \text{ ,}$$ which is impossible
since $m_1 \leq \frac{1}{2}$.

\item Suppose $Q_1 \in \tilde{E_2} \cap F_1$ then the log pair
$$K_{\tilde{X}_1 } + \lambda \tilde{D}_1  + \lambda a_2\tilde{E_2} + (\lambda(a_1 + a_2 + m_1) -1)F_1 $$ is not
log canonical at $Q_1$ and so are the pairs
$$
K_{\tilde{X}_1 } + \lambda \tilde{D}_1  + \tilde{E_2} +
(\lambda(a_1 + a_2 + m_1) - 1)F_1
$$
and
$$
K_{\tilde{X}_1 } + \lambda \tilde{D}_1  + \lambda a_2 \tilde{E_2} +
F_1 \text{ .}
$$
By adjunction it follows that
$$
2a_2 - a_1 - a_3 - m_1 = \tilde{D}_1  \cdot \tilde{E_2} \geq
 \text{mult}_{Q_1}\Bigl(\tilde{D}_1  \cdot \tilde{E_2} \Bigr) > \frac{12}{5} - a_1 - a_2 - m_1 $$
 and
   $$
m_1 = \tilde{D}_1  \cdot F_1 \geq
 \text{mult}_{Q_1}\Bigl(\tilde{D}_1  \cdot F_1 \Bigr) > \frac{1}{\lambda} - a_2 > \frac{6}{5} - a_2\text{ , }$$
 which imply that
 $$
 3a_2 - \frac{a_2}{2} \geq 3a_2 - a_3 > \frac{12}{5} \Rightarrow a_2 > \frac{6}{7} \cdot \frac{6}{5} \text{ .}
 $$

  \end{itemize}

  Consider now the blow-up $\sigma_2: \tilde{X}_2  \to \tilde{X}_1 $
of the surface $\tilde{X}_1$ at the point $Q_1$ that contracts the
$(-1)$-curve $F_2$ to the point $Q_1$. We then have

  \begin{eqnarray*}
   K_{\tilde{X}_2 } & \sim_{\mathbb{Q}} & \pi_2^*(K_X) + \tilde{F}_1 + 2F_2\\
  \tilde{D}_2  & \sim_{\mathbb{Q}} & \pi_2^*(D) -a_1 \tilde{E_1} -a_2 \tilde{E_2} - a_3 \tilde{E_3} -
a_4 \tilde{E_4} -(a_1 + a_2 + m_1) \tilde{F}_1-(a_1+2a_2+m_1+m_2) F_2\\
  \tilde{Z}_2  & \sim_{\mathbb{Q}} & \pi_2^*(Z) - \tilde{E_1} - \tilde{E_2} - \tilde{E_3} - \tilde{E_4} - 2\tilde{F}_1
- 3F_2 \text{ ,}
  \end{eqnarray*}
and
\begin{eqnarray*}
 0 \leq \tilde{D}_2  \cdot \tilde{Z}_2  & = & 1 - a_1 - a_4\\
0 \leq \tilde{E}_1 \cdot \tilde{D}_2  & = & 2a_1 - a_2 - m_1\\
0 \leq \tilde{E}_2 \cdot \tilde{D}_2  & = & 2a_2 - a_1 - a_3 - m_1 - m_2\\
0 \leq \tilde{E}_3 \cdot \tilde{D}_2  & = & 2a_3 - a_2 -a_4\\
0 \leq \tilde{E}_4 \cdot \tilde{D}_2  & = & 2a_4 - a_3\\
0 \leq \tilde{F}_1 \cdot \tilde{D}_2  & = & m_1 - m_2\\
0 \leq F_2 \cdot \tilde{D}_2  & = & m_2 \text{ ,}
\end{eqnarray*}
where $m_2 = \text{mult}_{Q_1} \tilde{D}_2 $.

Because of the equivalence

\begin{eqnarray*}
\lefteqn{\pi_2^*(K_X + \lambda D) \sim_{\mathbb{Q}} }\\
 & & K_{\tilde{X}_2 } + \lambda \tilde{D}_2  + \lambda a_1 \tilde{E}_1
+ \lambda a_2 \tilde{E}_2 + \lambda a_3 \tilde{E}_3 +  \lambda a_4
\tilde{E}_4 +\\
 & & (\lambda(a_1 + a_2 + m_1) -1)
 \tilde{F}_1 + (\lambda(a_1 + 2a_2 + m_1 + m_2) - 2) F_2
\end{eqnarray*}

there is a point $Q_2 \in F_2$ such that the pair
$$
K_{ \tilde{X}_2} + \lambda \tilde{D}_2  + \lambda a_1 \tilde{E}_1
+ \lambda a_2 \tilde{E}_2 + \lambda a_3 \tilde{E}_3 + \lambda a_4
\tilde{E}_4 + (\lambda(a_1 + a_2 + m_1) -1) \tilde{F}_1 +
(\lambda(a_1 + 2a_2 + m_1 + m_2) - 2) F_2
$$
is not log canonical at $Q_2$.
\begin{itemize}
\item Suppose $Q_2 \in F_2 \cap \tilde{F}_1$, then the log pair
$$
 \left( \tilde{X}_2 , \lambda \tilde{D}_2  + (\lambda(a_1 + a_2 + m_1) -1)
 \tilde{F}_1 + (\lambda(a_1 + 2a_2 + m_1 + m_2) - 2) F_2 \right)
$$
 is not log canonical at $Q_2$ and so are the log pairs
 $$
K_{\tilde{X}_2 } + \lambda \tilde{D}_2  + \tilde{F}_1 +
(\lambda(a_1 + 2a_2 + m_1 + m_2) - 2) F_2
 $$
and
$$
K_{\tilde{X}_2 } + \lambda \tilde{D}_2  + (\lambda(a_1 + a_2 +
m_1) -1) \tilde{F}_1 +  F_2 \text{ .}
$$
By adjunction it follows that
 $$
 m_1 - m_2 = \tilde{D}_2 \cdot \tilde{F}_1 \geq \text{mult}_{Q_2} \Bigl(\tilde{D}_2 \cdot
\tilde{F}_1 \Bigr) > \frac{18}{5} - ( a_1 + 2a_2 + m_1 + m_2 )
 $$
 and
 $$
 m_2 = \tilde{D}_2  \cdot F_2 \geq \text{mult}_{Q_2} \Bigl( \tilde{D}_2 \cdot F_2 \Bigr)>
\frac{12}{5} - ( a_1 + a_2 + m_1) \text{ ,}
$$ which is false.

\item Suppose $O \in F_2 \backslash (\tilde{F}_1 \cup
\tilde{E_2})$, then the log pair
$$
 K_{\tilde{X}_2 } + \lambda \tilde{D}_2  + (\lambda(a_1 + 2a_2 + m_1 + m_2) - 2) F_2
$$
is not log canonical at $Q_2$ and so is the log pair
$$
 K_{\tilde{X}_2 } + \lambda  \tilde{D}_2  +  F_2   \text{   , since   }\frac{5}{6} (a_1 + 2a_2 + m_1 + m_2) -
2 \leq 1 \text{ .}
$$
By adjunction it follows that
$$
 m_2 = \tilde{D}_2 \cdot F_2 \geq \text{mult}_{Q_2} \Bigl(\tilde{D}_2 \cdot F_2 \Bigr) >
\frac{6}{5}  \text{ ,}
 $$
 which is false.

\item Suppose $Q_2 \in F_2 \cap \tilde{E}_2$, then the log pair
$$
 K_{\tilde{X}_2 } + \lambda \tilde{D}_2  + (\lambda(a_1 + 2a_2 + m_1 + m_2) - 2) F_2 + \lambda a_2 \tilde{E}_2
$$
 is not log canonical at $Q_2$ and so are the log pairs
 $$
 K_{\tilde{X}_2 } + \lambda \tilde{D}_2  + (\lambda(a_1 + 2a_2 + m_1 + m_2) - 2) F_2 +  \tilde{E}_2
 $$
and
$$
K_{\tilde{X}_2 } + \lambda \tilde{D}_2  + F_2 + \lambda a_2
\tilde{E}_2 \text{ .}
$$
By adjunction it follows that
 $$
 2a_2 - a_1 - a_3 - m_1 - m_2 = \tilde{D}_2 \cdot \tilde{E}_2 \geq \text{mult}_{Q_2}
\Bigl(\tilde{D}_2 \cdot \tilde{E}_2\Bigr) > \frac{18}{5} - ( a_1 +
2a_2 + m_1 + m_2 )
 $$
 and
 $$
 m_2 = \tilde{D}_2 \cdot F_2 \geq \text{mult}_{Q_2} \Bigl(\tilde{D}_2 \cdot F_2 \Bigr) >
\frac{6}{5} - a_2 \text{, }
 $$
which implies that
$$
4a_2 - \frac{a_2}{2} \geq 4a_2 - a_3 > \frac{18}{5} \Rightarrow
a_2 > \frac{9}{10} \cdot \frac{6}{5} \text{ .}
$$
\end{itemize}

Consider now the blow-up $\sigma_k: \tilde{X}_k  \to
\tilde{X}_{k-1} $ of the surface $\tilde{X}_{k-1}$ at the point
$Q_{k-1}$ that contracts the $(-1)$-curve $F_k$ to the point
$Q_{k-1}$. We then have

\begin{eqnarray*}
K_{\tilde{X}_k } & \sim_{\mathbb{Q}} & \pi_k^*(K_X) + \tilde{F}_1
+ 2\tilde{F}_2 + 3 \tilde{F}_3
+...+ (k-1) \tilde{F}_{k-1}+ k F_k\\
\tilde{D}_k  & \sim_{\mathbb{Q}} & \pi_k^*(D) -a_1 \tilde{E_1}
-a_2 \tilde{E_2} - a_3 \tilde{E_3} - a_4 \tilde{E_4} -(a_1
+ a_2 + m_1) \tilde{F}_1-(a_1+2a_2+m_1+m_2) \tilde{F}_2 -...\\
& & - \left( a_1 + (k-1) a_2 + m_1 + m_2+ ..+ m_{k-1} \right)
\tilde{F}_{k-1}
- \left( a_1 + k a_2 + m_1 + m_2+ ..+ m_k \right) F_k  \\
\tilde{Z}_k  & \sim_{\mathbb{Q}} & \pi_k^*(Z) - \tilde{E_1} -
\tilde{E_2} - \tilde{E_3} - \tilde{E_4} - 2\tilde{F}_1 -
3\tilde{F}_2 -...- k \tilde{F}_{k-1} - (k+1) F_k \text{ ,}
\end{eqnarray*}
and
\begin{eqnarray*}
 0 \leq \tilde{D}_k  \cdot \tilde{Z}_k  & = & 1 - a_1 - a_4\\
0 \leq \tilde{E}_1 \cdot \tilde{D}_k  & = & 2a_1 - a_2 - m_1\\
0 \leq \tilde{E}_2 \cdot \tilde{D}_k  & = & 2a_2 - a_1 - a_3 - m_1 - m_2 -...- m_k\\
0 \leq \tilde{E}_3 \cdot \tilde{D}_k  & = & 2a_3 - a_2 -a_4\\
0 \leq \tilde{E}_4 \cdot \tilde{D}_k  & = & 2a_4 - a_3\\
0 \leq \tilde{F}_1 \cdot \tilde{D}_k  & = & m_1 - m_2\\
0 \leq \tilde{F}_2 \cdot \tilde{D}_k  & = & m_2 - m_3\\
 & . & \\
  & . & \\
   & . & \\
0 \leq \tilde{F}_{k-1} \cdot \tilde{D}_k  & = & m_{k-1} - m_k\\
0 \leq F_k \cdot \tilde{D}_k  & = & m_k \text{ ,}
\end{eqnarray*}
where $m_i = \text{mult}_{Q_{i-1}} \tilde{D}_i \text{ , for  }
i=1,...,k$.

Because of the equivalence
\begin{eqnarray*}
\lefteqn{\pi_k^*(K_X + \lambda D)  \sim_{\mathbb{Q}} }\\
 & & K_{\tilde{X}_k } + \lambda \tilde{D}_k  + \lambda a_1 \tilde{E}_1
+ \lambda a_2 \tilde{E}_2 + \lambda a_3 \tilde{E}_3 +  \lambda a_4
\tilde{E}_4 +\\
 & & (\lambda(a_1 + a_2 + m_1) -1)
 \tilde{F}_1 + (\lambda(a_1 + 2a_2 + m_1 + m_2) - 2) F_2 +...+\\
& & \left( \lambda ( a_1 + (k-1) a_2 + m_1 + m_2+...+ m_{k-1}) -
(k-1) \right) \tilde{F}_{k-1} +\\
& & \left( \lambda (a_1 + k a_2 + m_1 + m_2+...+ m_k ) -k \right)
F_k
\end{eqnarray*}

there is a point $Q_k \in F_k$ such that the pair
\begin{eqnarray*}
 K_{\tilde{X}_k } + \lambda \tilde{D}_k + \lambda a_2 \tilde{E}_2+ \left( \lambda (a_1 + (k-1) a_2 + m_1 +
m_2+...+ m_{k-1}) - (k-1) \right) \tilde{F}_{k-1}\\
 + \left( \lambda ( a_1 + k a_2 + m_1 + m_2+...+ m_k) -k \right) F_k
\end{eqnarray*}
is not log canonical at $Q_k$.
\begin{itemize}
\item Suppose $Q_k \in F_k \cap \tilde{F}_{k-1}$, then the log
pair
\begin{eqnarray*}
 K_{\tilde{X}_k } + \lambda \tilde{D}_k + \left( \lambda (a_1 + (k-1) a_2 + m_1 + m_2+...+ m_{k-1}) - (k-1)
\right) \tilde{F}_{k-1}\\
  + \left( \lambda ( a_1 + k a_2 + m_1 + m_2+...+ m_k) -k \right) F_k
\end{eqnarray*}
 is not log canonical at $Q_2$ and so are the log pairs
 $$
K_{\tilde{X}_k } + \lambda \tilde{D}_k  + \tilde{F}_{k-1} + \left(
\lambda ( a_1 + k a_2 + m_1 + m_2+...+ m_k) -k \right) F_k
 $$
and
$$
K_{\tilde{X}_k } + \lambda \tilde{D}_k  + \left( \lambda (a_1 +
(k-1) a_2 + m_1 + m_2+...+ m_{k-1}) - (k-1) \right)
\tilde{F}_{k-1} +  F_k \text{ .}
$$
By adjunction it follows that
 $$
 m_{k-1} - m_k = \tilde{D}_k \cdot \tilde{F}_{k-1} \geq \text{mult}_{Q_k}
\Bigl(\tilde{D}_k \cdot \tilde{F}_{k-1} \Bigr) > \frac{1}{\lambda}  (k+1) - (
a_1 + ka_2 + m_1 + m_2 +...+m_k ) \text{, }
 $$
  which is a contradiction. Indeed from the inequality above we have that
$$
 a_1 + ka_2 + m_1 + m_2 +...+m_{k-2} + 2 m_{k-1} > \frac{1}{\lambda}  (k+1)
 $$
 but since $m_1 \geq m_2 \geq ...\geq m_k $ we get that
 $$
 a_1 + ka_2 + k m_1  > \frac{1}{\lambda}  (k+1) \text{ .}
 $$
 However the inequality $0 \leq \tilde{E}_1 \cdot \tilde{D}_k   =  2a_1 - a_2 - m_1$
 finally gives us
 $$(2k+1) a_1 > \lambda (k+1) \Rightarrow a_4 < \frac{4k-1}{5(2k+1)} \Rightarrow a_2 \leq 3a_4 <
\frac{4k-1}{4k+2} \cdot \frac{6}{5} \text{ , }$$ which is false
since after the $(k-1)$-th blow up (and before the $k$-th blow up)
we have
$$
a_2 > \frac{6}{5} \cdot \frac{3k}{3k+1} > \frac{4k-1}{4k+2} \cdot
\frac{6}{5} \text{ .}
$$

\medskip

\item Suppose $Q_k \in F_k \backslash (\tilde{F}_{k-1} \cup
\tilde{E_2})$, then the log pair
$$
 K_{\tilde{X}_k } + \lambda \tilde{D}_k  + \left( \lambda (a_1 + k a_2 + m_1 + m_2+...+ m_k) -k \right) F_k
$$
is not log canonical at $Q_k$ and so is the log pair
$$
 K_{\tilde{X}_k } + \lambda \tilde{D}_k  +  F_k   \text{   , since   } \left( \lambda (a_1 + k a_2 + m_1
+ m_2+...+ m_k) -k \right) \leq 1 \text{ .}
$$
By adjunction it follows that
$$
 m_k = \tilde{D}_k \cdot F_k \geq \text{mult}_{Q_k} \Bigl(\tilde{D}_k \cdot F_k \Bigr) >
\frac{6}{5}  \text{ ,}
 $$
 which is false, since $ \frac{1}{2} \geq m_1 \geq m_2 \geq ... \geq m_k$.

\medskip

\item Suppose $Q_k \in F_k \cap \tilde{E}_2$, then the log pair
$$
 K_{\tilde{X}_k } + \lambda \tilde{D}_k  + \left( \lambda (a_1 + k a_2 + m_1 + m_2+...+ m_k) -k \right)
F_k + \lambda a_2 \tilde{E}_2
$$
 is not log canonical at $Q_k$ and so are the log pairs
 $$
 K_{\tilde{X}_k } + \lambda \tilde{D}_k  + \left( \lambda (a_1 + k a_2 + m_1 + m_2+...+ m_k) -k \right)
F_k +  \tilde{E}_2
 $$
and
$$
K_{\tilde{X}_k } + \lambda \tilde{D}_k  + F_k + \lambda a_2
\tilde{E}_2 \text{ .}
$$
By adjunction it follows that
 $$
 2a_2 - a_1 - a_3 - m_1 - m_2-...-m_k = \tilde{D}_k \cdot \tilde{E}_2 \geq
\text{mult}_{Q_k} \Bigl(\tilde{D}_k \cdot \tilde{E}_2\Bigr) >
\frac{6}{5} (k+1) - \left( a_1 + k a_2 + m_1 + m_2+...+ m_k
\right)
 $$
 and
 $$
 m_k = \tilde{D}_k \cdot F_k \geq \text{mult}_{Q_k} \Bigl(\tilde{D}_k \cdot F_k \Bigr) >
\frac{6}{5} - a_2 \text{, }
 $$
which implies that
$$
(k+2) a_2 - \frac{2}{3} a_2\geq (k+2) a_2 - a_3 > \frac{6}{5}
(k+1) \Rightarrow a_2
> \frac{6}{5} \cdot \frac{3(k+1)}{3(k+1)+1} \text{ .}
$$
\end{itemize}

\begin{remark}
 It remains to be shown that after the $k$-th blow up we have
 $$
 \left( \lambda (a_1 + k a_2 + m_1 + m_2+...+ m_k) -k \right) \leq 1 \text{, }
 $$
 and for this it is enough to show that
$$
 \left( \frac{5}{6} (a_1 + k a_2 + m_1 + m_2+...+ m_k) -k \right) \leq 1 \text{ . }
 $$
 Suppose that we have blown up $k-1$ times, then $a_2 > \frac{6}{5} \cdot \frac{3k}{3k+1}$.
 Let us assume on the contrary that
 \begin{eqnarray*}
 \frac{5}{6} (a_1 + k a_2 + m_1 + m_2+...+ m_k) -k  > 1 & \Rightarrow & a_1 + k a_2 + m_1 + m_2+...+
m_k > \frac{6}{5} (k + 1)\\
  a_1 + 2k a_1  \geq a_1 + k a_2 + k m_1  > \frac{6}{5} (k + 1) & \Rightarrow &  a_1 >
\frac{6}{5} \cdot \frac{k+1}{2k+1}\\
  a_4 \leq 1 - a_1 < 1 - \frac{6}{5} \cdot \frac{k+1}{2k+1} & \Rightarrow & a_2 \leq 3 a_4 < \frac{6}{5} \cdot \frac{4k-1}{4k+2}
\text{, }
 \end{eqnarray*}
 which is a contradiction.

\end{remark}

\begin{itemize}

\item If $Q\in E_2 \cap E_3$ then the log pair
$$
K_{\tilde{X}} +\lambda \tilde{D} + \lambda a_2E_2 + \lambda a_3E_3
$$
is not log canonical at the point $Q$ and so are the log pairs
$$
K_{\tilde{X}} + \lambda \tilde{D} + E_2 + \lambda a_3E_3 \text{
and } K_{\tilde{X}} + \lambda \tilde{D} + \lambda a_2E_2 + E_3
\text{ .}
$$
By adjunction it follows that
$$
2a_2- \frac{1}{2} a_2 - a_3  \geq 2a_2 - a_1 - a_3 = \tilde{D}
\cdot E_2 \geq
 \text{mult}_Q\Bigl(\tilde{D} \cdot E_2 \Bigr) > \frac{6}{5} - a_3 \text{ ,}
$$
and
$$
2a_3 - a_2 - a_4 = \tilde{D} \cdot E_3 \geq
\text{mult}_Q\Bigl(\tilde{D} \cdot E_3 \Bigr)  > \frac{6}{5} - a_2
\text{ .}
$$
These imply that $a_2 > \frac{4}{5}$ and $a_3 > \frac{4}{5}$.

Consider now the blow-up $\pi_2 : \tilde{\tilde{X}} \to \tilde{X}$
of the surface $\tilde{X}$ at the point $Q$ that contracts the
$(-1)$-curve $E$ to the point $Q$. Then for the strict transforms
of the exceptional divisors $E_1, E_2, E_3, E_4$ we have
\begin{eqnarray*}
\tilde{E}_1 & \sim_{\mathbb{Q}} & \pi_2^*(E_1)\\
\tilde{E}_2 & \sim_{\mathbb{Q}} & \pi_2^*(E_2) - E\\
\tilde{E}_3 & \sim_{\mathbb{Q}} & \pi_2^*(E_3) - E\\
\tilde{E}_4 & \sim_{\mathbb{Q}} & \pi_2^*(E_4)\\
\end{eqnarray*}
Let now
$$\pi : \tilde{\tilde{X}} \stackrel{\pi_2}{\to} \tilde{X} \stackrel{\pi_1}{\to} X$$ be
the composition $\pi = \pi_1 \circ \pi_2$. We have
\begin{eqnarray*}
K_{\tilde{\tilde{X}}} & = & \pi_2^*(K_{\tilde{X}}) + E
  \sim_{\mathbb{Q}} \pi_2^*(\pi_1^*(K_X)) + E
  \sim_{\mathbb{Q}} \pi^*(K_X) + E
\end{eqnarray*}
and
\begin{eqnarray*}
\tilde{\tilde{D}} & = & \pi_2^*(\tilde{D}) - m E\\
 & \sim_{\mathbb{Q}} & \pi_2^*(\pi_1^*(D)-a_1E_1-a_2E_2-a_3E_3 - a_4 E_4 - b_1F_1 - b_2F_2 - b_3F_3 -
b_4F_4) - m E\\
 & \sim_{\mathbb{Q}} & \pi^*(D) - a_1 \tilde{E}_1 - a_2 \tilde{E}_2 - a_3 \tilde{E}_3 - a_4 \tilde{E}_4 -
b_1 \tilde{F}_1 - b_2 \tilde{F}_2 - b_3 \tilde{F}_3 - b_4
\tilde{F}_4 - (a_2 + a_3 + m) E \text{ ,}
\end{eqnarray*}
where $m = \text{mult}_{Q} \tilde{D}$. Also the strict transform
of the anticanonical curve $Z$ is
\begin{eqnarray*}
\tilde{\tilde{Z}} & \sim_{\mathbb{Q}} & \pi_2^*(\tilde{Z})\\
 & \sim_{\mathbb{Q}} & \pi_2^*(\pi_1^*(Z) - E_1 - E_2 - E_3 - E_4)\\
 & \sim_{\mathbb{Q}} & \pi^*(Z) - \tilde{E}_1 - \tilde{E}_2 - \tilde{E}_3 - \tilde{E}_4 - 2E \text{ .}
\end{eqnarray*}
From the inequalities
\begin{eqnarray*}
0  \leq  \tilde{\tilde{D}} \cdot \tilde{\tilde{Z}} & = & 1 - a_1 - a_4\\
0  \leq  \tilde{E}_1 \cdot \tilde{\tilde{D}} & = & 2a_1 - a_2\\
0  \leq  \tilde{E}_2 \cdot \tilde{\tilde{D}} & = & 2a_2 - a_1 - a_3 - m\\
0  \leq  \tilde{E}_3 \cdot \tilde{\tilde{D}} & = & 2a_3 - a_2 - a_4 - m\\
0  \leq  \tilde{E}_4 \cdot \tilde{\tilde{D}} & = & 2a_4 - a_3\\
0  \leq  E \cdot \tilde{\tilde{D}} & = & m
\end{eqnarray*}
we get that $m = \text{mult}_{Q} \tilde{D} \leq \frac{1}{2}$.

The equivalence
$$
K_{\tilde{\tilde{X}}} + \lambda \tilde{\tilde{D}} + \lambda
a_2\tilde{E}_2 + \lambda a_3\tilde{E}_3 + \Bigl( \lambda (a_2 +
a_3 + \text{mult}_{Q} \tilde{D}) - 1 \Bigr) E \sim_{\mathbb{Q}}
\pi^*(K_{\tilde{X}} + \lambda \tilde{D} + \lambda a_2E_2 + \lambda
a_3E_3 )
$$
implies that there is a point $R \in E$ such that the pair
$$
K_{\tilde{\tilde{X}}} + \lambda \tilde{\tilde{D}} + \lambda
a_2\tilde{E}_2 + \lambda a_3\tilde{E}_3 + \Bigl( \lambda ( a_2 +
a_3 + \text{mult}_{Q} \tilde{D}) - 1 \Bigr) E
$$
is not log canonical at $R$.

\begin{itemize}
\item[(i)]  If $Q\in \tilde{E}_2 \cap E$ then the log pair
$$
K_{\tilde{\tilde{X}}} + \lambda \tilde{\tilde{D}} + \lambda
a_2\tilde{E}_2 + \Bigl( \lambda ( a_2  + a_3 + \text{mult}_{Q}
\tilde{D}) - 1 \Bigr) E
$$
is not log canonical at the point $R$ and so is the log pair
$$
K_{\tilde{\tilde{X}}} + \lambda \tilde{\tilde{D}}  + \lambda
a_2\tilde{E}_2 + E \text{ .}
$$
By adjunction and inequality ~\ref{mult>2} it follows that
$$
2 - \frac{4}{5} - a_2 \geq 2 - a_2 - a_3 \geq \text{mult}_{Q}
\tilde{D} = \tilde{\tilde{D}} \cdot E \geq
 \text{mult}_R\Bigl(\tilde{\tilde{D}} \cdot E \Bigr) > \frac{6}{5} - a_2  \text{ ,}
$$
which is false.

\item[(ii)]  If $Q\in  E \backslash (E_2 \cup E_3)$ then the log
pair
$$
K_{\tilde{\tilde{X}}}  + \lambda \tilde{\tilde{D}} + \Bigl(
\lambda ( a_2  + a_3 + \text{mult}_{Q} \tilde{D}) - 1 \Bigr) E
$$
is not log canonical at the point $R$ and so is the log pair
$$
K_{\tilde{\tilde{X}}} + \lambda \tilde{\tilde{D}}  +  E \text{ .}
$$
By adjunction it follows that
$$
 \text{mult}_{Q} \tilde{D} = \tilde{\tilde{D}} \cdot E \geq
 \text{mult}_R\Bigl(\tilde{\tilde{D}} \cdot E \Bigr) > \frac{6}{5} \text{ ,}
$$
which is false.

\item[(iii)]  If $Q\in \tilde{E}_3 \cap E$ then the log pair
$$
K_{\tilde{\tilde{X}}} + \lambda \tilde{\tilde{D}} + \lambda
a_3\tilde{E}_3 + \Bigl( \lambda ( a_2  + a_3 + \text{mult}_{Q}
\tilde{D}) - 1 \Bigr) E
$$
is not log canonical at the point $R$ and so is the log pair
$$
K_{\tilde{\tilde{X}}} + \lambda \tilde{\tilde{D}}  + \lambda a_3
\tilde{E}_3 + E \text{ .}
$$
By adjunction and inequality ~\ref{mult>2} it follows that
$$
 2 - \frac{4}{5} - a_2 \geq 2 - a_2 - a_3 \geq \text{mult}_{Q} \tilde{D} =
\tilde{\tilde{D}} \cdot E \geq
 \text{mult}_R\Bigl(\tilde{\tilde{D}} \cdot E \Bigr) > \frac{6}{5} - a_3 \text{ ,}
$$
which is false.

\end{itemize}

\end{itemize}

\end{proof}

\subsection{Del Pezzo surfaces of degree 1 with one
$\mathbb{A}_5$ type singularity}

In this section we will prove the following.

\begin{lemma}
\label{A5} Let $X$ be a del Pezzo surface with one Du Val
singularity of type $\mathbb{A}_5$ and $K_X^2=1$. Then the global
log canonical threshold of $X$ is
$$
\mathrm{lct} (X) = \frac{2}{3} \text{ .}
$$
\end{lemma}

\begin{proof}

Let $X$ be a del Pezzo surface with at most one Du Val singularity
of type $\mathbb{A}_5$ and $K_X^2=1$. Let $\pi_1: \tilde{X} \to X$
be the minimal resolution of $X$. The following diagram shows how
the exceptional curves intersect each other.
\bigskip

$\mathbb{A}_5$. \xymatrix{ {\bullet}^{E_1} \ar@{-}[r] &
{\bullet}^{E_2} \ar@{-}[r] & {\bullet}^{E_3} \ar@{-}[r] &
{\bullet}^{E_4} \ar@{-}[r] & {\bullet}^{E_5}}

\bigskip

Suppose that $\mathrm{lct}(X) < \frac{2}{3}$, then there exists an
effective $\mathbb{Q}$-divisor $D\in X$ such that $D
\sim_{\mathbb{Q}} -K_X$ and the log pair $(X, \lambda D)$ is  not
log canonical, where $\lambda < \frac{2}{3}$. We derive that the
pair $(X, \lambda D)$ is log canonical outside of a point $P\in X$
and not log canonical at $P$. Then
$$
\tilde{D} \sim_{\mathbb{Q}} \pi_1^*(D) - a_1E_1 - a_2E_2 - a_3E_3
- a_4E_4 - a_5E_5 \text{ and } \tilde{Z} \sim_{\mathbb{Q}}
\pi_1^*(Z) - E_1 - E_2 - E_3 - E_4 - E_5\text{ .}
$$
Let $Z$ be the curve in $|-K_X|$ that contains $P$. Since the
curve $Z$ is irreducible we may assume that the support of $D$
does not contain $Z$.

From the inequalities
\begin{eqnarray*}
0 \leq \tilde{D} \cdot \tilde{Z} & = & 1 - a_1 - a_5\\
0 \leq E_1 \cdot \tilde{D} & = & 2a_1 - a_2\\
0 \leq E_2 \cdot \tilde{D} & = & 2a_2 - a_1 - a_3\\
0 \leq E_3 \cdot \tilde{D} & = & 2a_3 - a_2 - a_4\\
0 \leq E_4 \cdot \tilde{D} & = & 2a_4 - a_3 - a_5\\
0 \leq E_5 \cdot \tilde{D} & = & 2a_5 - a_4
\end{eqnarray*}
we see that
$$a_1 \leq \frac{5}{6} \text{, } a_2 \leq \frac{4}{3} \text{, } a_3 \leq \frac{3}{2}
 \text{, } a_4 \leq \frac{4}{3} \text{, } a_5 \leq \frac{5}{6} $$
 and what is more
 $$
 2a_5 \geq a_4 \text{ , } \frac{3}{2} a_4 \geq a_3 \text{ , } \frac{4}{3} a_3 \geq a_2
\text{ , } \frac{5}{4} a_2 \geq a_1 \text{ .}
 $$
 Furthermore there exists a curve $ L_3 \in X$, that passes through the point
$P$, whose strict
 transform is a $(-1)$-curve that intersects the fundamental cycle as following.
 $$ \tilde{L}_3 \cdot E_3  = 1$$
 and
 $$
 L_3 \cdot E_j = 0 \text{ for  } j= 1,2,4,5 \text{ .}
 $$
 Then we easily get that
 \begin{eqnarray*}
 \tilde{L}_3 & \sim_{\mathbb{Q}} & \pi^*(L_3) - \frac{1}{2} E_1 - E_2 - \frac{3}{2} E_3 - E_4 -
\frac{1}{2} E_5 \text{ .}
 \end{eqnarray*}

 The image of $L_3$ under involution is either fixed or $L_3$ is mapped to another curve $L_3'$.
In either case we can assume that the irreducible line $L_3$ is
not contained in $\text{Supp}(D)$ and thus deduce the inequality
 $$
 0 \leq \tilde{L_3} \cdot \tilde{D} = 1 - a_3 \text{ , }
 $$

The equivalence
$$
K_{\tilde{X}} + \lambda \tilde{D} + a_1 \lambda E_1 + a_2 \lambda
E_2 + a_3 \lambda E_3 + a_4 \lambda E_4 + a_5 \lambda E_5
\sim_{\mathbb{Q}} \pi_1^*(K_X + D)
$$
implies that there is a point $Q\in E_1\cup E_2\cup E_3\cup E_4
\cup E_5$ such that the pair
$$K_{\tilde{X}} + \lambda \tilde{D} + a_1 \lambda E_1 + a_2 \lambda E_2 + a_3 \lambda E_3
+ a_4 \lambda E_4 + a_5 \lambda E_5$$ is not log canonical at $Q$.

\begin{itemize}
\item If the point $Q \in E_1$ and $Q\not \in E_2$ then
$$
K_{\tilde{X}} + \lambda \tilde{D} +  a_1 \lambda E_1
$$
is not log canonical at the point $Q$ and so is the pair
$$
K_{\tilde{X}} + \lambda \tilde{D} + E_1 \text{ , since  } a_1
\lambda \leq 1 \text{ .}
$$
By adjunction $(E_1, \lambda \tilde{D}|_{E_1})$ is not log
canonical at $Q$ and
$$1 \geq 2a_1 - \frac{4}{5} a_1 \geq 2a_1 - a_2 = \tilde{D} \cdot E_1 \geq
 \text{mult}_Q\Bigl(\tilde{D} \cdot E_1 \Bigr)  > \frac{1}{\lambda} > \frac{3}{2} \text{
,}$$ which is a contradiction.

\item If $Q\in E_1 \cap E_2$ then the log pair
$$
K_{\tilde{X}} + \lambda \tilde{D} +  a_1 \lambda E_1 + a_2 \lambda
E_2
$$
is not log canonical at the point $Q$ and so are the log pairs
$$
K_{\tilde{X}} + \lambda \tilde{D} + E_1 + a_2 \lambda E_2 \text{
and  } K_{\tilde{X}} + \lambda \tilde{D} + a_1 \lambda E_1 + E_2
\text{ .}
$$
By adjunction it follows that
$$
2a_2 - a_1 - a_3 = \tilde{D} \cdot E_2 \geq
\text{mult}_Q\Bigl(\tilde{D}|_{E_2} \Bigr) =
 \text{mult}_Q\Bigl(\tilde{D} \cdot E_2 \Bigr) > \frac{1}{\lambda} - a_1 > \frac{3}{2} -
a_1 \text{ ,}
$$
and
$$
2a_1 - a_2 = \tilde{D} \cdot E_1 \geq
\text{mult}_Q\Bigl(\tilde{D}|_{E_1} \Bigr) =
\text{mult}_Q\Bigl(\tilde{D} \cdot E_1 \Bigr)  > \frac{1}{\lambda}
- a_2 > \frac{3}{2} - a_2 \text{ .}
$$
From the first inequality we get $a_3 \geq \frac{9}{10}$ and then
we see that
$$
1 \geq a_1 + a_5 \geq a_1 + \frac{1}{2} \cdot \frac{2}{3} a_3 >
\frac{3}{4} + \frac{3}{10} > 1 \text{ , }
$$
which is a contradiction.

\item If $Q\in E_2$ but $Q\not \in E_1 \cup E_3$ then
$$
K_{\tilde{X}} + \lambda \tilde{D} + a_2 \lambda E_2
$$
is not log canonical at the point $Q$ and so is the pair
$$
K_{\tilde{X}} + \lambda \tilde{D} + E_2 \text{ , since  } a_2
\lambda \leq 1 \text{ .}
$$
By adjunction $(E_2, \lambda \tilde{D}|_{E_2})$ is not log
canonical at $Q$ and
$$1 \geq 2a_2 - \frac{1}{2} a_2 - \frac{3}{4} a_2 \geq 2a_2 - a_1 - a_3 = \tilde{D} \cdot
E_2 \geq
 \text{mult}_Q\Bigl(\tilde{D} \cdot E_2 \Bigr) > \frac{1}{\lambda} > \frac{3}{2}  \text{,
}$$ which is a contradiction.

\item If $Q\in E_2 \cap E_3$ then the log pair
$$
K_{\tilde{X}} + \lambda \tilde{D} +  a_2 \lambda E_2 + a_3 \lambda
E_3
$$
is not log canonical at the point $Q$ and so are the log pairs
$$
K_{\tilde{X}} + \lambda \tilde{D} + a_2 \lambda E_2 + E_3 \text{ ,
since } \lambda a_3 \leq 1 \text{ . }
$$
By adjunction it follows that
$$
2a_3 - a_2 -a_4= \tilde{D} \cdot E_3 \geq
\text{mult}_Q\Bigl(\tilde{D} \cdot E_3 \Bigr)  > \frac{1}{\lambda}
- a_2 > \frac{3}{2} - a_2 \text{ .}
$$
which, together with the inequality $a_4 \geq \frac{2}{3} a_3$,
implies that $a_3 > \frac{9}{8}$. However, this is impossible
since $a_3 \leq 1$.

\item If $Q\in E_3$ but $Q\not \in E_2 \cup E_4$ then
$$
K_{\tilde{X}} + \lambda \tilde{D} + a_3 \lambda E_3
$$
is not log canonical at the point $Q$ and so is the pair
$$
K_{\tilde{X}} + \lambda \tilde{D} + E_3 \text{ , since  } a_3
\lambda \leq 1 \text{ .}
$$
By adjunction $(E_3, \lambda \tilde{D}|_{E_3})$ is not log
canonical at $Q$ and
$$1 \geq \frac{2}{3} a_3 \geq 2a_3 - \frac{2}{3}a_3 - \frac{2}{3} a_3 \geq 2a_3 - a_2 -
a_4 = \tilde{D} \cdot E_3 \geq
 \text{mult}_Q\Bigl(\tilde{D} \cdot E_3 \Bigr) > \frac{3}{2}  \text{ ,}$$
which is false.

\bigskip

%In the graph that follows we can see the configuration of
%exceptional curves in the smooth surface $\tilde{X}$.

%\begin{figure}[h]
%\centering
%\includegraphics{A5_pict1.eps}
%\caption{Minimal resolution of an $\mathbb{A}_5$ singularity}
%\end{figure}

We will now show the existence of the curve $L_3$ which gave us
the inequality $a_3 \leq 1 $. If we now contract the curves $C,
E_5, E_4, E_3$ in this order we obtain two curves intersecting
each other as following.

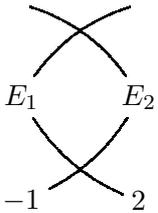
\xymatrix{ \ar @/^/ @{-}[dr]  &  \ar @/_/ @{-}[dl]\\
E_1 \ar @/_/ @{-} [dr] & E_2 \ar @/^/ @{-} [dl] \\
 -1 & 2  }

%\begin{figure}[h]
%\centering
%\includegraphics{A5_pict2.eps}
%\end{figure}

 However the resulting surface is isomorphic to $\mathbb{P}^2$ blown up at 4 points and
we know the configuration of all the -1 curves in this case.
Therefore there is always a -1 curve $L_2$ that intersects the
exceptional curve $E_2$ and not $E_3$. Indeed the resulting
surface is a smooth del Pezzo surface of degree 5 and $E_2+E_3$ is
its anticanonical divisor. Every -1 curve intersects the
anticanonical divisor only at one point by adjunction formula and
thus $L_2$ cannot intersect $E_1$.

If we now contract the curves $L_2, C, E_5, E_4$ we obtain a
smooth del Pezzo surface of degree 6 and we have the following
configuration of lines.

 
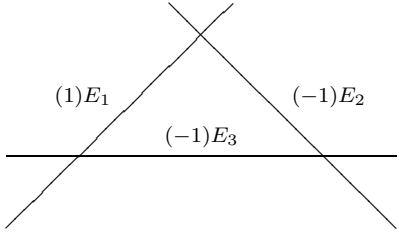
\xymatrix{
& & \ar@{-}[dddrrr]^{(-1)E_2} &  &  &\\
&  &  &   & &\\
\ar@{-}[rrrrr]^{(-1)E_3}& &  &  & &\\
\ar@{-}[uuurrr]^{(1)E_1}& & &  & & }

%\begin{figure}[h]
%\centering
%\includegraphics{A5_pict3.eps}
%\end{figure}

Therefore there exist -1 curves $L_2'$ and $L_3$ which intersect
the exceptional curves $E_2$ and $E_3$ transversaly and do not
intersect any other exceptional curve.

Therefore there are two -1 curves $L_2$ and $L_2'$ intersecting
the exceptional curve $E_2$ and -1 curve $L_3$ intersecting $E_3$,
such that $\tilde{L_2} \cdot \tilde{L_3} = \tilde{L_2'} \cdot
\tilde{L_3} = \tilde{L_2} \cdot \tilde{L_2'} = 0$. Involution
either fixes the curve $L_3$ or sends it to another line $L_3'$.
In any case we deduce the inequality $a_3 \leq 1$.

\end{itemize}
\end{proof}

\subsection{Del Pezzo surfaces of degree 1 with at most one
$\mathbb{A}_6$ type singularity}

In this section we will prove the following.

\begin{lemma}
\label{A6} Let $X$ be a del Pezzo surface with at most one Du Val
singularity of type $\mathbb{A}_6$ and $K_X^2=1$. Then the global
log canonical threshold of $X$ is
$$
\mathrm{lct} (X) = \frac{2}{3} \text{ .}
$$
\end{lemma}

\begin{proof}

Let $X$ be a del Pezzo surface with at most one Du Val singularity
of type $\mathbb{A}_6$ and $K_X^2=1$. Suppose that
$\mathrm{lct}(X) < \mathrm{lct}_2(X) \leq \frac{2}{3}$, then there
exists an effective $\mathbb{Q}$-divisor $D\in X$, such that $D
\sim_{\mathbb{Q}} -K_X$ and the log pair $(X, \lambda D)$ is  not
log canonical, where $\lambda < \frac{2}{3}$.

Let $Z$ be the unique curve in $|-K_X|$ that contains $P$. Since
the curve $Z$ is irreducible we may assume that the support of $D$
does not contain $Z$.

We derive that the pair $(X, \lambda D)$ is log canonical outside
of the singular point $P\in X$ and not log canonical at $P$. Let
$\pi_1: \tilde{X} \to X$ be the minimal resolution of $X$. The
following diagram shows how the exceptional curves intersect each
other.
\bigskip

$\mathbb{A}_6$. \xymatrix{ {\bullet}^{E_1} \ar@{-}[r] &
{\bullet}^{E_2} \ar@{-}[r] & {\bullet}^{E_3} \ar@{-}[r] &
{\bullet}^{E_4} \ar@{-}[r] & {\bullet}^{E_5} \ar@{-}[r] &
{\bullet}^{E_6}}
\bigskip

Then
$$
\tilde{D} \sim_{\mathbb{Q}} \pi_1^*(D) - a_1E_1 - a_2E_2 - a_3E_3
- a_4E_4 - a_5E_5 - a_6E_6  \text{ and } \tilde{Z}
\sim_{\mathbb{Q}} \pi_1^*(Z) - E_1 - E_2 - E_3 - E_4 - E_5 -
E_6\text{ .}
$$
From the inequalities
\begin{eqnarray*}
0 \leq \tilde{D} \cdot \tilde{Z} & = & 1 - a_1 - a_6\\
0 \leq E_1 \cdot \tilde{D} & = & 2a_1 - a_2\\
0 \leq E_2 \cdot \tilde{D} & = & 2a_2 - a_1 - a_3\\
0 \leq E_3 \cdot \tilde{D} & = & 2a_3 - a_2 - a_4\\
0 \leq E_4 \cdot \tilde{D} & = & 2a_4 - a_3 - a_5\\
0 \leq E_5 \cdot \tilde{D} & = & 2a_5 - a_4 - a_6\\
0 \leq E_6 \cdot \tilde{D} & = & 2a_6 - a_5
\end{eqnarray*}
we see that
$$
2a_6 \geq a_5 \text{ , } \frac{3}{2} a_5 \geq a_4 \text{ , }
\frac{4}{3} a_4 \geq a_3 \text{ , } \frac{5}{4} a_3 \geq a_2
\text{ , } \frac{6}{5} a_2 \geq a_1
$$
and
$$
2a_1 \geq a_2 \text{ , } \frac{3}{2} a_2 \geq a_3 \text{ , }
\frac{4}{3} a_3 \geq a_4 \text{ , } \frac{5}{4} a_4 \geq a_5
\text{ , } \frac{6}{5} a_5 \geq a_6 \text{ .}
$$
Moreover, for these non-negative numbers we get the bounds
$$
a_1 \leq \frac{6}{7} \text{, } a_2 \leq \frac{10}{7} \text{, } a_3
\leq \frac{12}{7}
 \text{, } a_4 \leq \frac{12}{7} \text{, } a_5 \leq \frac{10}{7} \text{, } a_6 \leq
\frac{6}{7} \text{ .}
 $$
 Furthermore there are four curves $ L_2, L_3, L_4,L_5 \in X$ that pass through the point
$P$, such that their strict transforms
 in $\tilde{X}$ are the $(-1)$-curves $\tilde{L_2}, \tilde{L_3}, \tilde{L_4},
\tilde{L_5}$ that intersect the fundamental cycle as following
 $$
 \tilde{L_2} \cdot E_2 = \tilde{L_3} \cdot E_3 = \tilde{L_4} \cdot E_4 =\tilde{L_5} \cdot
E_5 =  1
 $$
 and
 $$
 \tilde{L}_i \cdot E_j = 0 \text{ for all } i=2,3,4,5 \text{ and } j= 1,...,6 \text{ with
} i \not = j \text{ .}
 $$
 We can easily see that
 \begin{eqnarray*}
  \tilde{L_2} & \sim_{\mathbb{Q}} & \pi^* (L_2) - \frac{5}{7} E_1 - \frac{10}{7} E_2 - \frac{8}{7}
E_3 - \frac{6}{7} E_4 - \frac{4}{7}E_5- \frac{2}{7}E_6\text{ , }\\
 \tilde{L_3}  & \sim_{\mathbb{Q}} & \pi^* (L_3) - \frac{4}{7} E_1 - \frac{8}{7} E_2 - \frac{12}{7}
E_3 - \frac{9}{7} E_4 - \frac{6}{7}E_5- \frac{3}{7}E_6\text{ , }\\
 \tilde{L_4}  & \sim_{\mathbb{Q}} & \pi^* (L_4) - \frac{3}{7} E_1 - \frac{6}{7} E_2 - \frac{9}{7}
E_3 - \frac{12}{7} E_4 - \frac{8}{7}E_5- \frac{4}{7}E_6\text{ , }\\
 \tilde{L_5} & \sim_{\mathbb{Q}} & \pi^* (L_5) - \frac{2}{7} E_1 - \frac{4}{7} E_2 - \frac{6}{7} E_3
- \frac{8}{7} E_4 - \frac{10}{7}E_5- \frac{5}{7}E_6 \text{ .}
 \end{eqnarray*}
We have that $L_2 + L_5 \in |-2K_X|$ and $L_3 + L_4 \in |-2K_X|$
and we can assume that at least one member from each pair $L_2 +
L_5$ and $L_3 + L_4$ is not contained in the support of $D$. Thus
$0 \leq \tilde{L_3} \cdot \tilde{D} = 1 - a_3$ or $0 \leq
\tilde{L_4} \cdot \tilde{D} = 1 - a_4$. The equivalence
$$
K_{\tilde{X}} + \lambda \tilde{D} + a_1 \lambda E_1 + a_2 \lambda
E_2 + a_3 \lambda E_3 + a_4 \lambda E_4 + a_5 \lambda E_5 + a_6
\lambda E_6= \pi_1^*(K_X + \lambda D)
$$
implies that there is a point $Q\in E_1\cup E_2\cup E_3\cup E_4
\cup E_5 \cup E_6$ such that the pair
$$K_{\tilde{X}} + \lambda \tilde{D} + a_1 \lambda E_1 + a_2 \lambda E_2 + a_3 \lambda E_3
+ a_4 \lambda E_4 + a_5 \lambda E_5 + a_6 \lambda E_6$$ is not log
canonical at $Q$.

\begin{itemize}
\item If the point $Q \in E_1$ and $Q\not \in E_2$ then
$$
K_{\tilde{X}} + \lambda \tilde{D} +  a_1 \lambda E_1
$$
is not log canonical at the point $Q$ and so is the pair
$$
K_{\tilde{X}} + \lambda \tilde{D} + E_1 \text{ , since  } a_1
\lambda \leq 1 \text{ .}
$$
By adjunction $(E_1, \lambda \tilde{D}|_{E_1})$ is not log
canonical at $Q$ and
$$1 \geq \frac{7}{6} a_1 \geq 2a_1 - \frac{5}{6} a_1 \geq  2a_1 - a_2 = \tilde{D} \cdot
E_1 \geq
 \text{mult}_Q\Bigl(\tilde{D} \cdot E_1 \Bigr)  > \frac{1}{\lambda} > \frac{3}{2} \text{
,}$$ which is false.

\item If $Q\in E_1 \cap E_2$ then the log pair
$$
K_{\tilde{X}} + \lambda \tilde{D} +  a_1 \lambda E_1 + a_2 \lambda
E_2
$$
is not log canonical at the point $Q$ and so are the log pairs
$$
K_{\tilde{X}} + \lambda \tilde{D} + E_1 + a_2 \lambda E_2 \text{
and  } K_{\tilde{X}} + \lambda \tilde{D} + a_1 \lambda E_1 + E_2
\text{ .}
$$
By adjunction it follows that
$$
2 \frac{5}{4} a_3 - a_1 - a_3 \geq 2 a_2 - a_1 - a_3 = \tilde{D}
\cdot E_2 \geq
 \text{mult}_Q\Bigl(\tilde{D} \cdot E_2 \Bigr)  > \frac{3}{2} - a_1 \Rightarrow a_3 > 1
\text{ ,}
$$
and
$$
2a_1 - a_2 = \tilde{D} \cdot E_1 \geq
\text{mult}_Q\Bigl(\tilde{D}|_{E_1} \Bigr) =
\text{mult}_Q\Bigl(\tilde{D} \cdot E_1 \Bigr)  > \frac{3}{2} - a_2
\Rightarrow a_1 > \frac{3}{4} \text{ .}
$$
which leads to contradiction, as
$$
1 \geq a_1 + a_6 \geq a_1 + \frac{1}{2} \frac{2}{3} \frac{3}{4}
a_3 > \frac{3}{4} + \frac{1}{4} = 1 \text{ .}
$$

\item If $Q\in E_2$ but $Q\not \in E_1 \cup E_3$ then
$$
K_{\tilde{X}} + \lambda \tilde{D} + a_2 \lambda E_2
$$
is not log canonical at the point $Q$ and so is the pair
$$
K_{\tilde{X}} + \lambda \tilde{D} + E_2 \text{ , since  } a_2
\lambda \leq 1 \text{ .}
$$
By adjunction $(E_2, \lambda \tilde{D}|_{E_2})$ is not log
canonical at $Q$ and
$$1 \geq \frac{7}{10} a_2 \geq 2a_2 - \frac{4}{5} a_2 - \frac{1}{2} a_2  \geq 2a_2 - a_1
- a_3 = \tilde{D} \cdot E_2 \geq
 \text{mult}_Q\Bigl(\tilde{D} \cdot E_2 \Bigr) > \frac{3}{2} \text{ ,}$$
which is false.

\item If $Q\in E_3$ but $Q\not \in E_2 \cup E_4$ then
$$
K_{\tilde{X}} + \lambda \tilde{D} + a_3 \lambda E_3
$$
is not log canonical at the point $Q$ and so is the pair
$$
K_{\tilde{X}} + \lambda \tilde{D} + E_3 \text{ , since  } a_3
\lambda \leq 1 \text{ .}
$$
By adjunction $(E_3, \lambda \tilde{D}|_{E_3})$ is not log
canonical at $Q$ and
$$1 \geq \frac{7}{12} a_3 \geq 2a_3 - \frac{2}{3} a_3 - \frac{3}{4} a_3  = 2a_3 - a_2 -
a_4 = \tilde{D} \cdot E_3 \geq  \text{mult}_Q\Bigl(\tilde{D} \cdot
E_3 \Bigr) > \frac{1}{\lambda} > \frac{3}{2} \text{ ,}$$ which is
false.

\item If $Q\in E_3 \cap E_4$ then the log pair
$$
K_{\tilde{X}} + \lambda \tilde{D} +  a_3 \lambda E_3 + a_4 \lambda
E_4
$$
is not log canonical at the point $Q$ and so are the log pairs
$$
K_{\tilde{X}} + \lambda \tilde{D} + E_3 + a_4 \lambda E_4 \text{
and  } K_{\tilde{X}} + \lambda \tilde{D} + a_3 \lambda E_3 + E_4
\text{ ,  since  } \lambda a_4 \leq 1 \text{ .}
$$
By adjunction it follows that
$$
2a_3 - \frac{2}{3} a_3 - a_4 \geq 2a_3 - a_2 - a_4 = \tilde{D}
\cdot E_3 \geq
 \text{mult}_Q\Bigl(\tilde{D} \cdot E_3 \Bigr) > \frac{1}{\lambda} - a_4 > \frac{3}{2} -
a_4 \text{ ,}
$$
and
$$2a_4 - a_3 - \frac{2}{3} a_4 \geq
2a_4 - a_3 -a_5= \tilde{D} \cdot E_4 \geq
\text{mult}_Q\Bigl(\tilde{D} \cdot E_4 \Bigr)  > \frac{1}{\lambda}
- a_3 > \frac{3}{2} - a_3 \text{ .}
$$
This means that $a_3>1$ and $a_4>1$ which is impossible.

\item If $Q\in E_2 \cap E_3$ then the log pair
$$
K_{\tilde{X}} + \lambda \tilde{D} +  a_2 \lambda E_2 + a_3 \lambda
E_3
$$
is not log canonical at the point $Q$ and so are the log pairs
$$
K_{\tilde{X}} + \lambda \tilde{D} + E_2 + a_3 \lambda E_3 \text{
and  } K_{\tilde{X}} + \lambda \tilde{D} + a_2 \lambda E_2 + E_3
\text{ ,  since  } \lambda a_3 \leq 1 \text{ .}
$$
By adjunction it follows that
$$
2a_2 - a_1 - a_3 = \tilde{D} \cdot E_2 \geq
\text{mult}_Q\Bigl(\tilde{D}|_{E_2} \Bigr) > \frac{1}{\lambda} -
a_3 > \frac{3}{2} - a_3 \Rightarrow a_2>1 \text{ ,}
$$
and
$$
2a_3 - a_2 -a_4= \tilde{D} \cdot E_3 \geq
\text{mult}_Q\Bigl(\tilde{D}|_{E_3} \Bigr) > \frac{1}{\lambda} -
a_2 > \frac{3}{2} - a_2 \Rightarrow a_3 > \frac{6}{5} \text{ .}
$$

In $L_2 \not \in \text{Supp}D$ then  $0 \leq \tilde{D} \cdot L_3 =
1- a_3$ which contradicts the above inequalities and the same
holds when $L_2 \not \in \text{Supp}D$. Therefore we assume that
$D = a L_2 + c L_3 + \Omega$.

In the following graph we can see how the exceptional curves
intersect each other.

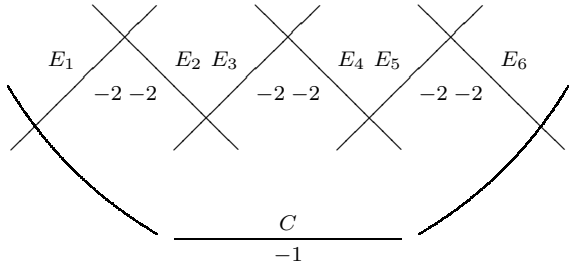
\begin{center}
\xymatrix{
 & \ar@{-}[ddrr]^{E_2}_{-2} &  & \ar@{-}[ddrr]^{E_4}_{-2} & & \ar@{-}[ddrr]^{E_6}_{-2} &
&\\
 \ar @/_/ @{-}[ddrr] &  &   & & & & &\\
\ar@{-}[uurr]^{E_1}_{-2} &  & \ar@{-}[uurr]^{E_3}_{-2} & &
\ar@{-}[uurr]^{E_5}_{-2} & &
&\\
& & \ar@{-}[rrr]_{-1}^{C} & & & \ar @/_/ @{-}[uurr] & & }
\end{center}

\bigskip
If we now contract the curves $C, E_6, E_5, E_4, E_3$ in this
order we get a del Pezzo surface of degree 6 and here is how the
remaining curves intersect.

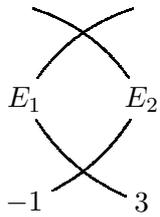
\xymatrix{ \ar @/^/ @{-}[dr]  &  \ar @/_/ @{-}[dl]\\
E_1 \ar @/_/ @{-} [dr] & E_2 \ar @/^/ @{-} [dl] \\
 -1 & 3  }

 \bigskip
However the resulting surface is isomorphic to $\mathbb{P}^2$
blown up at 3 points and we know the configuration of all the -1
curves in this case, they form a hexagon. Therefore there is
always a -1 curve $L_2$ that intersects the exceptional curve
$E_2$ and not $E_3$. Indeed the resulting surface is a smooth del
Pezzo surface of degree 6 and $E_2+E_3$ is its anticanonical
divisor. Every -1 curve intersects the anticanonical divisor only
at one point by adjunction formula and thus $L_2$ cannot intersect
$E_1$.

If we now contract the -1 curve $L_2$ in the original setting we
have a smooth surface of degree 2 with the following configuration
of curves.

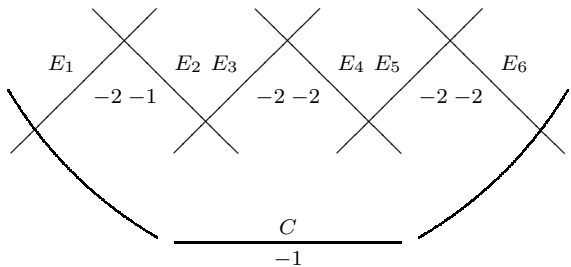
\begin{center}
\xymatrix{
 & \ar@{-}[ddrr]^{E_2}_{-1} &  & \ar@{-}[ddrr]^{E_4}_{-2} & & \ar@{-}[ddrr]^{E_6}_{-2} &
&\\
 \ar @/_/ @{-}[ddrr] &  &   & & & & &\\
\ar@{-}[uurr]^{E_1}_{-2} &  & \ar@{-}[uurr]^{E_3}_{-2} & &
\ar@{-}[uurr]^{E_5}_{-2} & &
&\\
& & \ar@{-}[rrr]_{-1}^{C} & & & \ar @/_/ @{-}[uurr] & & }
\end{center}

If after that we contract the curves $ C, E_6, E_5, E_4$ we obtain
a smooth del Pezzo surface of degree 6 and we have the following
configuration of lines. 
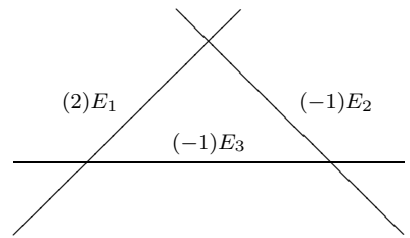
\xymatrix{
& & \ar@{-}[dddrrr]^{(-1)E_2} &  &  &\\
&  &  &   & &\\
\ar@{-}[rrrrr]^{(-1)E_3}& &  &  & &\\
\ar@{-}[uuurrr]^{(2)E_1}& & &  & & }

Therefore there exist -1 curves $L_2'$ and $L_3$ which intersect
the exceptional curves $E_2$ and $E_3$ transversaly and do not
intersect any other exceptional curve.

Therefore there are two -1 curves $L_2$ and $L_2'$ intersecting
the exceptional curve $E_2$ and -1 curve $L_3$ intersecting $E_3$,
such that $\tilde{L_2} \cdot \tilde{L_3} = \tilde{L_2'} \cdot
\tilde{L_3} = \tilde{L_2} \cdot \tilde{L_2'} = 0$.

Furthermore we also have the involusive images of the curves $L_2,
L_2', L_3$. In total we get six curves $ L_2, L_2', L_3,
L_4,L_5,L_5' \in X$ that pass through the point $P$, such that
their strict transforms
 in $\tilde{X}$ are the $(-1)$-curves
 $\tilde{L_2},\tilde{L_2'} , \tilde{L_3}, \tilde{L_4}, \tilde{L_5}, \tilde{L_5'}$ that
intersect the fundamental cycle as following
 $$
 \tilde{L_2} \cdot E_2 = \tilde{L_3} \cdot E_3 = \tilde{L_4} \cdot E_4 =\tilde{L_5} \cdot
E_5 =  1
 $$
 and
 $$
 L_i \cdot E_j = 0 \text{ for all } i=2,3,4,5 \text{ and } j= 1,...,6 \text{ with } i
\not = j \text{ .}
 $$
 We can easily see that
 \begin{eqnarray*}
  \tilde{L_2} & \sim_{\mathbb{Q}} & \pi^* (L_2) - \frac{5}{7} E_1 - \frac{10}{7} E_2 - \frac{8}{7} E_3 -
\frac{6}{7} E_4 - \frac{4}{7}E_5- \frac{2}{7}E_6\text{ , }\\
 \tilde{L_2'} & \sim_{\mathbb{Q}} & \pi^* (L_2') - \frac{5}{7} E_1 - \frac{10}{7} E_2 - \frac{8}{7} E_3 -
\frac{6}{7} E_4 - \frac{4}{7}E_5- \frac{2}{7}E_6\text{ , }\\
 \tilde{L_3}  & \sim_{\mathbb{Q}} & \pi^* (L_3) - \frac{4}{7} E_1 - \frac{8}{7} E_2 - \frac{12}{7} E_3 -
\frac{9}{7} E_4 - \frac{6}{7}E_5- \frac{3}{7}E_6\text{ , }\\
 \tilde{L_4}  & \sim_{\mathbb{Q}} & \pi^* (L_4) - \frac{3}{7} E_1 - \frac{6}{7} E_2 - \frac{9}{7} E_3 -
\frac{12}{7} E_4 - \frac{8}{7}E_5- \frac{4}{7}E_6\text{ , }\\
 \tilde{L_5} & \sim_{\mathbb{Q}} & \pi^* (L_5) - \frac{2}{7} E_1 - \frac{4}{7} E_2 - \frac{6}{7} E_3 -
\frac{8}{7} E_4 - \frac{10}{7}E_5- \frac{5}{7}E_6 \text{ ,}\\
 \tilde{L_5'} & \sim_{\mathbb{Q}} & \pi^* (L_5') - \frac{2}{7} E_1 - \frac{4}{7} E_2 - \frac{6}{7} E_3 -
\frac{8}{7} E_4 - \frac{10}{7}E_5- \frac{5}{7}E_6 \text{ .}
 \end{eqnarray*}

 We compute the intersection matrix for the curves $L_2, L_2', L_3$ and we see that these
three divisors are linearly independent. We know
 $\text{Pic}(\tilde{X})= \mathbb{Z}^9$ and we collapse six exceptional -2 curves,
therefore
 $\text{Pic}(X)= \mathbb{Z} \oplus \mathbb{Z} \oplus \mathbb{Z}$. Therefore this is a basis of the $\text{Pic}(X)= \mathbb{Z} \oplus \mathbb{Z} \oplus
\mathbb{Z}$.

 \begin{eqnarray*}
 L_2^2 & = & \pi^* (L_2) \cdot \pi^* (L_2) = \tilde{L_2} \cdot \pi^* (L_2) =
 \tilde{L_2}^2 +  \frac{10}{7}  \tilde{L_2} \cdot E_2 = -1 + \frac{10}{7} = \frac{3}{7}
\text{ ,}\\
 L_2'^2 & = & \pi^* (L_2') \cdot \pi^* (L_2') = \tilde{L_2'} \cdot \pi^* (L_2') =
 \tilde{L_2'}^2 +  \frac{10}{7}  \tilde{L_2'} \cdot E_2 = -1 + \frac{10}{7} = \frac{3}{7}
\text{ ,}\\
 L_3^2 & = & \pi^* (L_3) \cdot \pi^* (L_3) = \tilde{L_3} \cdot \pi^* (L_3) =
 \tilde{L_3}^2 +  \frac{12}{7}  \tilde{L_3} \cdot E_3 = -1 + \frac{12}{7} = \frac{5}{7}
\text{ ,}\\
 L_2' \cdot L_3 & = & \pi^* (L_2') \cdot \pi^* (L_3) = \tilde{L_2'} \cdot \pi^* (L_3) =
 \tilde{L_2'} \cdot \tilde{L_3} + \frac{8}{7}  \tilde{L_2'} \cdot E_2
 = \frac{8}{7}\\
 L_2 \cdot L_3 & = & \pi^* (L_2) \cdot \pi^* (L_3) = \tilde{L_2} \cdot \pi^* (L_3) =
 \tilde{L_2} \cdot \tilde{L_3} + \frac{8}{7}  \tilde{L_2} \cdot E_2 = \frac{8}{7}\\
 L_2 \cdot L_2' & = &  \pi^* (L_2) \cdot \pi^* (L_2') = \tilde{L_2} \cdot \pi^* (L_2') =
 \tilde{L_2} \cdot \tilde{L_2'} + \frac{10}{7}  \tilde{L_2} \cdot E_2 = \frac{10}{7}
 \end{eqnarray*}
 Now we would like to calculate $D= a L_2 + b L_2' + c L_3$. We have the following system
of equations.

 \begin{eqnarray*}
 \frac{3}{7} a +  \frac{10}{7}b+  \frac{8}{7}c & = & 1\\
  \frac{10}{7} a +  \frac{3}{7}b+  \frac{8}{7}c & = & 1\\
   \frac{8}{7} a +  \frac{8}{7}b+  \frac{5}{7}c & = & 1
 \end{eqnarray*}

 Therefore our effective divisor $D$ is
 $$
 D = \frac{1}{3} L_2 + \frac{1}{3} L_2' + \frac{1}{3} L_3 \text{ .}
 $$
 and
 $$
 \tilde{D} = \pi^* (D) - \frac{2}{3} E_1 - \frac{4}{3} E_2 - \frac{4}{3} E_3 -  E_4 -
\frac{2}{3}E_5- \frac{1}{3}E_6\text{ . }
 $$
We should note here that the divisor $\tilde{D}$ is a simple
normal crossings divisor and thus if we blow up more we do not
improve the log canonical threshold. There is no need to blow up
further and $(X, \lambda D)$  is log canonical which is a
contradiction.

\end{itemize}
\end{proof}

\subsection{Del Pezzo surfaces of degree 1 with exactly one
$\mathbb{A}_7$ type singularity}

In this section we will prove the following.

\begin{lemma}
\label{A7} Let $X$ be a del Pezzo surface with at most one Du Val
singularity of type $\mathbb{A}_7$ and $K_X^2=1$. Then the global
log canonical threshold of $X$ is
$$
\mathrm{lct} (X) = \frac{1}{2} \text{  or  } \frac{8}{15} \text{
.}
$$
\end{lemma}

\bigskip

The del Pezzo surface $X$ can be realised as the double cover
$$
X  \stackrel{2:1}{\longrightarrow} \mathbb{P}(1,1,2) \text{ ,}
$$
which is ramified along a sextic curve $R\in \mathbb{P}(1,1,2)$.
Let $\pi_1: \tilde{X} \to X$ be the minimal resolution of $X$. The
following diagram shows how the exceptional curves intersect each
other.
\bigskip

$\mathbb{A}_7$. \xymatrix{ {\bullet}^{E_1} \ar@{-}[r] &
{\bullet}^{E_2} \ar@{-}[r] & {\bullet}^{E_3} \ar@{-}[r] &
{\bullet}^{E_4} \ar@{-}[r] & {\bullet}^{E_5} \ar@{-}[r] &
{\bullet}^{E_6} \ar@{-}[r] & {\bullet}^{E_7}}
\bigskip
If the ramification divisor $R$ is irreducible then this implies
the existence of a -1 curve $\tilde{L}_4$ which intersects the
fundamental cycle only at the central exceptional curve $E_4$ and
this intersection is transversal. In the case the ramification
divisor $R$ is reducible no such line exists. Therefore we should
consider two cases depending on the existence or not of the -1
curve $\tilde{L}_4$.

\subsection{The ramification divisor is irreducible.}
\begin{proof}
Suppose $\mathrm{lct}(X) < \frac{1}{2}$. Then there exists an
effective $\mathbb{Q}$-divisor $D\in X$ and a positive rational
number $\lambda < \frac{1}{2}$, such that the log pair $(X,
\lambda D)$ is  not log canonical and $D \sim_{\mathbb{Q}} -K_X$,
where $\lambda < \frac{1}{2}$. Therefore the log pair $(X, \lambda
D)$ is also not log canonical.

Let $Z$ be the curve in $|-K_X|$ that contains $P$. Since the
curve $Z$ is irreducible we may assume that the support of $D$
does not contain $Z$.

We derive that the pair $(X, \lambda D)$ is log canonical outside
of a point $P\in X$ and not log canonical at $P$.

Then
\begin{eqnarray*}
\tilde{D} & \sim_{\mathbb{Q}} & \pi_1^*(D) - a_1E_1 - a_2E_2 -
a_3E_3 - a_4E_4 - a_5E_5 -a_6E_6 -a_7E_7
\text{ and }\\
\tilde{Z} & \sim_{\mathbb{Q}} & \pi_1^*(Z) - E_1 - E_2 - E_3 - E_4
- E_5 - E_6 - E_7 \text{ .}
\end{eqnarray*}
From the inequalities
\begin{eqnarray*}
0 \leq \tilde{D} \cdot \tilde{Z} & = & 1 - a_1 - a_7\\
0 \leq E_1 \cdot \tilde{D} & = & 2a_1 - a_2\\
0 \leq E_2 \cdot \tilde{D} & = & 2a_2 - a_1 - a_3\\
0 \leq E_3 \cdot \tilde{D} & = & 2a_3 - a_2 - a_4\\
0 \leq E_4 \cdot \tilde{D} & = & 2a_4 - a_3 - a_5\\
0 \leq E_5 \cdot \tilde{D} & = & 2a_5 - a_4 - a_6\\
0 \leq E_6 \cdot \tilde{D} & = & 2a_6 - a_5 - a_7\\
0 \leq E_7 \cdot \tilde{D} & = & 2a_7 - a_6
\end{eqnarray*}
we get
$$
2a_7 \geq a_6 \text{ , } \frac{3}{2} a_6 \geq a_5 \text{ , }
\frac{4}{3} a_5 \geq a_4 \text{ , } \frac{5}{4} a_4 \geq a_3
\text{ , } \frac{6}{5} a_3 \geq a_2 \text{ , } \frac{7}{6} a_2
\geq a_1
$$
and moreover
$$
a_1 \leq \frac{7}{8} \text{ , } a_2 \leq \frac{12}{8} \text{ , }
a_3 \leq \frac{15}{8} \text{ , } a_4 \leq 2 \text{ , } a_5 \leq
\frac{15}{8} \text{ , } a_6 \leq \frac{12}{8} \text{ , } a_7 \leq
\frac{7}{8} \text{ . }
$$
Furthermore there are lines $L_2, L_4, L_6 \in X$ that pass
through the point $P$ whose strict
 transforms are $(-1)$-curves that intersect the fundamental cycle as following.
 $$ L_2 \cdot E_2 = L_4 \cdot E_4 = L_6 \cdot E_6 = 1$$
 and
 $$
 L_i \cdot E_j = 0 \text{ for all } i,j= 2,4,6 \text{ with } i \not = j \text{ .}
 $$
Then we easily get that
\begin{eqnarray*}
 \tilde{L_2} & \sim_{\mathbb{Q}} & \pi^*(L_2) - \frac{3}{4} E_1 - \frac{3}{2} E_2 - \frac{5}{4} E_3 - E_4
- \frac{3}{4} E_5 -\frac{1}{2} E_6 - \frac{1}{4} E_7\\
\tilde{L_4} & \sim_{\mathbb{Q}} & \pi^*(L_4) - \frac{1}{2} E_1 -
E_2 - \frac{3}{2} E_3 - 2E_4 -
\frac{3}{2} E_5 - E_6 - \frac{1}{2} E_7\\
\tilde{L_6} & \sim_{\mathbb{Q}} & \pi^*(L_6) - \frac{1}{4} E_1 -
\frac{1}{2} E_2 - \frac{3}{4} E_3 - E_4 -\frac{5}{4} E_5
-\frac{3}{2} E_6 -\frac{3}{4} E_7 \text{ .}
\end{eqnarray*}
We observe that $2L_4$ is a Cartier divisor in the bianticanonical
linear system $|-2K_X|$. We will show that $a_4 \leq 1$ which is a
key inequality for what will follow. Indeed, since $L_4$ is
irreducible and $L_4 \sim_{\mathbb{Q}} -K_X$, we can assume that
$L_4 \not \in \text{Supp}(D)$. Then
$$
0 \leq \tilde{L_4} \cdot \tilde{D} = 1 - a_4 \text{ .}
$$
The equivalence
$$
K_{\tilde{X}} + \lambda \tilde{D} + \lambda a_1E_1 + \lambda
a_2E_2 + \lambda a_3E_3 + \lambda a_4E_4 + \lambda a_5E_5 +
\lambda a_6E_6 + \lambda a_7E_7 \sim_{\mathbb{Q}} \pi_1^*(K_X + D)
$$
implies that there is a point $Q\in E_1\cup E_2\cup E_3\cup E_4
\cup E_5 \cup E_6 \cup E_7$, such that the pair
$$K_{\tilde{X}} + \lambda \tilde{D} + \lambda a_1E_1 + \lambda a_2E_2 + \lambda a_3E_3 +
\lambda a_4E_4 + \lambda a_5E_5 + \lambda a_6E_6 + \lambda
a_7E_7$$ is not log canonical at $Q$.

\begin{itemize}
\item If the point $Q \in E_1$ and $Q\not \in E_2$ then
$$
K_{\tilde{X}} + \lambda \tilde{D} +  a_1 \lambda E_1
$$
is not log canonical at the point $Q$ and so is the pair
$$
K_{\tilde{X}} + \lambda \tilde{D} + E_1 \text{ , since  } a_1
\lambda \leq 1 \text{ .}
$$
By adjunction $(E_1, \lambda \tilde{D}|_{E_1})$ is not log
canonical at $Q$ and
$$
2 \cdot \frac{7}{6} a_2 -a_2 \geq 2a_1 - a_2 = \tilde{D} \cdot E_1
\geq
 \text{mult}_Q\Bigl(\tilde{D} \cdot E_1 \Bigr)  > \frac{1}{\lambda} > 2 \text{ ,}$$
which is false, since $a_2 \leq \frac{12}{8}$.

\item If $Q\in E_1 \cap E_2$ then the log pair
$$
K_{\tilde{X}} + \lambda \tilde{D} +  \lambda a_1 E_1 +  \lambda
a_2E_2
$$
is not log canonical at the point $Q$ and so is the log pair
$$
K_{\tilde{X}} + \lambda \tilde{D} + E_1 + \lambda a_2 E_2  \text{
.}
$$
By adjunction it follows that
$$
2a_1 - a_2 = \tilde{D} \cdot E_1 \geq
\text{mult}_Q\Bigl(\tilde{D}|_{E_1} \Bigr) =
\text{mult}_Q\Bigl(\tilde{D} \cdot E_1 \Bigr)  > \frac{1}{\lambda}
- a_2 > 2 - a2 \text{ ,}
$$
which is false, since $a_1 \leq \frac{7}{8}$.

\item If $Q\in E_2$ but $Q\not \in E_1 \cup E_3$ then
$$
K_{\tilde{X}} + \lambda \tilde{D} +  \lambda a_2 E_2
$$
is not log canonical at the point $Q$ and so is the pair
$$
K_{\tilde{X}} + \lambda \tilde{D} + E_2 \text{ , since  }  \lambda
a_2 \leq 1 \text{ .}
$$
By adjunction $(E_2, \lambda \tilde{D}|_{E_2})$ is not log
canonical at $Q$ and
$$2a_2 - \frac{5}{6} a2 \geq 2a_2 - a_1 - a_3 = \tilde{D} \cdot E_2 \geq
 \text{mult}_Q\Bigl(\tilde{D} \cdot E_2 \Bigr) > \frac{1}{\lambda} > 2  \text{ ,}$$
which is false, since $a_2 \leq \frac{12}{8}$.

\item If $Q\in E_2 \cap E_3$ then the log pair
$$
K_{\tilde{X}} + \lambda \tilde{D} +  \lambda a_2 E_2 +  \lambda
a_3E_3
$$
is not log canonical at the point $Q$ and so is the log pair
$$
K_{\tilde{X}} + \lambda \tilde{D} +  \lambda a_2E_2 + E_3 \text{ ,
since } \lambda a_3 < 1 \text{ .}
$$
By adjunction it follows that
$$
2a_3 - a_2 -a_4= \tilde{D} \cdot E_3 \geq
\text{mult}_Q\Bigl(\tilde{D}|_{E_3} \Bigr) =
\text{mult}_Q\Bigl(\tilde{D} \cdot E_3 \Bigr)  > \frac{1}{\lambda}
- a_2 > 2 - a_2 \text{ ,}
$$
which, along with the inequality $a_4 \geq \frac{4}{5} a_3$,
implies that $a_4 > 1$, which is impossible.

\item If $Q\in E_3$ but $Q\not \in E_2 \cup E_4$ then
$$
K_{\tilde{X}} + \lambda \tilde{D} +  \lambda a_3E_3
$$
is not log canonical at the point $Q$ and so is the pair
$$
K_{\tilde{X}} + \lambda \tilde{D} + E_3 \text{ , since  }  \lambda
a_3\leq 1 \text{ .}
$$
By adjunction $(E_3, \lambda \tilde{D}|_{E_3})$ is not log
canonical at $Q$ and
$$2a_3 - a_2 - a_4 = \tilde{D} \cdot E_3 \geq \text{mult}_Q\Bigl(\tilde{D}|_{E_3}
\Bigr) =
 \text{mult}_Q\Bigl(\tilde{D} \cdot E_3 \Bigr) > \frac{1}{\lambda} > 2  \text{ .}$$
This inequality together with $a_4 \geq \frac{4}{5} a_3$  implies
that $a_4 > 1$, which is impossible.

\item If $Q\in E_3 \cap E_4$ then the log pair
$$
K_{\tilde{X}} + \lambda \tilde{D} +  \lambda a_3 E_3 +  \lambda
a_4 E_4
$$
is not log canonical at the point $Q$ and so is the log pair
$$
K_{\tilde{X}} + \lambda \tilde{D} +  \lambda a_3E_3 + E_4 \text{ ,
since } \lambda a_4 \leq 1 \text{ .}
$$
By adjunction it follows that

$$
2a_4 - a_3 -a_5= \tilde{D} \cdot E_4 \geq
\text{mult}_Q\Bigl(\tilde{D}|_{E_4} \Bigr) =
\text{mult}_Q\Bigl(\tilde{D} \cdot E_4 \Bigr)  > \frac{1}{\lambda}
- a_3 > 2 - a_3 \text{ ,}
$$
which contradicts $a_4 \leq 1$.

\item If $Q\in E_4$ but $Q\not \in E_3 \cup E_5$ then
$$
K_{\tilde{X}} + \lambda \tilde{D} +  \lambda a_4E_4
$$
is not log canonical at the point $Q$ and so is the pair
$$
K_{\tilde{X}} + \lambda \tilde{D} + E_4 \text{ , since  }  \lambda
a_4\leq 1 \text{ .}
$$
By adjunction $(E_4, \lambda \tilde{D}|_{E_4})$ is not log
canonical at $Q$ and
$$2a_4 - a_3 - a_5 = \tilde{D} \cdot E_4 \geq \text{mult}_Q\Bigl(\tilde{D}|_{E_4}
\Bigr) =
 \text{mult}_Q\Bigl(\tilde{D} \cdot E_4 \Bigr) > \frac{1}{\lambda} > 2  \text{ ,}$$
which is false since $a_4 \leq 1$.

\end{itemize}
\end{proof}

\subsection{The ramification divisor $R$ is reducible.}
\begin{proof} Let $X$ be a del Pezzo surface with exactly one Du
Val singularity of type $\mathbb{A}_7$ and $K_X^2=1$. Suppose
$\mathrm{lct}(X) < \frac{8}{15}$. Then there exists an effective
$\mathbb{Q}$-divisor $D\in X$ and a positive rational number
$\lambda < \frac{8}{15}$, such that the log pair $(X, \lambda D)$
is  not log canonical and $D \sim_{\mathbb{Q}} -K_X$, where
$\lambda < \frac{8}{15}$. Therefore the log pair $(X, \lambda D)$
is also not log canonical.

Let $Z$ be the curve in $|-K_X|$ that contains $P$. Since the
curve $Z$ is irreducible we may assume that the support of $D$
does not contain $Z$.

We derive that the pair $(X, \lambda D)$ is log canonical outside
of a point $P\in X$ and not log canonical at $P$. Then
\begin{eqnarray*}
\tilde{D} & \sim_{\mathbb{Q}} & \pi_1^*(D) - a_1E_1 - a_2E_2 -
a_3E_3 - a_4E_4 - a_5E_5 -a_6E_6 -a_7E_7
\text{ and }\\
\tilde{Z} & \sim_{\mathbb{Q}} & \pi_1^*(Z) - E_1 - E_2 - E_3 - E_4
- E_5 - E_6 - E_7 \text{ .}
\end{eqnarray*}
From the inequalities
\begin{eqnarray*}
0 \leq \tilde{D} \cdot \tilde{Z} & = & 1 - a_1 - a_7\\
0 \leq E_1 \cdot \tilde{D} & = & 2a_1 - a_2\\
0 \leq E_2 \cdot \tilde{D} & = & 2a_2 - a_1 - a_3\\
0 \leq E_3 \cdot \tilde{D} & = & 2a_3 - a_2 - a_4\\
0 \leq E_4 \cdot \tilde{D} & = & 2a_4 - a_3 - a_5\\
0 \leq E_5 \cdot \tilde{D} & = & 2a_5 - a_4 - a_6\\
0 \leq E_6 \cdot \tilde{D} & = & 2a_6 - a_5 - a_7\\
0 \leq E_7 \cdot \tilde{D} & = & 2a_7 - a_6
\end{eqnarray*}
we get
$$
2a_7 \geq a_6 \text{ , } \frac{3}{2} a_6 \geq a_5 \text{ , }
\frac{4}{3} a_5 \geq a_4 \text{ , } \frac{5}{4} a_4 \geq a_3
\text{ , } \frac{6}{5} a_3 \geq a_2 \text{ , } \frac{7}{6} a_2
\geq a_1
$$
and moreover
$$
a_1 \leq \frac{7}{8} \text{ , } a_2 \leq \frac{12}{8} \text{ , }
a_3 \leq \frac{15}{8} \text{ , } a_4 \leq 2 \text{ , } a_5 \leq
\frac{15}{8} \text{ , } a_6 \leq \frac{12}{8} \text{ , } a_7 \leq
\frac{7}{8} \text{ . }
$$
Furthermore there are lines $L_2,L_5, L_6 \in X$ that pass through
the point $P$ whose strict
 transforms are $(-1)$-curves that intersect the fundamental cycle as following.
 $$ L_2 \cdot E_2 = L_5 \cdot E_5 = L_6 \cdot E_6 = 1$$
 and
 $$
 L_i \cdot E_j = 0 \text{ for all } i,j= 2,5,6 \text{ with } i \not = j \text{ .}
 $$
Then we easily get that
\begin{eqnarray*}
 \tilde{L}_2 & \sim_{\mathbb{Q}} & \pi^*(L_2) - \frac{3}{4} E_1 - \frac{3}{2} E_2 - \frac{5}{4} E_3 - E_4
- \frac{3}{4} E_5 -\frac{1}{2} E_6 - \frac{1}{4} E_7\\
 \tilde{L}_3 & \sim_{\mathbb{Q}} & \pi^*(L_3) - \frac{5}{8} E_1 - \frac{5}{4} E_2 - \frac{15}{8} E_3 -
\frac{3}{2} E_4 - \frac{9}{8} E_5 -\frac{3}{4} E_6 - \frac{3}{8} E_7\\
 \tilde{L}_5 & \sim_{\mathbb{Q}} & \pi^*(L_5) - \frac{3}{8} E_1 - \frac{3}{4} E_2 - \frac{9}{8} E_3 -
\frac{3}{2} E_4 - \frac{15}{8} E_5 -\frac{5}{4} E_6 - \frac{5}{8} E_7\\
\tilde{L}_6 & \sim_{\mathbb{Q}} & \pi^*(L_6) - \frac{1}{4} E_1 -
\frac{1}{2} E_2 - \frac{3}{4} E_3 - E_4 -\frac{5}{4} E_5
-\frac{3}{2} E_6 -\frac{3}{4} E_7 \text{ .}
\end{eqnarray*}
Since $L_2 + 2L_3 \in |-3K_X|$ we can assume that $L_2 \not \in
\text{Supp}D$ or $L_3 \not \in \text{Supp}D$ and then
\begin{eqnarray*}
 0 \leq \tilde{L}_2 \cdot \tilde{D} = 1 - a_2 \Rightarrow a_2 \leq 1 \text{  or }
 0 \leq \tilde{L}_3 \cdot \tilde{D} = 1 - a_3 \Rightarrow a_3 \leq 1 \text{ .}
\end{eqnarray*}
In the same way we obtain that
\begin{eqnarray*}
 0 \leq \tilde{L}_3 \cdot \tilde{D} = 1 - a_2 \Rightarrow a_3 \leq 1 \text{  or }
 0 \leq \tilde{L}_5 \cdot \tilde{D} = 1 - a_3 \Rightarrow a_5 \leq 1 \text{ , }
\end{eqnarray*}
since $L_3 + L_5 \in |-2 K_X|$ together with Lemma implies that
either $L_3 \not \in D$ or $L_5 \not \in D$.

The equivalence
$$
K_{\tilde{X}} + \lambda \tilde{D} + \lambda a_1E_1 + \lambda
a_2E_2 + \lambda a_3E_3 + \lambda a_4E_4 + \lambda a_5E_5 +
\lambda a_6E_6 + \lambda a_7E_7 \sim_{\mathbb{Q}} \pi_1^*(K_X +
\lambda D)
$$
implies that there is a point $Q\in E_1\cup E_2\cup E_3\cup E_4
\cup E_5 \cup E_6 \cup E_7$, such that the pair
$$K_{\tilde{X}} + \lambda \tilde{D} + \lambda a_1E_1 + \lambda a_2E_2 + \lambda a_3E_3 +
\lambda a_4E_4 + \lambda a_5E_5 + \lambda a_6E_6 + \lambda
a_7E_7$$ is not log canonical at $Q$.

\begin{itemize}
\item If the point $Q \in E_1$ and $Q\not \in E_2$ then
$$
K_{\tilde{X}} + \lambda \tilde{D} +  a_1 \lambda E_1
$$
is not log canonical at the point $Q$ and so is the pair
$$
K_{\tilde{X}} + \lambda \tilde{D} + E_1 \text{ , since  } a_1
\lambda \leq 1 \text{ .}
$$
By adjunction $(E_1, \lambda \tilde{D}|_{E_1})$ is not log
canonical at $Q$ and
$$
1 \geq 2 a_1 -\frac{6}{7} a_1 \geq 2a_1 - a_2 = \tilde{D} \cdot
E_1 \geq
 \text{mult}_Q\Bigl(\tilde{D} \cdot E_1 \Bigr)  > \frac{1}{\lambda} > \frac{15}{8} \text{
,}$$ which is false.

\item If $Q\in E_1 \cap E_2$ then the log pair
$$
K_{\tilde{X}} + \lambda \tilde{D} +  \lambda a_1 E_1 +  \lambda
a_2E_2
$$
is not log canonical at the point $Q$ and so is the log pair
$$
K_{\tilde{X}} + \lambda \tilde{D} + E_2 + \lambda a_1 E_1  \text{
.}
$$
By adjunction it follows that
$$
2a_2 - \frac{5}{6} a_2 -a_1 \geq 2a_2 - a_1 - a_3= \tilde{D} \cdot
E_2 \geq \text{mult}_Q\Bigl(\tilde{D} \cdot E_2 \Bigr)  >
\frac{1}{\lambda} - a_1 > \frac{15}{8} - a_1 \text{ ,}
$$
which is false.

\item If $Q\in E_2$ but $Q\not \in E_1 \cup E_3$ then
$$
K_{\tilde{X}} + \lambda \tilde{D} +  \lambda a_2 E_2
$$
is not log canonical at the point $Q$ and so is the pair
$$
K_{\tilde{X}} + \lambda \tilde{D} + E_2 \text{ , since  }  \lambda
a_2 \leq 1 \text{ .}
$$
By adjunction $(E_2, \lambda \tilde{D}|_{E_2})$ is not log
canonical at $Q$ and
$$1\geq \frac{2}{3} a_2 \geq 2a_2 - \frac{5}{6} a_2 \geq 2a_2 - a_1 - a_3 = \tilde{D}
\cdot E_2 \geq
 \text{mult}_Q\Bigl(\tilde{D} \cdot E_2 \Bigr) > \frac{1}{\lambda} > \frac{15}{8}  \text{
,}$$ which is false.

\item If $Q\in E_2 \cap E_3$ then the log pair
$$
K_{\tilde{X}} + \lambda \tilde{D} +  \lambda a_2 E_2 +  \lambda
a_3E_3
$$
is not log canonical at the point $Q$ and so are the log pairs
$$
K_{\tilde{X}} + \lambda \tilde{D} +  \lambda a_3E_3 + E_2 \text{ ,
since } \lambda a_2 < 1
$$
and
$$
K_{\tilde{X}} + \lambda \tilde{D} + E_3 + \lambda a_2 E_2 \text{ ,
since } \lambda a_3 < 1 \text{ .}
$$
By adjunction it follows that
$$
2a_2 - a_3 -a_1= \tilde{D} \cdot E_2 \geq
\text{mult}_Q\Bigl(\tilde{D} \cdot E_2 \Bigr)  > \frac{1}{\lambda}
- a_3 > \frac{15}{8} - a_3
$$
and
$$
2a_3 - a_2 -\frac{4}{5}a_3 \geq 2a_3 - a_2 -a_4= \tilde{D} \cdot
E_3 \geq \text{mult}_Q\Bigl(\tilde{D} \cdot E_3 \Bigr)  >
\frac{1}{\lambda} - a_2 > \frac{15}{8} - a_2 \text{ .}
$$
This implies that
$$
\frac{3}{2} a_2 \geq 2a_2 - \frac{a_2}{2} \geq 2a_2 - a_1 >
\frac{15}{8} \Rightarrow a_2 \frac{5}{4}
$$ and
$$
\frac{6}{5} a_3 \geq 2a_3 - \frac{4}{5}a_3 \geq 2a_3 - a_4 >
\frac{15}{8} \Rightarrow a_3 \frac{25}{16}
$$ which is false, since either $a_2 \leq 1$ or $a_3 \leq 1$.

\item If $Q\in E_3$ but $Q\not \in E_2 \cup E_4$ then
$$
K_{\tilde{X}} + \lambda \tilde{D} +  \lambda a_3E_3
$$
is not log canonical at the point $Q$ and so is the pair
$$
K_{\tilde{X}} + \lambda \tilde{D} + E_3 \text{ , since  }  \lambda
a_3\leq 1 \text{ .}
$$
By adjunction $(E_3, \lambda \tilde{D}|_{E_3})$ is not log
canonical at $Q$ and
$$\frac{8}{15}a_3 \geq 2a_3 - \frac{2}{3}a_3 -\frac{4}{5}a_3 \geq 2a_3 - a_2 - a_4
= \tilde{D} \cdot E_3 \geq
 \text{mult}_Q\Bigl(\tilde{D} \cdot E_3 \Bigr) > \frac{1}{\lambda} > \frac{15}{8}
\text{, }$$
 which is impossible.

\item If $Q\in E_3 \cap E_4$ then the log pair
$$
K_{\tilde{X}} + \lambda \tilde{D} +  \lambda a_3 E_3 +  \lambda
a_4 E_4
$$
is not log canonical at the point $Q$ and so are the log pairs
$$
K_{\tilde{X}} + \lambda \tilde{D} +  \lambda a_3E_3 + E_4 \text{ ,
since } \lambda a_4 \leq 1
$$
and
$$
K_{\tilde{X}} + \lambda \tilde{D} + E_3 + \lambda a_4 E_4 \text{ ,
since } \lambda a_3 \leq 1 \text{ .}
$$
By adjunction it follows that
$$
2a_4 - a_3 -a_5= \tilde{D} \cdot E_4 \geq
\text{mult}_Q\Bigl(\tilde{D}|_{E_4} \Bigr) =
\text{mult}_Q\Bigl(\tilde{D} \cdot E_4 \Bigr) > \frac{15}{8} - a_3
\text{ ,}
$$
and this implies that
$$
2a_4 - \frac{3}{4} a_4 \geq 2a_4 - a_5 > \frac{15}{8} \Rightarrow
a_4 > \frac{3}{2}
$$
We have that either $a_5 \leq 1$ or $a_3 \leq 1$ and this implies
that $a_4 \leq \frac{4}{3} a_3 \leq \frac{4}{3}$.

\item If $Q\in E_4$ but $Q\not \in E_3 \cup E_5$ then
$$
K_{\tilde{X}} + \lambda \tilde{D} +  \lambda a_4E_4
$$
is not log canonical at the point $Q$ and so is the pair
$$
K_{\tilde{X}} + \lambda \tilde{D} + E_4 \text{ , since  }  \lambda
a_4\leq 1 \text{ .}
$$
By adjunction $(E_4, \lambda \tilde{D}|_{E_4})$ is not log
canonical at $Q$ and
$$1\geq \frac{a_4}{2} \geq 2a_4 - \frac{3}{4} a_4 -\frac{3}{4}a_4 \geq 2a_4 - a_3 - a_5 =
\tilde{D} \cdot E_4 \geq
 \text{mult}_Q\Bigl(\tilde{D} \cdot E_4 \Bigr) > \frac{15}{8} \text{ ,}$$
which is false.

\end{itemize}

\end{proof}

\subsection{Del Pezzo surfaces of degree 1 with exactly one
$\mathbb{A}_8$ type singularity}

In this section we will prove the following.

\begin{lemma}
\label{A8} Let $X$ be a del Pezzo surface with at most one Du Val
singularity of type $\mathbb{A}_8$ and $K_X^2=1$. Then the global
log canonical threshold of $X$ is
$$
\mathrm{lct} (X) = \frac{1}{2} \text{ .}
$$
\end{lemma}

\begin{proof}

Let $X$ be a del Pezzo surface with exactly one Du Val singularity
of type $\mathbb{A}_8$ and $K_X^2=1$. Suppose $\mathrm{lct}(X) <
\frac{1}{2}$. Then there exists an effective $\mathbb{Q}$-divisor
$D\in X$ and a positive rational number $\lambda < \frac{1}{2}$,
such that the log pair $(X, \lambda D)$ is  not log canonical and
$D=-K_X$, where $\lambda < \frac{1}{2}$. Therefore the log pair
$(X, \lambda D)$ is also not log canonical.

Let $Z$ be the curve in $|-K_X|$ that contains $P$. Since the
curve $Z$ is irreducible we may assume that the support of $D$
does not contain $Z$.

We derive that the pair $(X, \lambda D)$ is log canonical outside
of a point $P\in X$ and not log canonical at $P$. Let $\pi_1:
\tilde{X} \to X$ be the minimal resolution of $X$. The following
diagram shows how the exceptional curves intersect each other.
\bigskip

$\mathbb{A}_8$. \xymatrix{ {\bullet}^{E_1} \ar@{-}[r] &
{\bullet}^{E_2} \ar@{-}[r] & {\bullet}^{E_3} \ar@{-}[r] &
{\bullet}^{E_4} \ar@{-}[r] & {\bullet}^{E_5} \ar@{-}[r] &
{\bullet}^{E_6} \ar@{-}[r] & {\bullet}^{E_7} \ar@{-}[r] &
{\bullet}^{E_8}}
\bigskip

Then
\begin{eqnarray*}
\tilde{D} & \sim_{\mathbb{Q}} & \pi_1^*(D) - a_1E_1 - a_2E_2 -
a_3E_3 - a_4E_4 - a_5E_5 -a_6E_6
-a_7E_7-a_8E_8 \text{ and }\\
\tilde{Z} & \sim_{\mathbb{Q}} & \pi_1^*(Z) - E_1 - E_2 - E_3 - E_4
- E_5 - E_6 - E_7- E_8 \text{ .}
\end{eqnarray*}
From the inequalities
\begin{eqnarray*}
0 \leq \tilde{D} \cdot \tilde{Z} & = & 1 - a_1 - a_8\\
0 \leq E_1 \cdot \tilde{D} & = & 2a_1 - a_2\\
0 \leq E_2 \cdot \tilde{D} & = & 2a_2 - a_1 - a_3\\
0 \leq E_3 \cdot \tilde{D} & = & 2a_3 - a_2 - a_4\\
0 \leq E_4 \cdot \tilde{D} & = & 2a_4 - a_3 - a_5\\
0 \leq E_5 \cdot \tilde{D} & = & 2a_5 - a_4 - a_6\\
0 \leq E_6 \cdot \tilde{D} & = & 2a_6 - a_5 - a_7\\
0 \leq E_7 \cdot \tilde{D} & = & 2a_7 - a_6 -a_8\\
0 \leq E_8 \cdot \tilde{D} & = & 2a_8 - a_7\\
\end{eqnarray*}
we get
$$
2a_8 \geq a_7 \text{ , } \frac{3}{2} a_7 \geq a_6 \text{ , }
\frac{4}{3} a_6 \geq a_5 \text{ , } \frac{5}{4} a_5 \geq a_4
\text{ , } \frac{6}{5} a_4 \geq a_3 \text{ , } \frac{7}{6} a_3
\geq a_2  \text{ , } \frac{8}{7} a_2 \geq a_1
$$
and moreover
$$
a_1 \leq \frac{8}{9} \text{ , } a_2 \leq \frac{14}{9} \text{ , }
a_3 \leq 2 \text{ , } a_4 \leq \frac{20}{9} \text{ , } a_5 \leq
\frac{20}{9} \text{ , } a_6 \leq 2 \text{ , } a_7 \leq
\frac{14}{9} \text{ , } a_8 \leq \frac{8}{9}  \text{ .}
$$
Furthermore there are lines $L_3, L_6 \in X$ that pass through the
point $P$ whose strict
 transforms are $(-1)$-curves that intersect the fundamental cycle as following.
 $$ L_3 \cdot E_3 = L_6 \cdot E_6 = 1$$
 and
 $$
 L_i \cdot E_j = 0 \text{ for all } i,j= 3,6 \text{ with } i \not = j \text{ .}
 $$
Then we easily get that
\begin{eqnarray*}
 \tilde{L_3} & \sim_{\mathbb{Q}} & \pi^*(L_3) -\frac{2}{3} E_1 -\frac{4}{3} E_2 - 2 E_3 -\frac{5}{3} E_4
-\frac{4}{3} E_5 - E_6 -\frac{2}{3} E_7-\frac{1}{3}E_8\\
\tilde{L_6} & \sim_{\mathbb{Q}} & \pi^*(L_6) -\frac{1}{3} E_1
-\frac{2}{3} E_2 -E_3 -\frac{4}{3} E_4 -\frac{5}{3} E_5-2 E_6
-\frac{4}{3} E_7-\frac{2}{3}E_8 \text{ .}
\end{eqnarray*}
We observe that $L_3+L_4$ is a Cartier divisor in the
bianticanonical linear system $|-2K_X|$. Since $L_3$ and $L_6$ are
irreducible and $L_3 \sim_{\mathbb{Q}} L_4 \sim_{\mathbb{Q}}
-K_X$, we can assume that $L_3 \not \in \text{Supp}(D)$ and $L_6
\not \in \text{Supp}(D)$. Then
$$
0 \leq \tilde{L_3} \cdot \tilde{D} = 1 - a_3 \text{ and } 0 \leq
\tilde{L_6} \cdot \tilde{D} = 1 - a_6  \text{  .}
$$
The equivalence
$$
K_{\tilde{X}} + \lambda \tilde{D} + \lambda a_1E_1 + \lambda
a_2E_2 + \lambda a_3E_3 + \lambda a_4E_4 + \lambda a_5E_5 +
\lambda a_6E_6 + \lambda a_7E_7+ \lambda a_8E_8 \sim_{\mathbb{Q}}
\pi_1^*(K_X + D)
$$
implies that there is a point $Q\in E_1\cup E_2\cup E_3\cup E_4
\cup E_5 \cup E_6 \cup E_7 \cup E_8$, such that the pair
$$K_{\tilde{X}} + \lambda \tilde{D} + \lambda a_1E_1 + \lambda a_2E_2 + \lambda a_3E_3 +
\lambda a_4E_4 + \lambda a_5E_5 + \lambda a_6E_6 + \lambda a_7E_7
+ \lambda a_8E_8$$ is not log canonical at $Q$.

\begin{itemize}
\item If the point $Q \in E_1$ and $Q\not \in E_2$ then
$$
K_{\tilde{X}} + \lambda \tilde{D} +  a_1 \lambda E_1
$$
is not log canonical at the point $Q$ and so is the pair
$$
K_{\tilde{X}} + \lambda \tilde{D} + E_1 \text{ , since  } a_1
\lambda \leq 1 \text{ .}
$$
By adjunction $(E_1, \lambda \tilde{D}|_{E_1})$ is not log
canonical at $Q$ and
$$ 2a_1 - \frac{7}{8} a_1 \geq
2a_1 - a_2 = \tilde{D} \cdot E_1 \geq
 \text{mult}_Q\Bigl(\tilde{D} \cdot E_1 \Bigr)  > \frac{1}{\lambda} > 2 \text{ ,}$$
which is false since $a_1 \leq \frac{8}{9}$.

\item If $Q\in E_1 \cap E_2$ then the log pair
$$
K_{\tilde{X}} + \lambda \tilde{D} +  \lambda a_1 E_1 +  \lambda
a_2E_2
$$
is not log canonical at the point $Q$ and so is the log pair
$$
K_{\tilde{X}} + \lambda \tilde{D} + E_1 + \lambda a_2 E_2  \text{
, since } \lambda a_1 \leq 1 \text{ .}
$$
By adjunction it follows that
$$  2a_1 - a_2 = \tilde{D} \cdot E_1 \geq
 \text{mult}_Q\Bigl(\tilde{D} \cdot E_1 \Bigr) > \frac{1}{\lambda} - a_2 > 2 - a_2 \text{
,}
$$
which is false since $a_1 \leq 1$.

\item If $Q\in E_2$ but $Q\not \in E_1 \cup E_3$ then
$$
K_{\tilde{X}} + \lambda \tilde{D} +  \lambda a_2 E_2
$$
is not log canonical at the point $Q$ and so is the pair
$$
K_{\tilde{X}} + \lambda \tilde{D} + E_2 \text{ , since  }  \lambda
a_2 \leq  1 \text{ .}
$$
By adjunction $(E_2, \lambda \tilde{D}|_{E_2})$ is not log
canonical at $Q$ and
$$ \frac{9}{14} a_2 \geq 2a_2 - \frac{1}{2} a_2 - \frac{6}{7} a_2 \geq 2a_2 - a_1 - a_3 =
\tilde{D} \cdot E_2 \geq
 \text{mult}_Q\Bigl(\tilde{D} \cdot E_2 \Bigr) > \frac{1}{\lambda} > 2  \text{ ,}$$
which is false, since $a_2 \leq \frac{14}{9}$.

\item If $Q\in E_2 \cap E_3$ then the log pair
$$
K_{\tilde{X}} + \lambda \tilde{D} +  \lambda a_2 E_2 +  \lambda
a_3E_3
$$
is not log canonical at the point $Q$ and so is the log pair
$$
K_{\tilde{X}} + \lambda \tilde{D} +  \lambda a_2E_2 + E_3 \text{ ,
since } \lambda a_3 < 1 \text{ .}
$$
By adjunction it follows that
$$
2a_3 - \frac{5}{6} a_3 -a_2 \geq  2a_3 - a_2 -a_4= \tilde{D} \cdot
E_3 \geq \text{mult}_Q\Bigl(\tilde{D}|_{E_3} \Bigr)  >
\frac{1}{\lambda} - a_2 > 2 - a_2 \text{ ,}
$$
which is impossible, since $a_3 \leq 1$.

\item If $Q\in E_3$ but $Q\not \in E_2 \cup E_4$ then
$$
K_{\tilde{X}} + \lambda \tilde{D} +  \lambda a_3E_3
$$
is not log canonical at the point $Q$ and so is the pair
$$
K_{\tilde{X}} + \lambda \tilde{D} + E_3 \text{ , since  }  \lambda
a_3\leq 1 \text{ .}
$$
By adjunction $(E_3, \lambda \tilde{D}|_{E_3})$ is not log
canonical at $Q$ and
$$2a_3 - \frac{2}{3} a_3 - \frac{5}{6} a_3 \geq 2a_3 - a_2 - a_4 = \tilde{D} \cdot E_3
\geq \text{mult}_Q\Bigl(\tilde{D}|_{E_3} \Bigr) >
\frac{1}{\lambda} > 2 \text{ .}$$ and this  implies that $a_3 >
4$, which is impossible.

\item If $Q\in E_3 \cap E_4$ then the log pair
$$
K_{\tilde{X}} + \lambda \tilde{D} +  \lambda a_3 E_3 +  \lambda
a_4 E_4
$$
is not log canonical at the point $Q$ and so is the log pair
$$
K_{\tilde{X}} + \lambda \tilde{D} +  E_3 + \lambda a_4 E_4 \text{
, since } \lambda a_3 \leq 1 \text{ .}
$$
By adjunction it follows that

$$
2a_3 - \frac{2}{3} a_3 - a_4\geq 2a_3 - a_2 -a_4= \tilde{D} \cdot
E_3 \geq \text{mult}_Q\Bigl(\tilde{D}|_{E_3} \Bigr)  >
\frac{1}{\lambda} - a_4 > 2 - a_4 \text{ ,}
$$
which contradicts $a_3 \leq 1$.

\item If $Q\in E_4$ but $Q\not \in E_3 \cup E_5$ then
$$
K_{\tilde{X}} + \lambda \tilde{D} +  \lambda a_4E_4
$$
is not log canonical at the point $Q$ and so is the pair
$$
K_{\tilde{X}} + \lambda \tilde{D} + E_4 \text{ , since  }  \lambda
a_4\leq \frac{4}{3} \lambda a_3 \leq \frac{8}{9} \text{ .}
$$
By adjunction $(E_4, \lambda \tilde{D}|_{E_4})$ is not log
canonical at $Q$ and
$$2a_4 - \frac{3}{4} a_4 - \frac{4}{5} a_4 \geq 2a_4 - a_3 - a_5 = \tilde{D} \cdot E_4
\geq \text{mult}_Q\Bigl(\tilde{D}|_{E_4} \Bigr) >
\frac{1}{\lambda} > 2 \text{ ,}$$ implies that $a_4 >
\frac{40}{9}$ which is false.

\item If $Q\in E_4 \cap E_5$ then the log pair
$$
K_{\tilde{X}} + \lambda \tilde{D} +  \lambda a_4 E_4 +  \lambda
a_5 E_5
$$
is not log canonical at the point $Q$ and so is the log pair
$$
K_{\tilde{X}} + \lambda \tilde{D}  +  E_4 + \lambda a_5 E_5 \text{
, since } \lambda a_4 \leq \frac{4}{3} \lambda a_3 < 1 \text{ .}
$$
By adjunction it follows that
$$
2a_4 - \frac{3}{4} a_4 - a_5\geq 2a_4 - a_3 -a_5= \tilde{D} \cdot
E_4 \geq \text{mult}_Q\Bigl(\tilde{D}|_{E_4} \Bigr)  >
\frac{1}{\lambda} - a_5 > 2 - a_5 \text{ ,}
$$
which implies $a_4 > \frac{8}{5}$ and contradicts $a_4 \leq
\frac{4}{3}a_3 \leq \frac{4}{3}$.
\end{itemize}
\end{proof}

\subsection{Del Pezzo surfaces of degree 1 with exactly one
$\mathbb{D}_4$ type singularity}

In this section we will prove the following.

\begin{lemma}
\label{degree1D4} Let $X$ be a del Pezzo surface with at most one
Du Val singularity of type $\mathbb{D}_4$ and $K_X^2 = 1$. Then
the global log canonical threshold of $X$ is
$$
\mathrm{lct} (X) =  \frac{1}{2} \text{ .}
$$
\end{lemma}

\begin{proof}

 Suppose $\mathrm{lct}(X) <
\frac{1}{2} $. Then there exists an effective $\mathbb{Q}$-divisor
$D\in X$ and a rational number $\lambda < \frac{1}{2}$, such that
the log pair $(X,\lambda D)$ is  not log canonical and $D
\sim_{\mathbb{Q}} -K_X$.

Let $Z$ be the unique curve in $|-K_X|$ that contains $P$. Since
the curve $Z$ is irreducible we may assume that the support of $D$
does not contain $Z$.

We derive that the pair $(X,\lambda D)$ is log canonical
everywhere outside of a singular point $P\in X$ and is not log
canonical at $P$. Let $\pi: \tilde{X} \to X$ be the minimal
resolution of $X$. The following diagram shows how the exceptional
curves intersect each other.
\bigskip

$\mathbb{D}_4$. 
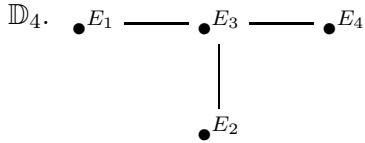
\xymatrix{ {\bullet}^{E_1} \ar@{-}[r] &
{\bullet}^{E_3} \ar@{-}[r] \ar@{-}[d] & {\bullet}^{E_4}\\ &
{\bullet}^{E_2} &}
\bigskip

Then
$$
\tilde{D} \sim_{\mathbb{Q}} \pi^*(D)-a_1E_1-a_2E_2-a_3E_3-a_4E_4
\text{ and } \tilde{Z} \sim_{\mathbb{Q}}
\pi^*(Z)-E_1-E_2-2E_3-E_4\text{ .}
$$
From the inequalities
\begin{eqnarray*}
0 \leq \tilde{D} \cdot \tilde{Z} & = & 1 - a_3\\
0 \leq E_1 \cdot \tilde{D} & = & 2a_1 - a_3\\
0 \leq E_2 \cdot \tilde{D} & = & 2a_2 - a_3\\
0 \leq E_3 \cdot \tilde{D} & = & 2a_3 - a_1 - a_2 - a_4\\
0 \leq E_4 \cdot \tilde{D} & = & 2a_4 - a_3
\end{eqnarray*}
we get the following upper bounds $a_1 \leq 1 \text{, } a_2 \leq 1
\text{, } a_3 \leq 1 \text{, } a_4 \leq 1$. The equivalence
$$
K_{\tilde{X}} + \lambda \tilde{D} + \lambda a_1E_1 + \lambda
a_2E_2 +  \lambda a_3E_3 + \lambda a_4E_4 \sim_{\mathbb{Q}}
\pi^*(K_X+\lambda D)
$$
implies that there is a point $Q\in E_1\cup E_2\cup E_3\cup E_4$
such that the pair $K_{\tilde{X}} + \lambda \tilde{D} + \lambda
a_1E_1 + \lambda a_2E_2 +  \lambda a_3E_3 + \lambda a_4E_4$ is not
log canonical at $Q$.

\begin{itemize}
\item If the point $Q \in E_1$ and $Q\not \in E_3$ then
$$
K_{\tilde{X}} + \lambda \tilde{D} + \lambda a_1E_1
$$
is not log canonical at the point $Q$ and so is the pair
$$
K_{\tilde{X}} + \lambda \tilde{D} + E_1 \text{ .}
$$
By adjunction $(E_1, \lambda \tilde{D}|_{E_1})$ is not log
canonical at $Q$ and
$$2a_1-a_3 = \tilde{D} \cdot E_1 \geq \text{mult}_Q\Bigl(\tilde{D}|_{E_1} \Bigr) =
 \text{mult}_Q\Bigl(\tilde{D} \cdot E_1 \Bigr)  > 2 \text{ ,}$$
which along with the inequalities $a_3\leq 2 a_1  \text{, } a_3
\leq 2 a_2  \text{, } a_1 + a_2 + a_4 \leq 2a_3 \text{, } a_3 \leq
2 a_4$ implies that $a_1>2$ which is false.

\item If $Q\in E_3$ but $Q\not \in E_1 \cup E_2 \cup E_4$ then
$$
K_{\tilde{X}} + \lambda \tilde{D} +  \lambda a_3E_3
$$
is not log canonical at the point $Q$ and so is the pair
$$
K_{\tilde{X}} + \lambda \tilde{D} + E_3 \text{, since } \lambda
a_3 \leq 1 \text{ .}
$$
By adjunction $(E_3, \lambda \tilde{D}|_{E_3})$ is not log
canonical at $Q$ and
$$2a_3 - a_1 - a_2 - a_4 = \tilde{D} \cdot E_1 \geq \text{mult}_Q\Bigl(\tilde{D}|_{E_3}
\Bigr) =
 \text{mult}_Q\Bigl(\tilde{D} \cdot E_3 \Bigr) > 2 \text{ ,}$$
which along with $a_3\leq 2 a_1  \text{, } a_3 \leq 2a_2  \text{,
} a_3 \leq 2a_4$ implies that $a_3 >4$ which is false.

\item If $Q\in E_1 \cap E_3$ then the log pair
$$
K_{\tilde{X}} + \lambda \tilde{D} + \lambda a_1E_1 +  \lambda
a_3E_3
$$
is not log canonical at the point $Q$ and so is the log pair
$$
K_{\tilde{X}} + \lambda \tilde{D} + E_1 +  \lambda a_3E_3 \text{
.}
$$
By adjunction it follows that
$$
2a_1-a_3 = \tilde{D} \cdot E_1 \geq
\text{mult}_Q\Bigl(\tilde{D}|_{E_1} \Bigr) =
 \text{mult}_Q\Bigl(\tilde{D} \cdot E_1 \Bigr) > 2 - a_3 \text{ .}
$$
and this implies that $a_1>1$ which is not possible.

\end{itemize}
\end{proof}

\subsection{Del Pezzo surfaces of degree 1 with exactly one
$\mathbb{D}_5$ singularity.}

In this section we will prove the following.

\begin{lemma}
\label{D5} Let $X$ be a del Pezzo surface with exactly one Du Val
singularity of type $\mathbb{D}_5$ and $K_X^2=1$. Then the global
log canonical threshold of $X$ is
$$
\mathrm{lct} (X) = \frac{1}{2} \text{ .}
$$
\end{lemma}

\begin{proof}

Suppose that $\mathrm{lct}(X)<\frac{1}{2}$,  then there exists a
$\mathbb{Q}$-divisor $D \in X$ and a rational number $\lambda <
\frac{1}{2}$, such that the log pair $(X,\lambda D)$ is  not log
canonical and $D \sim_{\mathbb{Q}} -K_X$. We derive that the pair
$(X,\lambda D)$ is log canonical outside of a point $P\in X$ and
not log canonical at $P$. Let $\pi: \tilde{X} \to X$ be the
minimal resolution of $X$. The configuration of the exceptional
curves is given by the following Dynkin diagram.
\bigskip

$\mathbb{D}_5$. 
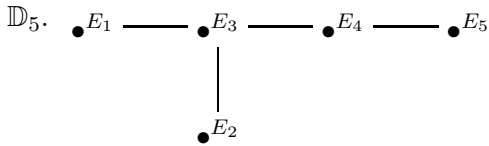
\xymatrix{ {\bullet}^{E_1} \ar@{-}[r] &
{\bullet}^{E_3} \ar@{-}[r] \ar@{-}[d] & {\bullet}^{E_4} \ar@{-}[r]
& {\bullet}^{E_5}\\ & {\bullet}^{E_2} &}
\bigskip

Then
$$
\tilde{D} \sim_{\mathbb{Q}}
\pi^*(D)-a_1E_1-a_2E_2-2a_3E_3-2a_4E_4-a_5E_5 \text{ and }
\tilde{Z} \sim_{\mathbb{Q}} \pi^*(Z)-E_1-E_2-2E_3-2E_4-E_5\text{
.}
$$
From the inequalities
\begin{eqnarray*}
0 \leq \tilde{D} \cdot \tilde{Z} & = & 1 - 2a_4\\
0 \leq E_1 \cdot \tilde{D} & = & 2a_1 - 2a_3\\
0 \leq E_2 \cdot \tilde{D} & = & 2a_2 - 2a_3\\
0 \leq E_3 \cdot \tilde{D} & = & 4a_3 - a_1 - a_2 - 2a_4\\
0 \leq E_4 \cdot \tilde{D} & = & 4a_4 - 2a_3 - a_5\\
0 \leq E_5 \cdot \tilde{D} & = & 2a_5 - 2a_4\\
\end{eqnarray*}
we see that $a_1 \leq \frac{5}{4} \text{, } a_2 \leq \frac{5}{4}
\text{, } a_3 \leq \frac{3}{4} \text{, } a_4 \leq \frac{1}{2}
\text{, } a_5 \leq 1$. The equivalence
$$
K_{\tilde{X}} + \lambda \tilde{D} + \lambda a_1E_1 + \lambda
a_2E_2 + 2 \lambda a_3E_3 + 2\lambda a_4E_4 + \lambda a_5E_5
\sim_{\mathbb{Q}} \pi^*(K_X+\lambda D)
$$
implies that there is a point $Q\in E_1\cup E_2\cup E_3\cup E_4
\cup E_5$ such that the pair
$$
K_{\tilde{X}} + \lambda \tilde{D} + \lambda a_1E_1 + \lambda
a_2E_2 + 2 \lambda a_3E_3 + 2\lambda a_4E_4 + \lambda a_5E_5
$$
is not log canonical at $Q$.
\begin{itemize}
\item If the point $Q \in E_1 \backslash E_3$ then
$$
K_{\tilde{X}} + \lambda \tilde{D} + \lambda a_1E_1
$$
is not log canonical at the point $Q$ and so is the pair
$$
K_{\tilde{X}} + \lambda \tilde{D} + E_1 \text{ .}
$$
By adjunction $(E_1, \lambda \tilde{D}|_{E_1})$ is not log
canonical at $Q$ and
$$2a_1-2a_3 = \tilde{D} \cdot E_1 \geq \text{mult}_Q\Bigl(\tilde{D}|_{E_1} \Bigr) =
 \text{mult}_Q\Bigl(\tilde{D} \cdot E_1 \Bigr)  > 2 \text{ ,}$$
which along with the inequalities $a_3\leq a_1  \text{, } a_3 \leq
a_2  \text{, } a_1 + a_2 + 2a_4 \leq 4a_3 \text{, } a_4 \leq a_5
\text{, } 4a_4 - 2a_3 - a_5$ implies that $a_1>\frac{5}{2}$ which
is false.

\item If $Q\in E_3$ but $Q\not \in E_1 \cup E_2 \cup E_4$ then
$$
K_{\tilde{X}} + \lambda \tilde{D} + 2 \lambda a_3E_3
$$
is not log canonical at the point $Q$ and so is the pair
$$
K_{\tilde{X}} + \lambda \tilde{D} + E_3 \text{, since } a_3 \leq 1
\text{ .}
$$
By adjunction $(E_3, \lambda \tilde{D}|_{E_3})$ is not log
canonical at $Q$ and
$$4a_3 - a_1 - a_2 - 2a_4 = \tilde{D} \cdot E_3 \geq \text{mult}_Q\Bigl(\tilde{D}|_{E_3}
\Bigr) =
 \text{mult}_Q\Bigl(\tilde{D} \cdot E_3 \Bigr) > 2 \text{ ,}$$
which along with the inequalities $a_3\leq a_1  \text{, } a_3 \leq
a_2 \text{, } a_1 + a_2 + 2a_4 \leq 4a_3 \text{, } a_4 \leq a_5
\text{, } 4a_4 - 2a_3 - a_5$ implies that $a_1>3$ which is false.

\item If $Q\in E_1 \cap E_3$ then the log pair
$$
K_{\tilde{X}} + \lambda \tilde{D} + \lambda a_1E_1 + 2 \lambda
a_3E_3
$$
is not log canonical at the point $Q$ and so are the log pairs
$$
K_{\tilde{X}} + \lambda \tilde{D} + E_1 + 2 \lambda a_3E_3 \text{
and  } K_{\tilde{X}} + \lambda \tilde{D} + \lambda a_1E_1 + E_3
\text{ .}
$$
By adjunction it follows that
$$
2a_1-2a_3 = \tilde{D} \cdot E_1 \geq
\text{mult}_Q\Bigl(\tilde{D}|_{E_1} \Bigr) =
 \text{mult}_Q\Bigl(\tilde{D} \cdot E_1 \Bigr) > 2 - 2a_3 \text{ .}
$$
and
$$
4a_3 - a_1 - a_2 - 2a_4 = \tilde{D} \cdot E_3 \geq
\text{mult}_Q\Bigl(\tilde{D}|_{E_3} \Bigr) =
 \text{mult}_Q\Bigl(\tilde{D} \cdot E_3 \Bigr) > 2 - a_1 \text{ .}
$$
which is not possible.

\item $Q\in E_5 \backslash E_4$ then the log pair
$$
K_{\tilde{X}} + \lambda \tilde{D} + \lambda a_5E_5
$$
is not log canonical at the point $Q$ and so is the pair
$$
K_{\tilde{X}} + \lambda \tilde{D} + E_5 \text{ .}
$$
By adjunction $(E_5, \lambda \tilde{D}|_{E_5})$ is not log
canonical at $Q$ and
$$2a_5-2a_4 = \tilde{D} \cdot E_5 \geq \text{mult}_Q\Bigl(\tilde{D}|_{E_5} \Bigr) =
 \text{mult}_Q\Bigl(\tilde{D} \cdot E_5 \Bigr) > 2 \text{ ,}$$
which along with the inequalities
$$
a_3 \leq a_1 \text{, } a_3 \leq a_2 \text{, } a_1+a_2+2a_4 \leq
4a_3 \text{, } 2a_3+a_5 \leq 4a_4 \text{, } a_4 \leq a_5
$$
implies that $a_5>2$ which is false.

\item $Q\in E_4 \backslash (E_3 \cap E_5)$ then the log pair
$$
K_{\tilde{X}} + \lambda \tilde{D} + \lambda a_4E_4
$$
is not log canonical at the point $Q$ and so is the pair
$$
K_{\tilde{X}} + \lambda \tilde{D} + E_4 \text{ .}
$$
By adjunction $(E_4, \lambda \tilde{D}|_{E_4})$ is not log
canonical at $Q$ and
$$4a_4 - 2a_3 - a_5 = \tilde{D} \cdot E_4 \geq \text{mult}_Q\Bigl(\tilde{D}|_{E_4} \Bigr)
=
 \text{mult}_Q\Bigl(\tilde{D} \cdot E_4 \Bigr) > 2 \text{ ,}$$
which along with the inequalities
$$
a_3 \leq a_1 \text{, } a_3 \leq a_2 \text{, } a_1+a_2+2a_4 \leq
4a_3 \text{, } 2a_3+a_5 \leq 4a_4 \text{, } a_4 \leq a_5
$$
implies that $a_4>2$ which is false.

\item $Q\in E_4 \cap E_5$ then the log pair
$$
K_{\tilde{X}} + \lambda \tilde{D} +2\lambda a_4E_4 + \lambda
a_5E_5
$$
is not log canonical at the point $Q$ and so is the log pair
$$
K_{\tilde{X}} + \lambda \tilde{D} + E_5 + 2\lambda a_4E_4 \text{
.}
$$
By adjunction it follows that
$$
2a_5-2a_4 = \tilde{D} \cdot E_5 \geq
\text{mult}_Q\Bigl(\tilde{D}|_{E_5} \Bigr) =
 \text{mult}_Q\Bigl(\tilde{D} \cdot E_5 \Bigr) > 2 - 2a_4 \text{ .}
$$
and we see then that $a_5>1$ which is not possible .
\end{itemize}
\end{proof}

\subsection{Del Pezzo surfaces of degree 1 with exactly one
$\mathbb{D}_6$ singularity.}

In this section we will prove the following.

\begin{lemma}
\label{D6} Let $X$ be a del Pezzo surface with exactly one Du Val
singularity of type $\mathbb{D}_6$ and $K_X^2=1$. Then the global
log canonical threshold of $X$ is
$$
\mathrm{lct} (X) = \frac{1}{2} \text{ .}
$$
\end{lemma}

\begin{proof}

Suppose that $\mathrm{lct}(X)<\frac{1}{2}$,  then there exists a
$\mathbb{Q}$-divisor $D \in X$ and a rational number $\lambda <
\frac{1}{2}$, such that the log pair $(X,\lambda D)$ is  not log
canonical and $D \sim_{\mathbb{Q}} -K_X$. We derive that the pair
$(X,\lambda D)$ is log canonical outside of a point $P\in X$ and
not log canonical at $P$. Let $\pi: \tilde{X} \to X$ be the
minimal resolution of $X$. The configuration of the exceptional
curves is given by the following Dynkin diagram.
\bigskip

$\mathbb{D}_6$. 
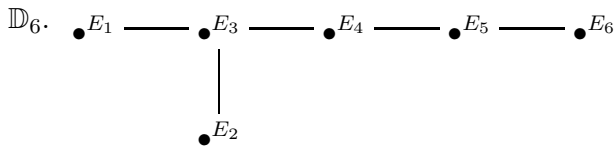
\xymatrix{ {\bullet}^{E_1} \ar@{-}[r] &
{\bullet}^{E_3} \ar@{-}[r] \ar@{-}[d] & {\bullet}^{E_4} \ar@{-}[r]
& {\bullet}^{E_5} \ar@{-}[r] & {\bullet}^{E_6}\\ & {\bullet}^{E_2}
&}
\bigskip

Then
$$
\tilde{D} \sim_{\mathbb{Q}}
\pi^*(D)-a_1E_1-a_2E_2-2a_3E_3-2a_4E_4-2a_5E_5-a_6E_6 \text{ and }
\tilde{Z} \sim_{\mathbb{Q}}
\pi^*(Z)-E_1-E_2-2E_3-2E_4-2E_5-E_6\text{ .}
$$
From the inequalities
\begin{eqnarray*}
0 \leq \tilde{D} \cdot \tilde{Z} & = & 1 - 2a_5\\
0 \leq E_1 \cdot \tilde{D} & = & 2a_1 - 2a_3\\
0 \leq E_2 \cdot \tilde{D} & = & 2a_2 - 2a_3\\
0 \leq E_3 \cdot \tilde{D} & = & 4a_3 - a_1 - a_2 - 2a_4\\
0 \leq E_4 \cdot \tilde{D} & = & 4a_4 - 2a_3 - 2a_5\\
0 \leq E_5 \cdot \tilde{D} & = & 4a_5 - 2a_4 - a_6\\
0 \leq E_6 \cdot \tilde{D} & = & 2a_6 - 2a_5
\end{eqnarray*}
we see that $a_1 \leq \frac{3}{2} \text{, } a_2 \leq \frac{3}{2}
\text{, } a_3 \leq 1 \text{, } a_4 \leq \frac{3}{4} \text{, } a_5
\leq \frac{1}{2} \text{, } a_6 \leq 1$. The equivalence
$$
K_{\tilde{X}} + \lambda \tilde{D} + \lambda a_1E_1 + \lambda
a_2E_2 + 2 \lambda a_3E_3 + 2\lambda a_4E_4 + 2 \lambda a_5E_5 +
\lambda a_6E_6 \sim_{\mathbb{Q}} \pi^*(K_X+\lambda D)
$$
implies that there is a point $Q\in E_1\cup E_2\cup E_3\cup E_4
\cup E_5 \cup E_6$ such that the pair
$$
K_{\tilde{X}} + \lambda \tilde{D} + \lambda a_1E_1 + \lambda
a_2E_2 + 2 \lambda a_3E_3 + 2\lambda a_4E_4 + 2\lambda a_5E_5 +
\lambda a_6E_6
$$
is not log canonical at $Q$.

\begin{itemize}
\item If the point $Q \in E_1 \backslash E_3$ then
$$
K_{\tilde{X}} + \lambda \tilde{D} + \lambda a_1E_1
$$
is not log canonical at the point $Q$ and so is the pair
$$
K_{\tilde{X}} + \lambda \tilde{D} + E_1 \text{ .}
$$
By adjunction $(E_1, \lambda \tilde{D}|_{E_1})$ is not log
canonical at $Q$ and
$$2a_1-2a_3 = \tilde{D} \cdot E_1 \geq \text{mult}_Q\Bigl(\tilde{D}|_{E_1} \Bigr) =
 \text{mult}_Q\Bigl(\tilde{D} \cdot E_1 \Bigr)  > 2 \text{ ,}$$
which is false.

\item If $Q\in E_3$ but $Q\not \in E_1 \cup E_2 \cup E_4$ then
$$
K_{\tilde{X}} + \lambda \tilde{D} + 2 \lambda a_3E_3
$$
is not log canonical at the point $Q$ and so is the pair
$$
K_{\tilde{X}} + \lambda \tilde{D} + E_3 \text{, since } a_3 \leq 1
\text{ .}
$$
By adjunction $(E_3, \lambda \tilde{D}|_{E_3})$ is not log
canonical at $Q$ and
$$4a_3 - a_1 - a_2 - 2a_4 = \tilde{D} \cdot E_3 \geq \text{mult}_Q\Bigl(\tilde{D}|_{E_3}
\Bigr) =
 \text{mult}_Q\Bigl(\tilde{D} \cdot E_3 \Bigr) > 2 \text{ ,}$$
which is false.

\item If $Q\in E_1 \cap E_3$ then the log pair
$$
K_{\tilde{X}} + \lambda \tilde{D} + \lambda a_1E_1 + 2 \lambda
a_3E_3
$$
is not log canonical at the point $Q$ and so are the log pairs
$$
K_{\tilde{X}} + \lambda \tilde{D} + E_1 + 2 \lambda a_3E_3 \text{
and } K_{\tilde{X}} + \lambda \tilde{D} + \lambda a_1E_1 + E_3
\text{ .}
$$
By adjunction it follows that
$$
2a_1-2a_3 = \tilde{D} \cdot E_1 \geq
\text{mult}_Q\Bigl(\tilde{D}|_{E_1} \Bigr) =
 \text{mult}_Q\Bigl(\tilde{D} \cdot E_1 \Bigr) > 2 - 2a_3 \text{ .}
$$
and
$$
4a_3 - a_1 - a_2 - 2a_4 = \tilde{D} \cdot E_3 \geq
\text{mult}_Q\Bigl(\tilde{D}|_{E_3} \Bigr) =
 \text{mult}_Q\Bigl(\tilde{D} \cdot E_3 \Bigr) > 2 - a_1
$$
which is false.

\item If $Q\in E_3 \cap E_4$ then the log pair
$$
K_{\tilde{X}} + \lambda \tilde{D} + 2 \lambda a_3E_3 + 2 \lambda
a_4E_4
$$
is not log canonical at the point $Q$ and so are the log pairs
$$
K_{\tilde{X}} + \lambda \tilde{D} + E_3 + 2 \lambda a_4E_4 \text{
and } K_{\tilde{X}} + \lambda \tilde{D} + 2\lambda a_3E_3 + E_4
\text{ .}
$$
By adjuction
$$
4a_3 - a_1 - a_2 - 2a_4 = \tilde{D} \cdot E_3 \geq
\text{mult}_Q\Bigl(\tilde{D}|_{E_3} \Bigr) =
 \text{mult}_Q\Bigl(\tilde{D} \cdot E_3 \Bigr) > 2 - 2a_4
$$
and
$$
4a_4 - 2a_3 - 2a_5 = \tilde{D} \cdot E_4 \geq
\text{mult}_Q\Bigl(\tilde{D}|_{E_4} \Bigr) =
 \text{mult}_Q\Bigl(\tilde{D} \cdot E_4 \Bigr) > 2 - 2a_3 \text{ .}
$$
which is false.

\item If the point $Q \in E_4 \backslash (E_3 \cup E_5)$ then
$$
K_{\tilde{X}} + \lambda \tilde{D} + 2\lambda a_4E_4
$$
is not log canonical at the point $Q$ and so is the pair
$$
K_{\tilde{X}} + \lambda \tilde{D} + E_4 \text{ .}
$$
By adjunction $(E_4, \lambda \tilde{D}|_{E_4})$ is not log
canonical at $Q$ and
$$4a_4 - 2a_3 - 2a_5 = \tilde{D} \cdot E_4 \geq \text{mult}_Q\Bigl(\tilde{D}|_{E_4}
\Bigr) =
 \text{mult}_Q\Bigl(\tilde{D} \cdot E_4 \Bigr)  > 2 \text{ ,}$$
which is false.

\item $Q\in E_5 \backslash (E_4 \cup E_6)$ then the log pair
$$
K_{\tilde{X}} + \lambda \tilde{D} + 2 \lambda a_5E_5
$$
is not log canonical at the point $Q$ and so is the pair
$$
K_{\tilde{X}} + \lambda \tilde{D} + E_5 \text{ .}
$$
By adjunction $(E_5, \lambda \tilde{D}|_{E_5})$ is not log
canonical at $Q$ and
$$4a_5 - 2a_4 - a_6 = \tilde{D} \cdot E_5 \geq \text{mult}_Q\Bigl(\tilde{D}|_{E_5} \Bigr)
=
 \text{mult}_Q\Bigl(\tilde{D} \cdot E_5 \Bigr) > 2 \text{ ,}$$
which is false.

\item If $Q\in E_4 \cap E_5$ then the log pair
$$
K_{\tilde{X}} + \lambda \tilde{D} + 2 \lambda a_4E_4 + 2 \lambda
a_5E_5
$$
is not log canonical at the point $Q$ and so are the log pairs
$$
K_{\tilde{X}} + \lambda \tilde{D} + E_4 + 2 \lambda a_5E_5 \text{
and } K_{\tilde{X}} + \lambda \tilde{D} + 2\lambda a_4E_4 + E_5
\text{ .}
$$
By adjuction
$$
4a_4 - 2a_3 - 2a_5 = \tilde{D} \cdot E_4 \geq
\text{mult}_Q\Bigl(\tilde{D}|_{E_4} \Bigr) =
 \text{mult}_Q\Bigl(\tilde{D} \cdot E_4 \Bigr) > 2 - 2a_5
$$
and
$$
4a_5 - 2a_4 - a_6 = \tilde{D} \cdot E_5 \geq
\text{mult}_Q\Bigl(\tilde{D}|_{E_5} \Bigr) =
 \text{mult}_Q\Bigl(\tilde{D} \cdot E_5 \Bigr) > 2 - 2a_4 \text{ .}
$$\item $Q\in E_5 \backslash (E_4 \cup E_6)$ then the log pair
$$
K_{\tilde{X}} + \lambda \tilde{D} + 2 \lambda a_5E_5
$$
is not log canonical at the point $Q$ and so is the pair
$$
K_{\tilde{X}} + \lambda \tilde{D} + E_5 \text{ .}
$$
By adjunction $(E_5, \lambda \tilde{D}|_{E_5})$ is not log
canonical at $Q$ and
$$4a_5 - 2a_4 - a_6 = \tilde{D} \cdot E_5 \geq \text{mult}_Q\Bigl(\tilde{D}|_{E_5} \Bigr)
=
 \text{mult}_Q\Bigl(\tilde{D} \cdot E_5 \Bigr) > 2 \text{ ,}$$
which is false.

\item If $Q\in E_5 \cap E_6$ then the log pair
$$
K_{\tilde{X}} + \lambda \tilde{D} + 2 \lambda a_5E_5 + \lambda
a_6E_6
$$
is not log canonical at the point $Q$ and so is the log pair
$$
K_{\tilde{X}} + \lambda \tilde{D} + 2 \lambda a_5E_5 + E_6 \text{
.}
$$
By adjuction
$$
2a_6 - 2a_5 = \tilde{D} \cdot E_6 \geq
\text{mult}_Q\Bigl(\tilde{D}|_{E_6} \Bigr) =
 \text{mult}_Q\Bigl(\tilde{D} \cdot E_6 \Bigr) > 2 - 2a_5 \text{ .}
$$
which is false.

\item If the point $Q \in E_6 \backslash E_5$ then
$$
K_{\tilde{X}} + \lambda \tilde{D} + \lambda a_6E_6
$$
is not log canonical at the point $Q$ and so is the pair
$$
K_{\tilde{X}} + \lambda \tilde{D} + E_6 \text{ .}
$$
By adjunction $(E_6, \lambda \tilde{D}|_{E_6})$ is not log
canonical at $Q$ and
$$2a_6-2a_5 = \tilde{D} \cdot E_6 \geq \text{mult}_Q\Bigl(\tilde{D}|_{E_6} \Bigr) =
 \text{mult}_Q\Bigl(\tilde{D} \cdot E_6 \Bigr)  > 2 \text{ ,}$$
which is false.

\end{itemize}
\end{proof}

\subsection{Del Pezzo surfaces of degree 1 with exactly one
$\mathbb{D}_7$ singularity.}

In this section we will prove the following.

\begin{lemma}
\label{D7} Let $X$ be a del Pezzo surface with exactly one Du Val
singularity of type $\mathbb{D}_7$ and $K_X^2=1$. Then the global
log canonical threshold of $X$ is
$$
\mathrm{lct} (X) = \frac{2}{5} \text{ .}
$$
\end{lemma}

\begin{proof}

Suppose that $\mathrm{lct}(X)<\frac{2}{5}$,  then there exist a
$\mathbb{Q}$-divisor $D \in X$ and a rational number $\lambda <
\frac{2}{5}$, such that the log pair $(X,\lambda D)$ is not log
canonical and $D \sim_{\mathbb{Q}} -K_X$. We derive that the pair
$(X,\lambda D)$ is log canonical outside of a point $P\in X$ and
not log canonical at $P$. Let $\pi: \tilde{X} \to X$ be the
minimal resolution of $X$. The configuration of the exceptional
curves is given by the following Dynkin diagram.
\bigskip

$\mathbb{D}_7$. 
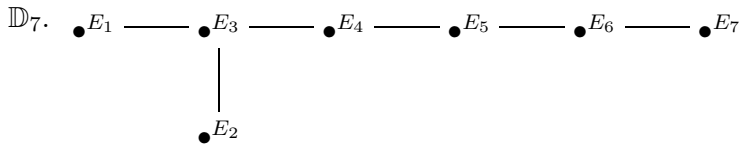
\xymatrix{ {\bullet}^{E_1} \ar@{-}[r] & {\bullet}^{E_3}
\ar@{-}[r] \ar@{-}[d] & {\bullet}^{E_4} \ar@{-}[r] & {\bullet}^{E_5}
\ar@{-}[r] & {\bullet}^{E_6} \ar@{-}[r] & {\bullet}^{E_7} \\ & {\bullet}^{E_2} &}
\bigskip

Then
\begin{eqnarray*}
\tilde{D} & \sim_{\mathbb{Q}} &
\pi^*(D)-a_1E_1-a_2E_2-a_3E_3-a_4E_4-a_5E_5-a_6E_6-a_7E_7 \text{
and
}\\
\tilde{Z} & \sim_{\mathbb{Q}} & \pi^*(Z)-E_1-E_2-2E_3-2E_4-2E_5-2E_6-E_7\text{ .}\\
\end{eqnarray*}

From the inequalities
\begin{eqnarray*}
0 \leq \tilde{D} \cdot \tilde{Z} & = & 1 - a_6\\
0 \leq E_1 \cdot \tilde{D} & = & 2a_1 - a_3\\
0 \leq E_2 \cdot \tilde{D} & = & 2a_2 - a_3\\
0 \leq E_3 \cdot \tilde{D} & = & 2a_3 - a_1 - a_2 - a_4\\
0 \leq E_4 \cdot \tilde{D} & = & 2a_4 - a_3 - a_5\\
0 \leq E_5 \cdot \tilde{D} & = & 2a_5 - a_4 - a_6\\
0 \leq E_6 \cdot \tilde{D} & = & 2a_6 - a_5 - a_7\\
0 \leq E_7 \cdot \tilde{D} & = & 2a_7 - a_6\\
\end{eqnarray*}
we see that $a_1 \leq \frac{7}{4} \text{, } a_2 \leq \frac{7}{4}
\text{, } a_3 \leq \frac{5}{2} \text{, } a_4 \leq \frac{8}{4}=2
\text{, } a_5 \leq \frac{6}{4} \text{, } a_6 \leq \frac{4}{4}=1
\text{, } a_7 \leq 1 \text{ .}$ Moreover we get the inequalities
$$
2a_7 \geq a_6 \text{, } \frac{3}{2} a_6 \geq a_5 \text{, }
\frac{4}{3} a_5 \geq a_4 \text{, } \frac{5}{4} a_4 \geq a_3
$$
and
$$
2a_1 \geq a_3 \text{, } 2a_2 \geq a_3 \text{, } a_3 \geq a_4 \geq
a_5 \geq a_6 \geq a_7 \text{ .}
$$
The equivalence
$$
K_{\tilde{X}} + \lambda \tilde{D} + \lambda a_1E_1 + \lambda
a_2E_2 + \lambda a_3E_3 + \lambda a_4E_4 + \lambda a_5E_5 +
\lambda a_6E_6 + \lambda a_7E_7 \sim_{\mathbb{Q}} \pi^*(K_X +
\lambda D)
$$
implies that there is a point $Q\in E_1\cup E_2\cup E_3\cup E_4
\cup E_5 \cup E_6 \cup E_7$ such that the pair
$$
K_{\tilde{X}} + \lambda \tilde{D} + \lambda a_1E_1 + \lambda
a_2E_2 + \lambda a_3E_3 + \lambda a_4E_4 + \lambda a_5E_5 +
\lambda a_6E_6 + \lambda a_7E_7
$$
is not log canonical at $Q$.

\begin{itemize}
\item If the point $Q \in E_1 \backslash E_3$ then
$$
K_{\tilde{X}} + \lambda \tilde{D} + \lambda a_1E_1
$$
is not log canonical at the point $Q$ and so is the pair
$$
K_{\tilde{X}} + \lambda \tilde{D} + E_1 \text{ , since } \lambda
a_1 \leq 1 \text{ .}
$$
By adjunction $(E_1, \lambda \tilde{D}|_{E_1})$ is not log
canonical at $Q$ and
$$2a_1 - a_3 = \tilde{D} \cdot E_1 \geq \text{mult}_Q\Bigl(\tilde{D}|_{E_1} \Bigr) =
 \text{mult}_Q\Bigl(\tilde{D} \cdot E_1 \Bigr)  > \frac{1}{\lambda} > \frac{5}{2} \text{
,}$$ which is false.

\item If $Q\in E_1 \cap E_3$ then the log pair
$$
K_{\tilde{X}} + \lambda \tilde{D} + \lambda a_1E_1 + \lambda
a_3E_3
$$
is not log canonical at the point $Q$ and so is the log pair
$$
K_{\tilde{X}} + \lambda \tilde{D} + \lambda a_1 E_1 + E_3 \text{,
since } \lambda a_3 = 1 \text{ .}
$$
By adjunction $(E_3, \lambda \tilde{D}|_{E_3})$ is not log
canonical at $Q$ and
$$ 2a_3 - a_1 - a_2 - a_4 = \tilde{D} \cdot E_3 \geq
 \text{mult}_Q\Bigl(\tilde{D} \cdot E_3 \Bigr)  > \frac{1}{\lambda} > \frac{5}{2} - a_1
\text{ ,}$$
 implies that
$$
\frac{7}{4} \geq \frac{7}{10} a_3 \geq 2a_3 - \frac{a_3}{2} -
\frac{4}{5} a_3 \geq 2a_3 - a_2 - a_4 > \frac{5}{2}
$$
which is a contradiction.

\item If $Q\in E_3$ but $Q\not \in E_1 \cup E_2 \cup E_4$ then
$$
K_{\tilde{X}} + \lambda \tilde{D} + \lambda a_3E_3
$$
is not log canonical at the point $Q$ and so is the pair
$$
K_{\tilde{X}} + \lambda \tilde{D} +  E_3 \text{, since } \lambda
a_3 \leq 1 \text{ .}
$$
By adjunction $(E_3, \lambda \tilde{D}|_{E_3})$ is not log
canonical at $Q$ and
$$ 2a_3 - a_1 - a_2 - a_4 = \tilde{D} \cdot E_3 \geq
 \text{mult}_Q\Bigl(\tilde{D} \cdot E_3 \Bigr)  > \frac{1}{\lambda} > 2  \text{ ,}$$
 implies that
$$
 \frac{1}{2} \geq \frac{a_3}{5} \geq 2a_3 - \frac{a_3}{2}- \frac{a_3}{2} - \frac{4}{5}
a_3 \geq 2a_3  - a_2 - a_4 > \frac{5}{2}
$$
which is a contradiction.

\item If $Q\in E_3 \cap E_4$ then the log pair
$$
K_{\tilde{X}} + \lambda \tilde{D} + \lambda a_3E_3 + \lambda
a_4E_4
$$
is not log canonical at the point $Q$ and so is the log pair
$$
K_{\tilde{X}} + \lambda \tilde{D} +  E_3 + \lambda a_4E_4 \text{,
since } \lambda a_3 \leq 1 \text{ .}
$$
By adjuction
$$2a_3 - \frac{a_3}{2} - \frac{a_3}{2} -a_4\geq 2a_3 - a_1 - a_2 - a_4 = \tilde{D} \cdot
E_3 \geq
 \text{mult}_Q\Bigl(\tilde{D} \cdot E_3 \Bigr)   > \frac{5}{2} - a_4 \text{ ,}$$
which is false.

\item If the point $Q \in E_4 \backslash (E_3 \cup E_5)$ then
$$
K_{\tilde{X}} + \lambda \tilde{D} + \lambda a_4E_4
$$
is not log canonical at the point $Q$ and so is the pair
$K_{\tilde{X}} + \lambda \tilde{D} + E_4$. By adjunction $(E_4,
\lambda \tilde{D}|_{E_4})$ is not log canonical at $Q$ and
$$2a_4 - a_3 - a_5 = \tilde{D} \cdot E_4 \geq \text{mult}_Q\Bigl(\tilde{D}|_{E_4}
\Bigr) =
 \text{mult}_Q\Bigl(\tilde{D} \cdot E_4 \Bigr)  > \frac{5}{2} \text{ ,}$$
which is false.

\item If $Q\in E_4 \cap E_5$ then the log pair
$$
K_{\tilde{X}} + \lambda \tilde{D} + \lambda a_4E_4 + \lambda
a_5E_5
$$
is not log canonical at the point $Q$ and so are the log pairs
$$
K_{\tilde{X}} + \lambda \tilde{D} + E_4 + \lambda a_5E_5 \text{
and } K_{\tilde{X}} + \lambda \tilde{D} + \lambda a_4E_4 + E_5
\text{ .}
$$
By adjuction
$$
2a_4 - a_3 - a_5 = \tilde{D} \cdot E_4 \geq
\text{mult}_Q\Bigl(\tilde{D}|_{E_4} \Bigr) =
 \text{mult}_Q\Bigl(\tilde{D} \cdot E_4 \Bigr) > \frac{5}{2} - a_5
$$
and
$$
2a_5 - a_4 - a_6 = \tilde{D} \cdot E_5 \geq
\text{mult}_Q\Bigl(\tilde{D}|_{E_5} \Bigr) =
 \text{mult}_Q\Bigl(\tilde{D} \cdot E_5 \Bigr) > \frac{5}{2} - a_4 \text{ .}
$$
which is false.

\item $Q\in E_5 \backslash (E_4 \cup E_6)$ then the log pair
$$
K_{\tilde{X}} + \lambda \tilde{D} + \lambda a_5E_5
$$
is not log canonical at the point $Q$ and so is the pair
$$
K_{\tilde{X}} + \lambda \tilde{D} + E_5 \text{ .}
$$
By adjunction $(E_5, \lambda \tilde{D}|_{E_5})$ is not log
canonical at $Q$ and
$$2a_5 - a_4 - a_6 = \tilde{D} \cdot E_5 \geq \text{mult}_Q\Bigl(\tilde{D}|_{E_5}
\Bigr) =
 \text{mult}_Q\Bigl(\tilde{D} \cdot E_5 \Bigr) > \frac{5}{2} \text{ ,}$$
which is false.

\item If $Q\in E_5 \cap E_6$ then the log pair
$$
K_{\tilde{X}} + \lambda \tilde{D} + \lambda a_5E_5 + \lambda
a_6E_6
$$
is not log canonical at the point $Q$ and so are the log pairs
$$
K_{\tilde{X}} + \lambda \tilde{D} + \lambda a_5E_5 + E_6 \text{
and } K_{\tilde{X}} + \lambda \tilde{D} + E_5 + \lambda a_6 E_6
\text{ .}
$$
By adjuction
$$
2a_6 - a_5 - a_7 = \tilde{D} \cdot E_6 \geq
\text{mult}_Q\Bigl(\tilde{D}|_{E_6} \Bigr) =
 \text{mult}_Q\Bigl(\tilde{D} \cdot E_6 \Bigr) > \frac{5}{2} - a_5
$$ and
$$
2a_5 - a_4 - a_6 = \tilde{D} \cdot E_5 \geq
\text{mult}_Q\Bigl(\tilde{D}|_{E_5} \Bigr) =
 \text{mult}_Q\Bigl(\tilde{D} \cdot E_5 \Bigr) > \frac{5}{2} - a_6 \text{ .}
$$
which is false.

\item If the point $Q \in E_6 \backslash (E_5 \cup E_7)$ then
$$
K_{\tilde{X}} + \lambda \tilde{D} + \lambda a_6E_6
$$
is not log canonical at the point $Q$ and so is the pair
$$
K_{\tilde{X}} +  \tilde{D} + E_6 \text{ .}
$$
By adjunction $(E_6, \lambda \tilde{D}|_{E_6})$ is not log
canonical at $Q$ and
$$2a_6 - a_5 - a_7 = \tilde{D} \cdot E_6 \geq \text{mult}_Q\Bigl(\tilde{D}|_{E_6} \Bigr) =
 \text{mult}_Q\Bigl(\tilde{D} \cdot E_6 \Bigr)  > \frac{5}{2} \text{ ,}$$
which is false.

\item If $Q\in E_6 \cap E_7$ then the log pair
$$
K_{\tilde{X}} + \lambda \tilde{D} + \lambda a_6E_6 + \lambda
a_7E_7
$$
is not log canonical at the point $Q$ and so is the log pair
$$
K_{\tilde{X}} + \lambda \tilde{D} + \lambda a_6E_6 + E_7 \text{ .}
$$
By adjuction
$$
2a_7 - a_6 = \tilde{D} \cdot E_7 \geq
\text{mult}_Q\Bigl(\tilde{D}|_{E_7} \Bigr) =
 \text{mult}_Q\Bigl(\tilde{D} \cdot E_7 \Bigr) > \frac{5}{2} - a_6 \text{ .}
$$
which is false.

\item If the point $Q \in E_7 \backslash E_6$ then
$$
K_{\tilde{X}} + \lambda \tilde{D} + \lambda a_7E_7
$$
is not log canonical at the point $Q$ and so is the pair
$$
K_{\tilde{X}} + \lambda \tilde{D} + E_7 \text{ .}
$$
By adjunction $(E_7, \lambda \tilde{D}|_{E_7})$ is not log
canonical at $Q$ and
$$2a_7 - a_6 = \tilde{D} \cdot E_7 \geq \text{mult}_Q\Bigl(\tilde{D}|_{E_7} \Bigr) =
 \text{mult}_Q\Bigl(\tilde{D} \cdot E_7 \Bigr)  > \frac{5}{2} \text{ ,}$$
which is false.

\end{itemize}
\end{proof}

\subsection{Del Pezzo surfaces of degree 1 with exactly one
$\mathbb{D}_8$ singularity.}

In this section we will prove the following.

\begin{lemma}
\label{D8} Let $X$ be a del Pezzo surface with exactly one Du Val
singularity of type $\mathbb{D}_8$ and $K_X^2=1$. Then the global
log canonical threshold of $X$ is
$$
\mathrm{lct} (X) = \frac{1}{3} \text{ .}
$$
\end{lemma}

\begin{proof}

Suppose that $\mathrm{lct}(X)<\frac{1}{3}$,  then there exists a
$\mathbb{Q}$-divisor $D \in X$ and a rational number $\lambda <
\frac{1}{3}$, such that the log pair $(X,\lambda D)$ is  not log
canonical and $D \sim_{\mathbb{Q}} -K_X$. We derive that the pair
$(X,\lambda D)$ is log canonical outside of a point $P\in X$ and
not log canonical at $P$. Let $\pi: \tilde{X} \to X$ be the
minimal resolution of $X$. The configuration of the exceptional
curves is given by the following Dynkin diagram.
\bigskip

$\mathbb{D}_8$.
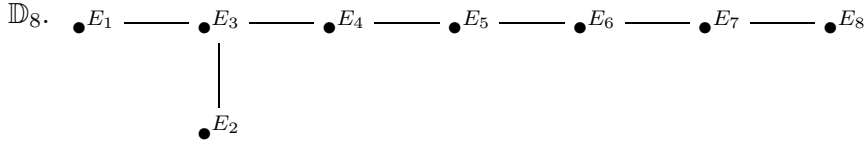
\xymatrix{ {\bullet}^{E_1} \ar@{-}[r] & {\bullet}^{E_3}
\ar@{-}[r] \ar@{-}[d] & {\bullet}^{E_4} \ar@{-}[r] & {\bullet}^{E_5}
\ar@{-}[r] & {\bullet}^{E_6} \ar@{-}[r] & {\bullet}^{E_7} \ar@{-}[r] & {\bullet}^{E_8}\\ & {\bullet}^{E_2} &}
\bigskip

Then
\begin{eqnarray*}
\tilde{D} & \sim_{\mathbb{Q}} &
\pi^*(D)-a_1E_1-a_2E_2-a_3E_3-a_4E_4-a_5E_5-a_6E_6-a_7E_7-a_8E_8
\text{  and }\\
\tilde{Z} & \sim_{\mathbb{Q}} & \pi^*(Z)-E_1-E_2-2E_3-2E_4-2E_5-2E_6-2E_7-a_8E_8\text{ .}\\
\end{eqnarray*}

From the inequalities
\begin{eqnarray*}
0 \leq \tilde{D} \cdot \tilde{Z} & = & 1 - a_7\\
0 \leq E_1 \cdot \tilde{D} & = & 2a_1 - a_3\\
0 \leq E_2 \cdot \tilde{D} & = & 2a_2 - a_3\\
0 \leq E_3 \cdot \tilde{D} & = & 2a_3 - a_1 - a_2 - a_4\\
0 \leq E_4 \cdot \tilde{D} & = & 2a_4 - a_3 - a_5\\
0 \leq E_5 \cdot \tilde{D} & = & 2a_5 - a_4 - a_6\\
0 \leq E_6 \cdot \tilde{D} & = & 2a_6 - a_5 - a_7\\
0 \leq E_7 \cdot \tilde{D} & = & 2a_7 - a_6 - a_8\\
0 \leq E_8 \cdot \tilde{D} & = & 2a_8 - a_7\\
\end{eqnarray*}
we see that $a_1 \leq \frac{8}{4} \text{, } a_2 \leq \frac{8}{4}
\text{, } a_3 \leq 3 \text{, } a_4 \leq \frac{5}{2}  \text{, } a_5
\leq 2   \text{, } a_6 \leq \frac{3}{2} \text{, } a_7 \leq 1
\text{, } a_8 \leq 1 \text{ .}$ Moreover we get the inequalities
$$
2a_1 \geq a_3 \text{, } 2a_2 \geq a_3 \geq a_4 \geq a_5 \geq a_6
\geq a_7 \geq a_8
$$
and
$$
2a_8 \geq a_7 \text{, } \frac{3}{2} a_7 \geq a_6 \text{, }
\frac{4}{3} a_6 \geq a_5 \text{, } \frac{5}{4} a_5 \geq a_4
\text{, } \frac{6}{5} a_4 \geq a_3 \text{ .}
$$
The equivalence
$$
K_{\tilde{X}} + \lambda \tilde{D} + \lambda a_1E_1 + \lambda
a_2E_2 + \lambda a_3E_3 +  \lambda a_4E_4 + \lambda a_5E_5 +
\lambda a_6E_6 + \lambda a_7E_7 + \lambda a_8E_8 \sim_{\mathbb{Q}}
\pi^*(K_X+ \lambda D)
$$
implies that there is a point $Q\in E_1\cup E_2\cup E_3\cup E_4
\cup E_5 \cup E_6 \cup E_7 \cup E_8$ such that the pair
$$
K_{\tilde{X}} + \lambda \tilde{D} + \lambda a_1E_1 + \lambda
a_2E_2 + \lambda a_3E_3 +  \lambda a_4E_4 +  \lambda a_5E_5 +
\lambda a_6E_6 + \lambda a_7E_7 + \lambda a_8E_8
$$
is not log canonical at $Q$.

\begin{itemize}
\item If the point $Q \in E_1 \backslash E_3$ then
$$
K_{\tilde{X}} + \lambda \tilde{D} + \lambda a_1E_1
$$
is not log canonical at the point $Q$ and so is the pair
$$
K_{\tilde{X}} + \lambda \tilde{D} + E_1 \text{ , since } \lambda
a_1 \leq 1 \text{ .}
$$
By adjunction $(E_1, \lambda \tilde{D}|_{E_1})$ is not log
canonical at $Q$ and
$$2a_1 - a_3 = \tilde{D} \cdot E_1 \geq \text{mult}_Q\Bigl(\tilde{D}|_{E_1} \Bigr) =
 \text{mult}_Q\Bigl(\tilde{D} \cdot E_1 \Bigr)  > 3 \text{ ,}$$
which is false.

\item If $Q\in E_1 \cap E_3$ then the log pair
$$
K_{\tilde{X}} + \lambda \tilde{D} + \lambda a_1E_1 +  \lambda
a_3E_3
$$
is not log canonical at the point $Q$ and so is the log pair
$$
K_{\tilde{X}} + \lambda \tilde{D} + \lambda a_1E_1 + E_3 \text{,
since } \lambda a_3 \leq 1 \text{ .}
$$
By adjunction $(E_3, \lambda \tilde{D}|_{E_3})$ is not log
canonical at $Q$ and
$$2a_3 - a_1 - a_2 - a_4 = \tilde{D} \cdot E_3 \geq \text{mult}_Q\Bigl(\tilde{D}|_{E_3}
\Bigr) =
 \text{mult}_Q\Bigl(\tilde{D} \cdot E_3 \Bigr)  > 3 - a_1 \text{ ,}$$
implies that
$$
2 \geq \frac{2}{3}a_3 \geq 2a_3 - \frac{a_3}{2} - \frac{5}{6}a_3
\geq 2a_3 - a_2 - a_4 > 3
$$
which is false.

\item If $Q\in E_3$ but $Q\not \in E_1 \cup E_2 \cup E_4$ then
$$
K_{\tilde{X}} + \lambda \tilde{D} +\lambda a_3E_3
$$
is not log canonical at the point $Q$ and so is the pair
$$
K_{\tilde{X}} + \lambda \tilde{D} +  E_3 \text{  , since  }
\lambda a_3 \geq 1 \text{ .}
$$
By adjunction $(E_3, \lambda \tilde{D}|_{E_3})$ is not log
canonical at $Q$ and
$$\frac{1}{2} \geq  \frac{a_3}{6} \geq 2a_3 - \frac{a_3}{2}- \frac{a_3}{2} -
\frac{5}{6}a_3 \geq 2a_3 - a_1 - a_2 - a_4 = \tilde{D} \cdot E_3
\geq
 \text{mult}_Q\Bigl(\tilde{D} \cdot E_3 \Bigr)  > 3  \text{ ,}$$
which is false.

\item If $Q\in E_3 \cap E_4$ then the log pair
$$
K_{\tilde{X}} + \lambda \tilde{D} + \lambda  a_3E_3 + \lambda
a_4E_4
$$
is not log canonical at the point $Q$ and so is the log pair
$$
K_{\tilde{X}} + \lambda  \tilde{D} + E_3 + \lambda a_4 E_4 \text{
.}
$$
By adjunction $(E_3, \lambda \tilde{D}|_{E_3})$ is not log
canonical at $Q$ and
$$ a_3 - a_4 \geq 2a_3 - \frac{a_3}{2}- \frac{a_3}{2} - a_4 \geq 2a_3 - a_1 - a_2 - a_4 =
\tilde{D} \cdot E_3 \geq
 \text{mult}_Q\Bigl(\tilde{D} \cdot E_3 \Bigr)  > \frac{1}{\lambda} - a_4 >3 -a_4  \text{
,}$$ which is false.

\item If the point $Q \in E_4 \backslash (E_3 \cup E_5)$ then
$$
K_{\tilde{X}} + \lambda \tilde{D} + \lambda a_4E_4
$$
is not log canonical at the point $Q$ and so is the pair
$$
K_{\tilde{X}} + \lambda \tilde{D} + E_4 \text{, since } \lambda
a_4 < 1 \text{ .}
$$
By adjunction $(E_4, \lambda \tilde{D}|_{E_4})$ is not log
canonical at $Q$ and
$$  \frac{1}{4} \geq \frac{1}{5}a_4 \geq 2a_4 - a_4 - \frac{4}{5}a_4 \geq 2a_4 - a_3 -
a_5 = \tilde{D} \cdot E_4 \geq
 \text{mult}_Q\Bigl(\tilde{D} \cdot E_4 \Bigr)  > \frac{1}{\lambda} > 3 \text{ ,}$$
which is false.

\item If $Q\in E_4 \cap E_5$ then the log pair
$$
K_{\tilde{X}} + \lambda \tilde{D} + \lambda a_4E_4 + \lambda
a_5E_5
$$
is not log canonical at the point $Q$ and so are the log pairs
$$
K_{\tilde{X}} + \lambda \tilde{D} + E_4 + \lambda a_5E_5 \text{
and } K_{\tilde{X}} + \lambda \tilde{D} + \lambda a_4E_4 + E_5
\text{ .}
$$
By adjuction
$$
 a_4 - a_5 \geq 2a_4 - a_3 -a_5 = \tilde{D} \cdot E_4 \geq
 \text{mult}_Q\Bigl(\tilde{D} \cdot E_4 \Bigr) > \frac{1}{\lambda} - a_5 > 3 - a_5 \text{
.}
$$
 which is a contradiction.

\item $Q\in E_5 \backslash (E_4 \cup E_6)$ then the log pair
$$
K_{\tilde{X}} + \lambda \tilde{D} + \lambda a_5E_5
$$
is not log canonical at the point $Q$ and so is the pair
$$
K_{\tilde{X}} + \lambda \tilde{D} + E_5 \text{ .}
$$
By adjunction $(E_5, \lambda \tilde{D}|_{E_5})$ is not log
canonical at $Q$ and
$$2a_5 - a_4 - a_6 = \tilde{D} \cdot E_5 \geq \text{mult}_Q\Bigl(\tilde{D}|_{E_5}
\Bigr) =
 \text{mult}_Q\Bigl(\tilde{D} \cdot E_5 \Bigr) > 3 \text{ ,}$$
which is false.

\item If $Q\in E_5 \cap E_6$ then the log pair
$$
K_{\tilde{X}} + \lambda \tilde{D} + \lambda a_5E_5 + \lambda
a_6E_6
$$
is not log canonical at the point $Q$ and so are the log pairs
$$
K_{\tilde{X}} + \lambda \tilde{D} + \lambda a_5E_5 + E_6 \text{
and } K_{\tilde{X}} + \lambda \tilde{D} + E_5 + \lambda a_6E_6
\text{ .}
$$
By adjuction
$$
2a_6 - a_5 - a_7 = \tilde{D} \cdot E_6 \geq
\text{mult}_Q\Bigl(\tilde{D}|_{E_6} \Bigr) =
 \text{mult}_Q\Bigl(\tilde{D} \cdot E_6 \Bigr) > 3 - a_5
$$ and
$$
2a_5 - a_4 - a_6 = \tilde{D} \cdot E_5 \geq
\text{mult}_Q\Bigl(\tilde{D}|_{E_5} \Bigr) =
 \text{mult}_Q\Bigl(\tilde{D} \cdot E_5 \Bigr) > 3 - a_6 \text{ .}
$$
which is false.

\item If the point $Q \in E_6 \backslash (E_5 \cup E_7)$ then
$$
K_{\tilde{X}} + \lambda \tilde{D} + \lambda a_6E_6
$$
is not log canonical at the point $Q$ and so is the pair
$$
K_{\tilde{X}} + \lambda \tilde{D} + E_6 \text{ .}
$$
By adjunction $(E_6, \lambda \tilde{D}|_{E_6})$ is not log
canonical at $Q$ and
$$2a_6-a_5-a_7 = \tilde{D} \cdot E_6 \geq \text{mult}_Q\Bigl(\tilde{D}|_{E_6} \Bigr) =
 \text{mult}_Q\Bigl(\tilde{D} \cdot E_6 \Bigr)  > 3 \text{ ,}$$
which is false.

\item If $Q\in E_6 \cap E_7$ then the log pair
$$
K_{\tilde{X}} + \lambda \tilde{D} + \lambda a_6E_6 + \lambda
a_7E_7
$$
is not log canonical at the point $Q$ and so is the log pair
$$
K_{\tilde{X}} + \lambda \tilde{D} + \lambda a_6E_6 + E_7 \text{ .}
$$
By adjuction
$$
2a_7 - a_6 - a_8 = \tilde{D} \cdot E_7 \geq
\text{mult}_Q\Bigl(\tilde{D}|_{E_7} \Bigr) =
 \text{mult}_Q\Bigl(\tilde{D} \cdot E_7 \Bigr) > 3 - a_6 \text{ .}
$$
which is false.

\item If the point $Q \in E_7 \backslash (E_6 \cup E_8)$ then
$$
K_{\tilde{X}} + \lambda \tilde{D} + \lambda a_7E_7
$$
is not log canonical at the point $Q$ and so is the pair
$$
K_{\tilde{X}} + \lambda \tilde{D} + E_7 \text{ .}
$$
By adjunction $(E_7, \lambda \tilde{D}|_{E_7})$ is not log
canonical at $Q$ and
$$2a_7 - a_6 - a_8 = \tilde{D} \cdot E_7 \geq \text{mult}_Q\Bigl(\tilde{D}|_{E_7} \Bigr) =
 \text{mult}_Q\Bigl(\tilde{D} \cdot E_7 \Bigr)  > 3 \text{ ,}$$
which is false.

\item If $Q\in E_7 \cap E_8$ then the log pair
$$
K_{\tilde{X}} + \lambda \tilde{D} + \lambda a_7E_7 + \lambda
a_8E_8
$$
is not log canonical at the point $Q$ and so is the log pair
$$
K_{\tilde{X}} + \lambda \tilde{D} + \lambda a_7E_7 + E_8 \text{ .}
$$
By adjuction
$$
2a_8 - a_7 = \tilde{D} \cdot E_8 \geq
\text{mult}_Q\Bigl(\tilde{D}|_{E_8} \Bigr) =
 \text{mult}_Q\Bigl(\tilde{D} \cdot E_8 \Bigr) > 3 - a_7 \text{ .}
$$
which is false.

\item If the point $Q \in E_8 \backslash E_7$ then
$$
K_{\tilde{X}} + \lambda \tilde{D} + \lambda a_8E_8
$$
is not log canonical at the point $Q$ and so is the pair
$$
K_{\tilde{X}} + \lambda \tilde{D} + E_8 \text{ .}
$$
By adjunction $(E_8, \lambda \tilde{D}|_{E_8})$ is not log
canonical at $Q$ and
$$2a_8 - a_7 = \tilde{D} \cdot E_8 \geq \text{mult}_Q\Bigl(\tilde{D}|_{E_8} \Bigr) =
 \text{mult}_Q\Bigl(\tilde{D} \cdot E_8 \Bigr)  > 3 \text{ ,}$$
which is false.

\end{itemize}
\end{proof}

\subsection{Del Pezzo surface of degree 1 with exactly an
$\mathbb{E}_6$ type singular point}

In this section we will prove the following.

\begin{lemma}
\label{E6} Let $X$ be a del Pezzo surface with exactly one Du Val
singularity of type $\mathbb{E}_6$ and $K_X^2=1$. Then the global
log canonical threshold of $X$ is
$$
\mathrm{lct} (X) = \frac{1}{3} \text{ .}
$$
\end{lemma}

\begin{proof}

Suppose that $\mathrm{lct}(X)<\frac{1}{3}$,  then there exists a
$\mathbb{Q}$-divisor $D \in X$ and a rational number $\lambda <
\frac{1}{3}$, such that the log pair $(X,\lambda D)$ is  not log
canonical and $D \sim_{\mathbb{Q}} -K_X$. We derive that the pair
$(X,\lambda D)$ is log canonical outside of a point $P\in X$ and
not log canonical at $P$. Let $\pi: \tilde{X} \to X$ be the
minimal resolution of $X$. The configuration of the exceptional
curves is given by the following Dynkin diagram.
\bigskip

$\mathbb{E}_6$. 
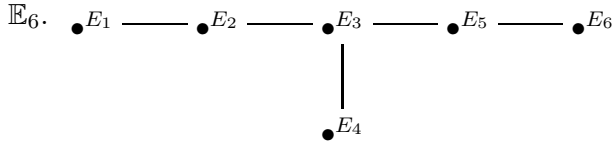
\xymatrix{ {\bullet}^{E_1} \ar@{-}[r] & {\bullet}^{E_2}
\ar@{-}[r]  & {\bullet}^{E_3} \ar@{-}[r] \ar@{-}[d] &
{\bullet}^{E_5} \ar@{-}[r] & {\bullet}^{E_6}\\
 & & {\bullet}^{E_4} & & }
\bigskip

\noindent Then
$$
\tilde{D} \sim_{\mathbb{Q}}
\pi^*(D)-a_1E_1-2a_2E_2-3a_3E_3-2a_4E_4-2a_5E_5-a_6E_6
$$
and
$$
\tilde{Z} \sim_{\mathbb{Q}}
\pi^*(Z)-E_1-2E_2-3E_3-2E_4-2E_5-E_6\text{ .}
$$
The inequalities
\begin{eqnarray*}
0 \leq \tilde{D} \cdot \tilde{Z} & = & 1 - 2a_4\\
0 \leq E_1 \cdot \tilde{D} & = & 2a_1 - 2a_2\\
0 \leq E_2 \cdot \tilde{D} & = & 4a_2 - a_1 - 3a_3\\
0 \leq E_3 \cdot \tilde{D} & = & 6a_3 - 2a_2 - 2a_4 - 2a_5\\
0 \leq E_4 \cdot \tilde{D} & = & 4a_4 - 3a_3\\
0 \leq E_5 \cdot \tilde{D} & = & 4a_5 - 3a_3 - a_6\\
0 \leq E_6 \cdot \tilde{D} & = & 2a_6 - 2a_5\\
\end{eqnarray*}
imply that $a_1=a_6 \leq \frac{4}{3} \text{, } a_2=a_5\leq
\frac{5}{6} \text{, } a_3 \leq \frac{2}{3} \text{, } a_4 \leq
\frac{1}{2} \text{ .}$ The equivalence
$$
K_{\tilde{X}} + \lambda \tilde{D} +
 \lambda a_1E_1 + 2 \lambda a_2E_2 + 3 \lambda a_3E_3 + 2 \lambda a_4E_4 +
2 \lambda a_5E_5 + \lambda a_6E_6 \sim_{\mathbb{Q}}
\pi^*(K_X+\lambda D)
$$
implies that there is a point $Q\in E_1\cup E_2\cup E_3\cup E_4
\cup E_5 \cup E_6$ such that the pair
$$
K_{\tilde{X}} + \lambda \tilde{D} + \lambda a_1E_1 + 2 \lambda
a_2E_2 + 3 \lambda a_3E_3 + 2 \lambda a_4E_4 + 2 \lambda a_5E_5 +
\lambda a_6E_6
$$
is not log canonical at $Q$.

\begin{itemize}
\item If the point $Q \in E_1$ and $Q\not \in E_2$ then
$$
K_{\tilde{X}} + \lambda \tilde{D} + \lambda a_1E_1
$$
is not log canonical at the point $Q$ and so is the pair
$$
K_{\tilde{X}} + \lambda \tilde{D} + E_1 \text{ since }\lambda a_1
\leq 1 \text{ .}
$$
By adjunction $(E_1, \lambda \tilde{D}|_{E_1})$ is not log
canonical at $Q$ and
$$2a_1 - 2a_2 = \tilde{D} \cdot E_1 \geq \text{mult}_Q\Bigl(\tilde{D}|_{E_1} \Bigr) =
 \text{mult}_Q\Bigl(\tilde{D} \cdot E_1 \Bigr)  > 3 \text{ ,}$$ which is false.

 \item If the point $Q \in E_1 \cap E_2$ then
$$
K_{\tilde{X}} + \lambda \tilde{D} + \lambda a_1E_1 + 2 \lambda
a_2E_2
$$
is not log canonical at the point $Q$ and so are the pairs
$$
K_{\tilde{X}} + \lambda \tilde{D} + E_1 + 2 \lambda a_2E_2 \text{
and } K_{\tilde{X}} + \lambda \tilde{D} + \lambda E_1 + E_2 \text{
.}
$$
By adjunction
$$2a_1 - 2a_2 = \tilde{D} \cdot E_1 \geq \text{mult}_Q\Bigl(\tilde{D}|_{E_1} \Bigr) =
 \text{mult}_Q\Bigl(\tilde{D} \cdot E_1 \Bigr)  > 3 - 2a_2 \text{ and }$$

$$4a_2 - a_1 -3a_3 = \tilde{D} \cdot E_2 \geq \text{mult}_Q\Bigl(\tilde{D}|_{E_2} \Bigr) =
 \text{mult}_Q\Bigl(\tilde{D} \cdot E_2 \Bigr)  > 3-a_1 \text{ ,}$$
which is false.

\item If the point $Q \in E_2 \backslash (E_1 \cup E_3)$ then
$$
K_{\tilde{X}} + \lambda \tilde{D} + 2 \lambda a_2E_2
$$
is not log canonical at the point $Q$ and so is the pair
$$
K_{\tilde{X}} + \lambda \tilde{D} + E_2 \text{ .}
$$
By adjunction
$$4a_2 - a_1 -3a_3 = \tilde{D} \cdot E_2 \geq \text{mult}_Q\Bigl(\tilde{D}|_{E_2} \Bigr) =
 \text{mult}_Q\Bigl(\tilde{D} \cdot E_2 \Bigr)  > 3 \text{ ,}$$
which is false.

\item If the point $Q \in E_2 \cap E_3$ then
$$
K_{\tilde{X}} + \lambda \tilde{D} + 2\lambda a_2E_2 + 3\lambda
a_3E_3
$$
is not log canonical at the point $Q$ and so are the pairs
$$
K_{\tilde{X}} + \lambda \tilde{D} + E_2 + 3 \lambda a_3E_3 \text{
and } K_{\tilde{X}} + \lambda \tilde{D} + 2\lambda a_2E_2 + E_3
\text{ .}
$$
By adjunction
$$4a_2 - 3a_3 - a_1 = \tilde{D} \cdot E_2 \geq \text{mult}_Q\Bigl(\tilde{D}|_{E_2} \Bigr)
=
 \text{mult}_Q\Bigl(\tilde{D} \cdot E_2 \Bigr)  > 3 - 3a_3 \text{ and }$$
and
$$6a_3 - 2a_2 -2a_4-2a_5 = \tilde{D} \cdot E_3 \geq \text{mult}_Q\Bigl(\tilde{D}|_{E_3}
\Bigr) =
 \text{mult}_Q\Bigl(\tilde{D} \cdot E_3 \Bigr)  > 3-2a_2 \text{ ,}$$
which is false.

\item If the point $Q \in E_3 \backslash (E_2 \cup E_4 \cup E_5)$
then
$$
K_{\tilde{X}} + \lambda \tilde{D} + 3\lambda a_3E_3
$$
is not log canonical at the point $Q$ and so is the pair
$$
K_{\tilde{X}} + \lambda \tilde{D} + E_3 \text{ .}
$$
By adjunction $(E_3, \lambda \tilde{D}|_{E_3})$ is not log
canonical at $Q$ and
$$6a_3 - 2a_2 - 2a_4 - 2a_5 = \tilde{D} \cdot E_3 \geq
\text{mult}_Q\Bigl(\tilde{D}|_{E_3} \Bigr) =
 \text{mult}_Q\Bigl(\tilde{D} \cdot E_3 \Bigr)  > 3 \text{ ,}$$ which is false.

\item If the point $Q \in E_3 \cap E_4$ then
$$
K_{\tilde{X}} + \lambda \tilde{D} + 3\lambda a_3E_3 + 2\lambda
a_4E_4
$$
is not log canonical at the point $Q$ and so are the pairs
$$
K_{\tilde{X}} + \lambda \tilde{D} + E_3 + 2 \lambda a_4E_4 \text{
and } K_{\tilde{X}} + \lambda \tilde{D} + 3\lambda a_3E_3 + E_4
\text{ .}
$$
By adjunction
$$6a_3 - 2a_2 - 2a_4 - 2a_5 = \tilde{D} \cdot E_3 \geq
\text{mult}_Q\Bigl(\tilde{D}|_{E_3} \Bigr) =
 \text{mult}_Q\Bigl(\tilde{D} \cdot E_3 \Bigr)  > 3 - 2a_4 \text{ and }$$

$$4a_4 - 3a_3 = \tilde{D} \cdot E_4 \geq \text{mult}_Q\Bigl(\tilde{D}|_{E_4} \Bigr) =
 \text{mult}_Q\Bigl(\tilde{D} \cdot E_4 \Bigr)  > 3 - 3a_3 \text{ ,}$$
which is false.

\item If the point $Q \in E_4 \backslash E_3$ then
$$
K_{\tilde{X}} + \lambda \tilde{D} + 2\lambda a_4E_4
$$
is not log canonical at the point $Q$ and so is the pair
$$
K_{\tilde{X}} + \lambda \tilde{D} + E_4 \text{ .}
$$
By adjunction $(E_4, \lambda \tilde{D}|_{E_4})$ is not log
canonical at $Q$ and
$$4a_4 - 3a_3 = \tilde{D} \cdot E_4 \geq \text{mult}_Q\Bigl(\tilde{D}|_{E_4} \Bigr) =
 \text{mult}_Q\Bigl(\tilde{D} \cdot E_4 \Bigr)  > 3 \text{ ,}$$ which is false.
\end{itemize}
\end{proof}

\subsection{Del Pezzo surface of degree 1 with exactly one $\mathbb{E}_7$
type singularity}

In this section we will prove the following.

\begin{lemma}
\label{E7} Let $X$ be a del Pezzo surface with exactly one Du Val
singularity of type $\mathbb{E}_7$ and $K_X^2=1$. Then the global
log canonical threshold of $X$ is
$$
\mathrm{lct} (X) = \frac{1}{4} \text{ .}
$$
\end{lemma}

\begin{proof}

Suppose that $\mathrm{lct}(X)<\frac{1}{4}$,  then there exists a
$\mathbb{Q}$-divisor $D \in X$ and a rational number $\lambda <
\frac{1}{4}$, such that the log pair $(X,\lambda D)$ is  not log
canonical and $D \sim_{\mathbb{Q}} -K_X$. We derive that the pair
$(X,\lambda D)$ is log canonical outside of a point $P\in X$ and
not log canonical at $P$. Let $\pi: \tilde{X} \to X$ be the
minimal resolution of $X$. The configuration of the exceptional
curves is given by the following Dynkin diagram.
\bigskip

$\mathbb{E}_7$.
\xymatrix{ {\bullet}^{E_1} \ar@{-}[r] & {\bullet}^{E_2}
\ar@{-}[r]  & {\bullet}^{E_3} \ar@{-}[r] \ar@{-}[d] &
{\bullet}^{E_5} \ar@{-}[r] & {\bullet}^{E_6} \ar@{-}[r] & {\bullet}^{E_7}\\
 & & {\bullet}^{E_4} & & }
\bigskip

\noindent Then
$$
\tilde{D} \sim_{\mathbb{Q}} \pi^*(D) - 2a_1E_1 - 3a_2E_2 - 4a_3E_3
- 2a_4E_4 - 3a_5E_5 - 2a_6E_6 - a_7E_7
$$
and
$$
\tilde{Z} \sim_{\mathbb{Q}} \pi^*(Z) - 2E_1 - 3E_2 - 4E_3 - 2E_4 -
3E_5 - 2E_6 - E_7\text{ .}
$$
The inequalities
\begin{eqnarray*}
0 \leq \tilde{D} \cdot \tilde{Z} & = & 1 - 2a_1\\
0 \leq E_1 \cdot \tilde{D} & = & 4a_1 - 3a_2\\
0 \leq E_2 \cdot \tilde{D} & = & 6a_2 - 2a_1 - 4a_3\\
0 \leq E_3 \cdot \tilde{D} & = & 8a_3 - 3a_2 - 3a_5 - 2a_4\\
0 \leq E_4 \cdot \tilde{D} & = & 4a_4 - 4a_3\\
0 \leq E_5 \cdot \tilde{D} & = & 6a_5 - 4a_3 - 2a_6\\
0 \leq E_6 \cdot \tilde{D} & = & 4a_6 - 3a_5 - a_7\\
0 \leq E_7 \cdot \tilde{D} & = & 2a_7 - 2a_6\\
\end{eqnarray*}
imply that $a_1 \leq \frac{1}{2} \text{, } a_2 \leq \frac{2}{3}
\text{, }
 a_3 \leq \frac{3}{4} \text{, } a_4 \leq \frac{7}{8} \text{, }
 a_5 \leq \frac{5}{6} \text{, } a_6 \leq 1 \text{, }
 a_7 \leq \frac{3}{2} \text{ .}
$ The equivalence
$$
K_{\tilde{X}} + \lambda \tilde{D} + 2 \lambda a_1E_1 + 3 \lambda
a_2E_2 + 4 \lambda a_3E_3 + 2 \lambda a_4E_4 + 3 \lambda a_5E_5 +
2 \lambda a_6E_6 + \lambda a_7E_7 = \pi^*(K_X+\lambda D)
$$
implies that there is a point $Q\in E_1\cup E_2\cup E_3\cup E_4
\cup E_5 \cup E_6 \cup E_7$ such that the pair
$$
K_{\tilde{X}} + \lambda \tilde{D} + 2 \lambda a_1E_1 + 3 \lambda
a_2E_2 + 4 \lambda a_3E_3 + 2 \lambda a_4E_4 + 3 \lambda a_5E_5 +
2 \lambda a_6E_6 + \lambda a_7E_7
$$
is not log canonical at $Q$.

\begin{itemize}
\item If the point $Q \in E_1$ and $Q\not \in E_2$ then
$$
K_{\tilde{X}} + \lambda \tilde{D} + 2 \lambda a_1E_1
$$
is not log canonical at the point $Q$ and so is the pair
$$
K_{\tilde{X}} + \lambda \tilde{D} + E_1 \text{ since } 2 \lambda
a_1 \leq 1 \text{ .}
$$
By adjunction $(E_1, \lambda \tilde{D}|_{E_1})$ is not log
canonical at $Q$ and
$$4a_1 - 3a_2 = \tilde{D} \cdot E_1 \geq \text{mult}_Q\Bigl(\tilde{D}|_{E_1} \Bigr) =
 \text{mult}_Q\Bigl(\tilde{D} \cdot E_1 \Bigr)  > 4 \text{ ,}$$ which is false.

 \item If the point $Q \in E_1 \cap E_2$ then
$$
K_{\tilde{X}} + \lambda \tilde{D} + 2 \lambda a_1E_1 + 3 \lambda
a_2E_2
$$
is not log canonical at the point $Q$ and so are the pairs
$$
K_{\tilde{X}} + \lambda \tilde{D} + E_1 + 3 \lambda a_2E_2 \text{
and } K_{\tilde{X}} + \lambda \tilde{D} + 2 \lambda a_1E_1 + E_2
\text{ .}
$$
By adjunction
$$4a_1 - 3a_2 = \tilde{D} \cdot E_1 \geq \text{mult}_Q\Bigl(\tilde{D}|_{E_1} \Bigr) =
 \text{mult}_Q\Bigl(\tilde{D} \cdot E_1 \Bigr)  > 4 - 3a_2 \text{ and }$$
and
$$6a_2 - 2a_1 -4a_3 = \tilde{D} \cdot E_2 \geq \text{mult}_Q\Bigl(\tilde{D}|_{E_2} \Bigr)
=
 \text{mult}_Q\Bigl(\tilde{D} \cdot E_2 \Bigr)  > 4-2a_1 \text{ ,}$$
which is false.

\item If the point $Q \in E_2 \backslash (E_1 \cup E_3)$ then
$$
K_{\tilde{X}} + \lambda \tilde{D} + 3 \lambda a_2E_2
$$
is not log canonical at the point $Q$ and so is the pair
$$
K_{\tilde{X}} + \lambda \tilde{D} + E_2 \text{ .}
$$
By adjunction
$$6a_2-2a_1-4a_3 = \tilde{D} \cdot E_2 \geq \text{mult}_Q\Bigl(\tilde{D}|_{E_2} \Bigr) =
 \text{mult}_Q\Bigl(\tilde{D} \cdot E_2 \Bigr)  > 4 \text{ ,}$$
which is false.

\item If the point $Q \in E_2 \cap E_3$ then
$$
K_{\tilde{X}} + \lambda \tilde{D} + 3 \lambda a_2E_2 + 4 \lambda
a_3E_3
$$
is not log canonical at the point $Q$ and so are the pairs
$$
K_{\tilde{X}} + \lambda \tilde{D} + E_2 + 4 \lambda a_3E_3 \text{
and } K_{\tilde{X}} + \lambda \tilde{D} + 3 \lambda a_2E_2 + E_3
\text{ .}
$$
By adjunction
$$6a_2-2a_1-4a_3 = \tilde{D} \cdot E_2 \geq \text{mult}_Q\Bigl(\tilde{D}|_{E_2} \Bigr)
=
 \text{mult}_Q\Bigl(\tilde{D} \cdot E_2 \Bigr)  > 4 - 4a_3 \text{ and }$$
and
$$8a_3 - 3a_2 -3a_5-2a_4 = \tilde{D} \cdot E_3 \geq \text{mult}_Q\Bigl(\tilde{D}|_{E_3}
\Bigr) =
 \text{mult}_Q\Bigl(\tilde{D} \cdot E_3 \Bigr)  > 4-3a_2 \text{ ,}$$
which is false.

\item If the point $Q \in E_3 \backslash (E_2 \cup E_4 \cup E_5)$
then
$$
K_{\tilde{X}} + \lambda \tilde{D} + 4 \lambda a_3E_3
$$
is not log canonical at the point $Q$ and so is the pair
$$
K_{\tilde{X}} + \lambda \tilde{D} + E_3 \text{ .}
$$
By adjunction $(E_3, \lambda \tilde{D}|_{E_3})$ is not log
canonical at $Q$ and
$$8a_3 - 3a_2 - 2a_4 - 3a_5 = \tilde{D} \cdot E_3 \geq
\text{mult}_Q\Bigl(\tilde{D}|_{E_3} \Bigr) =
 \text{mult}_Q\Bigl(\tilde{D} \cdot E_3 \Bigr)  > 4 \text{ ,}$$ which is false.

\item If the point $Q \in E_3 \cap E_4$ then
$$
K_{\tilde{X}} + \lambda \tilde{D} + 4 \lambda a_3E_3 + 2\lambda
a_4E_4
$$
is not log canonical at the point $Q$ and so are the pairs
$$
K_{\tilde{X}} + \lambda \tilde{D} + E_3 + 2 \lambda a_4E_4 \text{
and } K_{\tilde{X}} + \lambda \tilde{D} + 4 \lambda a_3E_3 + E_4
\text{ .}
$$
By adjunction
$$8a_3 - 3a_2 - 2a_4 - 3a_5 = \tilde{D} \cdot E_3 \geq
\text{mult}_Q\Bigl(\tilde{D}|_{E_3} \Bigr) =
 \text{mult}_Q\Bigl(\tilde{D} \cdot E_3 \Bigr)  > 4 - 2a_4 \text{ and }$$
and
$$4a_4 - 4a_3 = \tilde{D} \cdot E_4 \geq \text{mult}_Q\Bigl(\tilde{D}|_{E_4} \Bigr) =
 \text{mult}_Q\Bigl(\tilde{D} \cdot E_4 \Bigr)  > 4 - 4a_3 \text{ ,}$$
which is false.

\item If the point $Q \in E_4 \backslash E_3$ then
$$
K_{\tilde{X}} + \lambda \tilde{D} + 2\lambda a_4E_4
$$
is not log canonical at the point $Q$ and so is the pair
$$
K_{\tilde{X}} + \lambda \tilde{D} + E_4 \text{ .}
$$
By adjunction $(E_4, \lambda \tilde{D}|_{E_4})$ is not log
canonical at $Q$ and
$$4a_4 - 4a_3 = \tilde{D} \cdot E_4 \geq \text{mult}_Q\Bigl(\tilde{D}|_{E_4} \Bigr) =
 \text{mult}_Q\Bigl(\tilde{D} \cdot E_4 \Bigr)  > 4 \text{ ,}$$ which is false.

 \item If the point $Q \in E_3 \cap E_5$ then
$$
K_{\tilde{X}} + \lambda \tilde{D} + 4 \lambda a_3E_3 + 3 \lambda
a_5E_5
$$
is not log canonical at the point $Q$ and so are the pairs
$$
K_{\tilde{X}} + \lambda \tilde{D} + E_3 + 3 \lambda a_5E_5 \text{
and } K_{\tilde{X}} + \lambda \tilde{D} + 4 \lambda a_3E_3 + E_5
\text{ .}
$$
By adjunction
$$8a_3 - 3a_2 - 2a_4 - 3a_5 = \tilde{D} \cdot E_3 \geq
\text{mult}_Q\Bigl(\tilde{D}|_{E_3} \Bigr) =
 \text{mult}_Q\Bigl(\tilde{D} \cdot E_3 \Bigr)  > 4 - 3a_5 \text{ and }$$
and
$$6a_5 - 4a_3 - 2a_6 = \tilde{D} \cdot E_5 \geq \text{mult}_Q\Bigl(\tilde{D}|_{E_5}
\Bigr) =
 \text{mult}_Q\Bigl(\tilde{D} \cdot E_5 \Bigr)  > 4 - 4a_3 \text{ ,}$$
which is false.

\item If the point $Q \in E_5 \backslash (E_3 \cup E_6)$ then
$$
K_{\tilde{X}} + \lambda \tilde{D} + 3\lambda a_5E_5
$$
is not log canonical at the point $Q$ and so is the pair
$$
K_{\tilde{X}} + \lambda \tilde{D} + E_5 \text{ .}
$$
By adjunction $(E_5, \lambda \tilde{D}|_{E_5})$ is not log
canonical at $Q$ and
$$6a_5 - 4a_3 - 2a_6= \tilde{D} \cdot E_5 \geq \text{mult}_Q\Bigl(\tilde{D}|_{E_5} \Bigr)
=
 \text{mult}_Q\Bigl(\tilde{D} \cdot E_5 \Bigr)  > 4 \text{ ,}$$ which is false.

 \item If the point $Q \in E_5 \cap E_6$ then
$$
K_{\tilde{X}} + \lambda \tilde{D} + 3 \lambda a_5E_5 + 2 \lambda
a_6E_6
$$
is not log canonical at the point $Q$ and so are the pairs
$$
K_{\tilde{X}} + \lambda \tilde{D} + E_5 + 2 \lambda a_6E_6 \text{
and } K_{\tilde{X}} + \lambda \tilde{D} + 3 \lambda a_5E_5 + E_6
\text{ .}
$$
By adjunction
$$6a_5 - 4a_3 - 2a_6 = \tilde{D} \cdot E_5 \geq
\text{mult}_Q\Bigl(\tilde{D}|_{E_5} \Bigr) =
 \text{mult}_Q\Bigl(\tilde{D} \cdot E_5 \Bigr)  > 4 - 2a_6 \text{ and }$$
and
$$4a_6 - 3a_5 - a_7 = \tilde{D} \cdot E_6 \geq \text{mult}_Q\Bigl(\tilde{D}|_{E_6} \Bigr)
=
 \text{mult}_Q\Bigl(\tilde{D} \cdot E_6 \Bigr)  > 4 - 3a_5 \text{ ,}$$
which is false.

\item If the point $Q \in E_6 \backslash (E_5 \cup E_7)$ then
$$
K_{\tilde{X}} + \lambda \tilde{D} + 2\lambda a_6E_6
$$
is not log canonical at the point $Q$ and so is the pair
$$
K_{\tilde{X}} + \lambda \tilde{D} + E_6 \text{ .}
$$
By adjunction $(E_6, \lambda \tilde{D}|_{E_6})$ is not log
canonical at $Q$ and
$$4a_6 - 3a_5 - a_7= \tilde{D} \cdot E_6 \geq \text{mult}_Q\Bigl(\tilde{D}|_{E_6} \Bigr) =
 \text{mult}_Q\Bigl(\tilde{D} \cdot E_6 \Bigr)  > 4 \text{ ,}$$ which is false.

 \item If the point $Q \in E_6 \cap E_7$ then
$$
K_{\tilde{X}} + \lambda \tilde{D} + 2 \lambda a_6E_6 +  \lambda
a_7E_7
$$
is not log canonical at the point $Q$ and so are the pairs
$$
K_{\tilde{X}} + \lambda \tilde{D} + E_6 + \lambda a_7E_7 \text{
and } K_{\tilde{X}} + \lambda \tilde{D} + 2 \lambda a_6E_6 + E_7
\text{ .}
$$
By adjunction
$$4a_6 - 3a_5 - a_7 = \tilde{D} \cdot E_6 \geq
\text{mult}_Q\Bigl(\tilde{D}|_{E_6} \Bigr) =
 \text{mult}_Q\Bigl(\tilde{D} \cdot E_6 \Bigr)  > 4 - a_7 \text{ and }$$
and
$$2a_7 - 2a_6 = \tilde{D} \cdot E_7 \geq \text{mult}_Q\Bigl(\tilde{D}|_{E_7} \Bigr) =
 \text{mult}_Q\Bigl(\tilde{D} \cdot E_7 \Bigr)  > 4 - 2a_6 \text{ ,}$$
which is false.

\item If the point $Q \in E_7 \backslash E_6$ then
$$
K_{\tilde{X}} + \lambda \tilde{D} + \lambda a_7E_7
$$
is not log canonical at the point $Q$ and so is the pair
$$
K_{\tilde{X}} + \lambda \tilde{D} + E_7 \text{ .}
$$
By adjunction $(E_7, \lambda \tilde{D}|_{E_7})$ is not log
canonical at $Q$ and
$$2a_7 - 2a_6= \tilde{D} \cdot E_7 \geq \text{mult}_Q\Bigl(\tilde{D}|_{E_7} \Bigr) =
 \text{mult}_Q\Bigl(\tilde{D} \cdot E_7 \Bigr)  > 4 \text{ ,}$$ which is false.

\end{itemize}
\end{proof}

\subsection{Del Pezzo surface of degree 1 with exactly one $\mathbb{E}_8$
type singular point}

In this section we will prove the following.

\begin{lemma}
\label{E8} Let $X$ be a del Pezzo surface with exactly one Du Val
singularity of type $\mathbb{E}_8$ and $K_X^2=1$. Then the global
log canonical threshold of $X$ is
$$
\mathrm{lct} (X) = \frac{1}{6} \text{ .}
$$
\end{lemma}

\begin{proof}

Suppose that $\mathrm{lct}(X)<\frac{1}{6}$ , then there exists a
$\mathbb{Q}$-divisor $D \in X$, such that the log pair $(X,\lambda
D)$ is  not log canonical for a rational number $\lambda <
\frac{1}{6}$ and $D \sim_{\mathbb{Q}} -K_X$. We derive that the
pair $(X,\lambda D)$ is log canonical outside of a point $P\in X$
and not log canonical at $P$. Let $\pi: \tilde{X} \to X$ be the
minimal resolution of $X$. The configuration of the exceptional
curves is given by the following Dynkin diagram.
\bigskip

$\mathbb{E}_8$.
\xymatrix{ {\bullet}^{E_1} \ar@{-}[r] & {\bullet}^{E_2}
\ar@{-}[r]  & {\bullet}^{E_3} \ar@{-}[r] \ar@{-}[d] &
{\bullet}^{E_5} \ar@{-}[r] & {\bullet}^{E_6} \ar@{-}[r] & {\bullet}^{E_7} \ar@{-}[r] & {\bullet}^{E_8}\\
 & & {\bullet}^{E_4} & & }
\bigskip

\noindent Then
$$
\tilde{D} \sim_{\mathbb{Q}} \pi^*(D) - 2 a_1E_1 - 4 a_2E_2 - 6
a_3E_3 - 3 a_4E_4 - 5 a_5E_5 - 4 a_6E_6 - 3 a_7E_7 - 2 a_8E_8
$$
and
$$
\tilde{Z} \sim_{\mathbb{Q}} \pi^*(Z)- 2 E_1 - 4 E_2 - 6 E_3 - 3
E_4 - 5 E_5 - 4 E_6 - 3 E_7 - 2 E_8\text{ .}
$$
We have the inequalities
\begin{eqnarray*}
0 \leq \tilde{D} \cdot \tilde{Z} & = & 1 - 2a_8\\
0 \leq E_1 \cdot \tilde{D} & = & 4a_1 - 4a_2\\
0 \leq E_2 \cdot \tilde{D} & = & 8a_2 - 2a_1 - 6a_3\\
0 \leq E_3 \cdot \tilde{D} & = & 12a_3 - 4a_2 - 5a_5 - 3a_4\\
0 \leq E_4 \cdot \tilde{D} & = & 6a_4 - 6a_3\\
0 \leq E_5 \cdot \tilde{D} & = & 10a_5 - 6a_3 - 4a_6\\
0 \leq E_6 \cdot \tilde{D} & = & 8a_6 - 5a_5 - 3a_7\\
0 \leq E_7 \cdot \tilde{D} & = & 6a_7 - 4a_6 - 2a_8\\
0 \leq E_8 \cdot \tilde{D} & = & 4a_8 - 3a_7 \text{ .}\\
\end{eqnarray*}
The equivalence
$$
K_{\tilde{X}} + \lambda \tilde{D} + 2 \lambda a_1E_1 + 4 \lambda
a_2E_2 + 6 \lambda a_3E_3 + 3 \lambda a_4E_4 + 5 \lambda a_5E_5 +
4 \lambda a_6E_6 + 3 \lambda a_7E_7 + 2 \lambda a_8E_8
\sim_{\mathbb{Q}} \pi^*(K_X+ \lambda D)
$$
implies that there is a point $Q\in E_1\cup E_2\cup E_3\cup E_4
\cup E_5 \cup E_6 \cup E_7 \cup E_8$ such that the pair
$$
K_{\tilde{X}} + \lambda \tilde{D} + 2 \lambda a_1E_1 + 4 \lambda
a_2E_2 + 6 \lambda a_3E_3 + 3 \lambda a_4E_4 + 5 \lambda a_5E_5 +
4 \lambda a_6E_6 + 3 \lambda a_7E_7 + 2 \lambda a_8E_8
$$
is not log canonical at $Q$.

\begin{itemize}
\item If the point $Q \in E_1$ and $Q\not \in E_2$ then
$$
K_{\tilde{X}} + \lambda \tilde{D} + 2 \lambda a_1E_1
$$
is not log canonical at the point $Q$ and so is the pair
$$
K_{\tilde{X}} + \lambda \tilde{D} + E_1 \text{ since } 2 \lambda
a_1 \leq 1 \text{ .}
$$
By adjunction $(E_1, \lambda \tilde{D}|_{E_1})$ is not log
canonical at $Q$ and
$$4a_1 - 4a_2 = \tilde{D} \cdot E_1 \geq \text{mult}_Q\Bigl(\tilde{D}|_{E_1} \Bigr) =
 \text{mult}_Q\Bigl(\tilde{D} \cdot E_1 \Bigr)  > 6 \text{ ,}$$ which is false.

 \item If the point $Q \in E_1 \cap E_2$ then
$$
K_{\tilde{X}} + \lambda \tilde{D} + 2 \lambda a_1E_1 + 4 \lambda
a_2E_2
$$
is not log canonical at the point $Q$ and so are the pairs
$$
K_{\tilde{X}} + \lambda \tilde{D} + E_1 + 4 \lambda a_2E_2 \text{
and } K_{\tilde{X}} + \lambda \tilde{D} + 2 \lambda a_1E_1 + E_2
\text{ .}
$$
By adjunction
$$4a_1 - 4a_2 = \tilde{D} \cdot E_1 \geq \text{mult}_Q\Bigl(\tilde{D}|_{E_1} \Bigr) =
 \text{mult}_Q\Bigl(\tilde{D} \cdot E_1 \Bigr)  > 6 - 4a_2 \text{ and }$$
and
$$8a_2 - 2a_1 -6a_3 = \tilde{D} \cdot E_2 \geq \text{mult}_Q\Bigl(\tilde{D}|_{E_2} \Bigr)
=
 \text{mult}_Q\Bigl(\tilde{D} \cdot E_2 \Bigr)  > 6 - 2a_1 \text{ ,}$$
which is false.

\item If the point $Q \in E_2 \backslash (E_1 \cup E_3)$ then
$$
K_{\tilde{X}} + \lambda \tilde{D} + 4 \lambda a_2E_2
$$
is not log canonical at the point $Q$ and so is the pair
$$
K_{\tilde{X}} + \lambda \tilde{D} + E_2 \text{ .}
$$
By adjunction
$$8a_2-2a_1-6a_3 = \tilde{D} \cdot E_2 \geq \text{mult}_Q\Bigl(\tilde{D}|_{E_2} \Bigr) =
 \text{mult}_Q\Bigl(\tilde{D} \cdot E_2 \Bigr)  > 6 \text{ ,}$$
which is false.

\item If the point $Q \in E_2 \cap E_3$ then
$$
K_{\tilde{X}} + \lambda \tilde{D} + 4 \lambda a_2E_2 + 6 \lambda
a_3E_3
$$
is not log canonical at the point $Q$ and so are the pairs
$$
K_{\tilde{X}} + \lambda \tilde{D} + E_2 + 6 \lambda a_3E_3 \text{
and } K_{\tilde{X}} + \lambda \tilde{D} + 4 \lambda a_2E_2 + E_3
\text{ .}
$$
By adjunction
$$8a_2-2a_1-6a_3 = \tilde{D} \cdot E_2 \geq \text{mult}_Q\Bigl(\tilde{D}|_{E_2} \Bigr)
=
 \text{mult}_Q\Bigl(\tilde{D} \cdot E_2 \Bigr)  > 6 - 6a_3 \text{ and }$$
and
$$12a_3 - 4a_2 -5a_5-3a_4 = \tilde{D} \cdot E_3 \geq \text{mult}_Q\Bigl(\tilde{D}|_{E_3}
\Bigr) =
 \text{mult}_Q\Bigl(\tilde{D} \cdot E_3 \Bigr)  > 6 - 4a_2 \text{ ,}$$
which is false.

\item If the point $Q \in E_3 \backslash (E_2 \cup E_4 \cup E_5)$
then
$$
K_{\tilde{X}} + \lambda \tilde{D} + 6 \lambda a_3E_3
$$
is not log canonical at the point $Q$ and so is the pair
$$
K_{\tilde{X}} + \lambda \tilde{D} + E_3 \text{ .}
$$
By adjunction $(E_3, \lambda \tilde{D}|_{E_3})$ is not log
canonical at $Q$ and
$$12a_3 - 4a_2 - 3a_4 - 5a_5 = \tilde{D} \cdot E_3 \geq
\text{mult}_Q\Bigl(\tilde{D}|_{E_3} \Bigr) =
 \text{mult}_Q\Bigl(\tilde{D} \cdot E_3 \Bigr)  > 6 \text{ ,}$$ which is false.

\item If the point $Q \in E_3 \cap E_4$ then
$$
K_{\tilde{X}} + \lambda \tilde{D} + 6 \lambda a_3E_3 + 3 \lambda
a_4E_4
$$
is not log canonical at the point $Q$ and so are the pairs
$$
K_{\tilde{X}} + \lambda\tilde{D} + E_3 + 3 \lambda a_4E_4 \text{
and } K_{\tilde{X}} + \lambda \tilde{D} + 6 \lambda a_3E_3 + E_4
\text{ .}
$$
By adjunction
$$12a_3 - 4a_2 - 3a_4 - 5a_5 = \tilde{D} \cdot E_3 \geq
\text{mult}_Q\Bigl(\tilde{D}|_{E_3} \Bigr) =
 \text{mult}_Q\Bigl(\tilde{D} \cdot E_3 \Bigr)  > 6 - 3a_4 \text{ and }$$
and
$$6a_4 - 6a_3 = \tilde{D} \cdot E_4 \geq \text{mult}_Q\Bigl(\tilde{D}|_{E_4} \Bigr) =
 \text{mult}_Q\Bigl(\tilde{D} \cdot E_4 \Bigr)  > 6 - 6a_3 \text{ ,}$$
which is false.

\item If the point $Q \in E_4 \backslash E_3$ then
$$
K_{\tilde{X}} + \lambda \tilde{D} + 3 \lambda a_4E_4
$$
is not log canonical at the point $Q$ and so is the pair
$$
K_{\tilde{X}} + \lambda \tilde{D} + E_4 \text{ .}
$$
By adjunction $(E_4, \lambda \tilde{D}|_{E_4})$ is not log
canonical at $Q$ and
$$6a_4 - 6a_3 = \tilde{D} \cdot E_4 \geq \text{mult}_Q\Bigl(\tilde{D}|_{E_4} \Bigr) =
 \text{mult}_Q\Bigl(\tilde{D} \cdot E_4 \Bigr)  > 6 \text{ ,}$$ which is
false.

 \item If the point $Q \in E_3 \cap E_5$ then
$$
K_{\tilde{X}} + \lambda \tilde{D} + 6 \lambda a_3E_3 + 5 \lambda
a_5E_5
$$
is not log canonical at the point $Q$ and so are the pairs
$$
K_{\tilde{X}} + \lambda \tilde{D} + E_3 + 5 \lambda a_5E_5 \text{
and } K_{\tilde{X}} + \lambda \tilde{D} + 6 \lambda a_3E_3 + E_5
\text{ .}
$$
By adjunction
$$12a_3 - 4a_2 - 3a_4 - 5a_5 = \tilde{D} \cdot E_3 \geq
\text{mult}_Q\Bigl(\tilde{D}|_{E_3} \Bigr) =
 \text{mult}_Q\Bigl(\tilde{D} \cdot E_3 \Bigr)  > 6 - 5a_5 \text{ and }$$
and
$$10a_5 - 6a_3 - 4a_6 = \tilde{D} \cdot E_5 \geq \text{mult}_Q\Bigl(\tilde{D}|_{E_5}
\Bigr) =
 \text{mult}_Q\Bigl(\tilde{D} \cdot E_5 \Bigr)  > 6 - 6a_3 \text{ ,}$$
which is false.

\item If the point $Q \in E_5 \backslash (E_3 \cup E_6)$ then
$$
K_{\tilde{X}} + \lambda \tilde{D} + 5 \lambda a_5E_5
$$
is not log canonical at the point $Q$ and so is the pair
$$
K_{\tilde{X}} + \lambda \tilde{D} + E_5 \text{ .}
$$
By adjunction $(E_5, \lambda \tilde{D}|_{E_5})$ is not log
canonical at $Q$ and
$$10a_5 - 6a_3 - 4a_6= \tilde{D} \cdot E_5 \geq \text{mult}_Q\Bigl(\tilde{D}|_{E_5}
\Bigr) =
 \text{mult}_Q\Bigl(\tilde{D} \cdot E_5 \Bigr)  > 6 \text{ ,}$$ which is false.

 \item If the point $Q \in E_5 \cap E_6$ then
$$
K_{\tilde{X}} + \lambda \tilde{D} + 5 \lambda a_5E_5 + 4 \lambda
a_6E_6
$$
is not log canonical at the point $Q$ and so are the pairs
$$
K_{\tilde{X}} + \lambda \tilde{D} + E_5 + 4 \lambda a_6E_6 \text{
and } K_{\tilde{X}} + \lambda \tilde{D} + 5 \lambda a_5E_5 + E_6
\text{ .}
$$
By adjunction
$$10a_5 - 6a_3 - 4a_6 = \tilde{D} \cdot E_5 \geq
\text{mult}_Q\Bigl(\tilde{D}|_{E_5} \Bigr) =
 \text{mult}_Q\Bigl(\tilde{D} \cdot E_5 \Bigr)  > 6 - 4a_6 \text{ and }$$
and
$$8a_6 - 5a_5 - 3a_7 = \tilde{D} \cdot E_6 \geq \text{mult}_Q\Bigl(\tilde{D}|_{E_6}
\Bigr) =
 \text{mult}_Q\Bigl(\tilde{D} \cdot E_6 \Bigr)  > 6 - 5a_5 \text{ ,}$$
which is false.

\item If the point $Q \in E_6 \backslash (E_5 \cup E_7)$ then
$$
K_{\tilde{X}} + \lambda \tilde{D} + 4 \lambda a_6E_6
$$
is not log canonical at the point $Q$ and so is the pair
$$
K_{\tilde{X}} + \lambda \tilde{D} + E_6 \text{ .}
$$
By adjunction $(E_6, \lambda \tilde{D}|_{E_6})$ is not log
canonical at $Q$ and
$$8a_6 - 5a_5 - 3a_7= \tilde{D} \cdot E_6 \geq \text{mult}_Q\Bigl(\tilde{D}|_{E_6} \Bigr)
=
 \text{mult}_Q\Bigl(\tilde{D} \cdot E_6 \Bigr)  > 6 \text{ ,}$$ which is false.

 \item If the point $Q \in E_6 \cap E_7$ then
$$
K_{\tilde{X}} + \lambda \tilde{D} + 4 \lambda a_6E_6 + 3 \lambda
a_7E_7
$$
is not log canonical at the point $Q$ and so are the pairs
$$
K_{\tilde{X}} + \lambda \tilde{D} + E_6 + 3 \lambda a_7E_7 \text{
and } K_{\tilde{X}} + \lambda \tilde{D} + 4 \lambda a_6E_6 + E_7
\text{ .}
$$
By adjunction
$$8a_6 - 5a_5 - 3a_7 = \tilde{D} \cdot E_6 \geq
\text{mult}_Q\Bigl(\tilde{D}|_{E_6} \Bigr) =
 \text{mult}_Q\Bigl(\tilde{D} \cdot E_6 \Bigr)  > 6 - 3a_7 \text{ and }$$
and
$$6a_7 - 4a_6 - 2a_8 = \tilde{D} \cdot E_7 \geq \text{mult}_Q\Bigl(\tilde{D}|_{E_7}
\Bigr) =
 \text{mult}_Q\Bigl(\tilde{D} \cdot E_7 \Bigr)  > 6 - 4a_6 \text{ ,}$$
which is false.

\item If the point $Q \in E_7 \backslash (E_6 \cup E_8)$ then
$$
K_{\tilde{X}} + \lambda \tilde{D} + 3 \lambda a_7E_7
$$
is not log canonical at the point $Q$ and so is the pair
$$
K_{\tilde{X}} + \lambda \tilde{D} + E_7 \text{ .}
$$
By adjunction $(E_7, \lambda \tilde{D}|_{E_7})$ is not log
canonical at $Q$ and
$$6a_7 - 4a_6 - 2a_8= \tilde{D} \cdot E_7 \geq \text{mult}_Q\Bigl(\tilde{D}|_{E_7} \Bigr)
=
 \text{mult}_Q\Bigl(\tilde{D} \cdot E_7 \Bigr)  > 6 \text{ ,}$$ which is false.

\item If the point $Q \in E_7 \cap E_8$ then
$$
K_{\tilde{X}} + \lambda \tilde{D} + 3 \lambda a_7E_7 + 2 \lambda
a_8E_8
$$
is not log canonical at the point $Q$ and so are the pairs
$$
K_{\tilde{X}} + \lambda \tilde{D} + E_7 + 2 \lambda a_8E_8 \text{
and } K_{\tilde{X}} + \lambda \tilde{D} + 3 \lambda a_7E_7 + E_8
\text{ .}
$$
By adjunction
$$6a_7 - 4a_6 - 2a_8 = \tilde{D} \cdot E_7 \geq
\text{mult}_Q\Bigl(\tilde{D}|_{E_7} \Bigr) =
 \text{mult}_Q\Bigl(\tilde{D} \cdot E_7 \Bigr)  > 6 - 2a_8 \text{ and }$$
and
$$4a_8 - 3a_7 = \tilde{D} \cdot E_8 \geq \text{mult}_Q\Bigl(\tilde{D}|_{E_8} \Bigr) =
 \text{mult}_Q\Bigl(\tilde{D} \cdot E_8 \Bigr)  > 6 - 3a_7 \text{ ,}$$
which is false.

\item If the point $Q \in E_8 \backslash E_7$ then
$$
K_{\tilde{X}} + \lambda \tilde{D} + 2 \lambda a_8E_8
$$
is not log canonical at the point $Q$ and so is the pair
$$
K_{\tilde{X}} + \lambda \tilde{D} + E_8 \text{ .}
$$
By adjunction $(E_8, \lambda \tilde{D}|_{E_8})$ is not log
canonical at $Q$ and
$$4a_8 - 3a_7= \tilde{D} \cdot E_8 \geq \text{mult}_Q\Bigl(\tilde{D}|_{E_8} \Bigr) =
 \text{mult}_Q\Bigl(\tilde{D} \cdot E_8 \Bigr)  > 6 \text{ ,}$$ which is false.

\end{itemize}
\end{proof}

\section{Del Pezzo surfaces of degree 1 with at least two Du Val singularities}

Suppose now that $X$ is a del Pezzo surface of degree 1  having at
least two Du Val singular points. We have the following result.

\begin{lemma}
\label{exceptional} Suppose that the surface $X$ has at least one
singularity of type $\mathbb{D}_4 \text{, } \mathbb{D}_5 \text{, }
\mathbb{D}_6 \text{, } \mathbb{E}_6 \text{ .}$ Then the global log
canonical threshold of $X$ is
$$
\mathrm{lct}\big(X\big)=\left\{%
\aligned
&1/3\hspace{0.5cm} \mathrm{when}\ \mathbb{E}_6 \in \mathrm{Sing}(X)\\%
&1/2\hspace{0.5cm} \mathrm{otherwise .}\ \\%
\endaligned\right.%
$$
\end{lemma}

\begin{proof}
We will only treat the case $\mathbb{D}_4$, as the rest of the
cases are similar. Since the linear system $|-K_X|$ is
1-dimensional there is a unique element $Z \in |-K_X|$ that passes
through  the singular point $\mathbb{D}_4$. This curve $Z$ is
irreducible and does not pass through any other singular point of
$X$. Let $\pi: \tilde{X} \to X$ be the minimal resolution of $X$.
Then
$$
\tilde{Z} \sim_{\mathbb{Q}} \pi^*(Z)-E_1-E_2-2E_3-E_4 \text{ , }
$$
where $E_1, E_2, E_3, E_4$ are the exceptional curves of $\pi$
that are contracted to the Du Val singular point $\mathbb{D}_4$.
This means that the global log canonical threshold is
$$
\mathrm{lct}(X) \leq \frac{1}{2} \text{ .}
$$
Now we assume that $\mathrm{lct}(X) < \frac{1}{2}$. Then  there
exists a $\mathbb{Q}$-divisor $D \in X$ such that $D
\sim_{\mathbb{Q}} -K_X$ and the log pair $(X, \lambda D)$ is  not
log canonical, for some rational number $\lambda < \frac{1}{2}$.
According to Lemma~\ref{SingPoint} the pair $(X, \lambda D)$ is
not log canonical at a singular point of $X$. If the pair $(X,
\lambda D)$ is  not log canonical at $\mathbb{D}_4$, we proceed as
in Lemma~\ref{degree1D4}, otherwise we follow the proof of
Theorem~\ref{dpezzoA1A2}. In any case we obtain a contradiction,
thus
$$
\mathrm{lct}(X) =\frac{1}{2} \text{ .}
$$
\end{proof}

\begin{lemma}
\label{cyclicquotientA5A6} Suppose that the surface $X$ has at
least one singularity of type $\mathbb{A}_5 \text{, } \mathbb{A}_6
\text{ .}$ Then the global log canonical threshold of $X$ is
$$
\mathrm{lct}\big(X\big)=\frac{2}{3} \text{ .}
$$
\end{lemma}

\begin{proof}
Again we will consider only the case that $X$ has at least one
$\mathbb{A}_5$ type singular point, as $\mathbb{A}_6$ can be
treated in a similar fashion. Let $\pi: \tilde{X} \to X$ be the
minimal resolution of $X$ and let $E_1, E_2, E_3, E_4, E_5$ be the
exceptional curves of $\pi$ that are contracted to the Du Val
singular point $\mathbb{A}_5$. Then we can always find a -1 curve
$\tilde{L}_3$ in $\tilde{X}$ that only intersects $E_3$
transversally among the exceptional curves of the fundamental
cycle. Then we have that
$$
\tilde{L}_3 \sim_{\mathbb{Q}} \pi^*(L_3)- \frac{1}{2}
E_1-E_2-\frac{3}{2}E_3-E_4- \frac{1}{2} E_5 \text{ , }
$$
 This means that the global log canonical threshold
is
$$
\mathrm{lct}(X) \leq \frac{2}{3} \text{ .}
$$
Now we assume that $\mathrm{lct}(X) < \frac{2}{3}$. Then  there
exists a $\mathbb{Q}$-divisor $D \in X$ such that $D
\sim_{\mathbb{Q}} -K_X$ and the log pair $(X, \lambda D)$ is  not
log canonical, for some rational number $\lambda < \frac{2}{3}$.
According to Lemma~\ref{SingPoint} the pair $(X, \lambda D)$ is
not log canonical at a singular point of $X$. If the pair $(X,
\lambda D)$ is  not log canonical at $\mathbb{A}_5$, we proceed as
in Lemma~\ref{degree1D4}, otherwise we follow the proof of
Theorem~\ref{dpezzoA1A2}. In any case we obtain a contradiction,
thus
$$
\mathrm{lct}(X) =\frac{2}{3} \text{ .}
$$
\end{proof}

\begin{lemma}
\label{cyclicquotientA4} Suppose that the surface $X$ has at least
one singularity of type $\mathbb{A}_4 \text{ .}$ Then the global
log canonical threshold of $X$ is
$$
\mathrm{lct}\big(X\big)=\left\{%
\aligned
 &2/3\ \mathrm{when}\ |-K_{X}|\ \mathrm{has\ a\ cuspidal\
curve}\
C\ \mathrm{such\ that}\ \mathrm{Sing}(C)= \mathbb{A}_2,\\%
&3/4\ \mathrm{when}\ |-K_{X}|\ \mathrm{has\ a\ cuspidal\ curve}\
C\ \mathrm{such\ that}\ \mathrm{Sing}(C)=\mathbb{A}_1\\%
& \hspace{0.75cm} \mathrm{ and\ no\ cuspidal\ curve}\ C\ \mathrm{such\ that}\ \mathrm{Sing}(C)=\mathbb{A}_2,\\%
&4/5\ \mathrm{in\ the\ remaining\ cases}.\\%
\endaligned\right.%
$$
\end{lemma}

\begin{proof} Let $\pi: \tilde{X} \to X$ be the
minimal resolution of $X$ and let $E_1, E_2, E_3, E_4$ be the
exceptional curves of $\pi$ that are contracted to the Du Val
singular point $\mathbb{A}_4$. In all the cases
when we have at least an $\mathbb{A}_4$ type singularity there exists
a unique smooth irreducible element $C$ of the linear system $|-2K_X|$, 
which passes through the intersection point $E_1 \cap E_2$.
For the pull back of the irreducible curve $C$ we have 
$$
\tilde{C} + E_1 + 2 E_2 + 2E_3 + E_4 \in |-K_{\tilde{X}}| \text{  .}
$$
If we blow up once more in order to get transversal intersections 
we see that the global log canonical threshold is
$$
\mathrm{lct}(X) \leq \mathrm{lct}(X, C) = \frac{4}{5} \text{ .}
$$
Now we assume that $\mathrm{lct}(X) <
\frac{4}{5}$. Then there exists a $\mathbb{Q}$-divisor $D \in X$
such that $D \sim_{\mathbb{Q}} -K_X$ and the log pair $(X, \lambda
D)$ is not log canonical, for some rational number $\lambda <
\frac{4}{5}$. According to Lemma~\ref{SingPoint} the pair $(X,
\lambda D)$ is  not log canonical at a singular point of $X$. If
the pair $(X, \lambda D)$ is  not log canonical at $\mathbb{A}_4$,
we proceed as in Lemma~\ref{A4}, otherwise we follow the proof of
Theorem~\ref{dpezzoA1A2}. In any case we obtain a contradiction,
and the result follows.
\end{proof}

\begin{lemma}
\label{cyclicquotientA3} Suppose that the surface $X$ has at least
one singularity of type $\mathbb{A}_3$ and no singularity of type
$\mathbb{A}_4$, $\mathbb{D}_4$, $\mathbb{D}_5$. Then the global
log canonical threshold of $X$ is
$$
\mathrm{lct}\big(X\big)=\left\{%
\aligned
&1\hspace{0.5cm} \mathrm{when}\ |-K_{X}|\ \mathrm{does\ not\ have\ cuspidal\ curves},\\%
&2/3\ \mathrm{when}\ |-K_{X}|\ \mathrm{has\ a\ cuspidal\ curve}\
C\ \mathrm{such\ that}\ \mathrm{Sing}(C)= \mathbb{A}_2,\\%
&3/4\ \mathrm{when}\ |-K_{X}|\ \mathrm{has\ a\ cuspidal\ curve}\
C\ \mathrm{such\ that}\ \mathrm{Sing}(C)=\mathbb{A}_1\\%
& \hspace{0.75cm} \mathrm{ and\ no\ cuspidal\ curve}\ C\ \mathrm{such\ that}\ \mathrm{Sing}(C)=\mathbb{A}_2,\\%
&5/6\ \mathrm{in\ the\ remaining\ cases}.\\%
\endaligned\right.%
$$
\end{lemma}

\begin{proof} Let $\pi: \tilde{X} \to X$ be the
minimal resolution of $X$ and let $E_1, E_2, E_3$ be the
exceptional curves of $\pi$ that are contracted to the Du Val
singular point $\mathbb{A}_3$. One can show that in all the cases
when we have at least an $\mathbb{A}_3$ type singularity we must
have $\mathrm{lct}(X) \leq 1$.

 Now we assume that $\mathrm{lct}(X) <
1$. Then there exists a $\mathbb{Q}$-divisor $D \in X$ such that
$D \sim_{\mathbb{Q}} -K_X$ and the log pair $(X, \lambda D)$ is
not log canonical, for some rational number $\lambda < 1$.
According to Lemma~\ref{SingPoint} the pair $(X, \lambda D)$ is
not log canonical at a singular point of $X$. If the pair $(X,
\lambda D)$ is  not log canonical at $\mathbb{A}_3$, we proceed as
in Lemma~\ref{A3}, otherwise we follow the proof of
Theorem~\ref{dpezzoA1A2}. In any case we obtain a contradiction,
and the result follows.
\end{proof}

\begin{lemma}
\label{exactlyA3or2A3} Suppose that the surface $X$ has exactly
two singularities of type $\mathbb{A}_3$. Then the global log
canonical threshold of $X$ is
$$
\mathrm{lct}\big(X\big)=1
$$
\end{lemma}

\begin{proof}
See Lemma~\ref{A3}.
\end{proof}

\section{Del Pezzo surfaces with Picard group $\mathbb{Z}$}

All possible singular points on a del Pezzo surface $X$ that has
only Du Val singular points and $\text{Pic}(X) \cong \mathbb{Z}$
are listed in \cite{Zhang}.

\subsection{Del Pezzo surface of degree 1 with an $\mathbb{E}_7$
and an $\mathbb{A}_1$ type singularity}

In this section we will prove the following.

\begin{lemma}
\label{E7+A1} Let $X$ be a del Pezzo surface with  one Du Val
singularity of type $\mathbb{E}_7$, one of type $\mathbb{A}_1$ and
$K_X^2=1$. Then the global log canonical threshold of $X$ is
$$
\mathrm{lct} (X) = \frac{1}{4} \text{ .}
$$
\end{lemma}

\begin{proof}

Suppose that $\mathrm{lct}(X)<\frac{1}{4}$ , then there exists a
$\mathbb{Q}$-divisor $D \in X$ such that the log pair $(X,\lambda
D)$ is  not log canonical and $D \sim_{\mathbb{Q}} -K_X$ for some
rational number $\lambda < \frac{1}{4}$. We derive that the pair
$(X, \lambda D)$ is log canonical outside of a point $P\in X$ and
not log canonical at $P$. Let $\pi: \tilde{X} \to X$ be the
minimal resolution of $X$. The configuration of the exceptional
curves is given by the following Dynkin diagram.
\bigskip

$\mathbb{E}_7 + \mathbb{A}_1 \text{ .  }$
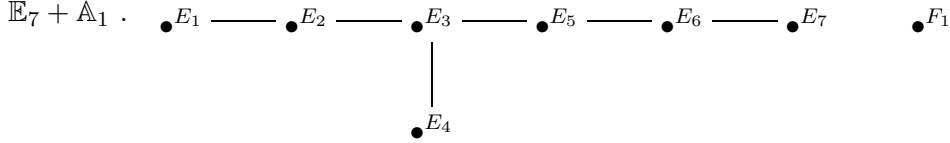
\xymatrix{ {\bullet}^{E_1} \ar@{-}[r] & {\bullet}^{E_2}
\ar@{-}[r]  & {\bullet}^{E_3} \ar@{-}[r] \ar@{-}[d] &
{\bullet}^{E_5} \ar@{-}[r] & {\bullet}^{E_6} \ar@{-}[r] & {\bullet}^{E_7} & {\bullet}^{F_1}\\
 & & {\bullet}^{E_4} & & & & }
\bigskip

\noindent Then
$$
\tilde{D} \sim_{\mathbb{Q}} \pi^*(D) - a_1E_1 - a_2E_2 - a_3E_3 -
a_4E_4 - a_5E_5 - a_6E_6 - a_7E_7 - b_1F_1 \text{ .}
$$
\noindent We should note here that there are two -1 curves $L_1,
L_7$ such that
$$
L_1 \cdot E_1 = L_7 \cdot E_7= L_7 \cdot F_1 = 1 \text{ .}
$$
Therefore we have
\begin{eqnarray*}
 \tilde{L}_1 & \sim_{\mathbb{Q}} & \pi^*(L_1) -2E_1 -3 E_2 - 4E_3 -2 E_4 -3 E_5-2 E_6-E_7\\
\tilde{L}_7 & \sim_{\mathbb{Q}} & \pi^*(L_7) -  E_1 - 2 E_2 - 3E_3
- \frac{3}{2} E_4 - \frac{5}{2} E_5 - 2E_6 - \frac{3}{2} E_7 -
\frac{1}{2} F_1\text{ .}
\end{eqnarray*}
and since  $L_7 \sim -K_X$ and $L_1 \sim -K_X$ we see that
$\mathrm{lct}(X) \leq \frac{1}{4}$.

The inequalities
\begin{eqnarray*}
0 \leq \tilde{D} \cdot \tilde{L}_7 & = & 1 - a_7 - b_1 \\
0 \leq \tilde{D} \cdot \tilde{L}_1 & = & 1 - a_1 \\
0 \leq E_1 \cdot \tilde{D} & = & 2a_1 - a_2\\
0 \leq E_2 \cdot \tilde{D} & = & 2a_2 - a_1 - a_3\\
0 \leq E_3 \cdot \tilde{D} & = & 2a_3 - a_2 - a_5 - a_4\\
0 \leq E_4 \cdot \tilde{D} & = & 2a_4 - a_3\\
0 \leq E_5 \cdot \tilde{D} & = & 2a_5 - a_3 - a_6\\
0 \leq E_6 \cdot \tilde{D} & = & 2a_6 - a_5 - a_7\\
0 \leq E_7 \cdot \tilde{D} & = & 2a_7 - a_6\\
0 \leq F_1 \cdot \tilde{D} & = & 2b_1
\end{eqnarray*}
imply that $a_1 \leq 1 \text{, } a_2 \leq 2 \text{, }
 a_3 \leq 3 \text{, } a_4 \leq \frac{7}{4} \text{, }
 a_5 \leq \frac{5}{2} \text{, } a_6 \leq 2 \text{, }
 a_7 \leq 1 \text{, } b_1 \leq 1 \text{ .}
$ Moreover we have
$$
2a_1 \geq a_2 \text{, } \frac{3}{2} a_2 \geq a_3 \text{, } 2 a_4
\geq a_3 \text{, } \frac{5}{6} a_3 \geq a_5 \text{, } \frac{4}{5}
a_5 \geq a_6
$$
and
$$
2a_7 \geq a_6 \text{, } \frac{3}{2} a_6 \geq a_5 \text{, }
\frac{4}{3} a_5 \geq a_3 \text{ .}
$$

\medskip

The equivalence
$$
K_{\tilde{X}} + \lambda \tilde{D} + \lambda a_1E_1 + \lambda
a_2E_2 +  \lambda a_3E_3 +  \lambda a_4E_4 + \lambda a_5E_5 +
\lambda a_6E_6 + \lambda a_7E_7 \sim_{\mathbb{Q}}
\pi^*(K_X+\lambda D)
$$
implies that there is a point $Q\in E_1\cup E_2\cup E_3\cup E_4
\cup E_5 \cup E_6 \cup E_7$ such that the pair
$$
K_{\tilde{X}} + \lambda \tilde{D} + \lambda a_1E_1 +  \lambda
a_2E_2 + \lambda a_3E_3 +  \lambda a_4E_4 +
 \lambda a_5E_5 +  \lambda a_6E_6 + \lambda a_7E_7
$$
is not log canonical at $Q$.

\begin{itemize}
\item If the point $Q \in E_1$ and $Q\not \in E_2$ then
$$
K_{\tilde{X}} + \lambda \tilde{D} + \lambda a_1E_1
$$
is not log canonical at the point $Q$ and so is the pair
$$
K_{\tilde{X}} + \lambda \tilde{D} + E_1 \text{ since } \lambda a_1
\leq 1 \text{ .}
$$
By adjunction $(E_1, \lambda \tilde{D}|_{E_1})$ is not log
canonical at $Q$ and
$$2 \geq 2a_1 \geq 2a_1 - a_2 = \tilde{D} \cdot E_1 \geq
 \text{mult}_Q\Bigl(\tilde{D} \cdot E_1 \Bigr)  > 4 \text{ ,}$$ which is false.

 \item If the point $Q \in E_1 \cap E_2$ then
$$
K_{\tilde{X}} + \lambda \tilde{D} +\lambda a_1E_1 + \lambda a_2E_2
$$
is not log canonical at the point $Q$ and so are the pairs
$$
K_{\tilde{X}} + \lambda \tilde{D} + E_1 + \lambda a_2E_2 \text{
and } K_{\tilde{X}} + \lambda \tilde{D} + \lambda a_1E_1 + E_2
\text{ .}
$$
By adjunction
$$2 - a_2 \geq 2a_1 - a_2 = \tilde{D} \cdot E_1 \geq \text{mult}_Q\Bigl(\tilde{D}|_{E_1}
\Bigr) =
 \text{mult}_Q\Bigl(\tilde{D} \cdot E_1 \Bigr)  > 4 - a_2 \text{ and }$$
and
$$2 - a_1 \geq 2a_2 - a_1 -a_3 = \tilde{D} \cdot E_2 \geq
\text{mult}_Q\Bigl(\tilde{D}|_{E_2} \Bigr) =
 \text{mult}_Q\Bigl(\tilde{D} \cdot E_2 \Bigr)  > 4 - a_1 \text{ ,}$$
which is false.

\item If the point $Q \in E_2 \backslash (E_1 \cup E_3)$ then
$$
K_{\tilde{X}} + \lambda \tilde{D} + \lambda a_2E_2
$$
is not log canonical at the point $Q$ and so is the pair
$$
K_{\tilde{X}} + \lambda \tilde{D} + E_2 \text{ .}
$$
By adjunction
$$2a_2 \geq 2a_2 -  a_1 - a_3 = \tilde{D} \cdot E_2 \geq
\text{mult}_Q\Bigl(\tilde{D}|_{E_2} \Bigr) =
 \text{mult}_Q\Bigl(\tilde{D} \cdot E_2 \Bigr)  > 4 \text{ ,}$$
implies that $a_2 > 2$ which is false.

\item If the point $Q \in E_2 \cap E_3$ then
$$
K_{\tilde{X}} + \lambda \tilde{D} + \lambda a_2E_2 + \lambda
a_3E_3
$$
is not log canonical at the point $Q$ and so are the pairs
$$
K_{\tilde{X}} + \lambda \tilde{D} + E_2 + \lambda a_3E_3 \text{
and } K_{\tilde{X}} + \lambda \tilde{D} + \lambda a_2E_2 + E_3
\text{ .}
$$
By adjunction
$$4 -a_3 \geq 2a_2 - a_3 \geq 2a_2-a_1-a_3 = \tilde{D} \cdot E_2 \geq
\text{mult}_Q\Bigl(\tilde{D}|_{E_2} \Bigr) =
 \text{mult}_Q\Bigl(\tilde{D} \cdot E_2 \Bigr)  > 4 - a_3 $$
which is a contradiction.

\item If the point $Q \in E_3 \backslash (E_2 \cup E_4 \cup E_5)$
then
$$
K_{\tilde{X}} + \lambda \tilde{D} + \lambda a_3E_3
$$
is not log canonical at the point $Q$ and so is the pair
$$
K_{\tilde{X}} + \lambda \tilde{D} + E_3 \text{ .}
$$
By adjunction $(E_3, \lambda \tilde{D}|_{E_3})$ is not log
canonical at $Q$ and
$$\frac{5}{2} \geq \frac{5}{6} a_3 \geq 2a_3 - \frac{2}{3} a_3 - \frac{a_3}{2}  \geq 2a_3
- a_2 - a_4 - a_5 = \tilde{D} \cdot E_3 \geq
 \text{mult}_Q\Bigl(\tilde{D} \cdot E_3 \Bigr)  > 4 \text{ ,}$$ which is false.

\item If the point $Q \in E_3 \cap E_4$ then
$$
K_{\tilde{X}} + \lambda \tilde{D} + \lambda a_3E_3 + \lambda
a_4E_4
$$
is not log canonical at the point $Q$ and so is the pair
$$
K_{\tilde{X}} + \lambda \tilde{D} + \lambda a_3E_3 + E_4 \text{ .}
$$
By adjunction
$$\frac{7}{2} - a_3 \geq 2a_4 - a_3 = \tilde{D} \cdot E_4 \geq
\text{mult}_Q\Bigl(\tilde{D}|_{E_4} \Bigr) =
 \text{mult}_Q\Bigl(\tilde{D} \cdot E_4 \Bigr)  > 4 - a_3 \text{ ,}$$
which is false.

\item If the point $Q \in E_4 \backslash E_3$ then
$$
K_{\tilde{X}} + \lambda \tilde{D} + \lambda a_4E_4
$$
is not log canonical at the point $Q$ and so is the pair
$$
K_{\tilde{X}} + \lambda \tilde{D} + E_4 \text{ .}
$$
By adjunction $(E_4, \lambda \tilde{D}|_{E_4})$ is not log
canonical at $Q$ and
$$\frac{7}{2} \geq 2a_4 \geq 2a_4 - a_3 = \tilde{D} \cdot E_4 \geq
 \text{mult}_Q\Bigl(\tilde{D} \cdot E_4 \Bigr)  > 4 \text{ ,}$$ which is false.

 \item If the point $Q \in E_3 \cap E_5$ then
$$
K_{\tilde{X}} + \lambda \tilde{D} + \lambda a_3E_3 + \lambda
a_5E_5
$$
is not log canonical at the point $Q$ and so are the pairs
$$
K_{\tilde{X}} + \lambda \tilde{D} +  \lambda a_3E_3 + E_5 \text{
.}
$$
By adjunction
$$\frac{10}{3} - a_3 \geq \frac{4}{3}a_5 -a_3 \geq 2a_5 - \frac{2}{3}a_5 - a_3 \geq 2a_5
- a_3 - a_6 = \tilde{D} \cdot E_5 \geq
 \text{mult}_Q\Bigl(\tilde{D} \cdot E_5 \Bigr)  > 4 - a_3 \text{ ,}$$
which is false.

\item If the point $Q \in E_5 \backslash (E_3 \cup E_6)$ then
$$
K_{\tilde{X}} + \lambda \tilde{D} + \lambda a_5E_5
$$
is not log canonical at the point $Q$ and so is the pair
$$
K_{\tilde{X}} + \lambda \tilde{D} + E_5 \text{ .}
$$
By adjunction $(E_5, \lambda \tilde{D}|_{E_5})$ is not log
canonical at $Q$ and
$$\frac{1}{3} \geq \frac{2}{15} a_5 \geq 2a_5 - \frac{6}{5} a_5 - \frac{2}{3} a_5 \geq
2a_5 - a_3 - a_6= \tilde{D} \cdot E_5 \geq
 \text{mult}_Q\Bigl(\tilde{D} \cdot E_5 \Bigr)  > 4 \text{ ,}$$ which is false.

 \item If the point $Q \in E_5 \cap E_6$ then
$$
K_{\tilde{X}} + \lambda \tilde{D} + \lambda a_5E_5 + \lambda
a_6E_6
$$
is not log canonical at the point $Q$ and so are the pairs
$$
K_{\tilde{X}} + \lambda \tilde{D} + E_5 + \lambda a_6E_6 \text{ .}
$$
By adjunction
$$2a_5 - \frac{6}{5} a_5 \geq 2a_5 - a_3 -a_6 = \tilde{D} \cdot E_5 \geq
\text{mult}_Q\Bigl(\tilde{D}|_{E_5} \Bigr) =
 \text{mult}_Q\Bigl(\tilde{D} \cdot E_5 \Bigr)  > 4 - a_6 $$
implies that $a_5 > 5$ which is false.

\item If the point $Q \in E_6 \backslash (E_5 \cup E_7)$ then
$$
K_{\tilde{X}} + \lambda \tilde{D} + \lambda a_6E_6
$$
is not log canonical at the point $Q$ and so is the pair
$$
K_{\tilde{X}} + \lambda \tilde{D} + E_6 \text{ .}
$$
By adjunction $(E_6, \lambda \tilde{D}|_{E_6})$ is not log
canonical at $Q$ and
$$2a_6 - \frac{5}{4}a_6 - \frac{a_6}{2} \geq 2a_6 - a_5 - a_7= \tilde{D} \cdot E_6 \geq
 \text{mult}_Q\Bigl(\tilde{D} \cdot E_6 \Bigr)  > 4$$
 implies that $a_6 > 16$ which is false.

 \item If the point $Q \in E_6 \cap E_7$ then
$$
K_{\tilde{X}} + \lambda \tilde{D} +  \lambda a_6E_6 + \lambda
a_7E_7
$$
is not canonical at the point $Q$ and so are the pairs
$$
K_{\tilde{X}} + \lambda \tilde{D} + E_6 + \lambda a_7E_7 \text{ .}
$$
By adjunction
$$2a_6 - a_7 \geq 2a_6 -a_5 - a_7 = \tilde{D} \cdot E_6 \geq
\text{mult}_Q\Bigl(\tilde{D}|_{E_6} \Bigr) =
 \text{mult}_Q\Bigl(\tilde{D} \cdot E_6 \Bigr)  > 4 - a_7 $$
implies that $a_6 >2$ which is false.

\item If the point $Q \in E_7 \backslash E_6$ then
$$
K_{\tilde{X}} + \lambda \tilde{D} + \lambda a_7E_7
$$
is not log canonical at the point $Q$ and so is the pair
$$
K_{\tilde{X}} + \lambda \tilde{D} + E_7 \text{ .}
$$
By adjunction $(E_7, \lambda \tilde{D}|_{E_7})$ is not log
canonical at $Q$ and
$$2a_7 - a_6= \tilde{D} \cdot E_7 \geq \text{mult}_Q\Bigl(\tilde{D}|_{E_7} \Bigr) =
 \text{mult}_Q\Bigl(\tilde{D} \cdot E_7 \Bigr)  > 4 \text{ ,}$$
 implies that $a_7 >2$ which is false.

 \item If the point $Q \in F_1$  then
$$
K_{\tilde{X}} + \lambda \tilde{D} + \lambda b_1F_1
$$
is not log canonical at the point $Q$ and so is the pair
$$
K_{\tilde{X}} + \lambda \tilde{D} + F_1 \text{ since } \lambda b_1
\leq 1 \text{ .}
$$
By adjunction $(F_1, \lambda \tilde{D}|_{F_1})$ is not log
canonical at $Q$ and
$$2 \geq 2b_1 = \tilde{D} \cdot F_1 \geq
 \text{mult}_Q\Bigl(\tilde{D} \cdot F_1 \Bigr)  > 4 \text{ ,}$$ which is false.

\end{itemize}
\end{proof}

\subsection{Del Pezzo surface of degree 1 with an $\mathbb{E}_6$
and an $\mathbb{A}_2$ type singularity}

In this section we will prove the following.

\begin{lemma}
\label{E6+A2} Let $X$ be a del Pezzo surface with one Du Val
singularity of type $\mathbb{E}_6$, one of type $\mathbb{A}_2$ and
$K_X^2=1$. Then the global log canonical threshold of $X$ is
$$
\mathrm{lct} (X) = \frac{1}{3} \text{ .}
$$
\end{lemma}

\begin{proof}

Suppose that $\mathrm{lct}(X)<\frac{1}{3}$,  then there exists a
$\mathbb{Q}$-divisor $D \in X$ such that the log pair $(X,\lambda
D)$ is  not log canonical, where $\lambda < \frac{1}{3}$  and $D
\sim_{\mathbb{Q}} -K_X$. We derive that the pair $(X,\lambda D)$
is log canonical everywhere outside of a point $P\in X$ and not
log canonical at $P$. Let $\pi: \tilde{X} \to X$ be the minimal
resolution of $X$. The configuration of the exceptional curves is
given by the following Dynkin diagram.
\bigskip

$\mathbb{E}_6 + \mathbb{A}_2 \text{ .  }$
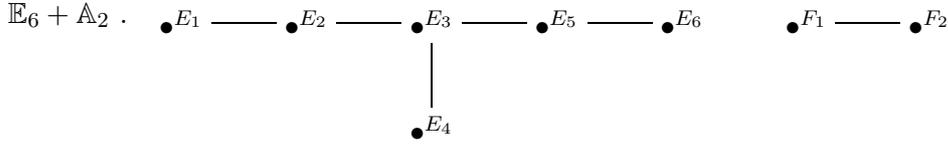
\xymatrix{ {\bullet}^{E_1} \ar@{-}[r] & {\bullet}^{E_2}
\ar@{-}[r]  & {\bullet}^{E_3} \ar@{-}[r] \ar@{-}[d] &
{\bullet}^{E_5} \ar@{-}[r] & {\bullet}^{E_6} & {\bullet}^{F_1} \ar@{-}[r] & {\bullet}^{F_2}\\
 & & {\bullet}^{E_4} & & & & }
\bigskip

\noindent Then
$$
\tilde{D} \sim_{\mathbb{Q}}
\pi^*(D)-a_1E_1-a_2E_2-a_3E_3-a_4E_4-a_5E_5-a_6E_6 - b_1F_1 -
b_2F_2 \text{ .}
$$

\noindent We should note here that there are two -1 curves
$\tilde{L}_4, \tilde{L}_6$ such that
$$
\tilde{L}_4 \cdot E_4 = \tilde{L}_6 \cdot E_6= \tilde{L}_6 \cdot
F_1 = 1 \text{ .}
$$
Therefore we have
\begin{eqnarray*}
\tilde{L}_4 & \sim_{\mathbb{Q}} & \pi^*(L_4) - E_1 - 2E_2 - 3E_3 - 2E_4 - 2E_5 - E_6\\
\tilde{L}_6 & \sim_{\mathbb{Q}} & \pi^*(L_6) - \frac{2}{3} E_1 -
\frac{4}{3} E_2 - 2E_3 - E_4 - \frac{5}{3} E_5 - \frac{4}{3} E_6 -
\frac{2}{3} F_1 - \frac{1}{3} F_2\text{ .}
\end{eqnarray*}
and since  $L_6 \sim_{\mathbb{Q}} -K_X$ and $L_4 \sim_{\mathbb{Q}}
-K_X$ we see that $\mathrm{lct}(X) \leq \frac{1}{3}$.

The inequalities
\begin{eqnarray*}
0 \leq \tilde{D} \cdot \tilde{L}_6 & = & 1 - a_6 - b_1\\
0 \leq \tilde{D} \cdot \tilde{L}_4 & = & 1 - a_4\\
0 \leq E_1 \cdot \tilde{D} & = & 2a_1 - a_2\\
0 \leq E_2 \cdot \tilde{D} & = & 2a_2 - a_1 - a_3\\
0 \leq E_3 \cdot \tilde{D} & = & 2a_3 - a_2 - a_4 - a_5\\
0 \leq E_4 \cdot \tilde{D} & = & 2a_4 - a_3\\
0 \leq E_5 \cdot \tilde{D} & = & 2a_5 - a_3 - a_6\\
0 \leq E_6 \cdot \tilde{D} & = & 2a_6 - a_5\\
0 \leq F_1 \cdot \tilde{D} & = & 2b_1 - b_2\\
0 \leq F_2 \cdot \tilde{D} & = & 2b_2 - b_1\\
\end{eqnarray*}
imply that $a_1 \leq \frac{4}{3} \text{, } a_2 \leq \frac{5}{3}
\text{, } a_3 \leq 2 \text{, } a_4 \leq 1 \text{, } a_5 \leq
\frac{5}{3} \text{, } a_6 \leq 1 \text{, } b_1 \leq 1 \text{, }
b_2 \leq 2 \text{ .}$ The equivalence
$$
K_{\tilde{X}} + \lambda \tilde{D} +
 \lambda a_1E_1 + \lambda  a_2E_2 + \lambda a_3E_3 +  \lambda a_4E_4 +
 \lambda a_5E_5 + \lambda a_6E_6 + \lambda b_1 F_1 + \lambda b_2 F_2
\sim_{\mathbb{Q}} \pi^*(K_X+ \lambda D)
$$
implies that there is a point $Q\in E_1\cup E_2\cup E_3\cup E_4
\cup E_5 \cup E_6 \cup F_1 \cup F_2$ such that the pair
$$
K_{\tilde{X}} + \lambda \tilde{D} +
 \lambda a_1E_1 + \lambda  a_2E_2 + \lambda a_3E_3 +  \lambda a_4E_4 +
 \lambda a_5E_5 + \lambda a_6E_6 + \lambda b_1 F_1 + \lambda b_2 F_2
$$
is not log canonical at $Q$.

\begin{itemize}
\item If the point $Q \in E_1$ and $Q\not \in E_2$ then
$$
K_{\tilde{X}} + \lambda \tilde{D} + \lambda a_1E_1
$$
is not log canonical at the point $Q$ and so is the pair
$$
K_{\tilde{X}} + \lambda \tilde{D} + E_1 \text{ since } \lambda a_1
\leq 1 \text{ .}
$$
By adjunction $(E_1, \lambda \tilde{D}|_{E_1})$ is not log
canonical at $Q$ and
$$2a_1 - \frac{5}{4} a_1 \geq 2a_1 - a_2 = \tilde{D} \cdot E_1 \geq
\text{mult}_Q\Bigl(\tilde{D}|_{E_1} \Bigr) =
 \text{mult}_Q\Bigl(\tilde{D} \cdot E_1 \Bigr)  > 3 $$
 implies that $a_1 > 4$ which is false.

 \item If the point $Q \in E_1 \cap E_2$ then
$$
K_{\tilde{X}} + \lambda \tilde{D} + \lambda a_1E_1 + \lambda
a_2E_2
$$
is not log canonical at the point $Q$ and so are the pairs
$$
K_{\tilde{X}} + \lambda \tilde{D} + E_1 + \lambda a_2E_2 \text{ .}
$$
By adjunction
$$2a_1 - a_2 = \tilde{D} \cdot E_1 \geq \text{mult}_Q\Bigl(\tilde{D}|_{E_1} \Bigr) =
 \text{mult}_Q\Bigl(\tilde{D} \cdot E_1 \Bigr)  > 3 - a_2 $$
implies that $a_1 > \frac{3}{2}$ which is false.

\item If the point $Q \in E_2 \backslash (E_1 \cup E_3)$ then
$$
K_{\tilde{X}} + \lambda \tilde{D} + \lambda a_2E_2
$$
is not log canonical at the point $Q$ and so is the pair
$$
K_{\tilde{X}} + \lambda \tilde{D} + E_2 \text{ .}
$$
By adjunction
$$2a_2 - \frac{1}{2}a_2 - \frac{6}{5}a_2 \geq2a_2 - a_1 - a_3 = \tilde{D} \cdot E_2 \geq
\text{mult}_Q\Bigl(\tilde{D}|_{E_2} \Bigr) =
 \text{mult}_Q\Bigl(\tilde{D} \cdot E_2 \Bigr)  > 3 \text{ ,}$$
implies that $a_2 > 10$ which is false.

\item If the point $Q \in E_2 \cap E_3$ then
$$
K_{\tilde{X}} + \lambda \tilde{D} + \lambda a_2E_2 + \lambda
a_3E_3
$$
is not log canonical at the point $Q$ and so are the pairs
$$
K_{\tilde{X}} + \lambda \tilde{D} + E_2 + \lambda a_3E_3 \text{ .}
$$
By adjunction
$$2a_2 - \frac{1}{2} a_2 - a_3 \geq 2a_2 - a_3 - a_1 = \tilde{D} \cdot E_2 \geq
\text{mult}_Q\Bigl(\tilde{D}|_{E_2} \Bigr) =
 \text{mult}_Q\Bigl(\tilde{D} \cdot E_2 \Bigr)  > 3 - a_3 $$
implies $a_2 > 2$ which is false.

\item If the point $Q \in E_3 \backslash (E_2 \cup E_4 \cup E_5)$
then
$$
K_{\tilde{X}} + \lambda \tilde{D} + \lambda a_3E_3
$$
is not log canonical at the point $Q$ and so is the pair
$$
K_{\tilde{X}} + \lambda \tilde{D} + E_3 \text{ .}
$$
By adjunction $(E_3, \lambda \tilde{D}|_{E_3})$ is not log
canonical at $Q$ and
$$2a_3 - \frac{2}{3}a_3 - \frac{1}{2}a_3 - \frac{2}{3} a_3 \geq 2a_3 - a_2 - a_4 - a_5 =
\tilde{D} \cdot E_3 \geq
 \text{mult}_Q\Bigl(\tilde{D} \cdot E_3 \Bigr)  > 3 $$
 implies $a_3 > 18$ which is false.

\item If the point $Q \in E_3 \cap E_4$ then
$$
K_{\tilde{X}} + \lambda \tilde{D} + \lambda a_3E_3 + \lambda
a_4E_4
$$
is not log canonical at the point $Q$ and so is the pair
$$
K_{\tilde{X}} + \lambda \tilde{D} + \lambda a_3E_3 + E_4 \text{ .}
$$
By adjunction
$$2a_4 - a_3 = \tilde{D} \cdot E_4 \geq \text{mult}_Q\Bigl(\tilde{D}|_{E_4} \Bigr) =
 \text{mult}_Q\Bigl(\tilde{D} \cdot E_4 \Bigr)  > 3 - a_3 \text{ ,}$$
which is false, since $a_4 \leq 1$.

\item If the point $Q \in E_4 \backslash E_3$ then
$$
K_{\tilde{X}} + \lambda \tilde{D} + \lambda a_4E_4
$$
is not log canonical at the point $Q$ and so is the pair
$$
K_{\tilde{X}} + \lambda \tilde{D} + E_4 \text{ .}
$$
By adjunction $(E_4, \lambda \tilde{D}|_{E_4})$ is not log
canonical at $Q$ and
$$2a_4 \geq 2a_4 - a_3 = \tilde{D} \cdot E_4 \geq \text{mult}_Q\Bigl(\tilde{D}|_{E_4}
\Bigr) =
 \text{mult}_Q\Bigl(\tilde{D} \cdot E_4 \Bigr)  > 3$$ which is false, since $a_4 \leq 1$.

\item If the point $Q \in F_1$ and $Q\not \in F_2$ then
$$
K_{\tilde{X}} + \lambda \tilde{D} + \lambda b_1F_1
$$
is not log canonical at the point $Q$ and so is the pair
$$
K_{\tilde{X}} + \lambda \tilde{D} + F_1 \text{ since } \lambda b_1
\leq 1 \text{ .}
$$
By adjunction $(F_1, \lambda \tilde{D}|_{F_1})$ is not log
canonical at $Q$ and
$$ 2b_1 - b_2 = \tilde{D} \cdot F_1 \geq \text{mult}_Q\Bigl(\tilde{D}|_{F_1} \Bigr) =
 \text{mult}_Q\Bigl(\tilde{D} \cdot F_1 \Bigr)  > 3 $$
 implies that $b_1 > \frac{3}{2}$ which is false.

 \item If the point $Q \in F_1 \cap F_2$ then
$$
K_{\tilde{X}} + \lambda \tilde{D} + \lambda b_1F_1 + \lambda
b_2F_2
$$
is not log canonical at the point $Q$ and so are the pairs
$$
K_{\tilde{X}} + \lambda \tilde{D} + F_1 +  b_2F_2 \text{ .}
$$
By adjunction
$$2b_1 - b_2 = \tilde{D} \cdot F_1 \geq \text{mult}_Q\Bigl(\tilde{D}|_{F_1} \Bigr) =
 \text{mult}_Q\Bigl(\tilde{D} \cdot F_1 \Bigr)  > 3 - b_2 $$
implies that $b_1 > \frac{3}{2}$ which is false.

\item If the point $Q \in F_2$ and $Q\not \in F_1$ then
$$
K_{\tilde{X}} + \lambda \tilde{D} + \lambda b_2F_2
$$
is not log canonical at the point $Q$ and so is the pair
$$
K_{\tilde{X}} + \lambda \tilde{D} + F_2 \text{ since }
\frac{1}{3}b_2 \leq 1 \text{ .}
$$
By adjunction $(F_2, \lambda \tilde{D}|_{F_2})$ is not log
canonical at $Q$ and
$$ \frac{3}{2}b_2 \geq 2b_2 - b_1 = \tilde{D} \cdot F_2 \geq
\text{mult}_Q\Bigl(\tilde{D}|_{F_2} \Bigr) =
 \text{mult}_Q\Bigl(\tilde{D} \cdot F_2 \Bigr)  > 3 $$
 implies that $b_2 > 2$ which is false.

\end{itemize}
\end{proof}

\subsection{Del Pezzo surface of degree 1 with an $\mathbb{A}_7$
and an $\mathbb{A}_1$ type singularity}

In this section we will prove the following.

\begin{lemma}
\label{A7+A1} Let $X$ be a del Pezzo surface with one Du Val
singularity of type $\mathbb{A}_7$, one of type $\mathbb{A}_1$ and
$K_X^2=1$. Then the global log canonical threshold of $X$ is
$$
\mathrm{lct} (X) = \frac{1}{2} \text{ .}
$$
\end{lemma}

\begin{proof}

Suppose that $\mathrm{lct}(X)<\frac{1}{2}$,  then there exists a
$\mathbb{Q}$-divisor $D \in X$ such that the log pair $(X, \lambda
D)$ is not log canonical, where $\lambda < \frac{1}{2}$  and $D
\sim_{\mathbb{Q}} -K_X$. We derive that the pair $(X, \lambda D)$
is log canonical everywhere outside of a point $P\in X$ and not
log canonical at $P$. Let $\pi: \tilde{X} \to X$ be the minimal
resolution of $X$. The configuration of the exceptional curves is
given by the following Dynkin diagram.
\bigskip

$\mathbb{A}_7 + \mathbb{A}_1$.
\xymatrix{ {\bullet}^{E_1} \ar@{-}[r] &
{\bullet}^{E_2} \ar@{-}[r] & {\bullet}^{E_3} \ar@{-}[r] &
{\bullet}^{E_4} \ar@{-}[r] & {\bullet}^{E_5} \ar@{-}[r] &
{\bullet}^{E_6} \ar@{-}[r] & {\bullet}^{E_7} & {\bullet}^{F_1}}
\bigskip

\noindent Then
$$
\tilde{D} \sim_{\mathbb{Q}}
\pi^*(D)-a_1E_1-a_2E_2-a_3E_3-a_4E_4-a_5E_5-a_6E_6 - a_7E_7 -
b_1F_1  \text{ .}
$$

\noindent We should note here that there are three -1 curves
$\tilde{L}_1, \tilde{L}_4, \tilde{L}_6$ such that
$$
\tilde{L}_1 \cdot E_1 = \tilde{L}_1 \cdot E_7 = \tilde{L}_4 \cdot
E_4 = \tilde{L}_6 \cdot E_6= \tilde{L}_6 \cdot F_1 = 1 \text{ .}
$$
Therefore we have
\begin{eqnarray*}
\tilde{L}_1 & \sim_{\mathbb{Q}} & \pi^*(L_1) - E_1 - E_2 - E_3 - E_4 - E_5 - E_6\\
\tilde{L}_4 & \sim_{\mathbb{Q}} & \pi^*(L_4) - \frac{1}{2} E_1 -
E_2 - \frac{3}{2} E_3 - 2 E_4 -
\frac{3}{2} E_5 - E_6 - \frac{1}{2} E_7\\
\tilde{L}_6 & \sim_{\mathbb{Q}} & \pi^*(L_6) - \frac{1}{4} E_1 -
\frac{1}{2} E_2 - \frac{3}{4} E_3 - E_4 - \frac{5}{4} E_5 -
\frac{3}{2} E_6 - \frac{3}{4} E_7 - \frac{1}{2} F_1\text{ .}
\end{eqnarray*}
and since  $L_6 \sim_{\mathbb{Q}} L_1 \sim_{\mathbb{Q}} L_4
\sim_{\mathbb{Q}} -K_X$ we see that $\mathrm{lct}(X) \leq
\frac{1}{2}$.

The inequalities
\begin{eqnarray*}
0 \leq \tilde{D} \cdot \tilde{L}_1 & = & 1 - a_1 - a_7\\
0 \leq \tilde{D} \cdot \tilde{L}_4 & = & 1 - a_4\\
0 \leq \tilde{D} \cdot \tilde{L}_6 & = & 1 - a_6 - b_1\\
0 \leq E_1 \cdot \tilde{D} & = & 2a_1 - a_2\\
0 \leq E_2 \cdot \tilde{D} & = & 2a_2 - a_1 - a_3\\
0 \leq E_3 \cdot \tilde{D} & = & 2a_3 - a_2 - a_4\\
0 \leq E_4 \cdot \tilde{D} & = & 2a_4 - a_3 - a_5\\
0 \leq E_5 \cdot \tilde{D} & = & 2a_5 - a_4 - a_6\\
0 \leq E_6 \cdot \tilde{D} & = & 2a_6 - a_5 - a_7\\
0 \leq E_7 \cdot \tilde{D} & = & 2a_7 - a_6\\
\end{eqnarray*}
imply that $a_1 \leq 1 \text{, } a_2 \leq \frac{3}{2} \text{, }
a_3 \leq \frac{5}{4} \text{, } a_4 \leq 1 \text{, } a_5 \leq
\frac{5}{4} \text{, } a_6 \leq 1 \text{, } a_7 \leq 1 \text{, }
b_1 \leq 1 \text{ .}$

Moreover we get
$$
2a_7 \geq a_6 \text{ , } \frac{3}{2} a_6 \geq a_5 \text{ , }
\frac{4}{3} a_5 \geq a_4 \text{ , } \frac{5}{4} a_4 \geq a_3
\text{ , } \frac{6}{5} a_3 \geq a_2 \text{ , } \frac{7}{6} a_2
\geq a_1 \text{ .}
$$

The equivalence
$$
K_{\tilde{X}} + \lambda \tilde{D} +
 \lambda a_1E_1 + \lambda a_2E_2 +  \lambda a_3E_3 + \lambda a_4E_4 +
\lambda a_5E_5 + \lambda a_6E_6 + \lambda a_7E_7 + \lambda b_1F_1
\sim_{\mathbb{Q}} \pi^*(K_X+ \lambda D)
$$
implies that there is a point $Q\in E_1\cup E_2\cup E_3\cup E_4
\cup E_5 \cup E_6 \cup E_7 \cup F_1$ such that the pair
$$
K_{\tilde{X}} + \lambda \tilde{D} +
 \lambda a_1E_1 + \lambda a_2E_2 +  \lambda a_3E_3 + \lambda a_4E_4 +
\lambda a_5E_5 + \lambda a_6E_6 + \lambda a_7E_7 + \lambda b_1F_1
$$
is not log canonical at $Q$.

\begin{itemize}
\item If the point $Q \in E_1$ and $Q\not \in E_2$ then
$$
K_{\tilde{X}} + \lambda \tilde{D} +  a_1 \lambda E_1
$$
is not log canonical at the point $Q$ and so is the pair
$$
K_{\tilde{X}} + \lambda \tilde{D} + E_1 \text{ , since  } a_1
\lambda \leq 1 \text{ .}
$$
By adjunction $(E_1, \lambda \tilde{D}|_{E_1})$ is not log
canonical at $Q$ and
$$
2a_1 - \frac{6}{7} a_1 \geq 2a_1 - a_2 = \tilde{D} \cdot E_1 \geq
 \text{mult}_Q\Bigl(\tilde{D} \cdot E_1 \Bigr)  > \frac{1}{\lambda} > 2 $$
implies that $a_1 > \frac{7}{4}$ which is false.

\item If $Q\in E_1 \cap E_2$ then the log pair
$$
K_{\tilde{X}} + \lambda \tilde{D} +  \lambda a_1 E_1 +  \lambda
a_2E_2
$$
is not log canonical at the point $Q$ and so is the log pair
$$
K_{\tilde{X}} + \lambda \tilde{D} + E_1 + \lambda a_2 E_2  \text{
.}
$$
By adjunction it follows that
$$
2a_1 - a_2 = \tilde{D} \cdot E_1 \geq
\text{mult}_Q\Bigl(\tilde{D}|_{E_1} \Bigr) =
\text{mult}_Q\Bigl(\tilde{D} \cdot E_1 \Bigr)  > \frac{1}{\lambda}
- a_2 > 2 - a2 \text{ ,}
$$
which is false, since $a_1 \leq 1$.

\item If $Q\in E_2$ but $Q\not \in E_1 \cup E_3$ then
$$
K_{\tilde{X}} + \lambda \tilde{D} +  \lambda a_2 E_2
$$
is not log canonical at the point $Q$ and so is the pair
$$
K_{\tilde{X}} + \lambda \tilde{D} + E_2 \text{ , since  }  \lambda
a_2 \leq 1 \text{ .}
$$
By adjunction $(E_2, \lambda \tilde{D}|_{E_2})$ is not log
canonical at $Q$ and
$$2a_2 - \frac{5}{6} a2 \geq 2a_2 - a_1 - a_3 = \tilde{D} \cdot E_2 \geq
 \text{mult}_Q\Bigl(\tilde{D} \cdot E_2 \Bigr) > \frac{1}{\lambda} > 2  \text{ ,}$$
which is false, since $a_2 \leq \frac{3}{2}$.

\item If $Q\in E_2 \cap E_3$ then the log pair
$$
K_{\tilde{X}} + \lambda \tilde{D} +  \lambda a_2 E_2 +  \lambda
a_3E_3
$$
is not log canonical at the point $Q$ and so is the log pair
$$
K_{\tilde{X}} + \lambda \tilde{D} +  \lambda a_2E_2 + E_3 \text{ ,
since } \lambda a_3 < 1 \text{ .}
$$
By adjunction it follows that
$$
2a_3 - a_2 -a_4= \tilde{D} \cdot E_3 \geq
\text{mult}_Q\Bigl(\tilde{D}|_{E_3} \Bigr) =
\text{mult}_Q\Bigl(\tilde{D} \cdot E_3 \Bigr)  > \frac{1}{\lambda}
- a_2 > 2 - a_2 \text{ ,}
$$
which, along with the inequality $a_4 \geq \frac{4}{5} a_3$,
implies that $a_4 > 1$, which is impossible.

\item If $Q\in E_3$ but $Q\not \in E_2 \cup E_4$ then
$$
K_{\tilde{X}} + \lambda \tilde{D} +  \lambda a_3E_3
$$
is not log canonical at the point $Q$ and so is the pair
$$
K_{\tilde{X}} + \lambda \tilde{D} + E_3 \text{ , since  }  \lambda
a_3\leq 1 \text{ .}
$$
By adjunction $(E_3, \lambda \tilde{D}|_{E_3})$ is not log
canonical at $Q$ and
$$2a_3 - a_2 - a_4 = \tilde{D} \cdot E_3 \geq \text{mult}_Q\Bigl(\tilde{D}|_{E_3}
\Bigr) =
 \text{mult}_Q\Bigl(\tilde{D} \cdot E_3 \Bigr) > \frac{1}{\lambda} > 2  \text{ .}$$
This inequality together with $a_4 \geq \frac{4}{5} a_3$  implies
that $a_4 > 1$, which is impossible.

\item If $Q\in E_3 \cap E_4$ then the log pair
$$
K_{\tilde{X}} + \lambda \tilde{D} +  \lambda a_3 E_3 +  \lambda
a_4 E_4
$$
is not log canonical at the point $Q$ and so is the log pair
$$
K_{\tilde{X}} + \lambda \tilde{D} +  \lambda a_3E_3 + E_4 \text{ ,
since } \lambda a_4 \leq 1 \text{ .}
$$
By adjunction it follows that

$$
2a_4 - a_3 -a_5= \tilde{D} \cdot E_4 \geq
\text{mult}_Q\Bigl(\tilde{D}|_{E_4} \Bigr) =
\text{mult}_Q\Bigl(\tilde{D} \cdot E_4 \Bigr)  > \frac{1}{\lambda}
- a_3 > 2 - a_3 \text{ ,}
$$
which contradicts $a_4 \leq 1$.

\item If $Q\in E_4$ but $Q\not \in E_3 \cup E_5$ then
$$
K_{\tilde{X}} + \lambda \tilde{D} +  \lambda a_4E_4
$$
is not log canonical at the point $Q$ and so is the pair
$$
K_{\tilde{X}} + \lambda \tilde{D} + E_4 \text{ , since  }  \lambda
a_4\leq 1 \text{ .}
$$
By adjunction $(E_4, \lambda \tilde{D}|_{E_4})$ is not log
canonical at $Q$ and
$$2a_4 - a_3 - a_5 = \tilde{D} \cdot E_4 \geq \text{mult}_Q\Bigl(\tilde{D}|_{E_4}
\Bigr) =
 \text{mult}_Q\Bigl(\tilde{D} \cdot E_4 \Bigr) > \frac{1}{\lambda} > 2  \text{ ,}$$
which is false since $a_4 \leq 1$.

\item If the point $Q \in F_1$  then
$$
K_{\tilde{X}} + \lambda \tilde{D} +  b_1 \lambda F_1
$$
is not log canonical at the point $Q$ and so is the pair
$$
K_{\tilde{X}} + \lambda \tilde{D} + F_1 \text{ , since  } b_1
\lambda \leq 1 \text{ .}
$$
By adjunction $(F_1, \lambda \tilde{D}|_{F_1})$ is not log
canonical at $Q$ and
$$
 2b_1 = \tilde{D} \cdot F_1 \geq
 \text{mult}_Q\Bigl(\tilde{D} \cdot F_1 \Bigr)  > \frac{1}{\lambda} > 2 \text{ ,}$$
which is false, since $b_1 \leq 1$.

\end{itemize}
\end{proof}

\subsection{Del Pezzo surface of degree 1 with an $\mathbb{D}_6$
and two $\mathbb{A}_1$ type singularities}

In this section we will prove the following.

\begin{lemma}
\label{D6+A1} Let $X$ be a del Pezzo surface with exactly one Du
Val singularity of type $\mathbb{D}_6$, two of type $\mathbb{A}_1$
and $K_X^2=1$. Then the global log canonical threshold of $X$ is
$$
\mathrm{lct} (X) = \frac{1}{2} \text{ .}
$$
\end{lemma}

\begin{proof}

Suppose that $\mathrm{lct}(X)<\frac{1}{2}$ , then there exists a
$\mathbb{Q}$-divisor $D \in X$ such that the log pair $(X,\lambda
D)$ is not log canonical, where $\lambda < \frac{1}{2}$ and $D
\sim_{\mathbb{Q}} -K_X$. We derive that the pair $(X,\lambda D)$
is log canonical everywhere outside of a point $P\in X$ and not
log canonical at $P$. Let $\pi: \tilde{X} \to X$ be the minimal
resolution of $X$. The configuration of the exceptional curves is
given by the following Dynkin diagram.
\bigskip

$\mathbb{D}_6 + 2 \mathbb{A}_1$.
 
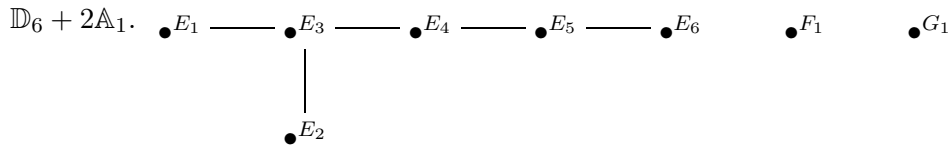
\xymatrix{ {\bullet}^{E_1} \ar@{-}[r] &
{\bullet}^{E_3} \ar@{-}[r] \ar@{-}[d] & {\bullet}^{E_4} \ar@{-}[r]
& {\bullet}^{E_5} \ar@{-}[r] & {\bullet}^{E_6} & {\bullet}^{F_1}  & {\bullet}^{G_1}  \\ & {\bullet}^{E_2}
& & & & &}
\bigskip

\noindent Then
$$
\tilde{D} \sim_{\mathbb{Q}}
\pi^*(D)-a_1E_1-a_2E_2-a_3E_3-a_4E_4-a_5E_5-a_6E_6 - b_1F_1 -
c_1G_1 \text{ .}
$$

\noindent We should note here that there are three -1 curves $L_1,
L_2, L_5$ such that
$$
L_1 \cdot E_1 = L_1 \cdot F_1 = L_2 \cdot E_2 = L_2 \cdot G_1 =
L_5 \cdot E_5 = 1 \text{ .}
$$
Therefore we have
\begin{eqnarray*}
\tilde{L}_1 & \sim_{\mathbb{Q}} & \pi^*(L_1) - \frac{3}{2} E_1 -
E_2 - 2E_3 - \frac{3}{2} E_4 -  E_5 -
\frac{1}{2} E_6 - \frac{1}{2} F_1 \\
\tilde{L}_2 & \sim_{\mathbb{Q}} & \pi^*(L_2) - E_1 - \frac{3}{2}
E_2 - 2E_3 - \frac{3}{2} E_4 - E_5 -
\frac{1}{2} E_6 - \frac{1}{2} G_1\\
\tilde{L}_5 & \sim_{\mathbb{Q}} & \pi^*(L_5) - E_1 - E_2 - 2E_3 -
2E_4 - 2E_5 -E_6 \text{ .}
\end{eqnarray*}
and since  $L_1 \sim_{\mathbb{Q}} L_2 \sim_{\mathbb{Q}} L_5
\sim_{\mathbb{Q}} -K_X$ we see that $\mathrm{lct}(X) \leq
\frac{1}{2}$.

The inequalities
\begin{eqnarray*}
0 \leq \tilde{D} \cdot \tilde{L}_1 & = & 1 - a_1 - b_1\\
0 \leq \tilde{D} \cdot \tilde{L}_2 & = & 1 - a_2 - c_1\\
0 \leq \tilde{D} \cdot \tilde{L}_5 & = & 1 - a_5\\
0 \leq E_1 \cdot \tilde{D} & = & 2a_1 - a_3\\
0 \leq E_2 \cdot \tilde{D} & = & 2a_2 - a_3\\
0 \leq E_3 \cdot \tilde{D} & = & 2a_3 - a_1 - a_2 - a_4\\
0 \leq E_4 \cdot \tilde{D} & = & 2a_4 - a_3 - a_5\\
0 \leq E_5 \cdot \tilde{D} & = & 2a_5 - a_4 - a_6\\
0 \leq E_6 \cdot \tilde{D} & = & 2a_6 - a_5\\
0 \leq F_1 \cdot \tilde{D} & = & 2b_1\\
0 \leq G_1 \cdot \tilde{D} & = & 2c_1
\end{eqnarray*}
imply that $a_1 \leq 1 \text{, } a_2 \leq 1 \text{, } a_3 \leq 2
\text{, } a_4 \leq \frac{3}{2} \text{, } a_5 \leq 1 \text{, } a_6
\leq 1 \text{, } b_1 \leq 1 \text{, } c_1 \leq 1 \text{ .}$
Moreover we have
$$
2a_1 \geq a_3 \text{, } 2a_2 \geq a_3 \text{, } a_3 \geq a_4
\text{, } a_4 \geq a_5 \text{, } a_5 \geq a_6
$$
and
$$
2a_6 \geq a_5 \text{, } \frac{3}{2} a_5 \geq a_4 \text{, }
\frac{4}{3} a_4 \geq a_3 \text{ .}
$$

The equivalence
$$
K_{\tilde{X}} + \lambda \tilde{D} +
 \lambda a_1E_1 +  \lambda a_2E_2 + \lambda a_3E_3 + \lambda a_4E_4 +
\lambda a_5E_5 + \lambda a_6E_6 + \lambda b_1F_1 + \lambda c_1G_1
\sim_{\mathbb{Q}} \pi^*(K_X+ \lambda D)
$$
implies that there is a point $Q$ such that $Q\in E_1\cup E_2\cup
E_3\cup E_4 \cup E_5 \cup E_6$ or $Q \in F_1$ or  $Q\in G_1$ where
the pair
$$
K_{\tilde{X}} + \lambda \tilde{D} +
 \lambda a_1E_1 +  \lambda a_2E_2 + \lambda a_3E_3 + \lambda a_4E_4 +
\lambda a_5E_5 + \lambda a_6E_6 + \lambda b_1F_1 + \lambda c_1G_1
$$
is not log canonical at $Q$.

\begin{itemize}
\item If the point $Q \in E_1 \backslash E_3$ then
$$
K_{\tilde{X}} + \lambda \tilde{D} + \lambda a_1E_1
$$
is not log canonical at the point $Q$ and so is the pair
$$
K_{\tilde{X}} + \lambda \tilde{D} + E_1 \text{ .}
$$
By adjunction $(E_1, \lambda \tilde{D}|_{E_1})$ is not log
canonical at $Q$ and
$$2a_1 \geq 2a_1-a_3 = \tilde{D} \cdot E_1 \geq \text{mult}_Q\Bigl(\tilde{D}|_{E_1}
\Bigr) =
 \text{mult}_Q\Bigl(\tilde{D} \cdot E_1 \Bigr)  > 2 \text{ ,}$$
implies that $a_1 > 1$ which is false.

\item If $Q\in E_3$ but $Q\not \in E_1 \cup E_2 \cup E_4$ then
$$
K_{\tilde{X}} + \lambda \tilde{D} + \lambda a_3E_3
$$
is not log canonical at the point $Q$ and so is the pair
$$
K_{\tilde{X}} + \lambda \tilde{D} + E_3 \text{, since } \lambda
a_3 \leq 1 \text{ .}
$$
By adjunction $(E_3, \lambda \tilde{D}|_{E_3})$ is not log
canonical at $Q$ and
$$2a_3 - \frac{1}{2} a_3 - \frac{1}{2} a_3 \geq 2a_3 - a_1 - a_2 - a_4 = \tilde{D} \cdot
E_3 \geq
 \text{mult}_Q\Bigl(\tilde{D} \cdot E_3 \Bigr) > 2 \text{ ,}$$
implies $a_3 >2$ which is false.

\item If $Q\in E_1 \cap E_3$ then the log pair
$$
K_{\tilde{X}} + \lambda \tilde{D} + \lambda a_1E_1 + \frac{1}{2}
a_3E_3
$$
is not log canonical at the point $Q$ and so are the log pairs
$$
K_{\tilde{X}} + \lambda \tilde{D} + E_1 + \lambda a_3E_3 \text{
and } K_{\tilde{X}} + \lambda \tilde{D} + \lambda a_1E_1 + E_3
\text{ .}
$$
By adjunction it follows that
$$
2a_1-a_3 = \tilde{D} \cdot E_1 \geq
\text{mult}_Q\Bigl(\tilde{D}|_{E_1} \Bigr) =
 \text{mult}_Q\Bigl(\tilde{D} \cdot E_1 \Bigr) > 2 - a_3 \text{ .}
$$
and
$$
2a_3 - a_1 - a_2 - a_4 = \tilde{D} \cdot E_3 \geq
\text{mult}_Q\Bigl(\tilde{D}|_{E_3} \Bigr) =
 \text{mult}_Q\Bigl(\tilde{D} \cdot E_3 \Bigr) > 2 - a_1 \text{ .}
$$
These inequalities imply that $a_1 >1$ and $a_3 >2$ which is a
contradiction.

\item If $Q\in E_3 \cap E_4$ then the log pair
$$
K_{\tilde{X}} + \lambda \tilde{D} + \lambda a_3E_3 + \lambda
a_4E_4
$$
is not log canonical at the point $Q$ and so is the log pair
$$
K_{\tilde{X}} + \lambda \tilde{D} + E_3 +  \lambda a_4E_4 \text{
.}
$$
By adjuction
$$
2a_3 - \frac{1}{2} a_3 - \frac{1}{2} a_3 - a_4 \geq 2a_3 - a_1 -
a_2 - a_4 = \tilde{D} \cdot E_3 \geq
 \text{mult}_Q\Bigl(\tilde{D} \cdot E_3 \Bigr) > 2 - a_4
$$
and this implies that $a_3 > 2$ which is false.

\item If the point $Q \in E_4 \backslash (E_3 \cup E_5)$ then
$$
K_{\tilde{X}} + \lambda \tilde{D} + \lambda a_4E_4
$$
is not log canonical at the point $Q$ and so is the pair
$$
K_{\tilde{X}} + \lambda \tilde{D} + E_4 \text{ .}
$$
By adjunction $(E_4, \lambda \tilde{D}|_{E_4})$ is not log
canonical at $Q$ and
$$2a_4  - a_4 - \frac{2}{3} a_4 \geq 2a_4 - a_3 - a_5 = \tilde{D} \cdot E_4 \geq
 \text{mult}_Q\Bigl(\tilde{D} \cdot E_4 \Bigr)  > 2 \text{ ,}$$
implies that $a_4 > 6$ which is false.

\item $Q\in E_5 \backslash (E_4 \cup E_6)$ then the log pair
$$
K_{\tilde{X}} + \lambda \tilde{D} + \lambda a_5E_5
$$
is not log canonical at the point $Q$ and so is the pair
$$
K_{\tilde{X}} + \lambda \tilde{D} + E_5 \text{ .}
$$
By adjunction $(E_5, \lambda \tilde{D}|_{E_5})$ is not log
canonical at $Q$ and
$$2a_5 - a_5 \ \frac{a_5}{2} \geq 2a_5 - a_4 - a_6 = \tilde{D} \cdot E_5 \geq
 \text{mult}_Q\Bigl(\tilde{D} \cdot E_5 \Bigr) > 2 \text{ ,}$$
implies that $a_5>4$ which is false.

\item If $Q\in E_4 \cap E_5$ then the log pair
$$
K_{\tilde{X}} + \lambda \tilde{D} + \lambda a_4E_4 + \lambda
a_5E_5
$$
is not log canonical at the point $Q$ and so are the log pairs
$$
K_{\tilde{X}} + \lambda \tilde{D} + E_4 + \lambda a_5E_5 \text{
and } K_{\tilde{X}} + \lambda \tilde{D} + \lambda a_4E_4 + E_5
\text{ .}
$$
By adjuction
$$
2a_4 - a_4 - a_5  \geq 2a_4 - a_3 - a_5 = \tilde{D} \cdot E_4 \geq
 \text{mult}_Q\Bigl(\tilde{D} \cdot E_4 \Bigr) > 2 - a_5
$$
and
$$
2a_5  - a_4 - \frac{a_5}{2} \geq 2a_5 - a_4 - a_6 = \tilde{D}
\cdot E_5 \geq
 \text{mult}_Q\Bigl(\tilde{D} \cdot E_5 \Bigr) > 2 - a_4 \text{ .}
$$
imply that $a_4>2$ and $a_5 > \frac{4}{3}$ which both lead to a
contradiction.

\item $Q\in E_5 \backslash (E_4 \cup E_6)$ then the log pair
$$
K_{\tilde{X}} + \lambda \tilde{D} + \lambda a_5E_5
$$
is not log canonical at the point $Q$ and so is the pair
$$
K_{\tilde{X}} + \lambda \tilde{D} + E_5 \text{ .}
$$
By adjunction $(E_5, \lambda \tilde{D}|_{E_5})$ is not log
canonical at $Q$ and
$$2a_5  - a_5 - \frac{a_5}{2} \geq 2a_5 - a_4 - a_6 = \tilde{D} \cdot E_5 \geq
\text{mult}_Q\Bigl(\tilde{D}|_{E_5} \Bigr) =
 \text{mult}_Q\Bigl(\tilde{D} \cdot E_5 \Bigr) > 2 \text{ ,}$$
implies $a_5 >4$ which is false.

\item If $Q\in E_5 \cap E_6$ then the log pair
$$
K_{\tilde{X}} + \lambda \tilde{D} + \lambda a_5E_5 + \lambda
a_6E_6
$$
is not log canonical at the point $Q$ and so is the log pair
$$
K_{\tilde{X}} + \lambda \tilde{D} + \lambda a_5E_5 + E_6 \text{ .}
$$
By adjuction
$$
2a_6 - a_5 = \tilde{D} \cdot E_6 \geq
\text{mult}_Q\Bigl(\tilde{D}|_{E_6} \Bigr) =
 \text{mult}_Q\Bigl(\tilde{D} \cdot E_6 \Bigr) > 2 - a_5 \text{ .}
$$
which is false.

\item If the point $Q \in E_6 \backslash E_5$ then
$$
K_{\tilde{X}} + \lambda \tilde{D} + \lambda a_6E_6
$$
is not log canonical at the point $Q$ and so is the pair
$$
K_{\tilde{X}} + \lambda \tilde{D} + E_6 \text{ .}
$$
By adjunction $(E_6, \lambda \tilde{D}|_{E_6})$ is not log
canonical at $Q$ and
$$2a_6 \geq 2a_6-a_5 = \tilde{D} \cdot E_6 \geq \text{mult}_Q\Bigl(\tilde{D}|_{E_6}
\Bigr) =
 \text{mult}_Q\Bigl(\tilde{D} \cdot E_6 \Bigr)  > 2 \text{ ,}$$
which is false.

\item If the point $Q \in F_1$ then
$$
K_{\tilde{X}} + \lambda \tilde{D} + \lambda b_1F_1
$$
is not log canonical at the point $Q$ and so is the pair
$$
K_{\tilde{X}} + \lambda \tilde{D} + F_1 \text{ .}
$$
By adjunction $(F_1, \lambda \tilde{D}|_{F_1})$ is not log
canonical at $Q$ and
$$2b_1 = \tilde{D} \cdot F_1 \geq \text{mult}_Q\Bigl(\tilde{D}|_{F_1} \Bigr) =
 \text{mult}_Q\Bigl(\tilde{D} \cdot F_1 \Bigr)  > 2 \text{ ,}$$
implies that $b_1 > 1$ which is false.

\end{itemize}
\end{proof}

\subsection{Del Pezzo surface of degree 1 with an $\mathbb{D}_5$
and an $\mathbb{A}_3$ type singularity}

In this section we will prove the following.

\begin{lemma}
\label{D5+A3} Let $X$ be a del Pezzo surface with one Du Val
singularity of type $\mathbb{D}_5$, one of type $\mathbb{A}_3$ and
$K_X^2=1$. Then the global log canonical threshold of $X$ is
$$
\mathrm{lct} (X) = \frac{1}{2} \text{ .}
$$
\end{lemma}

\begin{proof}

Suppose that $\mathrm{lct}(X)<\frac{1}{2}$,  then there exists a
$\mathbb{Q}$-divisor $D \in X$ such that $D \sim_{\mathbb{Q}}
-K_X$ and the log pair $(X,\lambda D)$ is not log canonical for
some rational number $\lambda< \frac{1}{2} $ . We derive that the
pair $(X, \lambda D)$ is log canonical everywhere outside of a
point $P\in X$ and not log canonical at $P$. Let $\pi: \tilde{X}
\to X$ be the minimal resolution of $X$. The configuration of the
exceptional curves is given by the following Dynkin diagram.
\bigskip

$\mathbb{D}_5 + \mathbb{A}_3$.
\xymatrix{ {\bullet}^{E_1} \ar@{-}[r] &
{\bullet}^{E_3} \ar@{-}[r] \ar@{-}[d] & {\bullet}^{E_4} \ar@{-}[r]
& {\bullet}^{E_5} & {\bullet}^{F_1} \ar@{-}[r] & {\bullet}^{F_2} \ar@{-}[r] 
& {\bullet}^{F_3} \\ & {\bullet}^{E_2}
& & & & & &}
\bigskip

\noindent Then
$$
\tilde{D} \sim_{\mathbb{Q}}
\pi^*(D)-a_1E_1-a_2E_2-a_3E_3-a_4E_4-a_5E_5 - b_1F_1 - b_2F_2 -
b_3F_3 \text{ .}
$$

\noindent We should note here that there are three -1 curves $L_1,
L_2, L_5$ such that
$$
L_1 \cdot E_1 = L_1 \cdot F_1 = L_2 \cdot E_2 = L_2 \cdot F_3 =
L_4 \cdot E_4 = 1 \text{ .}
$$
Therefore we have
\begin{eqnarray*}
\tilde{L}_1 & \sim_{\mathbb{Q}} & \pi^*(L_1) - \frac{5}{4} E_1 -
\frac{3}{4} E_2 - \frac{3}{2} E_3 - E_4
- \frac{1}{2} E_5 - \frac{3}{4} F_1 - \frac{1}{2} F_2 - \frac{1}{4} F_3 \\
\tilde{L}_2 & \sim_{\mathbb{Q}} & \pi^*(L_2) - \frac{3}{4} E_1 -
\frac{5}{4} E_2 - \frac{3}{2} E_3 - E_4
- \frac{1}{2} E_5 - \frac{1}{4} F_1 - \frac{1}{2} F_2 - \frac{3}{4} F_3 \\
\tilde{L}_4 & \sim_{\mathbb{Q}} & \pi^*(L_4) -  E_1 - E_2 - 2E_3 -
2E_4 - E_5 \text{ .}
\end{eqnarray*}
and since  $L_1 \sim_{\mathbb{Q}} L_2 \sim_{\mathbb{Q}} L_5
\sim_{\mathbb{Q}} -K_X$ we see that $\mathrm{lct}(X) \leq
\frac{1}{2}$.

The inequalities
\begin{eqnarray*}
0 \leq \tilde{D} \cdot \tilde{L}_1 & = & 1 - a_1 - b_1\\
0 \leq \tilde{D} \cdot \tilde{L}_2 & = & 1 - a_2 - b_3\\
0 \leq \tilde{D} \cdot \tilde{L}_5 & = & 1 - a_4\\
0 \leq E_1 \cdot \tilde{D} & = & 2a_1 - a_3\\
0 \leq E_2 \cdot \tilde{D} & = & 2a_2 - a_3\\
0 \leq E_3 \cdot \tilde{D} & = & 2a_3 - a_1 - a_2 - a_4\\
0 \leq E_4 \cdot \tilde{D} & = & 2a_4 - a_3 - a_5\\
0 \leq E_5 \cdot \tilde{D} & = & 2a_5 - a_4 \\
0 \leq F_1 \cdot \tilde{D} & = & 2b-1 - b_2\\
0 \leq F_2 \cdot \tilde{D} & = & 2b_2 - b_1 - b_3\\
0 \leq F_3 \cdot \tilde{D} & = & 2b_3 - b_2
\end{eqnarray*}
imply that $a_1 \leq 1 \text{, } a_2 \leq 1 \text{, } a_3 \leq
\frac{3}{2} \text{, } a_4 \leq 1 \text{, } a_5 \leq 1 \text{, }
b_1 \leq 1 \text{, }  b_2 \leq 2 \text{, } b_3 \leq 1 \text{ .}$
The equivalence
$$
K_{\tilde{X}} + \lambda \tilde{D} + \lambda a_1E_1 +  \lambda
a_2E_2 + \lambda a_3E_3 +  \lambda a_4E_4 + \lambda a_5E_5 +
\lambda b_1F_1 + \lambda b_2F_2 + \lambda b_3F_3 \sim_{\mathbb{Q}}
\pi^*(K_X + \lambda D)
$$
implies that there is a point $Q\in E_1\cup E_2\cup E_3\cup E_4
\cup E_5 \cup F_1 \cup F_2 \cup F_3$ such that the pair
$$
K_{\tilde{X}} + \lambda \tilde{D} + \lambda a_1E_1 +  \lambda
a_2E_2 + \lambda a_3E_3 +  \lambda a_4E_4 + \lambda a_5E_5 +
\lambda b_1F_1 + \lambda b_2F_2 + \lambda b_3F_3
$$
is not log canonical at $Q$.

\begin{itemize}
\item If the point $Q \in E_1 \backslash E_3$ then
$$
K_{\tilde{X}} + \lambda \tilde{D} + \lambda a_1E_1
$$
is not log canonical at the point $Q$ and so is the pair
$$
K_{\tilde{X}} + \lambda \tilde{D} + E_1 \text{ .}
$$
By adjunction $(E_1, \lambda \tilde{D}|_{E_1})$ is not log
canonical at $Q$ and
$$
2a_1 \geq 2a_1 - a_3 = \tilde{D} \cdot E_1 \geq
\text{mult}_Q\Bigl(\tilde{D}|_{E_1} \Bigr) =
 \text{mult}_Q\Bigl(\tilde{D} \cdot E_1 \Bigr)  > 2 \text{ ,}
$$
implies that $a_1>1$ which is false.

\item If $Q\in E_1 \cap E_3$ then the log pair
$$
K_{\tilde{X}} + \lambda \tilde{D} + \lambda a_1E_1 + \frac{1}{2}
a_3E_3
$$
is not log canonical at the point $Q$ and so is the log pair
$$
K_{\tilde{X}} + \lambda \tilde{D} + E_1 + \lambda a_3E_3 \text{ .}
$$
By adjunction it follows that
$$
2a_1-a_3 = \tilde{D} \cdot E_1 \geq
\text{mult}_Q\Bigl(\tilde{D}|_{E_1} \Bigr) =
 \text{mult}_Q\Bigl(\tilde{D} \cdot E_1 \Bigr) > 2 - a_3 \text{ .}
$$
which is not possible.

\item If $Q\in E_3$ but $Q\not \in E_1 \cup E_2 \cup E_4$ then
$$
K_{\tilde{X}} + \lambda \tilde{D} + \lambda a_3E_3
$$
is not log canonical at the point $Q$ and so is the pair
$$
K_{\tilde{X}} + \lambda \tilde{D} + E_3 \text{, since } a_3 \leq
\frac{3}{2} \text{ .}
$$
By adjunction $(E_3, \lambda \tilde{D}|_{E_3})$ is not log
canonical at $Q$ and
$$2a_3 - \frac{a_3}{2} - \frac{a_3}{2} \geq 2a_3 - a_1 - a_2 - a_4 = \tilde{D} \cdot E_3
\geq
 \text{mult}_Q\Bigl(\tilde{D} \cdot E_3 \Bigr) > 2 \text{ ,}$$
implies that $a_3>2$ which is false.

\item If $Q\in E_3 \cap E_4$ then the log pair
$$
K_{\tilde{X}} + \lambda \tilde{D} + \lambda a_3E_3 + \frac{1}{2}
a_4E_4
$$
is not log canonical at the point $Q$ and so is the log pair
$$
K_{\tilde{X}} + \lambda \tilde{D} + E_3 + \lambda a_4E_4 \text{ .}
$$
By adjunction it follows that
$$
2a_3 - \frac{a_3}{2} - \frac{a_3}{2} - a_4 \geq 2a_3 - a_1 - a_2 -
a_4 = \tilde{D} \cdot E_1 \geq
 \text{mult}_Q\Bigl(\tilde{D} \cdot E_1 \Bigr) > 2 - a_4 \text{ .}
$$
implies that $a_3>2$ which is not possible.

\item $Q\in E_4 \backslash (E_3 \cap E_5)$ then the log pair
$$
K_{\tilde{X}} + \lambda \tilde{D} + \lambda a_4E_4
$$
is not log canonical at the point $Q$ and so is the pair
$$
K_{\tilde{X}} + \lambda \tilde{D} + E_4 \text{ .}
$$
By adjunction $(E_4, \lambda \tilde{D}|_{E_4})$ is not log
canonical at $Q$ and
$$2a_4 \geq 2a_4 - a_3 - a_5 = \tilde{D} \cdot E_4 \geq
\text{mult}_Q\Bigl(\tilde{D}|_{E_4} \Bigr) =
 \text{mult}_Q\Bigl(\tilde{D} \cdot E_4 \Bigr) > 2 \text{ ,}$$
implies that $a_4 > 1$ which is false.

\item $Q\in E_4 \cap E_5$ then the log pair
$$
K_{\tilde{X}} + \lambda \tilde{D} +\lambda a_4E_4 + \lambda a_5E_5
$$
is not log canonical at the point $Q$ and so is the log pair
$$
K_{\tilde{X}} + \lambda \tilde{D} + E_5 + \lambda a_4E_4 \text{ .}
$$
By adjunction it follows that
$$
2a_5 - a_4 = \tilde{D} \cdot E_5 \geq
\text{mult}_Q\Bigl(\tilde{D}|_{E_5} \Bigr) =
 \text{mult}_Q\Bigl(\tilde{D} \cdot E_5 \Bigr) > 2 - a_4 \text{ .}
$$
and we see then that $a_5 > 1$ which is not possible.

\item $Q\in E_5 \backslash E_4$ then the log pair
$$
K_{\tilde{X}} + \lambda \tilde{D} + \lambda a_5E_5
$$
is not log canonical at the point $Q$ and so is the pair
$$
K_{\tilde{X}} + \lambda \tilde{D} + E_5 \text{ .}
$$
By adjunction $(E_5, \lambda \tilde{D}|_{E_5})$ is not log
canonical at $Q$ and
$$2a_5 \geq 2a_5 - a_4 = \tilde{D} \cdot E_5 \geq \text{mult}_Q\Bigl(\tilde{D}|_{E_5}
\Bigr) =
 \text{mult}_Q\Bigl(\tilde{D} \cdot E_5 \Bigr) > 2 \text{ ,}$$
implies that $a_5 > 1$ which is false.

\item If the point $Q \in F_1 \backslash F_3$ then
$$
K_{\tilde{X}} + \lambda \tilde{D} + \lambda b_1F_1
$$
is not log canonical at the point $Q$ and so is the pair
$$
K_{\tilde{X}} + \lambda \tilde{D} + F_1 \text{ .}
$$
By adjunction $(F_1, \lambda \tilde{D}|_{F_1})$ is not log
canonical at $Q$ and
$$
2b_1 \geq 2b_1 - b_2 = \tilde{D} \cdot F_1 \geq
\text{mult}_Q\Bigl(\tilde{D}|_{F_1} \Bigr) =
 \text{mult}_Q\Bigl(\tilde{D} \cdot F_1 \Bigr)  > 2 \text{ ,}
$$
implies that $b_1>1$ which is false.

\item If the point $Q \in F_1 \cap F_2$ then
$$
K_{\tilde{X}} + \lambda \tilde{D} + \lambda b_1F_1 + \lambda
b_2F_2
$$
is not log canonical at the point $Q$ and so is the pair
$$
K_{\tilde{X}} + \lambda \tilde{D} + F_1 + \lambda b_2F_2 \text{ .}
$$
By adjunction $(F_1, \lambda \tilde{D}|_{F_1})$ is not log
canonical at $Q$ and
$$
2b_1 - b_2 = \tilde{D} \cdot F_1 \geq
\text{mult}_Q\Bigl(\tilde{D}|_{F_1} \Bigr) =
 \text{mult}_Q\Bigl(\tilde{D} \cdot F_1 \Bigr)  > 2 - b_2 \text{ ,}
$$
implies that $b_1>1$ which is false.

\item If the point $Q \in F_2 \backslash F_1 \backslash F_3$ then
$$
K_{\tilde{X}} + \lambda \tilde{D} + \lambda b_2F_2
$$
is not log canonical at the point $Q$ and so is the pair
$$
K_{\tilde{X}} + \lambda \tilde{D} + F_2 \text{ .}
$$
By adjunction $(F_2, \lambda \tilde{D}|_{F_2})$ is not log
canonical at $Q$ and
$$
2b_2 - \frac{b_2}{2} - \frac{b_2}{2} \geq 2b_2 - b_1 - b_3 =
\tilde{D} \cdot F_2 \geq \text{mult}_Q\Bigl(\tilde{D}|_{F_2}
\Bigr) =
 \text{mult}_Q\Bigl(\tilde{D} \cdot F_2 \Bigr)  > 2 \text{ ,}
$$
implies that $b_2 > 2$ which is false.

\end{itemize}
\end{proof}

\subsection{Del Pezzo surface of degree 1 with two $\mathbb{D}_4$
type singularities}

In this section we will prove the following.

\begin{lemma}
\label{2D4} Let $X$ be a del Pezzo surface with two Du Val
singularity of type $\mathbb{D}_4$ and $K_X^2=1$. Then the global
log canonical threshold of $X$ is
$$
\mathrm{lct} (X) = \frac{1}{2} \text{ .}
$$
\end{lemma}

\begin{proof}

Suppose that  $\mathrm{lct}(X)<\frac{1}{2}$  then there exists a
$\mathbb{Q}$-divisor $D \in X$, such that the log pair $(X,\lambda
D)$ is  not log canonical, where $\lambda < \frac{1}{2}$ and $D
\sim_{\mathbb{Q}} -K_X$. We derive that the pair $(X, \lambda D)$
is log canonical outside of a point $P\in X$ and not log canonical
at $P$. Let $\pi: \tilde{X} \to X$ be the minimal resolution of
$X$. The configuration of the exceptional curves is given by the
following Dynkin diagram.
\bigskip

$\mathbb{D}_4 + \mathbb{D}_4$.
\xymatrix{ {\bullet}^{E_1} \ar@{-}[r] &
{\bullet}^{E_3} \ar@{-}[r] \ar@{-}[d] & {\bullet}^{E_4} & {\bullet}^{F_1} \ar@{-}[r] &
{\bullet}^{F_3} \ar@{-}[r] \ar@{-}[d] & {\bullet}^{F_4} \\ &
{\bullet}^{E_2} & & & {\bullet}^{F_2} &}
\bigskip

\noindent Then
$$
\tilde{D} \sim_{\mathbb{Q}} \pi^*(D)-a_1E_1-a_2E_2-a_3E_3-a_4E_4 -
b_1F_1 - b_2F_2 - b_3F_3 -b_4F_4 \text{ .}
$$

\noindent We should note here that there are four -1 curves
$\tilde{L}_1, \tilde{L}_2, \tilde{L}_4, \tilde{L}_3$ such that
$$
\tilde{L}_1 \cdot E_1 = \tilde{L}_1 \cdot F_1 = \tilde{L}_2 \cdot
E_2 = \tilde{L}_2 \cdot F_2 = \tilde{L}_4 \cdot E_4 = \tilde{L}_4
\cdot F_4 = \tilde{L}_3 \cdot E_3 = 1 \text{ .}
$$
Therefore we have
\begin{eqnarray*}
\tilde{L}_1 & \sim_{\mathbb{Q}} & \pi^*(L_1) -  E_1 - \frac{1}{2}
E_2 - E_3 - \frac{1}{2} E_4 - F_1 -
\frac{1}{2} F_2 - F_3 - \frac{1}{2} F_4 \\
\tilde{L}_2 & \sim_{\mathbb{Q}} & \pi^*(L_2) - \frac{1}{2} E_1 -
E_2 -  E_3 - \frac{1}{2}E_4 -
\frac{1}{2} F_1 - F_2 - F_3 - \frac{1}{2} F_4 \\
\tilde{L}_4 & \sim_{\mathbb{Q}} & \pi^*(L_4) - \frac{1}{2} E_1 -
\frac{1}{2} E_2 - E_3 - E_4
- \frac{1}{2} F_1 - \frac{1}{2} F_2 - F_3 - F_4 \\
\tilde{L}_3 & \sim_{\mathbb{Q}} & \pi^*(L_3) -  E_1 - E_2 - 2E_3 -
E_4 \text{ .}
\end{eqnarray*}
and since  $L_1 \sim_{\mathbb{Q}} L_2 \sim_{\mathbb{Q}} L_3
\sim_{\mathbb{Q}} L_4 \sim_{\mathbb{Q}} -K_X$ we see that
$\mathrm{lct}(X) \leq \frac{1}{2}$.

The inequalities
\begin{eqnarray*}
0 \leq \tilde{D} \cdot \tilde{L}_1 & = & 1 - a_1 - b_1\\
0 \leq \tilde{D} \cdot \tilde{L}_2 & = & 1 - a_2 - b_2\\
0 \leq \tilde{D} \cdot \tilde{L}_4 & = & 1 - a_4 - b_4\\
0 \leq \tilde{D} \cdot \tilde{L}_3 & = & 1 - a_3\\
0 \leq E_1 \cdot \tilde{D} & = & 2a_1 - a_3\\
0 \leq E_2 \cdot \tilde{D} & = & 2a_2 - a_3\\
0 \leq E_3 \cdot \tilde{D} & = & 2a_3 - a_1 - a_2 - a_4\\
0 \leq E_4 \cdot \tilde{D} & = & 2a_4 - a_3\\
0 \leq F_1 \cdot \tilde{D} & = & 2b_1 - b_3\\
0 \leq F_2 \cdot \tilde{D} & = & 2b_2 - b_3\\
0 \leq F_3 \cdot \tilde{D} & = & 2b_3 - b_1 - b_2 - b_4\\
0 \leq F_4 \cdot \tilde{D} & = & 2b_4 - b_3
\end{eqnarray*}
imply that $a_1 \leq 1 \text{, } a_2 \leq 1 \text{, } a_3 \leq 1
\text{, } a_4 \leq 1 \text{, }  b_1 \leq 1 \text{, }  b_2 \leq 1
\text{, } b_3 \leq \frac{3}{2} \text{, } b_4 \leq 1 \text{ .}$ The
equivalence
$$
K_{\tilde{X}} + \lambda \tilde{D} +
 \lambda a_1E_1 + \lambda a_2E_2 + \lambda a_3E_3 + \lambda a_4E_4  +
\lambda b_1F_1 + \lambda b_2F_2 + \lambda b_3F_3 + \lambda b_4F_4
\sim_{\mathbb{Q}} \pi^*(K_X + \lambda D)
$$
implies that there is a point $Q$ such that $Q\in E_1\cup E_2\cup
E_3\cup E_4$ or $Q\in F_1 \cup F_2\cup F_3\cup F_4$ where the pair
$$
K_{\tilde{X}} + \lambda \tilde{D} +
 \lambda a_1E_1 + \lambda a_2E_2 + \lambda a_3E_3 + \lambda a_4E_4  +
\lambda b_1F_1 + \lambda b_2F_2 + \lambda b_3F_3 + \lambda b_4F_4
$$
is not log canonical at $Q$.

\begin{itemize}
\item If the point $Q \in E_1$ and $Q\not \in E_3$ then
$$
K_{\tilde{X}} + \lambda \tilde{D} + \lambda a_1E_1
$$
is not log canonical at the point $Q$ and so is the pair
$$
K_{\tilde{X}} + \lambda \tilde{D} + E_1 \text{ .}
$$
By adjunction $(E_1, \lambda \tilde{D}|_{E_1})$ is not log
canonical at $Q$ and
$$
2a_1 - a_3 = \tilde{D} \cdot E_1 \geq
\text{mult}_Q\Bigl(\tilde{D}|_{E_1} \Bigr) =
\text{mult}_Q\Bigl(\tilde{D} \cdot E_1 \Bigr)  > 2 \text{ ,}
$$
implies that $a_1 > 1$ which is false.

\item If $Q\in E_3$ but $Q\not \in E_1 \cup E_2 \cup E_4$ then
$$
K_{\tilde{X}} + \lambda \tilde{D} + \lambda a_3E_3
$$
is not log canonical at the point $Q$ and so is the pair
$$
K_{\tilde{X}} + \lambda \tilde{D} + E_3 \text{, since } \lambda
a_3 \leq 1 \text{ .}
$$
By adjunction $(E_3, \lambda \tilde{D}|_{E_3})$ is not log
canonical at $Q$ and
$$2a_3 - \frac{a_3}{2} - \frac{a_3}{2} \geq 2a_3 - a_1 - a_2 - a_4 = \tilde{D} \cdot E_1
\geq
 \text{mult}_Q\Bigl(\tilde{D} \cdot E_3 \Bigr) > 2 \text{ ,}$$
implies that $a_3 > 2$ which is false.

\item If $Q\in E_1 \cap E_3$ then the log pair
$$
K_{\tilde{X}} + \lambda \tilde{D} + \lambda a_1E_1 + \lambda
a_3E_3
$$
is not log canonical at the point $Q$ and so is the log pair
$$
K_{\tilde{X}} + \lambda \tilde{D} + E_1 + \lambda a_3E_3 \text{ .}
$$
By adjunction it follows that
$$
2a_1 - a_3 = \tilde{D} \cdot E_1 \geq
\text{mult}_Q\Bigl(\tilde{D}|_{E_1} \Bigr) =
 \text{mult}_Q\Bigl(\tilde{D} \cdot E_1 \Bigr) > 2 - a_3 \text{ .}
$$
and we see then that $a_1 > 1$ which is not possible.

\end{itemize}
\end{proof}

\subsection{Del Pezzo surface of degree 1 with an $\mathbb{A}_5$,
an $\mathbb{A}_2$ and an $\mathbb{A}_1$ type singularity}

In this section we will prove the following.

\begin{lemma}
\label{A5+A2+A1} Let $X$ be a del Pezzo surface with an
$\mathbb{A}_5$, an $\mathbb{A}_2$ and an $\mathbb{A}_1$ type
singularity such that $K_X^2=1$. Then the global log canonical
threshold of $X$ is
$$
\mathrm{lct} (X) = \frac{2}{3} \text{ .}
$$
\end{lemma}

\begin{proof}

Suppose that $\mathrm{lct}(X)<\frac{2}{3}$ , then there exists a
$\mathbb{Q}$-divisor $D \in X$ such that the log pair
$(X,\frac{2}{3} D)$ is  not log canonical, where $\lambda <
\frac{2}{3}$  and $D \sim_{\mathbb{Q}} -K_X$. We derive that the
pair $(X,\lambda D)$ is log canonical outside of a point $P\in X$
and not log canonical at $P$. Let $\pi: \tilde{X} \to X$ be the
minimal resolution of $X$. The configuration of the exceptional
curves is given by the following Dynkin diagram.
\bigskip

$\mathbb{A}_5 + \mathbb{A}_2 + \mathbb{A}_1$.
\xymatrix{ {\bullet}^{E_1} \ar@{-}[r] &
{\bullet}^{E_2} \ar@{-}[r] & {\bullet}^{E_3} \ar@{-}[r] &
{\bullet}^{E_4} \ar@{-}[r] & {\bullet}^{E_5}  & {\bullet}^{F_1} \ar@{-}[r] &
{\bullet}^{F_2}  & {\bullet}^{G_1} }
\bigskip

\noindent Then
$$
\tilde{D} \sim_{\mathbb{Q}}
\pi^*(D)-a_1E_1-a_2E_2-a_3E_3-a_4E_4-a_5E_5 - b_1F_1 - b_2 F_2 -
c_1 G_1 \text{ .}
$$

\noindent We should note here that there are three -1 curves
$\tilde{L}_1, \tilde{L}_2, \tilde{L}_3$ such that
$$
L_1 \cdot E_1 = L_1 \cdot E_5 = L_2 \cdot E_2 = L_2 \cdot F_1 =
L_3 \cdot E_3= L_3 \cdot G_1 = 1 \text{ .}
$$
Therefore we have
\begin{eqnarray*}
\tilde{L}_1 & \sim_{\mathbb{Q}} & \pi^*(L_1) - E_1 - E_2 - E_3 - E_4 - E_5\\
\tilde{L}_2 & \sim_{\mathbb{Q}} & \pi^*(L_2) - \frac{2}{3} E_1 -
\frac{4}{3} E_2 - E_3 - \frac{2}{3} E_4
- \frac{1}{3} E_5 - \frac{2}{3} F_1 - \frac{1}{3} F_2 \\
\tilde{L}_3 & \sim_{\mathbb{Q}} & \pi^*(L_3) - \frac{1}{2} E_1 -
E_2 - \frac{3}{2} E_3 - E_4 - \frac{1}{2} E_5 - \frac{1}{2}
G_1\text{ .}
\end{eqnarray*}
and since  $L_1 \sim_{\mathbb{Q}} L_2 \sim_{\mathbb{Q}} L_3
\sim_{\mathbb{Q}} -K_X$ we see that $\mathrm{lct}(X) \leq
\frac{2}{3}$. Moreover there are two -1 curves $\tilde{L}_4$ and
$\tilde{L}_5$ that intersect the fundamental cycle as following
$$
\tilde{L}_4 \cdot F_1 = \tilde{L}_4 \cdot F_2 = 1 \text{ and }
\tilde{L}_4 \cdot G_1 = \tilde{L}_4 \cdot E_i = 0 \text{ for all }
i = 1,...,5
$$
and
$$
\tilde{L}_5 \cdot G_1 = 2 \text{ and } \tilde{L}_5 \cdot E_i =
\tilde{L}_5 \cdot F_j = 0 \text{  for all  } i=1 ,..., 5 \text{
and } j = 1, 2 \text{ .}
$$
The inequalities
\begin{eqnarray*}
0 \leq \tilde{D} \cdot \tilde{L}_1 & = & 1 - a_1 - a_5\\
0 \leq \tilde{D} \cdot \tilde{L}_2 & = & 1 - a_2 - b_1\\
0 \leq \tilde{D} \cdot \tilde{L}_3 & = & 1 - a_3 - c_1\\
0 \leq \tilde{D} \cdot \tilde{L}_4 & = & 1 - b_1 - b_2\\
0 \leq \tilde{D} \cdot \tilde{L}_5 & = & 1 - 2 c_1\\
0 \leq E_1 \cdot \tilde{D} & = & 2a_1 - a_2\\
0 \leq E_2 \cdot \tilde{D} & = & 2a_2 - a_1 - a_3\\
0 \leq E_3 \cdot \tilde{D} & = & 2a_3 - a_2 - a_4\\
0 \leq E_4 \cdot \tilde{D} & = & 2a_4 - a_3 - a_5\\
0 \leq E_5 \cdot \tilde{D} & = & 2a_5 - a_4\\
0 \leq F_1 \cdot \tilde{D} & = & 2b_1 - b_2\\
0 \leq F_2 \cdot \tilde{D} & = & 2b_2 - b_1\\
0 \leq G_1 \cdot \tilde{D} & = & 2c_1
\end{eqnarray*}
imply that
$$a_1 \leq \frac{5}{6} \text{, } a_2 \leq 1 \text{, } a_3 \leq 1 \text{, } a_4 \leq
\frac{4}{3} \text{, } a_5 \leq \frac{5}{6} \text{, } b_1 \leq
\frac{2}{3}  \text{, } b_2 \leq \frac{2}{3} \text{, } c_1 \leq
\frac{1}{2} \text{ , }$$ and what is more
 $$
 2a_5 \geq a_4 \text{ , } \frac{3}{2} a_4 \geq a_3 \text{ , } \frac{4}{3} a_3 \geq a_2
\text{ , } \frac{5}{4} a_2 \geq a_1 \text{ .}
 $$

The equivalence
$$
K_{\tilde{X}} + \lambda \tilde{D} +
 \lambda a_1E_1 + \lambda a_2E_2 + \lambda a_3E_3 + \lambda a_4E_4 +
\lambda a_5E_5 + \lambda b_1F_1 + \lambda b_2F_2 + \lambda c_1G_1
\sim_{\mathbb{Q}} \pi^*(K_X + \lambda D)
$$
implies that there is a point $Q\in E_1\cup E_2\cup E_3\cup E_4
\cup E_5  \cup F_1 \cup F_2 \cup G_1$ such that the pair
$$
K_{\tilde{X}} + \lambda \tilde{D} +
 \lambda a_1E_1 + \lambda a_2E_2 + \lambda a_3E_3 + \lambda a_4E_4 +
\lambda a_5E_5 + \lambda b_1F_1 + \lambda b_2F_2 + \lambda c_1G_1
$$
is not log canonical at $Q$.

\begin{itemize}
\item If the point $Q \in E_1$ and $Q\not \in E_2$ then
$$
K_{\tilde{X}} + \lambda \tilde{D} +  a_1 \lambda E_1
$$
is not log canonical at the point $Q$ and so is the pair
$$
K_{\tilde{X}} + \lambda \tilde{D} + E_1 \text{ , since  } a_1
\lambda \leq 1 \text{ .}
$$
By adjunction $(E_1, \lambda \tilde{D}|_{E_1})$ is not log
canonical at $Q$ and
$$2a_1 - \frac{4}{5} a_1 \geq 2a_1 - a_2 = \tilde{D} \cdot E_1 \geq
 \text{mult}_Q\Bigl(\tilde{D} \cdot E_1 \Bigr)  > \frac{1}{\lambda} > \frac{3}{2} \text{
,}$$
 implies that $a_1 > \frac{5}{4}$ which is a contradiction.

 \item If $Q\in E_1 \cap E_2$ then the log pair
$$
K_{\tilde{X}} + \lambda \tilde{D} +  a_1 \lambda E_1 + a_2 \lambda
E_2
$$
is not log canonical at the point $Q$ and so are the log pairs
$$
K_{\tilde{X}} + \lambda \tilde{D} + E_1 + a_2 \lambda E_2 \text{
and  } K_{\tilde{X}} + \lambda \tilde{D} + a_1 \lambda E_1 + E_2
\text{ .}
$$
By adjunction it follows that
$$
2a_2 - a_1 - a_3 = \tilde{D} \cdot E_2 \geq
\text{mult}_Q\Bigl(\tilde{D}|_{E_2} \Bigr) =
 \text{mult}_Q\Bigl(\tilde{D} \cdot E_2 \Bigr) > \frac{1}{\lambda} - a_1 > \frac{3}{2} -
a_1 \text{ ,}
$$
and
$$
2a_1 - a_2 = \tilde{D} \cdot E_1 \geq
\text{mult}_Q\Bigl(\tilde{D}|_{E_1} \Bigr) =
\text{mult}_Q\Bigl(\tilde{D} \cdot E_1 \Bigr)  > \frac{1}{\lambda}
- a_2 > \frac{3}{2} - a_2 \text{ .}
$$
From the first inequality we get $a_3 \geq \frac{9}{10}$ and then
we see that
$$
1 \geq a_1 + a_5 \geq a_1 + \frac{1}{2} \cdot \frac{2}{3} a_3 >
\frac{3}{4} + \frac{3}{10} > 1 \text{ , }
$$
which is a contradiction.

\item If $Q\in E_2$ but $Q\not \in E_1 \cup E_3$ then
$$
K_{\tilde{X}} + \lambda \tilde{D} + a_2 \lambda E_2
$$
is not log canonical at the point $Q$ and so is the pair
$$
K_{\tilde{X}} + \lambda \tilde{D} + E_2 \text{ , since  } a_2
\lambda \leq 1 \text{ .}
$$
By adjunction $(E_2, \lambda \tilde{D}|_{E_2})$ is not log
canonical at $Q$ and
$$2a_2 - \frac{3}{4} a_2 - \frac{a_2}{2} \geq 2a_2 - a_1 - a_3 = \tilde{D} \cdot E_2 \geq
\text{mult}_Q\Bigl(\tilde{D}|_{E_2} \Bigr) =
 \text{mult}_Q\Bigl(\tilde{D} \cdot E_2 \Bigr) > \frac{1}{\lambda} > \frac{3}{2}  \text{
.}$$ This implies that $a_2 > 2$ which is a contradiction.

\item If $Q\in E_2 \cap E_3$ then the log pair
$$
K_{\tilde{X}} + \lambda \tilde{D} +  a_2 \lambda E_2 + a_3 \lambda
E_3
$$
is not log canonical at the point $Q$ and so are the log pairs
$$
K_{\tilde{X}} + \lambda \tilde{D} + a_2 \lambda E_2 + E_3 \text{ ,
since } \lambda a_3 \leq 1 \text{ . }
$$
By adjunction it follows that
$$
2a_3 - a_2 -a_4= \tilde{D} \cdot E_3 \geq
\text{mult}_Q\Bigl(\tilde{D} \cdot E_3 \Bigr)  > \frac{1}{\lambda}
- a_2 > \frac{3}{2} - a2 \text{ .}
$$
which, together with the inequality $a_4 \geq \frac{2}{3} a_3$,
implies that $a_3 > \frac{9}{8}$. However, this is impossible
since $a_3 \leq 1$.

\item If $Q\in E_3$ but $Q\not \in E_2 \cup E_4$ then
$$
K_{\tilde{X}} + \lambda \tilde{D} + a_3 \lambda E_3
$$
is not log canonical at the point $Q$ and so is the pair
$$
K_{\tilde{X}} + \lambda \tilde{D} + E_3 \text{ , since  } a_3
\lambda \leq 1 \text{ .}
$$
By adjunction $(E_3, \lambda \tilde{D}|_{E_3})$ is not log
canonical at $Q$ and
$$2a_3 - a_2 - a_4 = \tilde{D} \cdot E_3 \geq \text{mult}_Q\Bigl(\tilde{D}|_{E_3}
\Bigr) =
 \text{mult}_Q\Bigl(\tilde{D} \cdot E_3 \Bigr) > \frac{1}{\lambda} > \frac{3}{2}  \text{
,}$$ which is false as we saw in the previous case.

\item If the point $Q \in F_1$ and $Q\not \in F_2$ then
$$
K_{\tilde{X}} + \lambda \tilde{D} +  b_1 \lambda F_1
$$
is not log canonical at the point $Q$ and so is the pair
$$
K_{\tilde{X}} + \lambda \tilde{D} + F_1 \text{ , since  } b_1
\lambda \leq 1 \text{ .}
$$
By adjunction $(F_1, \lambda \tilde{D}|_{F_1})$ is not log
canonical at $Q$ and
$$2b_1 - \frac{b_1}{2} \geq 2b_1 - b_2 = \tilde{D} \cdot F_1 \geq
 \text{mult}_Q\Bigl(\tilde{D} \cdot F_1 \Bigr)  > \frac{1}{\lambda} > \frac{3}{2} \text{
,}$$
 implies that $b_1 > 1$ which is a contradiction.

\item If $Q\in F_1 \cap F_2$ then the log pair
$$
K_{\tilde{X}} + \lambda \tilde{D} +  b_1 \lambda F_1 + b_2 \lambda
F_2
$$
is not log canonical at the point $Q$ and so are the log pairs
$$
K_{\tilde{X}} + \lambda \tilde{D} + b_1 \lambda F_1 +  F_2
$$
and
$$
K_{\tilde{X}} + \lambda \tilde{D} + F_1 + b_2 \lambda F_2  \text{
.}
$$
By adjunction it follows that
$$
2b_1 - b_2 = \tilde{D} \cdot F_1 \geq
\text{mult}_Q\Bigl(\tilde{D}|_{F_1} \Bigr) =
\text{mult}_Q\Bigl(\tilde{D} \cdot F_1 \Bigr)  > \frac{1}{\lambda}
- b_2 > \frac{3}{2} - b_2
$$
and
$$
2b_2 - b_1 = \tilde{D} \cdot F_2 \geq
\text{mult}_Q\Bigl(\tilde{D}|_{F_2} \Bigr) =
\text{mult}_Q\Bigl(\tilde{D} \cdot F_2 \Bigr)  > \frac{1}{\lambda}
- b_1 > \frac{3}{2} - b_1 \text{ .}
$$
This implies that
$$
1 \geq b_1 + b_2 > \frac{3}{4} + \frac{3}{4}
$$
which is a contradiction.

\item If the point $Q \in G_1$ then
$$
K_{\tilde{X}} + \lambda \tilde{D} +  c_1 \lambda G_1
$$
is not log canonical at the point $Q$ and so is the pair
$$
K_{\tilde{X}} + \lambda \tilde{D} + G_1 \text{ , since  } c_1
\lambda \leq 1 \text{ .}
$$
By adjunction $(G_1, \lambda \tilde{D}|_{G_1})$ is not log
canonical at $Q$ and
$$ 1 \geq 2c_1  = \tilde{D} \cdot G_1 \geq
 \text{mult}_Q\Bigl(\tilde{D} \cdot G_1 \Bigr)  > \frac{1}{\lambda} > \frac{3}{2} \text{
,}$$
 which is a contradiction.

\end{itemize}
\end{proof}

\subsection{Del Pezzo surface of degree 2 with an $\mathbb{E}_7$
type singularity}

In this section we will prove the following.

\begin{lemma}
\label{degree2E7} Let $X$ be a del Pezzo surface with  one Du Val
singularity of type $\mathbb{E}_7$ and $K_X^2=2$. Then the global
log canonical threshold of $X$ is
$$
\mathrm{lct} (X) = \frac{1}{6} \text{ .}
$$
\end{lemma}

\begin{proof}

Suppose that $\mathrm{lct}(X)<\frac{1}{6}$,  then there exists an
effective $\mathbb{Q}$-divisor $D \sim_{\mathbb{Q}} -K_X$ such
that the log pair $(X,\lambda D)$ is  not log canonical, where
$\lambda < \frac{1}{6}$. We derive that the pair $(X,\lambda D)$
is log canonical everywhere except for a point $P\in X$ at which
it is not log canonical. Let $\pi: \tilde{X} \to X$ be the minimal
resolution of $X$. The configuration of the exceptional curves is
given by the following Dynkin diagram.
\bigskip

$\mathbb{E}_7$.
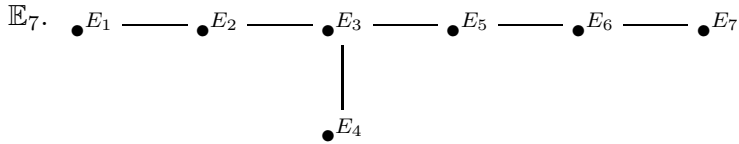
\xymatrix{ {\bullet}^{E_1} \ar@{-}[r] & {\bullet}^{E_2}
\ar@{-}[r]  & {\bullet}^{E_3} \ar@{-}[r] \ar@{-}[d] &
{\bullet}^{E_5} \ar@{-}[r] & {\bullet}^{E_6} \ar@{-}[r] & {\bullet}^{E_7}\\
 & & {\bullet}^{E_4} & & }
\bigskip

\noindent Then
$$
\tilde{D} \sim_{\mathbb{Q}} \pi^*(D) - a_1E_1 - a_2E_2 - a_3E_3 -
a_4E_4 - a_5E_5 - a_6E_6 - a_7E_7 \text{ .}
$$
By the way we obtain $\tilde{X}$ as the blow up of $\mathbb{P}^2$
at seven points we can see that there is a $-1$ curve $\tilde{L}$
that intersects the exceptional divisor $E_7$. In fact we have
$$
\tilde{L} \sim_{\mathbb{Q}} \pi^*(L) - E_1 - 2 E_2 - 3E_3 -
\frac{3}{2} E_4 - \frac{5}{2}E_5 - 2E_6 - \frac{3}{2} E_7 \text{
,}
$$
and since $2L \in |-K_{X}|$ we get that $\mathrm{lct} (X) \leq
\frac{1}{6}$.

The inequalities
\begin{eqnarray*}
0 \leq \tilde{D} \cdot \tilde{L} & = & 1 - a_7\\
0 \leq E_1 \cdot \tilde{D} & = & 2a_1 - 2a_2\\
0 \leq E_2 \cdot \tilde{D} & = & 2a_2 - a_1 - a_3\\
0 \leq E_3 \cdot \tilde{D} & = & 2a_3 - a_2 - a_5 - a_4\\
0 \leq E_4 \cdot \tilde{D} & = & 2a_4 - a_3\\
0 \leq E_5 \cdot \tilde{D} & = & 2a_5 - a_3 - a_6\\
0 \leq E_6 \cdot \tilde{D} & = & 2a_6 - a_5 - a_7\\
0 \leq E_7 \cdot \tilde{D} & = & 2a_7 - a_6\\
\end{eqnarray*}
imply that $a_1 \leq 2 \text{, } a_2 \leq 3 \text{, }
 a_3 \leq 4 \text{, } a_4 \leq \frac{7}{3} \text{, }
 a_5 \leq 3 \text{, } a_6 \leq 2 \text{, }
 a_7 \leq 1 \text{ .}
$ The equivalence
$$
K_{\tilde{X}} + \lambda \tilde{D} + \lambda  a_1E_1 + \lambda
a_2E_2 + \lambda  a_3E_3 + \lambda  a_4E_4 + \lambda  a_5E_5 +
\lambda  a_6E_6 + \lambda  a_7E_7 \sim_{\mathbb{Q}} \pi^*(K_X+
\lambda  D)
$$
implies that there is a point $Q\in E_1\cup E_2\cup E_3\cup E_4
\cup E_5 \cup E_6 \cup E_7$ such that the pair
$$
K_{\tilde{X}} + \lambda  \tilde{D} + \lambda  a_1E_1 + \lambda
a_2E_2 + \lambda  a_3E_3 + \lambda  a_4E_4 + \lambda  a_5E_5 +
\lambda  a_6E_6 + \lambda  a_7E_7
$$
is not log canonical at $Q$.

\begin{itemize}
\item If the point $Q \in E_1$ and $Q\not \in E_2$ then
$$
K_{\tilde{X}} + \lambda \tilde{D} + \lambda a_1E_1
$$
is not log canonical at the point $Q$ and so is the pair
$$
K_{\tilde{X}} + \lambda \tilde{D} + E_1 \text{ since } \lambda a_1
\leq 1 \text{ .}
$$
By adjunction $(E_1, \lambda \tilde{D}|_{E_1})$ is not log
canonical at $Q$ and
$$2a_1-\frac{3}{2} a_1 \geq 2a_1 - a_2 = \tilde{D} \cdot E_1 \geq
\text{mult}_Q\Bigl(\tilde{D} \cdot E_1 \Bigr)  > 6 \text{ ,}$$
which is false.

 \item If the point $Q \in E_1 \cap E_2$ then
$$
K_{\tilde{X}} + \lambda \tilde{D} + \lambda a_1E_1 + \lambda
a_2E_2
$$
is not log canonical at the point $Q$ and so are the pairs
$$
K_{\tilde{X}} + \lambda \tilde{D} + E_1 + \lambda a_2E_2 \text{
and } K_{\tilde{X}} + \lambda \tilde{D} + \lambda a_1E_1 + E_2
\text{ .}
$$
By adjunction
$$2a_1 - a_2 = \tilde{D} \cdot E_1 \geq \text{mult}_Q\Bigl(\tilde{D}|_{E_1} \Bigr) =
 \text{mult}_Q\Bigl(\tilde{D} \cdot E_1 \Bigr)  > 6 - a_2 \text{ and }$$
and
$$2a_2 - a_1 - a_3 = \tilde{D} \cdot E_2 \geq \text{mult}_Q\Bigl(\tilde{D}|_{E_2} \Bigr)
=
 \text{mult}_Q\Bigl(\tilde{D} \cdot E_2 \Bigr)  > 6 - a_1 \text{ ,}$$
which is false.

\item If the point $Q \in E_2 \backslash (E_1 \cup E_3)$ then
$$
K_{\tilde{X}} + \lambda \tilde{D} + \lambda a_2E_2
$$
is not log canonical at the point $Q$ and so is the pair
$$
K_{\tilde{X}} + \lambda \tilde{D} + E_2 \text{ .}
$$
By adjunction
$$2 a_2 - a_1 - a_3 = \tilde{D} \cdot E_2 \geq \text{mult}_Q\Bigl(\tilde{D}|_{E_2} \Bigr)
=
 \text{mult}_Q\Bigl(\tilde{D} \cdot E_2 \Bigr)  > 6 \text{ ,}$$
which is false.

\item If the point $Q \in E_2 \cap E_3$ then
$$
K_{\tilde{X}} + \lambda \tilde{D} + \lambda a_2E_2 + \lambda
a_3E_3
$$
is not log canonical at the point $Q$ and so are the pairs
$$
K_{\tilde{X}} + \lambda \tilde{D} + E_2 + \lambda a_3E_3 \text{
and } K_{\tilde{X}} + \lambda \tilde{D} + \lambda a_2E_2 + E_3
\text{ .}
$$
By adjunction
$$2a_2 - a_1 - a_3 = \tilde{D} \cdot E_2 \geq \text{mult}_Q\Bigl(\tilde{D}|_{E_2} \Bigr)
=
 \text{mult}_Q\Bigl(\tilde{D} \cdot E_2 \Bigr)  > 6 - a_3 \text{ and }$$
and
$$2a_3 - a_2 - a_5 - a_4 = \tilde{D} \cdot E_3 \geq \text{mult}_Q\Bigl(\tilde{D}|_{E_3}
\Bigr) =
 \text{mult}_Q\Bigl(\tilde{D} \cdot E_3 \Bigr)  > 6 - a_2 \text{ ,}$$
which is false.

\item If the point $Q \in E_3 \backslash (E_2 \cup E_4 \cup E_5)$
then
$$
K_{\tilde{X}} + \lambda \tilde{D} + \lambda a_3E_3
$$
is not log canonical at the point $Q$ and so is the pair
$$
K_{\tilde{X}} + \lambda \tilde{D} + E_3 \text{ .}
$$
By adjunction $(E_3, \lambda \tilde{D}|_{E_3})$ is not log
canonical at $Q$ and
$$2a_3 - a_2 - a_4 - a_5 = \tilde{D} \cdot E_3 \geq
\text{mult}_Q\Bigl(\tilde{D}|_{E_3} \Bigr) =
 \text{mult}_Q\Bigl(\tilde{D} \cdot E_3 \Bigr)  > 6 \text{ ,}$$ which is false.

\item If the point $Q \in E_3 \cap E_4$ then
$$
K_{\tilde{X}} + \lambda \tilde{D} + \lambda a_3E_3 + \lambda
a_4E_4
$$
is not log canonical at the point $Q$ and so are the pairs
$$
K_{\tilde{X}} + \lambda \tilde{D} + E_3 + \lambda a_4E_4 \text{
and } K_{\tilde{X}} + \lambda \tilde{D} + \lambda a_3E_3 + E_4
\text{ .}
$$
By adjunction
$$2a_3 - a_2 - a_4 - a_5 = \tilde{D} \cdot E_3 \geq
\text{mult}_Q\Bigl(\tilde{D}|_{E_3} \Bigr) =
 \text{mult}_Q\Bigl(\tilde{D} \cdot E_3 \Bigr)  > 6 - a_4 \text{ and }$$
and
$$2a_4 - a_3 = \tilde{D} \cdot E_4 \geq \text{mult}_Q\Bigl(\tilde{D}|_{E_4} \Bigr) =
 \text{mult}_Q\Bigl(\tilde{D} \cdot E_4 \Bigr)  > 6 - a_3 \text{ ,}$$
which is false.

\item If the point $Q \in E_4 \backslash E_3$ then
$$
K_{\tilde{X}} + \lambda \tilde{D} + \lambda a_4E_4
$$
is not log canonical at the point $Q$ and so is the pair
$$
K_{\tilde{X}} + \lambda \tilde{D} + E_4 \text{ .}
$$
By adjunction $(E_4, \lambda \tilde{D}|_{E_4})$ is not log
canonical at $Q$ and
$$2a_4 - a_3 = \tilde{D} \cdot E_4 \geq \text{mult}_Q\Bigl(\tilde{D}|_{E_4} \Bigr) =
 \text{mult}_Q\Bigl(\tilde{D} \cdot E_4 \Bigr)  > 6 \text{ ,}$$ which is false.

 \item If the point $Q \in E_3 \cap E_5$ then
$$
K_{\tilde{X}} + \lambda \tilde{D} + \lambda a_3E_3 + \lambda
a_5E_5
$$
is not log canonical at the point $Q$ and so are the pairs
$$
K_{\tilde{X}} + \lambda \tilde{D} + E_3 + \lambda a_5E_5 \text{
and } K_{\tilde{X}} + \lambda \tilde{D} + \lambda a_3E_3 + E_5
\text{ .}
$$
By adjunction
$$2a_3 - a_2 - a_4 - a_5 = \tilde{D} \cdot E_3 \geq
\text{mult}_Q\Bigl(\tilde{D}|_{E_3} \Bigr) =
 \text{mult}_Q\Bigl(\tilde{D} \cdot E_3 \Bigr)  > 6 - a_5 \text{ and }$$
and
$$2a_5 - a_3 - a_6 = \tilde{D} \cdot E_5 \geq \text{mult}_Q\Bigl(\tilde{D}|_{E_5}
\Bigr) =
 \text{mult}_Q\Bigl(\tilde{D} \cdot E_5 \Bigr)  > 6 - a_3 \text{ ,}$$
which is false.

\item If the point $Q \in E_5 \backslash (E_3 \cup E_6)$ then
$$
K_{\tilde{X}} + \lambda \tilde{D} + \lambda a_5E_5
$$
is not log canonical at the point $Q$ and so is the pair
$$
K_{\tilde{X}} + \lambda \tilde{D} + E_5 \text{ .}
$$
By adjunction $(E_5, \lambda \tilde{D}|_{E_5})$ is not log
canonical at $Q$ and
$$2a_5 - a_3 - a_6= \tilde{D} \cdot E_5 \geq \text{mult}_Q\Bigl(\tilde{D}|_{E_5} \Bigr)
=
 \text{mult}_Q\Bigl(\tilde{D} \cdot E_5 \Bigr)  > 6 \text{ ,}$$ which is false.

 \item If the point $Q \in E_5 \cap E_6$ then
$$
K_{\tilde{X}} + \lambda \tilde{D} + \lambda a_5E_5 + \lambda
a_6E_6
$$
is not log canonical at the point $Q$ and so are the pairs
$$
K_{\tilde{X}} + \lambda \tilde{D} + E_5 + \lambda a_6E_6 \text{
and } K_{\tilde{X}} + \lambda \tilde{D} + \lambda a_5E_5 + E_6
\text{ .}
$$
By adjunction
$$2a_5 - a_3 - a_6 = \tilde{D} \cdot E_5 \geq
\text{mult}_Q\Bigl(\tilde{D}|_{E_5} \Bigr) =
 \text{mult}_Q\Bigl(\tilde{D} \cdot E_5 \Bigr)  > 6 - a_6 \text{ and }$$
and
$$2a_6 - a_5 - a_7 = \tilde{D} \cdot E_6 \geq \text{mult}_Q\Bigl(\tilde{D}|_{E_6} \Bigr)
=
 \text{mult}_Q\Bigl(\tilde{D} \cdot E_6 \Bigr)  > 6 - a_5 \text{ ,}$$
which is false.

\item If the point $Q \in E_6 \backslash (E_5 \cup E_7)$ then
$$
K_{\tilde{X}} + \lambda \tilde{D} + \lambda a_6E_6
$$
is not log canonical at the point $Q$ and so is the pair
$$
K_{\tilde{X}} + \lambda \tilde{D} + E_6 \text{ .}
$$
By adjunction $(E_6,  \tilde{D}|_{E_6})$ is not log canonical at
$Q$ and
$$2a_6 - a_5 - a_7= \tilde{D} \cdot E_6 \geq \text{mult}_Q\Bigl(\tilde{D}|_{E_6} \Bigr) =
 \text{mult}_Q\Bigl(\tilde{D} \cdot E_6 \Bigr)  > 6 \text{ ,}$$ which is false.

 \item If the point $Q \in E_6 \cap E_7$ then
$$
K_{\tilde{X}} + \lambda \tilde{D} + \lambda a_6E_6 + \lambda
a_7E_7
$$
is not log canonical at the point $Q$ and so are the pairs
$$
K_{\tilde{X}} + \lambda \tilde{D} + E_6 + \lambda a_7E_7 \text{
and } K_{\tilde{X}} + \lambda \tilde{D} + \lambda a_6E_6 + E_7
\text{ .}
$$
By adjunction
$$2a_6 - a_5 - a_7 = \tilde{D} \cdot E_6 \geq
\text{mult}_Q\Bigl(\tilde{D}|_{E_6} \Bigr) =
 \text{mult}_Q\Bigl(\tilde{D} \cdot E_6 \Bigr)  > 6 - a_7 \text{ and }$$
and
$$2a_7 - a_6 = \tilde{D} \cdot E_7 \geq \text{mult}_Q\Bigl(\tilde{D}|_{E_7} \Bigr) =
 \text{mult}_Q\Bigl(\tilde{D} \cdot E_7 \Bigr)  > 6 - a_6 \text{ ,}$$
which is false.

\item If the point $Q \in E_7 \backslash E_6$ then
$$
K_{\tilde{X}} + \lambda \tilde{D} + \lambda a_7E_7
$$
is not log canonical at the point $Q$ and so is the pair
$$
K_{\tilde{X}} + \lambda \tilde{D} + E_7 \text{ .}
$$
By adjunction $(E_7, \lambda \tilde{D}|_{E_7})$ is not log
canonical at $Q$ and
$$2a_7 - a_6= \tilde{D} \cdot E_7 \geq \text{mult}_Q\Bigl(\tilde{D}|_{E_7} \Bigr) =
 \text{mult}_Q\Bigl(\tilde{D} \cdot E_7 \Bigr)  > 6 \text{ ,}$$ which is false.

\end{itemize}

\end{proof}

\subsection{Del Pezzo surfaces of degree 2 with one $\mathbb{D}_6$
and one $\mathbb{A}_1$ singularity.}

In this section we will prove the following.

\begin{lemma}
\label{degree2D6+A1} Let $X$ be a del Pezzo surface with  one Du
Val singularity of type $\mathbb{D}_6$, one of type $\mathbb{A}_1$
and $K_X^2=2$. Then the global log canonical threshold of $X$ is
$$
\mathrm{lct} (X) = \frac{1}{4} \text{ .}
$$
\end{lemma}

\begin{proof}

Suppose that $\mathrm{lct}(X)<\frac{1}{4}$,  then there exists a
$\mathbb{Q}$-divisor $D \sim_{\mathbb{Q}} -K_X$, such that the log
pair $(X,\lambda D)$ is  not log canonical for some rational
number $\lambda < \frac{1}{4}$. We derive that the pair
$(X,\lambda D)$ is log canonical outside of a point $P\in X$ and
not log canonical at $P$. Let $\pi: \tilde{X} \to X$ be the
minimal resolution of $X$. The configuration of the exceptional
curves is given by the following Dynkin diagram.
\bigskip

$\mathbb{D}_6 + \mathbb{A}_1$.
\xymatrix{ {\bullet}^{E_1} \ar@{-}[r] &
{\bullet}^{E_3} \ar@{-}[r] \ar@{-}[d] & {\bullet}^{E_4} \ar@{-}[r]
& {\bullet}^{E_5} \ar@{-}[r] & {\bullet}^{E_6}  & {\bullet}^{F_1}\\ & {\bullet}^{E_2}
&}
\bigskip

Then
$$
\tilde{D} \sim_{\mathbb{Q}}
\pi^*(D)-a_1E_1-a_2E_2-2a_3E_3-2a_4E_4-2a_5E_5-a_6E_6 - b_1 F_1
\text{ .}
$$
From the way we blow up $\mathbb{P}^2$ to obtain $\tilde{X}$ we
can see that there exist -1 curves $\tilde{L}_1, \tilde{L}_6$ such
that
$$
\tilde{L}_1 \cdot E_1 = \tilde{L}_6 \cdot E_6 = \tilde{L}_6 \cdot
F_1 = 1
$$
and therefore
$$
\tilde{L_1} \sim_{\mathbb{Q}} \pi^*(L_1)-\frac{3}{2}
E_1-E_2-2E_3-\frac{3}{2} E_4-E_5-\frac{1}{2} E_6
$$
and
$$
\tilde{L_6} \sim_{\mathbb{Q}} \pi^*(L_6)-\frac{1}{2}
E_1-\frac{1}{2} E_2-E_3- E_4-E_5- E_6 - \frac{1}{2} F_1 \text{ .}
$$
Since $2L_1 \sim_{\mathbb{Q}} 2L_6 \sim_{\mathbb{Q}} -K_X$ we get
that $\mathrm{lct}(X) \leq \frac{1}{4}$. From the inequalities
\begin{eqnarray*}
0 \leq \tilde{D} \cdot \tilde{L}_1 & = & 1 - a_1\\
0 \leq \tilde{D} \cdot \tilde{L}_6 & = & 1 - 2a_5 - b_1\\
0 \leq E_1 \cdot \tilde{D} & = & 2a_1 - 2a_3\\
0 \leq E_2 \cdot \tilde{D} & = & 2a_2 - 2a_3\\
0 \leq E_3 \cdot \tilde{D} & = & 4a_3 - a_1 - a_2 - 2a_4\\
0 \leq E_4 \cdot \tilde{D} & = & 4a_4 - 2a_3 - 2a_5\\
0 \leq E_5 \cdot \tilde{D} & = & 4a_5 - 2a_4 - a_6\\
0 \leq E_6 \cdot \tilde{D} & = & 2a_6 - 2a_5\\
0 \leq F_1 \cdot \tilde{D} & = & 2b_1\\
\end{eqnarray*}
we see that
$$
a_3 \leq a_1 \text{, } a_3 \leq a_2 \text{, } a_4 \leq a_3 \text{,
} a_5 \leq a_4 \text{, } a_6 \leq 2a_5
$$
and
$$
a_5 \leq a_6 \text{, } a_4 \leq \frac{3}{2} a_5 \text{, } a_3 \leq
\frac{4}{3} a_4 \text{, } a_1  \leq \frac{3}{2} a_3 \text{, } a_2
\leq \frac{3}{2}a_3 \text{, } a_1 + a_2 \leq \frac{5}{2}a_3
$$
In particular we get the following upper bounds
$$a_1 \leq 1 \text{, } a_2 \leq 1 \text{, } a_3 \leq \frac{2}{3}
\text{, } a_4 \leq \frac{3}{4} \text{, } a_5 \leq \frac{1}{2}
\text{, } a_6 \leq 1 \text{, } b_1 \leq 1 \text{ .}$$ The
equivalence
$$
K_{\tilde{X}} + \lambda \tilde{D} + \lambda a_1E_1 + \lambda
a_2E_2 + 2 \lambda a_3E_3 + 2 \lambda a_4E_4 + 2 \lambda a_5E_5 +
\lambda a_6E_6 \sim_{\mathbb{Q}} \pi^*(K_X+ \lambda D)
$$
implies that there is a point $Q\in E_1\cup E_2\cup E_3\cup E_4
\cup E_5 \cup E_6$ such that the pair
$$
K_{\tilde{X}} + \lambda \tilde{D} + \lambda a_1E_1 + \lambda
a_2E_2 + 2 \lambda a_3E_3 + 2 \lambda a_4E_4 + 2 \lambda a_5E_5 +
\lambda a_6E_6
$$
is not log canonical at $Q$.

\begin{itemize}
\item If the point $Q \in E_1 \backslash E_3$ then
$$
K_{\tilde{X}} + \lambda \tilde{D} + \lambda a_1E_1
$$
is not log canonical at the point $Q$ and so is the pair
$$
K_{\tilde{X}} + \lambda \tilde{D} + E_1 \text{, since } \lambda
a_1 \leq 1 \text{ .}
$$
By adjunction $(E_1, \lambda \tilde{D}|_{E_1})$ is not log
canonical at $Q$ and
$$2a_1 - \frac{4}{3} a_1 \geq 2a_1-2a_3 = \tilde{D} \cdot E_1 \geq
\text{mult}_Q\Bigl(\tilde{D}|_{E_1} \Bigr) =
 \text{mult}_Q\Bigl(\tilde{D} \cdot E_1 \Bigr)  > 4 \text{ ,}$$
implies that $a_1 \geq 6$ which is false.

\item If $Q\in E_3$ but $Q\not \in E_1 \cup E_2 \cup E_4$ then
$$
K_{\tilde{X}} + \lambda \tilde{D} + 2 \lambda a_3E_3
$$
is not log canonical at the point $Q$ and so is the pair
$$
K_{\tilde{X}} + \lambda \tilde{D} + E_3 \text{, since } 2 \lambda
a_3 \leq 1 \text{ .}
$$
By adjunction $(E_3, \lambda \tilde{D}|_{E_3})$ is not log
canonical at $Q$ and
$$4a_3-a_3-a_3- \frac{3}{2} a_3 \geq 4a_3 - a_1 - a_2 - 2a_4 = \tilde{D} \cdot E_3 \geq
\text{mult}_Q\Bigl(\tilde{D}|_{E_3} \Bigr) =
 \text{mult}_Q\Bigl(\tilde{D} \cdot E_3 \Bigr) > 4 \text{ ,}$$
implies that $a_3 \geq 8$ which is false.

\item If $Q\in E_1 \cap E_3$ then the log pair
$$
K_{\tilde{X}} + \lambda \tilde{D} + \lambda a_1E_1 + 2 \lambda
a_3E_3
$$
is not log canonical at the point $Q$ and so is the log pair
$$
K_{\tilde{X}} + \lambda \tilde{D} + \lambda a_1E_1 + E_3 \text{ .}
$$
By adjunction it follows that
$$
\frac{16}{3} a_4 - a_4 -2a_4 -a_1 \geq 4a_3 - a_1 - a_2 - 2a_4 =
\tilde{D} \cdot E_3 \geq \text{mult}_Q\Bigl(\tilde{D}|_{E_3}
\Bigr) =
 \text{mult}_Q\Bigl(\tilde{D} \cdot E_3 \Bigr) > 4 - a_1
$$
and this implies that $a_4 > \frac{12}{7}$ which is false.

\item If $Q\in E_3 \cap E_4$ then the log pair
$$
K_{\tilde{X}} + \lambda \tilde{D} + 2 \lambda a_3E_3 + 2 \lambda
a_4E_4
$$
is not log canonical at the point $Q$ and so is the log pair
$$
K_{\tilde{X}} + \lambda \tilde{D} + 2 \lambda a_3E_3 + E_4 \text{,
since } 2 \lambda a_4 \leq 1 \text{ .}
$$
By adjuction
$$
6a_5 -2a_5 -2a_3 \geq 4a_4 - 2a_3 - 2a_5 = \tilde{D} \cdot E_4
\geq \text{mult}_Q\Bigl(\tilde{D}|_{E_4} \Bigr) =
 \text{mult}_Q\Bigl(\tilde{D} \cdot E_4 \Bigr) > 4 - 2a_3 \text{ .}
$$
and this implies that $a_5 > 1$ which is false.

\item If the point $Q \in E_4 \backslash (E_3 \cup E_5)$ then
$$
K_{\tilde{X}} + \lambda \tilde{D} + 2\lambda a_4E_4
$$
is not log canonical at the point $Q$ and so is the pair
$$
K_{\tilde{X}} + \lambda \tilde{D} + E_4 \text{ .}
$$
By adjunction $(E_4, \lambda \tilde{D}|_{E_4})$ is not log
canonical at $Q$ and
$$4a_4 - 2a_4 - \frac{4}{3} a_4 \geq 4a_4 - 2a_3 - 2a_5 = \tilde{D} \cdot E_4 \geq
\text{mult}_Q\Bigl(\tilde{D}|_{E_4} \Bigr) =
 \text{mult}_Q\Bigl(\tilde{D} \cdot E_4 \Bigr)  > 4 \text{ ,}$$
and this implies that $a_4 > 6$ which is false.

\item $Q\in E_5 \backslash (E_4 \cup E_6)$ then the log pair
$$
K_{\tilde{X}} + \lambda \tilde{D} + 2\lambda a_5E_5
$$
is not log canonical at the point $Q$ and so is the pair
$$
K_{\tilde{X}} + \lambda \tilde{D} + E_5 \text{, since } 2 \lambda
a_5 \leq 1 \text{ .}
$$
By adjunction $(E_5, \lambda \tilde{D}|_{E_5})$ is not log
canonical at $Q$ and
$$4a_5-2a_5-a_5 \geq 4a_5 - 2a_4 - a_6 = \tilde{D} \cdot E_5 \geq
\text{mult}_Q\Bigl(\tilde{D}|_{E_5} \Bigr) =
 \text{mult}_Q\Bigl(\tilde{D} \cdot E_5 \Bigr) > 4 \text{ ,}$$
and this implies that $a_5 > 4$ which is false.

\item If $Q\in E_4 \cap E_5$ then the log pair
$$
K_{\tilde{X}} + \lambda \tilde{D} + 2 \lambda a_4E_4 + 2
\frac{1}{4} a_5E_5
$$
is not log canonical at the point $Q$ and so is the log pair
$$
K_{\tilde{X}} + \lambda \tilde{D} + 2\lambda a_4E_4 + E_5 \text{
.}
$$
By adjuction
$$
4a_5-a_5- 2a_4 \geq 4a_5 - 2a_4 - a_6 = \tilde{D} \cdot E_5 \geq
\text{mult}_Q\Bigl(\tilde{D}|_{E_5} \Bigr) =
 \text{mult}_Q\Bigl(\tilde{D} \cdot E_5 \Bigr) > 4 - 2a_4 \text{ .}
$$
implies that $a_5 > \frac{4}{3}$ which is false.

\item If $Q\in E_5 \cap E_6$ then the log pair
$$
K_{\tilde{X}} + \lambda \tilde{D} + 2 \lambda a_5E_5 + \lambda
a_6E_6
$$
is not log canonical at the point $Q$ and so is the log pair
$$
K_{\tilde{X}} + \lambda \tilde{D} + 2 \lambda a_5E_5 + E_6 \text{,
since } \lambda a_6 \leq 1
$$
By adjuction
$$
2a_6 - 2a_5 = \tilde{D} \cdot E_6 \geq
\text{mult}_Q\Bigl(\tilde{D}|_{E_6} \Bigr) =
 \text{mult}_Q\Bigl(\tilde{D} \cdot E_6 \Bigr) > 4 - 2a_5 \text{ .}
$$
implies that $a_6 >2$ which is false.

\item If the point $Q \in E_6 \backslash E_5$ then
$$
K_{\tilde{X}} + \lambda \tilde{D} + \lambda a_6E_6
$$
is not log canonical at the point $Q$ and so is the pair
$$
K_{\tilde{X}} + \lambda \tilde{D} + E_6 \text{ .}
$$
By adjunction $(E_6, \lambda \tilde{D}|_{E_6})$ is not log
canonical at $Q$ and
$$2a_6 - a_6 \geq 2a_6-2a_5 = \tilde{D} \cdot E_6 \geq
\text{mult}_Q\Bigl(\tilde{D}|_{E_6} \Bigr) =
 \text{mult}_Q\Bigl(\tilde{D} \cdot E_6 \Bigr)  > 4 \text{ ,}$$
implies that $a_6 > 4$ which is false.

\item If the point $Q \in F_1$ then
$$
K_{\tilde{X}} + \lambda \tilde{D} + \lambda b_1F_1
$$
is not log canonical at the point $Q$ and so is the pair
$$
K_{\tilde{X}} + \lambda \tilde{D} + F_1 \text{, since } \lambda
b_1 \leq 1 \text{ .}
$$
By adjunction $(F_1, \lambda \tilde{D}|_{F_1})$ is not log
canonical at $Q$ and
$$ 2b_1 = \tilde{D} \cdot F_1 \geq \text{mult}_Q\Bigl(\tilde{D}|_{F_1} \Bigr) =
 \text{mult}_Q\Bigl(\tilde{D} \cdot F_1 \Bigr)  > 4 \text{ ,}$$
 which is false.
\end{itemize}
\end{proof}

\subsection{Del Pezzo surfaces of degree 2 with one $\mathbb{D}_4$
and three $\mathbb{A}_1$ type singularities}

In this section we will prove the following.

\begin{lemma}
\label{D4+3A1} Let $X$ be a del Pezzo surface with  one Du Val
singularity of type $\mathbb{D}_4$, three of type $\mathbb{A}_1$
and $K_X^2=2$. Then the global log canonical threshold of $X$ is
$$
\mathrm{lct} (X) = \frac{1}{2} \text{ .}
$$
\end{lemma}

\begin{proof}

Let $X$ be a del Pezzo surface with  one Du Val singularity of
type $\mathbb{D}_4$, three $\mathbb{A}_1$ type singularities  and
$K_X^2=2$. Suppose $\mathrm{lct}(X) < \frac{1}{2} $. Then there
exists an effective $\mathbb{Q}$-divisor $D\in X$ such that the
log pair $(X,\lambda D)$ is  not log canonical for some rational
number $\lambda < \frac{1}{2}$ and $D \sim_{\mathbb{Q}} -K_X$.

We derive that the pair $(X, \lambda D)$ is log canonical outside
of a point $P\in X$ and not log canonical at $P$. Let $\pi:
\tilde{X} \to X$ be the minimal resolution of $X$. The following
diagram shows how the exceptional curves intersect each other.
\bigskip

$\mathbb{D}_4$. \xymatrix{ {\bullet}^{E_1} \ar@{-}[r] & {\bullet}^{E_3}
\ar@{-}[r] \ar@{-}[d] & {\bullet}^{E_4} & {\bullet}^{F_1} & {\bullet}^{F_2} &
{\bullet}^{F_4} \\ & {\bullet}^{E_2} &}
\bigskip

Then
$$
\tilde{D} \sim_{\mathbb{Q}}
\pi^*(D)-a_1E_1-a_2E_2-a_3E_3-a_4E_4-b_1F_1 - b_2F_2 - b_4F_4
\text{ .}
$$
From the inequalities
\begin{eqnarray*}
0 \leq \tilde{D} \cdot \tilde{L}_1 & = & 1 - a_1 - b_1\\
0 \leq \tilde{D} \cdot \tilde{L}_2 & = & 1 - a_2 - b_2\\
0 \leq \tilde{D} \cdot \tilde{L}_4 & = & 1 - a_4 - b_4\\
0 \leq E_1 \cdot \tilde{D} & = & 2a_1 - a_3\\
0 \leq E_2 \cdot \tilde{D} & = & 2a_2 - a_3\\
0 \leq E_3 \cdot \tilde{D} & = & 2a_3 - a_1 - a_2 - a_4\\
0 \leq E_4 \cdot \tilde{D} & = & 2a_4 - a_3\\
0 \leq F_1 \cdot \tilde{D} & = & 2b_1\\
0 \leq F_2 \cdot \tilde{D} & = & 2b_2\\
0 \leq F_4 \cdot \tilde{D} & = & 2b_4
\end{eqnarray*}
we see that
$$a_1 \leq 1 \text{, } a_2 \leq 1 \text{, } a_3 \leq 2 \text{, } a_4 \leq 1 \text{, }
b_1 \leq 1 \text{, } b_2 \leq 1 \text{, } b_4 \leq 1 \text{ .}$$
We should note here that there are three -1 curves $ \tilde{L}_1,
\tilde{L}_2, \tilde{L}_4$ such that
$$
\tilde{L}_1 \cdot E_1 = \tilde{L}_1 \cdot F_1 = \tilde{L}_2 \cdot
E_2 = L_2 \cdot F_2 = \tilde{L}_4 \cdot E_4= \tilde{L}_4 \cdot F_4
=1 \text{ .}
$$
Therefore we have
\begin{eqnarray*}
 \tilde{L_1} & \sim_{\mathbb{Q}} & \pi^*(L_1) -  E_1 - \frac{1}{2} E_2 - E_3 - \frac{1}{2} E_4 -
\frac{1}{2} F_1\\
\tilde{L_2} & \sim_{\mathbb{Q}} & \pi^*(L_2) - \frac{1}{2}  E_1 -
E_2 - E_3 - \frac{1}{2} E_4 -
\frac{1}{2} F_2\\
\tilde{L_4} & \sim_{\mathbb{Q}} & \pi^*(L_4) - \frac{1}{2} E_1 -
\frac{1}{2} E_2 - E_3 -  E_4 - \frac{1}{2} F_4 \text{ .}
\end{eqnarray*}

The equivalence
$$
K_{\tilde{X}} + \lambda \tilde{D} + \lambda a_1E_1 + \lambda
a_2E_2 + \lambda a_3E_3 + \lambda a_4E_4 + \lambda b_1F_1+ \lambda
b_2F_2 + \lambda b_4F_4 \sim_{\mathbb{Q}} \pi^*(K_X+ \lambda D)
$$
implies that there is a point $Q\in E_1\cup E_2\cup E_3\cup E_4
\cup F_1 \cup F_2 \cup F_4$ such that the pair
$$
K_{\tilde{X}} + \lambda \tilde{D} + \lambda a_1E_1 + \lambda
a_2E_2 + \lambda a_3E_3 + \lambda a_4E_4 + \lambda b_1F_1+ \lambda
b_2F_2 + \lambda b_4F_4
$$
is not log canonical at $Q$.

\begin{itemize}
\item If the point $Q \in E_1$ and $Q\not \in E_3$ then
$$
K_{\tilde{X}} + \lambda \tilde{D} + \lambda a_1E_1
$$
is not log canonical at the point $Q$ and so is the pair
$$
K_{\tilde{X}} + \lambda \tilde{D} + E_1 \text{ .}
$$
By adjunction $(E_1, \lambda \tilde{D}|_{E_1})$ is not log
canonical at $Q$ and
$$ 2a_1 - a_3 = \tilde{D} \cdot E_1 \geq \text{mult}_Q\Bigl(\tilde{D}|_{E_1} \Bigr) =
 \text{mult}_Q\Bigl(\tilde{D} \cdot E_1 \Bigr)  > 2 \text{ ,}$$
which implies that $a_1 > 1$ which is false.

\item If $Q\in E_3$ but $Q\not \in E_1 \cup E_2 \cup E_4$ then
$$
K_{\tilde{X}} + \lambda \tilde{D} +  \lambda a_3E_3
$$
is not log canonical at the point $Q$ and so is the pair
$$
K_{\tilde{X}} + \lambda \tilde{D} + E_3 \text{, since } \lambda
a_3 \leq 1 \text{ .}
$$
By adjunction $(E_3, \lambda \tilde{D}|_{E_3})$ is not log
canonical at $Q$ and
$$a_3 \geq 2a_3 - a_1 - a_2 - a_4 = \tilde{D} \cdot E_1 \geq
\text{mult}_Q\Bigl(\tilde{D}|_{E_3} \Bigr) =
 \text{mult}_Q\Bigl(\tilde{D} \cdot E_3 \Bigr) > 2 \text{ ,}$$
 which is false.

\item If $Q\in E_1 \cap E_3$ then the log pair
$$
K_{\tilde{X}} + \lambda \tilde{D} + \lambda a_1E_1 + \frac{1}{2}
a_3E_3
$$
is not log canonical at the point $Q$ and so is the log pair
$$
K_{\tilde{X}} + \lambda \tilde{D} + \lambda a_1 E_1 + E_3 \text{
.}
$$
By adjunction it follows that
$$
a_3 - a_1 \geq 2a_3 - a_1 - a_2 - a_4 = \tilde{D} \cdot E_3 \geq
\text{mult}_Q\Bigl(\tilde{D}|_{E_3} \Bigr) =
 \text{mult}_Q\Bigl(\tilde{D} \cdot E_3 \Bigr) > 2 -a_1 \text{ .}
$$
and we see then that $a_3>1$ which is not possible.

\item If $Q\in F_1$ then the log pair
$$
K_{\tilde{X}} + \lambda \tilde{D} + \lambda b_1F_1
$$
is not log canonical at the point $Q$ and so is the log pair
$$
K_{\tilde{X}} + \lambda \tilde{D} + F_1 \text{ .}
$$
By adjunction it follows that
$$
 2b_1 = \tilde{D} \cdot F_1 \geq \text{mult}_Q\Bigl(\tilde{D}|_{F_1} \Bigr) =
 \text{mult}_Q\Bigl(\tilde{D} \cdot F_1 \Bigr) > 2  \text{ .}
$$
and we see then that $b_1>1$ which is not possible.

\end{itemize}

\end{proof}

\subsection{Del Pezzo surfaces of degree 2 with two $\mathbb{A}_3$
and one $\mathbb{A}_1$ type singularity}

In this section we will prove the following.

\begin{lemma}
\label{2A3+A1} Let $X$ be a del Pezzo surface with two Du Val
singularities of type $\mathbb{A}_3$, one $\mathbb{A}_1$ type
singularity and $K_X^2=2$. Then the global log canonical threshold
of $X$ is
$$
\mathrm{lct} (X) = \frac{1}{2} \text{ .}
$$
\end{lemma}

\begin{proof}

Let $X$ be a del Pezzo surface with two Du Val singularities of
type $\mathbb{A}_3$, one $\mathbb{A}_1$ type singularity and
$K_X^2=2$. Suppose that $\mathrm{lct}(X) < \frac{1}{2} $, then
there exists an effective $\mathbb{Q}$-divisor $D\in X$ such that
the log pair $(X,\lambda D)$ is  not log canonical for some
rational number $\lambda < \frac{1}{2}$ and $D \sim_{\mathbb{Q}}
-K_X$.

Let $Z$ be the curve in $|-K_X|$ that contains $P$. Since the
curve $Z$ is irreducible we may assume that the support of $D$
does not contain $Z$.

We derive that the pair $(X, D)$ is log canonical outside of a
point $P\in X$ and not log canonical at $P$. Let $\pi_1: \tilde{X}
\to X$ be the minimal resolution of $X$. The following diagram
shows how the exceptional curves intersect each other.
\bigskip

$\mathbb{A}_3+\mathbb{A}_3+\mathbb{A}_1$ \xymatrix{ {\bullet}^{E_1}
\ar@{-}[r] & {\bullet}^{E_2} \ar@{-}[r] & {\bullet}^{E_3} & {\bullet}^{F_1}
\ar@{-}[r] & {\bullet}^{F_2} \ar@{-}[r] & {\bullet}^{F_3} & {\bullet}^{G_1} }
\bigskip

Then
$$
\tilde{D} \sim_{\mathbb{Q}}
\pi_1^*(D)-a_1E_1-a_2E_2-a_3E_3-b_1F_1-b_2F_2-b_3F_3 -c_1G_1
\text{ .}$$

From the inequalities
\begin{eqnarray*}
0 \leq \tilde{D} \cdot \tilde{L}_1 & = & 1 - a_1 - b_1\\
0 \leq \tilde{D} \cdot \tilde{L}_2 & = & 1 - a_2 - c_1 \\
0 \leq \tilde{D} \cdot \tilde{L}_3 & = & 1 - a_3 - b_3\\
0 \leq E_1 \cdot \tilde{D} & = & 2a_1 - a_2\\
0 \leq E_2 \cdot \tilde{D} & = & 2a_2 - a_1 - a_3\\
0 \leq E_3 \cdot \tilde{D} & = & 2a_3 - a_2\\
0 \leq F_1 \cdot \tilde{D} & = & 2b_1 - b_2\\
0 \leq F_2 \cdot \tilde{D} & = & 2b_2 - b_1 - b_3\\
0 \leq F_3 \cdot \tilde{D} & = & 2b_3 - b_2\\
0 \leq G_1 \cdot \tilde{D} & = & 2c_1
\end{eqnarray*}
we see that $a_1 \leq 1 \text{, } a_2 \leq 1 \text{, } a_3 \leq 1
\text{, } b_1 \leq 1 \text{, } b_2 \leq 2 \text{, } b_3 \leq 1
\text{, } c_1 \leq 1 \text{ .}$

We have three lines $\tilde{L}_1, \tilde{L}_2, \tilde{L}_3$
intersecting the fundamental cycle as following
\begin{eqnarray*}
\tilde{L}_1 \cdot E_1 = \tilde{L}_1 \cdot F_1 = 1 \text{, }\\
\tilde{L}_3 \cdot E_3 = \tilde{L}_3 \cdot F_3 = 1\text{, }\\
\tilde{L}_2 \cdot E_2 = \tilde{L}_2 \cdot G_1 =1 \text{, }
\end{eqnarray*}
and in particular we have
\begin{eqnarray*}
\tilde{L_1} & \sim_{\mathbb{Q}} &
\pi_1^*(L_1)-\frac{3}{4}E_1-\frac{1}{2}E_2-\frac{1}{4}E_3-\frac{3}{4}F_1-\frac{1}{2}F_2-
\frac{1}{4}F_3\\
\tilde{L_2} & \sim_{\mathbb{Q}} & \pi_1^*(L_2)-\frac{1}{2}E_1-E_2-\frac{1}{2}E_3-\frac{1}{2}G_1\\
\tilde{L_3} & \sim_{\mathbb{Q}} &
\pi_1^*(L_3)-\frac{1}{4}E_1-\frac{1}{2}E_2-\frac{3}{4}E_3-\frac{1}{4}F_1-\frac{1}{2}F_2-
\frac{3}{4}F_3 \text{ .}
\end{eqnarray*}

The equivalence
$$
K_{\tilde{X}} + \lambda \tilde{D} + \lambda a_1E_1 + \lambda
a_2E_2 + \lambda a_3E_3 + \lambda \tilde{D} + \lambda b_1F_1 +
\lambda b_2F_2 + \lambda b_3F_3 + \lambda c_1G_1 \sim_{\mathbb{Q}}
\pi_1^*(K_X+ \lambda D)
$$
implies that there is a point $Q\in E_1\cup E_2\cup E_3 \cup
F_1\cup F_2\cup F_3\cup G_1$ such that the pair
$$
K_{\tilde{X}} + \lambda \tilde{D} + \lambda a_1E_1 + \lambda
a_2E_2 + \lambda a_3E_3 + \lambda \tilde{D} + \lambda b_1F_1 +
\lambda b_2F_2 + \lambda b_3F_3 + \lambda c_1G_1
$$ is not log canonical at $Q$.

\begin{itemize}
\item If the point $Q \in E_1$ and $Q\not \in E_2$ then
$$
K_{\tilde{X}} + \lambda \tilde{D} + \lambda a_1E_1
$$
is not log canonical at the point $Q$ and so is the pair
$$
K_{\tilde{X}} + \lambda \tilde{D} + E_1 \text{ .}
$$
By adjunction $(E_1, \lambda \tilde{D}|_{E_1})$ is not log
canonical at $Q$ and
$$\frac{4}{3} a_1 \geq 2a_1 - a_2 = \tilde{D} \cdot E_1 \geq
\text{mult}_Q\Bigl(\tilde{D}|_{E_1} \Bigr) =
 \text{mult}_Q\Bigl(\tilde{D} \cdot E_1 \Bigr)  > 2 \text{ ,}$$
implies that $a_1 > \frac{3}{2}$ which is false.

\item If $Q\in E_2$ but $Q\not \in E_1 \cup E_3$ then
$$
K_{\tilde{X}} + \lambda \tilde{D} + \lambda a_2E_2
$$
is not log canonical at the point $Q$ and so is the pair
$$
K_{\tilde{X}} + \lambda \tilde{D} + E_2 \text{, since } \lambda
a_2 \leq 1 \text{ .}
$$
By adjunction $(E_2, \lambda \tilde{D}|_{E_2})$ is not log
canonical at $Q$ and
$$a_2 \geq 2 a_2 - a_1 - a_3 = \tilde{D} \cdot E_2 \geq
\text{mult}_Q\Bigl(\tilde{D}|_{E_2} \Bigr) =
 \text{mult}_Q\Bigl(\tilde{D} \cdot E_2 \Bigr) > 2 \text{ ,}$$
which is false.

\item If $Q\in E_1 \cap E_2$ then the log pair
$$
K_{\tilde{X}} + \lambda \tilde{D} + \lambda a_1E_1 + \lambda
a_2E_2
$$
is not log canonical at the point $Q$ and so is the log pair
$$
K_{\tilde{X}} + \lambda \tilde{D} + E_1 + \lambda a_2E_2 \text{ .}
$$
By adjunction it follows that
$$
2a_1 - a_2 = \tilde{D} \cdot E_1 \geq
\text{mult}_Q\Bigl(\tilde{D}|_{E_1} \Bigr) =
 \text{mult}_Q\Bigl(\tilde{D} \cdot E_1 \Bigr) > 2 - a_2
$$
and this implies that $a_1>1$ which is false .

\item If $Q\in G_1$ then the log pair
$$
K_{\tilde{X}} + \lambda \tilde{D} + \lambda c_1G_1
$$
is not log canonical at the point $Q$ and so is the log pair
$$
K_{\tilde{X}} + \lambda \tilde{D} + G_1 \text{ .}
$$
By adjunction it follows that
$$
2c_1 = \tilde{D} \cdot G_1 \geq
\text{mult}_Q\Bigl(\tilde{D}|_{G_1} \Bigr) =
 \text{mult}_Q\Bigl(\tilde{D} \cdot G_1 \Bigr) > 2
$$
which is false.

\end{itemize}
\end{proof}

\subsection{Del Pezzo surfaces of degree 2 with  one $\mathbb{A}_5$
and one $\mathbb{A}_2$ type singularity}

In this section we will prove the following.

\begin{lemma}
\label{A5+A2} Let $X$ be a del Pezzo surface with one Du Val
singularity of type $\mathbb{A}_5$, one of type $\mathbb{A}_2$ and
$K_X^2=2$.  Then the global log canonical threshold of $X$ is
$$
\mathrm{lct} (X) = \frac{1}{3} \text{ .}
$$
\end{lemma}

\begin{proof}

Suppose $\mathrm{lct}(X) < \frac{1}{3}$. Then there exists an
effective $\mathbb{Q}$-divisor $D\in X$ such that the log pair
$(X, \lambda D)$ is  not log canonical and $D \sim_{\mathbb{Q}}
-K_X$, where $\lambda < \frac{1}{3}$. Therefore the log pair $(X,
\lambda D)$ is also not log canonical.

Let $Z$ be the curve in $|-K_X|$ that contains $P$. Since the
curve $Z$ is irreducible we may assume that the support of $D$
does not contain $Z$.

We derive that the pair $(X, \lambda D)$ is log canonical outside
of a point $P\in X$ and not log canonical at $P$. Let $\pi_1:
\tilde{X} \to X$ be the minimal resolution of $X$. The following
diagram shows how the exceptional curves intersect each other.
\bigskip

$\mathbb{A}_5 + \mathbb{A}_2$. \xymatrix{ {\bullet}^{E_1} \ar@{-}[r] &
{\bullet}^{E_2} \ar@{-}[r] & {\bullet}^{E_3} \ar@{-}[r] & {\bullet}^{E_4} \ar@{-}[r]
& {\bullet}^{E_5} & {\bullet}^{F_1} \ar@{-}[r]  & {\bullet}^{F_2}}
\bigskip

Then
$$
\tilde{D} \sim_{\mathbb{Q}} \pi_1^*(D) - a_1E_1 - a_2E_2 - a_3E_3
- a_4E_4 - a_5E_5 - b_1F_1 - b_2F_2 \text{ .}
$$

We have three lines $\tilde{L}_1,\tilde{L}_3,\tilde{L}_5$
intersecting the fundamental cycle as following
$$
\tilde{L}_1 \cdot E_1= \tilde{L}_1 \cdot F_1 = \tilde{L}_3 \cdot
E_3  = \tilde{L}_5 \cdot E_5 = \tilde{L}_5 \cdot F_2 =  1
$$
Therefore
\begin{eqnarray*}
\tilde{L_1} & \sim_{\mathbb{Q}} & \pi_1^*(L_1) -
\frac{1}{6}E_1-\frac{1}{3}E_2-\frac{1}{2}E_3-\frac{2}{3}E_4-\frac{5}{6}E_5-
\frac{1}{3}F_1 -\frac{2}{3} F_2\\
\tilde{L_3} & \sim_{\mathbb{Q}} & \pi_1^*(L_3) - \frac{5}{6} E_1 -
\frac{2}{3}E_2-\frac{1}{2}E_3-\frac{1}{3}E_4-\frac{1}{6}E_5-
\frac{2}{3}F_1 -\frac{1}{3} F_2\\
\tilde{L_5} & \sim_{\mathbb{Q}} &
\pi_1^*(L_5)-\frac{1}{2}E_1-E_2-\frac{3}{2}E_3-E_4 -\frac{1}{2}E_5
\end{eqnarray*}
Since $2L_1 \sim_{\mathbb{Q}} 2L_3 \sim_{\mathbb{Q}} 2L_5
\sim_{\mathbb{Q}} -K_X$ we see that $\mathrm{lct}(X) \leq
\frac{1}{3}$.

From the inequalities
\begin{eqnarray*}
0 \leq \tilde{D} \cdot \tilde{L}_1 & = & 1 - a_1 - b_1\\
0 \leq \tilde{D} \cdot \tilde{L}_3 & = & 1 - a_3\\
0 \leq \tilde{D} \cdot \tilde{L}_5 & = & 1 - a_5 - b_2\\
0 \leq E_1 \cdot \tilde{D} & = & 2a_1 - a_2\\
0 \leq E_2 \cdot \tilde{D} & = & 2a_2 - a_1 - a_3\\
0 \leq E_3 \cdot \tilde{D} & = & 2a_3 - a_2 - a_4\\
0 \leq E_4 \cdot \tilde{D} & = & 2a_4 - a_3 - a_5\\
0 \leq E_5 \cdot \tilde{D} & = & 2a_5 - a_4\\
0 \leq F_1 \cdot \tilde{D} & = & 2b_1 - b_2\\
0 \leq F_2 \cdot \tilde{D} & = & 2b_2 - b_1
\end{eqnarray*}
we see that
$$a_1 \leq 1 \text{, } a_2 \leq \frac{4}{3} \text{, } a_3 \leq 1
 \text{, } a_4 \leq \frac{4}{3} \text{, } a_5 \leq 1 \text{, } b_1 \leq 1 \text{, } b_2
\leq 1$$
 and what is more
 $$
 2a_5 \geq a_4 \text{ , } \frac{3}{2} a_4 \geq a_3 \text{ , } \frac{4}{3} a_3 \geq a_2
\text{ , } \frac{5}{4} a_2 \geq a_1 \text{ .}
 $$

The equivalence
$$
K_{\tilde{X}} + \lambda \tilde{D} + \lambda a_1  E_1 + \lambda a_2
E_2 + \lambda a_3 E_3 +\lambda a_4  E_4 +\lambda a_5  E_5 +
 \lambda b_1F_1 + \lambda b_2 F_2 \sim_{\mathbb{Q}}
\pi_1^*(K_X + D)
$$
implies that there is a point $Q\in E_1\cup E_2\cup E_3\cup E_4
\cup E_5 \cup F_1 \cup F_2$ such that the pair
$$
K_{\tilde{X}} + \lambda \tilde{D} + \lambda a_1  E_1 + \lambda a_2
E_2 + \lambda a_3 E_3 +\lambda a_4  E_4 +\lambda a_5  E_5 +
 \lambda b_1F_1 + \lambda b_2 F_2
$$
is not log canonical at $Q$.

\begin{itemize}
\item If the point $Q \in E_1$ and $Q\not \in E_2$ then
$$
K_{\tilde{X}} + \lambda \tilde{D} +  a_1 \lambda E_1
$$
is not log canonical at the point $Q$ and so is the pair
$$
K_{\tilde{X}} + \lambda \tilde{D} + E_1 \text{ , since  } \lambda
a_1 \leq 1 \text{ .}
$$
By adjunction $(E_1, \lambda \tilde{D}|_{E_1})$ is not log
canonical at $Q$ and
$$2a_1 - \frac{4}{5} a_1 \geq 2a_1 - a_2 = \tilde{D} \cdot E_1 \geq
 \text{mult}_Q\Bigl(\tilde{D} \cdot E_1 \Bigr)  > \frac{1}{\lambda} > 3 \text{ ,}$$
implies that $a_1 > \frac{5}{2}$ which is a contradiction.

\item If $Q\in E_1 \cap E_2$ then the log pair
$$
K_{\tilde{X}} + \lambda \tilde{D} +  a_1 \lambda E_1 + a_2 \lambda
E_2
$$
is not log canonical at the point $Q$ and so are the log pairs
$$
K_{\tilde{X}} + \lambda \tilde{D} + E_1 + a_2 \lambda E_2 \text{
and  } K_{\tilde{X}} + \lambda \tilde{D} + a_1 \lambda E_1 + E_2
\text{ .}
$$
By adjunction it follows that
$$
2a_2 - a_1 - a_3 = \tilde{D} \cdot E_2 \geq
\text{mult}_Q\Bigl(\tilde{D}|_{E_2} \Bigr) =
 \text{mult}_Q\Bigl(\tilde{D} \cdot E_2 \Bigr) > \frac{1}{\lambda} - a_1 > 3 - a_1 \text{
,}
$$
and
$$
2a_1 - a_2 = \tilde{D} \cdot E_1 \geq
\text{mult}_Q\Bigl(\tilde{D}|_{E_1} \Bigr) =
\text{mult}_Q\Bigl(\tilde{D} \cdot E_1 \Bigr)  > \frac{1}{\lambda}
- a_2 > 3 - a_2 \text{ .}
$$
From the second inequality we get that $a_1 \geq \frac{3}{2}$
which is a contradiction.

\item If $Q\in E_2$ but $Q\not \in E_1 \cup E_3$ then
$$
K_{\tilde{X}} + \lambda \tilde{D} + a_2 \lambda E_2
$$
is not log canonical at the point $Q$ and so is the pair
$$
K_{\tilde{X}} + \lambda \tilde{D} + E_2 \text{ , since  } a_2
\lambda \leq 1 \text{ .}
$$
By adjunction $(E_2, \lambda \tilde{D}|_{E_2})$ is not log
canonical at $Q$ and
$$2a_2 - \frac{a_2}{2} - \frac{3}{4}a_2 \geq 2a_2 - a_1 - a_3 = \tilde{D} \cdot E_2 \geq
\text{mult}_Q\Bigl(\tilde{D}|_{E_2} \Bigr) =
 \text{mult}_Q\Bigl(\tilde{D} \cdot E_2 \Bigr) > \frac{1}{\lambda} > 3  \text{ .}$$
Then we get $a_2 > 4 $ which is a contradiction.

\item If $Q\in E_2 \cap E_3$ then the log pair
$$
K_{\tilde{X}} + \lambda \tilde{D} +  a_2 \lambda E_2 + a_3 \lambda
E_3
$$
is not log canonical at the point $Q$ and so are the log pairs
$$
K_{\tilde{X}} + \lambda \tilde{D} + a_2 \lambda E_2 + E_3 \text{ ,
since } \lambda a_3 \leq 1 \text{ . }
$$
By adjunction it follows that
$$
 2a_3 - a_2 -a_4= \tilde{D} \cdot E_3 \geq
\text{mult}_Q\Bigl(\tilde{D} \cdot E_3 \Bigr)  > \frac{1}{\lambda}
- a_2 > 3 - a_2 \text{ .}
$$
which, together with the inequality $a_4 \geq \frac{2}{3} a_3$,
implies that $a_3 > \frac{9}{4}$. However, this is impossible
since $a_3 \leq 1$.

\item If $Q\in E_3$ but $Q\not \in E_2 \cup E_4$ then
$$
K_{\tilde{X}} + \lambda \tilde{D} + a_3 \lambda E_3
$$
is not log canonical at the point $Q$ and so is the pair
$$
K_{\tilde{X}} + \lambda \tilde{D} + E_3 \text{ , since  } a_3
\lambda \leq 1 \text{ .}
$$
By adjunction $(E_3, \lambda \tilde{D}|_{E_3})$ is not log
canonical at $Q$ and
$$2a_3 - \frac{2}{3} a_3 - \frac{2}{3} a_3 \geq 2a_3 - a_2 - a_4 = \tilde{D} \cdot E_3
\geq \text{mult}_Q\Bigl(\tilde{D}|_{E_3} \Bigr) =
 \text{mult}_Q\Bigl(\tilde{D} \cdot E_3 \Bigr) > \frac{1}{\lambda} > 3  \text{ ,}$$
implies that $a_3 > \frac{9}{2}$ which is false.

\item If $Q\in F_1$ then the log pair
$$
K_{\tilde{X}} + \lambda \tilde{D} + \lambda b_1F_1
$$
is not log canonical at the point $Q$ and so is the log pair
$$
K_{\tilde{X}} + \lambda \tilde{D} + F_1 \text{ .}
$$
By adjunction it follows that
$$
 \frac{3}{2} b_1 \geq 2b_1 -b_2 = \tilde{D} \cdot F_1 \geq
\text{mult}_Q\Bigl(\tilde{D}|_{F_1} \Bigr) =
 \text{mult}_Q\Bigl(\tilde{D} \cdot F_1 \Bigr) > 3  \text{ .}
$$
and we see then that $b_1>2$ which is not possible.

\item If $Q\in F_1 \cap F_2$ then the log pair
$$
K_{\tilde{X}} + \lambda \tilde{D} + \lambda b_1F_1 +\lambda b_2F_2
$$
is not log canonical at the point $Q$ and so is the log pair
$$
K_{\tilde{X}} + \lambda \tilde{D} + F_1 + \lambda b_1F_1 \text{ .}
$$
By adjunction it follows that
$$
 2b_1 -b_2 = \tilde{D} \cdot F_1 \geq \text{mult}_Q\Bigl(\tilde{D}|_{F_1} \Bigr) =
 \text{mult}_Q\Bigl(\tilde{D} \cdot F_1 \Bigr) > 3 - b_2  \text{ .}
$$
and we see then that $b_1> \frac{3}{2}$ which is false.

\item If $Q\in F_2$ then the log pair
$$
K_{\tilde{X}} + \lambda \tilde{D} + \lambda b_2F_2
$$
is not log canonical at the point $Q$ and so is the log pair
$$
K_{\tilde{X}} + \lambda \tilde{D} + F_2 \text{ .}
$$
By adjunction it follows that
$$
 \frac{3}{2} b_2 \geq 2b_2 -b_1 = \tilde{D} \cdot F_2 \geq
\text{mult}_Q\Bigl(\tilde{D}|_{F_2} \Bigr) =
 \text{mult}_Q\Bigl(\tilde{D} \cdot F_2 \Bigr) > 3  \text{ .}
$$
and we see then that $b_2 > 2$ which is not possible.

\end{itemize}
\end{proof}

\subsection{Del Pezzo surfaces of degree 2 with exactly one
$\mathbb{A}_7$ type singularity}

In this section we will prove the following.

\begin{lemma}
\label{degree2A7} Let $X$ be a del Pezzo surface with at most one
Du Val singularity of type $\mathbb{A}_7$ and $K_X^2=1$. Then the
global log canonical threshold of $X$ is
$$
\mathrm{lct} (X) = \frac{1}{3} \text{ .}
$$
\end{lemma}

\begin{proof}

 Suppose that $\mathrm{lct}(X)
< \frac{1}{3}$, then there exists an effective
$\mathbb{Q}$-divisor $D\in X$ and a positive rational number
$\lambda < \frac{1}{3}$, such that the log pair $(X, \lambda D)$
is  not log canonical and $D \sim_{\mathbb{Q}} -K_X$, where
$\lambda < \frac{1}{3}$.

We derive that the pair $(X, \lambda D)$ is log canonical outside
of a point $P\in X$ and not log canonical at $P$. Let $\pi_1:
\tilde{X} \to X$ be the minimal resolution of $X$. The following
diagram shows how the exceptional curves intersect each other.
\bigskip

$\mathbb{A}_7$. \xymatrix{ {\bullet}^{E_1} \ar@{-}[r] & {\bullet}^{E_2}
\ar@{-}[r] & {\bullet}^{E_3} \ar@{-}[r] & {\bullet}^{E_4} \ar@{-}[r] &
{\bullet}^{E_5} \ar@{-}[r] & {\bullet}^{E_6} \ar@{-}[r] & {\bullet}^{E_7}}
\bigskip

Then
\begin{eqnarray*}
\tilde{D} & \sim_{\mathbb{Q}} & \pi_1^*(D) - a_1E_1 - a_2E_2 -
a_3E_3 - a_4E_4 - a_5E_5 -a_6E_6 -a_7E_7 \text{ .}
\end{eqnarray*}

Furthermore there are lines $\tilde{L}_2,  \tilde{L}_6 \in X$ that
pass through the point $P$ whose strict
 transforms are $(-1)$-curves that intersect the fundamental cycle as following.
 $$ \tilde{L}_2 \cdot E_2 =  \tilde{L}_6 \cdot E_6 = 1$$
 and
 $$
 \tilde{L}_i \cdot E_j = 0 \text{ for all } i,j= 2,6 \text{ with } i \not = j \text{ .}
 $$
Then we easily get that
\begin{eqnarray*}
 \tilde{L_2} & = & \pi^*(L_2) - \frac{3}{4} E_1 - \frac{3}{2} E_2 - \frac{5}{4} E_3 - E_4
- \frac{3}{4} E_5 -\frac{1}{2} E_6 - \frac{1}{4} E_7\\
\tilde{L_6} & = & \pi^*(L_6) - \frac{1}{4} E_1 - \frac{1}{2} E_2 -
\frac{3}{4} E_3 - E_4 -\frac{5}{4} E_5 -\frac{3}{2} E_6
-\frac{3}{4} E_7 \text{ .}
\end{eqnarray*}
Because $2L_2 \sim_{\mathbb{Q}} 2L_6 \sim_{\mathbb{Q}} -K_X$ we
have that $\mathrm{lct}(X) \leq \frac{1}{3}$.

From the inequalities
\begin{eqnarray*}
0 \leq \tilde{D} \cdot \tilde{L}_2 & = & 1 - a_2\\
0 \leq \tilde{D} \cdot \tilde{L}_6 & = & 1 - a_6\\
0 \leq E_1 \cdot \tilde{D} & = & 2a_1 - a_2\\
0 \leq E_2 \cdot \tilde{D} & = & 2a_2 - a_1 - a_3\\
0 \leq E_3 \cdot \tilde{D} & = & 2a_3 - a_2 - a_4\\
0 \leq E_4 \cdot \tilde{D} & = & 2a_4 - a_3 - a_5\\
0 \leq E_5 \cdot \tilde{D} & = & 2a_5 - a_4 - a_6\\
0 \leq E_6 \cdot \tilde{D} & = & 2a_6 - a_5 - a_7\\
0 \leq E_7 \cdot \tilde{D} & = & 2a_7 - a_6
\end{eqnarray*}
we get
$$
2a_7 \geq a_6 \text{ , } \frac{3}{2} a_6 \geq a_5 \text{ , }
\frac{4}{3} a_5 \geq a_4 \text{ , } \frac{5}{4} a_4 \geq a_3
\text{ , } \frac{6}{5} a_3 \geq a_2 \text{ , } \frac{7}{6} a_2
\geq a_1
$$
and
$$
2a_1 \geq a_2 \text{ , } \frac{3}{2} a_2 \geq a_3 \text{ , }
\frac{4}{3} a_3 \geq a_4 \text{ , } \frac{5}{4} a_4 \geq a_5
\text{ , } \frac{6}{5} a_5 \geq a_6 \text{ , } \frac{7}{6} a_6
\geq a_7 \text{ .}
$$
Therefore
$$
a_1 \leq \frac{7}{6} \text{ , } a_2 \leq 1 \text{ , } a_3 \leq
\frac{3}{2} \text{ , } a_4 \leq 2 \text{ , } a_5 \leq \frac{3}{2}
\text{ , } a_6 \leq 1 \text{ , } a_7 \leq \frac{7}{6} \text{ . }
$$
 The equivalence
$$
K_{\tilde{X}} + \lambda \tilde{D} + \lambda a_1E_1 + \lambda
a_2E_2 + \lambda a_3E_3 + \lambda a_4E_4 + \lambda a_5E_5 +
\lambda a_6E_6 + \lambda a_7E_7 \sim_{\mathbb{Q}} \pi_1^*(K_X + D)
$$
implies that there is a point $Q\in E_1\cup E_2\cup E_3\cup E_4
\cup E_5 \cup E_6 \cup E_7$, such that the pair
$$K_{\tilde{X}} + \lambda \tilde{D} + \lambda a_1E_1 + \lambda a_2E_2 + \lambda a_3E_3 +
\lambda a_4E_4 + \lambda a_5E_5 + \lambda a_6E_6 + \lambda
a_7E_7$$ is not log canonical at $Q$.

\begin{itemize}
\item If the point $Q \in E_1$ and $Q\not \in E_2$ then
$$
K_{\tilde{X}} + \lambda \tilde{D} +  a_1 \lambda E_1
$$
is not log canonical at the point $Q$ and so is the pair
$$
K_{\tilde{X}} + \lambda \tilde{D} + E_1 \text{ , since  } a_1
\lambda \leq 1 \text{ .}
$$
By adjunction $(E_1, \lambda \tilde{D}|_{E_1})$ is not log
canonical at $Q$ and
$$
 \frac{8}{7}a_1 \geq  2a_1 -\frac{6}{7}a_1 \geq 2a_1 - a_2 = \tilde{D} \cdot E_1
\geq
 \text{mult}_Q\Bigl(\tilde{D} \cdot E_1 \Bigr)  > \frac{1}{\lambda} > 3 \text{ ,}$$
which is false, since $a_1 \leq \frac{7}{6}$.

\item If $Q\in E_1 \cap E_2$ then the log pair
$$
K_{\tilde{X}} + \lambda \tilde{D} +  \lambda a_1 E_1 +  \lambda
a_2E_2
$$
is not log canonical at the point $Q$ and so is the log pair
$$
K_{\tilde{X}} + \lambda \tilde{D} + E_1 + \lambda a_2 E_2  \text{
.}
$$
By adjunction it follows that
$$
 2a_2 -\frac{5}{6}a_2 -a_1 \geq 2a_2 - a_1 - a_3 = \tilde{D} \cdot E_2 \geq
\text{mult}_Q\Bigl(\tilde{D}|_{E_2} \Bigr) =
\text{mult}_Q\Bigl(\tilde{D} \cdot E_2 \Bigr)  > \frac{1}{\lambda}
- a_2 > 3 - a_1 \text{ ,}
$$
which is false, since $a_2 \leq 1$.

\item If $Q\in E_2$ but $Q\not \in E_1 \cup E_3$ then
$$
K_{\tilde{X}} + \lambda \tilde{D} +  \lambda a_2 E_2
$$
is not log canonical at the point $Q$ and so is the pair
$$
K_{\tilde{X}} + \lambda \tilde{D} + E_2 \text{ , since  }  \lambda
a_2 \leq 1 \text{ .}
$$
By adjunction $(E_2, \lambda \tilde{D}|_{E_2})$ is not log
canonical at $Q$ and
$$2a_2 - \frac{5}{6} a_2 \geq 2a_2 - a_1 - a_3 = \tilde{D} \cdot E_2 \geq
 \text{mult}_Q\Bigl(\tilde{D} \cdot E_2 \Bigr) > \frac{1}{\lambda} > 3  \text{ ,}$$
which is false, since $a_2 \leq 1$.

\item If $Q\in E_2 \cap E_3$ then the log pair
$$
K_{\tilde{X}} + \lambda \tilde{D} +  \lambda a_2 E_2 +  \lambda
a_3E_3
$$
is not log canonical at the point $Q$ and so is the log pair
$$
K_{\tilde{X}} + \lambda \tilde{D} +  \lambda a_2E_2 + E_3 \text{ ,
since } \lambda a_3 < 1 \text{ .}
$$
By adjunction it follows that
$$
2a_3 - a_2 - \frac{4}{5} a_3 \geq 2a_3 - a_2 -a_4= \tilde{D} \cdot
E_3 \geq \text{mult}_Q\Bigl(\tilde{D}|_{E_3} \Bigr) =
\text{mult}_Q\Bigl(\tilde{D} \cdot E_3 \Bigr)  > \frac{1}{\lambda}
- a_2 > 3 - a_2 \text{ ,}
$$
which implies that $a_3 > \frac{5}{2}$, which is impossible.

\item If $Q\in E_3$ but $Q\not \in E_2 \cup E_4$ then
$$
K_{\tilde{X}} + \lambda \tilde{D} +  \lambda a_3E_3
$$
is not log canonical at the point $Q$ and so is the pair
$$
K_{\tilde{X}} + \lambda \tilde{D} + E_3 \text{ , since  }  \lambda
a_3\leq 1 \text{ .}
$$
By adjunction $(E_3, \lambda \tilde{D}|_{E_3})$ is not log
canonical at $Q$ and
$$2a_3  - \frac{4}{5} a_3 \geq 2a_3 - a_2 - a_4 = \tilde{D} \cdot E_3 \geq
 \text{mult}_Q\Bigl(\tilde{D} \cdot E_3 \Bigr) > \frac{1}{\lambda} > 3  \text{ .}$$
This inequality  implies that $a_3 > \frac{5}{2}$, which is
impossible.

\item If $Q\in E_3 \cap E_4$ then the log pair
$$
K_{\tilde{X}} + \lambda \tilde{D} +  \lambda a_3 E_3 +  \lambda
a_4 E_4
$$
is not log canonical at the point $Q$ and so is the log pair
$$
K_{\tilde{X}} + \lambda \tilde{D} +  \lambda a_3E_3 + E_4 \text{ ,
since } \lambda a_4 \leq 1 \text{ .}
$$
By adjunction it follows that

$$
2a_4 - a_3 - \frac{3}{4} a_4 \geq 2a_4 - a_3 -a_5= \tilde{D} \cdot
E_4 \geq \text{mult}_Q\Bigl(\tilde{D}|_{E_4} \Bigr) =
\text{mult}_Q\Bigl(\tilde{D} \cdot E_4 \Bigr)  > \frac{1}{\lambda}
- a_3 > 3 - a_3 \text{ ,}
$$
which contradicts $a_4 \leq 2$.

\item If $Q\in E_4$ but $Q\not \in E_3 \cup E_5$ then
$$
K_{\tilde{X}} + \lambda \tilde{D} +  \lambda a_4E_4
$$
is not log canonical at the point $Q$ and so is the pair
$$
K_{\tilde{X}} + \lambda \tilde{D} + E_4 \text{ , since  }  \lambda
a_4\leq 1 \text{ .}
$$
By adjunction $(E_4, \lambda \tilde{D}|_{E_4})$ is not log
canonical at $Q$ and
$$ 2a_4  - \frac{3}{4} a_4 \geq 2a_4 - a_3 - a_5 = \tilde{D} \cdot E_4 \geq \text{mult}_Q\Bigl(\tilde{D}|_{E_4}
\Bigr) =
 \text{mult}_Q\Bigl(\tilde{D} \cdot E_4 \Bigr) > \frac{1}{\lambda} > 3  \text{ ,}$$
which is false since $a_4 \leq 2$.

\end{itemize}
\end{proof}

\subsection{Del Pezzo surfaces of degree 6 with one $\mathbb{A}_2$
and one $\mathbb{A}_1$ type singularity}

In this section we will prove the following.

\begin{lemma}
\label{A2} Let $X$ be a del Pezzo surface with one Du Val
singularity of type $\mathbb{A}_2$, one of type $\mathbb{A}_1$ and
$K_X^2=6$. Then the global log canonical threshold of $X$ is
$$
\mathrm{lct} (X) = \frac{1}{6} \text{ .}
$$
\end{lemma}

\begin{proof}

 Suppose
that $\mathrm{lct}(X) < \frac{1}{6}$. Then there exists an
effective $\mathbb{Q}$-divisor $D\sim_{\mathbb{Q}} -K_X$ such that
the log pair $(X, \lambda D)$ is  not log canonical for a rational
number $\lambda < \frac{1}{6}$.

We derive that the pair $(X, \lambda D)$ is log canonical
everywhere except for a singular point $P$, at which point $P$ it
is not log canonical. Let $\pi_1: \tilde{X} \to X$ be the minimal
resolution of $X$. The following diagram shows how the exceptional
curves intersect each other.
\bigskip

$\mathbb{A}_2+\mathbb{A}_1$ \xymatrix{ {\bullet}^{E_1} \ar@{-}[r] &
{\bullet}^{E_2}  & {\bullet}^{E_3} }
\bigskip

Then
$$
\tilde{D} \sim_{\mathbb{Q}} \pi_1^*(D)-a_1E_1-a_2E_2-a_3E_3 \text{
.} $$

We have a $-1$ curve $\tilde{L}_1$ intersecting the fundamental
cycle as following
$$
\tilde{L}_1  \cdot E_2 = \tilde{L}_1  \cdot E_3 =1 \text{ , }
\tilde{L}_1  \cdot E_1=0
$$
and
\begin{eqnarray*}
\tilde{L}_1 & \sim_{\mathbb{Q}} &
\pi_1^*(L_1)-\frac{1}{3}E_1-\frac{2}{3}E_2-\frac{1}{2}E_3 \text{
.}
\end{eqnarray*}
Since $6L_1 \sim_{\mathbb{Q}} -K_X$ we see that $\mathrm{lct}(X)
\leq \frac{1}{6}$.

 From the
inequalities
\begin{eqnarray*}
0 \leq \tilde{D} \cdot \tilde{L}_1 & = & 1 - a_2 -  a_3 \\
0 \leq E_1 \cdot \tilde{D} & = & 2a_1 - a_2\\
0 \leq E_2 \cdot \tilde{D} & = & 2a_2 - a_1\\
0 \leq E_3 \cdot \tilde{D} & = & 2a_3
\end{eqnarray*}
we see that $a_1 \leq 2 \text{, } a_2 \leq 1 \text{, } a_3 \leq 1
\text{ .}$

The equivalence
$$
K_{\tilde{X}} + \lambda \tilde{D} + \lambda a_1E_1 + \lambda
a_2E_2 + \lambda a_3E_3 \sim_{\mathbb{Q}} \pi_1^*(K_X+ \lambda D)
$$
implies that there is a point $Q\in E_1\cup E_2\cup E_3 $ such
that the pair
$$K_{\tilde{X}} + \lambda \tilde{D} + \lambda a_1E_1 + \lambda a_2E_2 +
\lambda a_3E_3 $$ is not log canonical at $Q$.

\begin{itemize}
\item If the point $Q \in E_1$ and $Q\not \in E_2$ then
$$
K_{\tilde{X}} + \lambda \tilde{D} + \lambda a_1E_1
$$
is not log canonical at the point $Q$ and so is the pair
$$
K_{\tilde{X}} + \lambda \tilde{D} + E_1 \text{ , since  } \lambda
a_1 \leq 1 \text{ .}
$$
By adjunction $(E_1, \lambda \tilde{D}|_{E_1})$ is not log
canonical at $Q$ and
$$2a_1 - \frac{a_1}{2} \geq 2a_1 - a_2 = \tilde{D} \cdot E_1 \geq
\text{mult}_Q\Bigl(\tilde{D}|_{E_1} \Bigr) =
 \text{mult}_Q\Bigl(\tilde{D} \cdot E_1 \Bigr)  > 6 \text{ ,}$$
implies that $a_1 > 4$ which is false.

\item If $Q\in E_1 \cap E_2$ then the log pair
$$
K_{\tilde{X}} + \lambda \tilde{D} + \lambda a_1E_1 + \lambda
a_2E_2
$$
is not log canonical at the point $Q$ and so are the log pairs
$$
K_{\tilde{X}} + \lambda \tilde{D} + E_1 + \lambda a_2E_2 \text{
and } K_{\tilde{X}} + \lambda \tilde{D} + \lambda a_1E_1 + E_2
\text{ ,}
$$
   since   $\lambda a_2 \leq 1$. By adjunction it follows that
$$
2a_1 - a_2 = \tilde{D} \cdot E_1 \geq
\text{mult}_Q\Bigl(\tilde{D}|_{E_1} \Bigr) =
 \text{mult}_Q\Bigl(\tilde{D} \cdot E_1 \Bigr) > 6 - a_2
$$
and
$$
2a_2 - a_1 -a_3 = \tilde{D} \cdot E_2 \geq
\text{mult}_Q\Bigl(\tilde{D}|_{E_2} \Bigr) =
 \text{mult}_Q\Bigl(\tilde{D} \cdot E_2 \Bigr) > 6 - a_1 \text{ . }
$$
This implies that $a_1 >3 \text{, } a_2 > 3$ which is false.

\item If $Q\in E_3$ then the log pair
$$
K_{\tilde{X}} + \lambda \tilde{D} + \lambda a_3E_3
$$
is not log canonical at the point $Q$ and so is the log pair
$$
K_{\tilde{X}} + \lambda \tilde{D} + E_3 \text{ .}
$$
By adjunction it follows that
$$
2 \geq 2a_3 = \tilde{D} \cdot E_3 \geq
\text{mult}_Q\Bigl(\tilde{D}|_{E_3} \Bigr) =
 \text{mult}_Q\Bigl(\tilde{D} \cdot E_3 \Bigr) > 6
$$
which is a contradiction.

\end{itemize}
\end{proof}

\subsection{Del Pezzo surfaces of degree 5 with exactly one
$\mathbb{A}_4$ type singularity}

In this section we will prove the following.

\begin{lemma}
\label{degree5A4} Let $X$ be a del Pezzo surface with exactly one
Du Val singularity of type $\mathbb{A}_4$ and $K_X^2=5$. Then the
global log canonical threshold of $X$ is
$$
\mathrm{lct} (X) = \frac{1}{6} \text{ .}
$$
\end{lemma}

\begin{proof}

 Suppose $\mathrm{lct}(X) <
\frac{1}{6}$. Then there exist an effective $\mathbb{Q}$-divisor
$D\sim_{\mathbb{Q}} -K_X$ and a positive rational number $\lambda
< \frac{1}{6}$, such that the log pair $(X, \lambda D)$ is  not
log canonical.

We derive that the pair $(X, \lambda D)$ is log canonical
everywhere except for a Du Val point $P$, at which point the pair
is not log canonical. Let $\pi_1: \tilde{X} \to X$ be the minimal
resolution of $X$. The following diagram shows how the exceptional
curves intersect each other.
\bigskip

$\mathbb{A}_4$. \xymatrix{ {\bullet}^{E_1} \ar@{-}[r] & {\bullet}^{E_2}
\ar@{-}[r] & {\bullet}^{E_3} \ar@{-}[r] & {\bullet}^{E_4}   }
\bigskip

Then
\begin{eqnarray*}
\tilde{D} & \sim_{\mathbb{Q}} & \pi_1^*(D) - a_1E_1 - a_2E_2 -
a_3E_3 - a_4E_4 \text{ .}
\end{eqnarray*}

Furthermore there is a line $L_1 \in X$ that passes through the
point $P$, whose strict
 transform intersects the fundamental cycle as following
 $$ \tilde{L}_1 \cdot E_2 = 1 \text {  and  }
 \tilde{L}_1 \cdot E_j = 0 \text{ for all }  j = 1, 3, 4 \text{ .}
 $$
Then we easily get that
\begin{eqnarray*}
 \tilde{L_1} & \sim_{\mathbb{Q}} & \pi^*(L_1) - \frac{3}{5} E_1 - \frac{6}{5} E_2 - \frac{4}{5} E_3 -
\frac{2}{5} E_4 \text{ .}
\end{eqnarray*}
Because $5L_1 \sim_{\mathbb{Q}} -K_X$ we have that
$\mathrm{lct}(X) \leq \frac{1}{6}$.

Since $L_1$ is irreducible we can assume that $L_1 \not \in
\text{Supp}D$. Then from the inequalities
\begin{eqnarray*}
0 \leq \tilde{D} \cdot \tilde{L}_1 & = & 1 - a_2\\
0 \leq E_1 \cdot \tilde{D} & = & 2a_1 - a_2\\
0 \leq E_2 \cdot \tilde{D} & = & 2a_2 - a_1 - a_3\\
0 \leq E_3 \cdot \tilde{D} & = & 2a_3 - a_2 - a_4\\
0 \leq E_4 \cdot \tilde{D} & = & 2a_4 - a_3
\end{eqnarray*}
we get
$$
2 a_4 \geq a_3 \text{ , } \frac{3}{2} a_3 \geq a_2 \text{ , }
\frac{4}{3} a_2 \geq a_1
$$
and
$$
2a_1 \geq a_2 \text{ , } \frac{3}{2} a_2 \geq a_3 \text{ , }
\frac{4}{3} a_3 \geq a_4 \text{ .}
$$
Therefore
$$
a_1 \leq \frac{4}{3} \text{ , } a_2 \leq 1 \text{ , } a_3 \leq
\frac{3}{2} \text{ , } a_4 \leq 2  \text{ . }
$$
 The equivalence
$$
K_{\tilde{X}} + \lambda \tilde{D} + \lambda a_1E_1 + \lambda
a_2E_2 + \lambda a_3E_3 + \lambda a_4E_4 \sim_{\mathbb{Q}}
\pi_1^*(K_X + \lambda D)
$$
implies that there is a point $Q\in E_1\cup E_2\cup E_3\cup E_4 $,
such that the pair
$$K_{\tilde{X}} + \lambda \tilde{D} + \lambda a_1E_1 + \lambda a_2E_2 + \lambda a_3E_3 +
\lambda a_4E_4$$ is not canonical at $Q$.

\begin{itemize}
\item If the point $Q \in E_1$ and $Q\not \in E_2$ then
$$
K_{\tilde{X}} + \lambda \tilde{D} +  a_1 \lambda E_1
$$
is not log canonical at the point $Q$ and so is the pair
$$
K_{\tilde{X}} + \lambda \tilde{D} + E_1 \text{ , since  } a_1
\lambda \leq 1 \text{ .}
$$
By adjunction $(E_1, \lambda \tilde{D}|_{E_1})$ is not log
canonical at $Q$ and
$$
\frac{5}{3} \geq 2 a_1 - \frac{3}{4} a_1 \geq 2a_1 - a_2 =
\tilde{D} \cdot E_1 \geq
 \text{mult}_Q\Bigl(\tilde{D} \cdot E_1 \Bigr)  > \frac{1}{\lambda} > 6 \text{ ,}$$
which is a contradiction.

\item If $Q\in E_1 \cap E_2$ then the log pair
$$
K_{\tilde{X}} + \lambda \tilde{D} +  \lambda a_1 E_1 +  \lambda
a_2E_2
$$
is not log canonical at the point $Q$ and so is the log pair
$$
K_{\tilde{X}} + \lambda \tilde{D} + E_1 + \lambda a_2 E_2  \text{
.}
$$
By adjunction it follows that
$$
2a_1 - a_2 = \tilde{D} \cdot E_1 \geq
\text{mult}_Q\Bigl(\tilde{D}|_{E_1} \Bigr) =
\text{mult}_Q\Bigl(\tilde{D} \cdot E_1 \Bigr)  > \frac{1}{\lambda}
- a_2 > 6 - a_2 \text{ ,}
$$
and
$$ 2a_2 - a_1 - a_3 = \tilde{D} \cdot E_2 \geq
 \text{mult}_Q\Bigl(\tilde{D} \cdot E_2 \Bigr) > \frac{1}{\lambda} > 6 - a_1  \text{ .}$$
This implies that $ a_1>3 $ which is false.

\item If $Q\in E_2$ but $Q\not \in E_1 \cup E_3$ then
$$
K_{\tilde{X}} + \lambda \tilde{D} +  \lambda a_2 E_2
$$
is not log canonical at the point $Q$ and so is the pair
$$
K_{\tilde{X}} + \lambda \tilde{D} + E_2 \text{ , since  }  \lambda
a_2 \leq 1 \text{ .}
$$
By adjunction $(E_2, \lambda \tilde{D}|_{E_2})$ is not log
canonical at $Q$ and
$$2a_2 -\frac{a_2}{2} - \frac{2}{3}a_2 \geq 2a_2 - a_1 - a_3 = \tilde{D} \cdot E_2 \geq
 \text{mult}_Q\Bigl(\tilde{D} \cdot E_2 \Bigr) > \frac{1}{\lambda} >6 \text{ ,}$$
which is false, since $a_2 \leq 1$.

\item If $Q\in E_2 \cap E_3$ then the log pair
$$
K_{\tilde{X}} + \lambda \tilde{D} +  \lambda a_2 E_2 +  \lambda
a_3E_3
$$
is not log canonical at the point $Q$ and so is the log pair
$$
K_{\tilde{X}} + \lambda \tilde{D} +  \lambda a_2E_2 + E_3 \text{ ,
since } \lambda a_3 < 1 \text{ .}
$$
By adjunction it follows that
$$2a_2 -\frac{a_2}{2} - a_3 \geq 2a_2 - a_1 - a_3 = \tilde{D} \cdot E_2 \geq
 \text{mult}_Q\Bigl(\tilde{D} \cdot E_2 \Bigr) > \frac{1}{\lambda} >6 -a_3 $$
and
$$
2a_3 - a_2 -a_4= \tilde{D} \cdot E_3 \geq
\text{mult}_Q\Bigl(\tilde{D}|_{E_3} \Bigr) =
\text{mult}_Q\Bigl(\tilde{D} \cdot E_3 \Bigr)  > \frac{1}{\lambda}
- a_2 > 6 - a_2 \text{ . }
$$
This implies $a_2>4 \text{, } a_3>4 \text{, }$ which is false.
\end{itemize}
\end{proof}

\subsection{Del Pezzo surfaces of degree 4 with one $\mathbb{A}_3$
and two $\mathbb{A}_1$ type singularities}

In this section we will prove the following.

\begin{lemma}
\label{A3+A1} Let $X$ be a del Pezzo surface with two Du Val
singular points of type $\mathbb{A}_3$, two $\mathbb{A}_1$ type
singular points and $K_X^2=4$. Then the global log canonical
threshold of $X$ is
$$
\mathrm{lct} (X) = \frac{1}{4} \text{ .}
$$
\end{lemma}

\begin{proof}

 Suppose that $ \mathrm{lct}(X) < \frac{1}{4}$, then there
exists an effective $\mathbb{Q}$-divisor $D$ such that
$D\sim_{\mathbb{Q}} -K_X$ and the log pair $(X,\lambda D)$ is not
log canonical for some rational number $\lambda < \frac{1}{4}$.

We derive that the pair $(X, \lambda D)$ is log canonical
everywhere except for a singular point $P\in X$, where $(X,
\lambda D)$ is not log canonical. Let $\pi_1: \tilde{X} \to X$ be
the minimal resolution of $X$. The following diagram shows how the
exceptional curves intersect each other.
\bigskip

$\mathbb{A}_3+2 \mathbb{A}_1$ \xymatrix{ {\bullet}^{E_1} \ar@{-}[r] &
{\bullet}^{E_2} \ar@{-}[r] & {\bullet}^{E_3}  & {\bullet}^{F_1} & {\bullet}^{G_1}}
\bigskip

Then
$$
\tilde{D} \sim_{\mathbb{Q}} \pi_1^*(D)-a_1E_1-a_2E_2-a_3E_3
-b_1F_1 - c_1 G_1 \text{ .}$$

We have two lines $L_1,L_3$ intersecting the fundamental cycle as
following
\begin{eqnarray*}
\tilde{L_1} & \sim_{\mathbb{Q}} & \pi_1^*(L_1) - \frac{3}{4}E_1 -
\frac{1}{2}E_2 - \frac{1}{4}E_3 -
\frac{1}{2} F_1\\
\tilde{L_3} & \sim_{\mathbb{Q}} & \pi_1^*(L_3) - \frac{1}{4}E_1 -
\frac{1}{2}E_2 - \frac{3}{4}E_3 - \frac{1}{2} G_1
\end{eqnarray*}

Since $4L_1 \sim_{\mathbb{Q}} 4L_3 \sim_{\mathbb{Q}} -K_X$ we see
that $\mathrm{lct}(X) \leq \frac{1}{4}$. Moreover we can assume
that $L_1 \not \in \text{Supp}D$ and $L_3 \not \in \text{Supp}D$.
From the inequalities
\begin{eqnarray*}
0 \leq \tilde{D} \cdot \tilde{L}_1 & = & 1 - a_1 - b_1\\
0 \leq \tilde{D} \cdot \tilde{L}_3 & = & 1 - a_3 - c_1\\
0 \leq E_1 \cdot \tilde{D} & = & 2a_1 - a_2\\
0 \leq E_2 \cdot \tilde{D} & = & 2a_2 - a_1 -a_3\\
0 \leq E_3 \cdot \tilde{D} & = & 2a_3-a_2\\
0 \leq F_1 \cdot \tilde{D} & = & 2b_1\\
0 \leq G_1 \cdot \tilde{D} & = & 2c_1\\
\end{eqnarray*}
we see that
$$
a_2 \leq 2a_1 \text{, } a_3 \leq \frac{3}{2} a_2 \text{  and  }
a_2 \leq 2a_3 \text{, } a_1 \leq \frac{3}{2} a_2 \text{ .}
$$
Therefore we get the bounds
$$a_1 \leq 1 \text{, } a_2 \leq 2 \text{, } a_3 \leq 1 \text{, } b_1 \leq 1 \text{, } c_1
\leq 1 \text{ .}$$

%CHECK WHY $\mathrm{lct}(X) \leq \frac{1}{4}$!!!!!!!!!!!!!!!!!!!!!!!!!

The equivalence
$$
K_{\tilde{X}} + \lambda \tilde{D} + \lambda a_1E_1 + \lambda
a_2E_2 + \lambda a_3E_3 \sim_{\mathbb{Q}} \pi_1^*(K_X+ \lambda D)
$$
implies that there is a point $Q\in E_1\cup E_2\cup E_3 $ such
that the pair
$$K_{\tilde{X}} + \lambda \tilde{D} + \lambda a_1E_1 + \lambda a_2E_2 +
\lambda a_3E_3 $$ is not log canonical at $Q$.

\begin{itemize}
\item If the point $Q \in E_1$ and $Q\not \in E_2$ then
$$
K_{\tilde{X}} + \lambda \tilde{D} + \lambda a_1E_1
$$
is not log canonical at the point $Q$ and so is the pair
$$
K_{\tilde{X}} + \lambda \tilde{D} + E_1 \text{, since } \lambda
a_1 \leq 1 \text{ .}
$$
By adjunction $(E_1, \lambda \tilde{D}|_{E_1})$ is not log
canonical at $Q$ and
$$\frac{4}{3} a_1 \geq 2a_1 - \frac{2}{3} a_1 \geq  2a_1 - a_2 = \tilde{D} \cdot E_1 \geq
 \text{mult}_Q\Bigl(\tilde{D} \cdot E_1 \Bigr)  > 4 \text{ ,}$$
implies that $a_1 > 3$ which is false.

\item If $Q\in E_1 \cap E_2$ then the log pair
$$
K_{\tilde{X}} + \lambda \tilde{D} + \lambda a_1E_1 + \lambda
a_2E_2
$$
is not log canonical at the point $Q$ and so are the log pairs
$$
K_{\tilde{X}} + \lambda \tilde{D} + E_1 + \lambda a_2E_2 \text{
and } K_{\tilde{X}} + \lambda \tilde{D} + \lambda a_1E_1 + E_2
\text{ .}
$$
By adjunction it follows that
$$
2a_1 - a_2 = \tilde{D} \cdot E_1 \geq
 \text{mult}_Q\Bigl(\tilde{D} \cdot E_1 \Bigr) > 4 - a_2
$$
and
$$
\frac{3}{2} a_2 -a_1 \geq 2a_2-\frac{a_2}{2} - a_1 \geq 2a_2 - a_1
-a_3 = \tilde{D} \cdot E_2 \geq
 \text{mult}_Q\Bigl(\tilde{D} \cdot E_2 \Bigr) >4 - a_1 \text{ . }
$$
The first inequality implies that $a_1>2$, which is false.

\item If $Q\in E_2$ and $Q\not \in E_1 \cup E_3$ then the log pair
$$
K_{\tilde{X}} + \lambda \tilde{D} + \lambda a_2E_2
$$
is not log canonical at the point $Q$ and so is the log pair
$$
K_{\tilde{X}} + \lambda \tilde{D} + E_2 \text{ .}
$$
By adjunction it follows that
$$
2a_2-\frac{a_2}{2} -\frac{a_2}{2} \geq 2a_2-a_1-a_3 = \tilde{D}
\cdot E_3 \geq
 \text{mult}_Q\Bigl(\tilde{D} \cdot E_3 \Bigr) > 4
$$
and this implies that $a_2 > 4$, which is false.

\item If the point $Q \in F_1$  then
$$
K_{\tilde{X}} + \lambda \tilde{D} + \lambda b_1F_1
$$
is not log canonical at the point $Q$ and so is the pair
$$
K_{\tilde{X}} + \lambda \tilde{D} + F_1 \text{, since } \lambda
b_1 \leq 1 \text{ .}
$$
By adjunction $(F_1, \lambda \tilde{D}|_{F_1})$ is not log
canonical at $Q$ and
$$ 2b_1 = \tilde{D} \cdot F_1 \geq \text{mult}_Q\Bigl(\tilde{D}|_{F_1} \Bigr) =
 \text{mult}_Q\Bigl(\tilde{D} \cdot F_1 \Bigr)  > 4 \text{ ,}$$
implies that $b_1 > 2$ which is false.

\end{itemize}
\end{proof}

\subsection{Del Pezzo surfaces of degree 4 with one $\mathbb{D}_5$
type singularity.}

In this section we will prove the following.

\begin{lemma}
\label{degree4D5} Let $X$ be a del Pezzo surface with  one Du Val
singularity of type $\mathbb{D}_5$ and $K_X^2=4$. Then the global
log canonical threshold of $X$ is
$$
\mathrm{lct} (X) = \frac{1}{6} \text{ .}
$$
\end{lemma}

\begin{proof}

Suppose that $\mathrm{lct}(X)<\frac{1}{6}$,  then there exists a
$\mathbb{Q}$-divisor $D \in X$ such that $D \sim_{\mathbb{Q}}
-K_X$ and the log pair $(X, \lambda D)$ is  not log canonical, for
some rational number $\lambda < \frac{1}{6}$. We derive that the
pair $(X, \lambda D)$ is log canonical everywhere except for a
singular point $P\in X$, where $(X, \lambda D)$ is not log
canonical. Let $\pi: \tilde{X} \to X$ be the minimal resolution of
$X$. The configuration of the exceptional curves is given by the
following Dynkin diagram.
\bigskip

$\mathbb{D}_5$. 
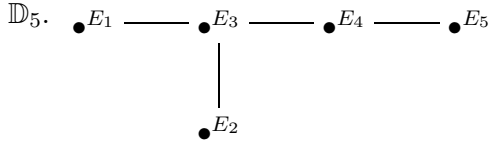
\xymatrix{ {\bullet}^{E_1} \ar@{-}[r] & {\bullet}^{E_3}
\ar@{-}[r] \ar@{-}[d] & {\bullet}^{E_4} \ar@{-}[r] & {\bullet}^{E_5}\\ &
{\bullet}^{E_2} &}
\bigskip

Then
$$
\tilde{D} \sim_{\mathbb{Q}}
\pi^*(D)-a_1E_1-a_2E_2-a_3E_3-a_4E_4-a_5E_5 \text{ .}
$$

We have a line $L_1$ intersecting the fundamental cycle as
following
\begin{eqnarray*}
\tilde{L_1} & \sim_{\mathbb{Q}} & \pi_1^*(L_1) - \frac{5}{4}E_1 -
\frac{3}{4}E_2 - \frac{3}{2}E_3 -
E_4 - \frac{1}{2} E_5\\
\end{eqnarray*}

Since $4L_1 \sim_{\mathbb{Q}} -K_X$ we see that $\mathrm{lct}(X)
\leq \frac{1}{6}$ and moreover we can assume that $L_1 \not \in
\text{Supp}D$.

From the inequalities
\begin{eqnarray*}
0 \leq \tilde{D} \cdot \tilde{L}_1 & = & 1 - a_1\\
0 \leq E_1 \cdot \tilde{D} & = & 2a_1 - a_3\\
0 \leq E_2 \cdot \tilde{D} & = & 2a_2 - a_3\\
0 \leq E_3 \cdot \tilde{D} & = & 2a_3 - a_1 - a_2 - a_4\\
0 \leq E_4 \cdot \tilde{D} & = & 2a_4 - a_3 - a_5\\
0 \leq E_5 \cdot \tilde{D} & = & 2a_5 - a_4\\
\end{eqnarray*}
we see that
$$
a_3 \leq 2 a_1 \text{, } a_3 \leq 2 a_2 \text{, } a_4 \leq a_3
\text{, } a_5 \leq a_4
$$
and
$$
a_4 \leq 2 a_5 \text{, } a_3 \leq \frac{3}{2} a_4 \text{, }  a_2
\leq \frac{5}{6} a_3 \text{, } a_1 \leq \frac{5}{6} a_3 \text{ .}
$$

In particular we get the following upper bounds
$$a_1 \leq 1 \text{, } a_2 \leq \frac{5}{3} \text{, }
a_5 \leq a_4 \leq a_3 \leq 2 \text{ .}$$

The equivalence
$$
K_{\tilde{X}} + \lambda \tilde{D} + \lambda a_1E_1 + \lambda
a_2E_2 + \lambda a_3E_3 + \lambda a_4E_4 + \lambda a_5E_5
\sim_{\mathbb{Q}} \pi^*(K_X + \lambda D)
$$
implies that there is a point $Q\in E_1\cup E_2\cup E_3\cup E_4
\cup E_5$ such that the pair
$$
K_{\tilde{X}} + \lambda \tilde{D} + \lambda a_1E_1 + \lambda
a_2E_2 + \lambda a_3E_3 + \lambda a_4E_4 + \lambda a_5E_5
$$
is not log canonical at $Q$.
\begin{itemize}
\item If the point $Q \in E_1 \backslash E_3$ then
$$
K_{\tilde{X}} + \lambda \tilde{D} + \lambda a_1E_1
$$
is not log canonical at the point $Q$ and so is the pair
$$
K_{\tilde{X}} + \lambda \tilde{D} + E_1 \text{ .}
$$
By adjunction $(E_1, \lambda \tilde{D}|_{E_1})$ is not log
canonical at $Q$ and
$$\frac{4}{5} \geq \frac{4}{5} a_1 \geq 2a_1 - \frac{6}{5}a_1 \geq 2a_1 - a_3 = \tilde{D} \cdot E_1
\geq
 \text{mult}_Q\Bigl(\tilde{D} \cdot E_1 \Bigr)  > 6 \text{ ,}$$
which is contradiction.

\item If $Q\in E_1 \cap E_3$ then the log pair
$$
K_{\tilde{X}} + \lambda \tilde{D} + \lambda a_1E_1 + \lambda
a_3E_3
$$
is not log canonical at the point $Q$ and so is the log pair
$$
K_{\tilde{X}} + \lambda \tilde{D} + E_1 + \lambda a_3E_3 \text{ .}
$$
By adjunction it follows that
$$
2a_1 - a_3 = \tilde{D} \cdot E_1 \geq
\text{mult}_Q\Bigl(\tilde{D}|_{E_1} \Bigr) =
 \text{mult}_Q\Bigl(\tilde{D} \cdot E_1 \Bigr) > 6 - a_3 \text{ .}
$$
and this implies that $a_3 > 3$, which is false.

\item If $Q\in E_3$ but $Q\not \in E_1 \cup E_2 \cup E_4$ then
$$
K_{\tilde{X}} + \lambda \tilde{D} + \lambda a_3E_3
$$
is not log canonical at the point $Q$ and so is the pair
$$
K_{\tilde{X}} + \lambda \tilde{D} + E_3 \text{, since } \lambda
a_3 \leq 1 \text{ .}
$$
By adjunction $(E_3, \lambda \tilde{D}|_{E_3})$ is not log
canonical at $Q$ and
$$2a_3- \frac{a_3}{2} - \frac{a_3}{2}-\frac{2}{3} a_3 \geq 2a_3 - a_1 - a_2 - a_4 = \tilde{D} \cdot E_3 \geq
 \text{mult}_Q\Bigl(\tilde{D} \cdot E_3 \Bigr) > 6 \text{ ,}$$
implies that $a_3>18$ which is false.

\item If $Q\in E_3 \cap E_4$ then the log pair
$$
K_{\tilde{X}} + \lambda \tilde{D} + \lambda a_3E_3 + \lambda
a_4E_4
$$
is not log canonical at the point $Q$ and so are the log pairs
$$
K_{\tilde{X}} + \lambda \tilde{D} + E_3 + \lambda a_4E_4 \text{
and  } K_{\tilde{X}} + \lambda \tilde{D} + \lambda a_3E_3 + E_4
\text{ .}
$$
By adjunction it follows that

$$
2a_4 - a_3 - a_5= \tilde{D} \cdot E_4 \geq
\text{mult}_Q\Bigl(\tilde{D}|_{E_4} \Bigr) =
 \text{mult}_Q\Bigl(\tilde{D} \cdot E_4 \Bigr) > 6 - a_3
$$
and
$$
2a_3 - a_2 - a_1 - a_4 = \tilde{D} \cdot E_3 \geq
\text{mult}_Q\Bigl(\tilde{D}|_{E_3} \Bigr) =
 \text{mult}_Q\Bigl(\tilde{D} \cdot E_3 \Bigr) > 6 - a_4 \text{ .}
$$

The last inequality implies that $a_3 > 6$ which is false.

\item $Q\in E_4 \backslash (E_3 \cap E_5)$ then the log pair
$$
K_{\tilde{X}} + \lambda \tilde{D} + \lambda a_4E_4
$$
is not log canonical at the point $Q$ and so is the pair
$$
K_{\tilde{X}} + \lambda \tilde{D} + E_4 \text{ .}
$$
By adjunction $(E_4, \lambda \tilde{D}|_{E_4})$ is not log
canonical at $Q$ and
$$2a_4 - a_3 - a_5 = \tilde{D} \cdot E_4 \geq \text{mult}_Q\Bigl(\tilde{D}|_{E_4} \Bigr)
=
 \text{mult}_Q\Bigl(\tilde{D} \cdot E_4 \Bigr) > 6 \text{ ,}$$
implies that $a_4>12$ which is false.

\item $Q\in E_5 \backslash E_4$ then the log pair
$$
K_{\tilde{X}} + \lambda \tilde{D} + \lambda a_5E_5
$$
is not log canonical at the point $Q$ and so is the pair
$$
K_{\tilde{X}} + \lambda \tilde{D} + E_5 \text{ .}
$$
By adjunction $(E_5, \lambda \tilde{D}|_{E_5})$ is not log
canonical at $Q$ and
$$2a_5 - a_4 = \tilde{D} \cdot E_5 \geq \text{mult}_Q\Bigl(\tilde{D}|_{E_5} \Bigr) =
 \text{mult}_Q\Bigl(\tilde{D} \cdot E_5 \Bigr) > 6 \text{ ,}$$
implies that $a_5 > 6$ which is false.

\item $Q\in E_4 \cap E_5$ then the log pair
$$
K_{\tilde{X}} + \lambda \tilde{D} + \lambda a_4E_4 + \lambda
a_5E_5
$$
is not log canonical at the point $Q$ and so is the log pair
$$
K_{\tilde{X}} + \lambda \tilde{D} + E_5 + \lambda a_4E_4 \text{ .}
$$
By adjunction it follows that
$$
2a_5 - a_4 = \tilde{D} \cdot E_5 \geq
\text{mult}_Q\Bigl(\tilde{D}|_{E_5} \Bigr) =
 \text{mult}_Q\Bigl(\tilde{D} \cdot E_5 \Bigr) > 6 - a_4 \text{ .}
$$
and we see then that $a_5 > 3$ which is not possible.
\end{itemize}
\end{proof}

\newpage

\section{Tables of Global Log Canonical Thresholds}

\begin{table}[htbp]
\label{table:degree1PicZ} \caption[Table 1]{Del Pezzo surfaces of
degree 1 and $\text{Pic}(X) \cong \mathbb{Z}$}
 \centerline{
\begin{tabular}{| l | c | }
\hline  {\bf \phantom0Singularity Type} & $\mathrm{lct}(X)$\\
\hline \phantom0   & \\
 \phantom{00}$\mathbb{E}_8$ & $\frac{1}{6}$\\
\phantom0  & \\
 \hline \phantom0   & \\
 \phantom{00}$\mathbb{E}_7 + \mathbb{A}_1 $ & $\frac{1}{4}$\\
 \phantom0 & \\
 \hline \phantom0  & \\
 \phantom{00} $\mathbb{E}_6+ \mathbb{A}_2 \text{, } \mathbb{D}_8$ & $\frac{1}{3}$\\
 \phantom0  & \\
 \hline \phantom0   & \\
 \phantom{00} $\mathbb{A}_8 \text{, } \mathbb{A}_7+ \mathbb{A}_1 \text{, } \mathbb{D}_6 + 2 \mathbb{A}_1
   \text{, } \mathbb{D}_5+ \mathbb{A}_3 \text{, } \mathbb{D}_4+ \mathbb{D}_4 $ & $\frac{1}{2}$\\
  \phantom0 & \\
\hline \phantom0   & \\
 \phantom{00}$2 \mathbb{A}_4 $ & $\frac{4}{5}$\\
 \phantom0 & \\
 \hline \phantom0   & \\
 \phantom{00}$ \mathbb{A}_5 + \mathbb{A}_2 + \mathbb{A}_1 $ & $\frac{2}{3}$\\
 \phantom0 & \\
 \hline \phantom0  & \\
 \phantom{00}$4 \mathbb{A}_2 $ and $|-K_X|$ has no cuspidal curves &  1 \\
\phantom0 & \\
 \hline \phantom0  & \\
\phantom{00}$4 \mathbb{A}_2$ and $|-K_X|$ has a cuspidal curve,&  $\frac{5}{6}$\\
\phantom{00}but no  cuspidal curve $C$ such that $\text{Sing}(C)= \mathbb{A}_2$ & \\
 \phantom0 & \\
 \hline \phantom0  & \\
 \phantom{00} $4 \mathbb{A}_2 $ and $|-K_X|$ has a cuspidal curve such that $\text{Sing}(C) = \mathbb{A}_2 $& $\frac{2}{3}$\\
 \phantom0 & \\
 \hline \phantom0  & \\
 \phantom{00}$2 \mathbb{A}_3 + 2 \mathbb{A}_1$ and $|-K_X|$ has no cuspidal curves & 1\\
 \phantom0 & \\
 \hline \phantom0  & \\
 \phantom{00}$2 \mathbb{A}_3 + 2 \mathbb{A}_1$ and $|-K_X|$ has a cuspidal curve , &  $\frac{5}{6}$\\
\phantom{00}but no  cuspidal curve $C$ such that $\text{Sing}(C)= \mathbb{A}_1 $ & \\
\phantom0 & \\
 \hline \phantom0  & \\
 \phantom{00}$2 \mathbb{A}_3 + 2 \mathbb{A}_1$ and $|-K_X|$ has a cuspidal curve with $\text{Sing}(C) = \mathbb{A}_1 $ & $\frac{3}{4}$\\
 \phantom0 & \\
\hline
\end{tabular}}
\end{table}

\newpage

\begin{table}[htbp]
\label{table:degree2PicZ} \caption[Table 1]{Del Pezzo surfaces of
degree 2 and $\text{Pic}(X) \cong \mathbb{Z}$}
 \centerline{
\begin{tabular}{| l | c | }
\hline  {\bf \phantom0Singularity Type} & $\mathrm{lct}(X)$\\
\hline \phantom0   & \\
 $\mathbb{E}_7$ & $\frac{1}{6}$\\
\phantom0  & \\
 \hline \phantom0   & \\
 $\mathbb{D}_6 + \mathbb{A}_1$ & $\frac{1}{4}$\\
 \phantom0 & \\
 \hline \phantom0  & \\
  $\mathbb{A}_7$ & $\frac{1}{3}$\\
 \phantom0  & \\
 \hline \phantom0   & \\
  $\mathbb{D}_4 + 3 \mathbb{A}_1$ & $\frac{1}{2}$\\
  \phantom0 & \\
  \hline \phantom0 & \\
 $\mathbb{A}_5 + \mathbb{A}_2$ & $\frac{1}{3}$\\
 \phantom0  & \\
 \hline \phantom0  & \\
 $2 \mathbb{A}_3 + \mathbb{A}_1$ & $\frac{1}{2}$\\
  \phantom0  & \\
\hline
\end{tabular}}
\end{table}

\begin{table}[htbp]
\label{table:degree3PicZ} \caption[Table 1]{Del Pezzo surfaces of
degree 3 and $\text{Pic}(X) \cong \mathbb{Z}$}
 \centerline{
\begin{tabular}{| l | c | }
\hline  {\bf \phantom0Singularity Type} & $\mathrm{lct}(X)$\\
\hline \phantom0   & \\
 $\mathbb{E}_6$ & $\frac{1}{6}$\\
\phantom0  & \\
 \hline \phantom0   & \\
 $\mathbb{A}_5 + \mathbb{A}_1$ & $\frac{1}{4}$\\
 \phantom0 & \\
 \hline \phantom0  & \\
  $3\mathbb{A}_2$ & $\frac{1}{3}$\\
 \phantom0  & \\
\hline
\end{tabular}}
\end{table}

\newpage

\begin{table}[htbp]
\label{table:degree4PicZ} \caption[Table 1]{Del Pezzo surfaces of
degree 4 and $\text{Pic}(X) \cong \mathbb{Z}$}
 \centerline{
\begin{tabular}{| l | c | }
\hline  {\bf \phantom0Singularity Type} & $\mathrm{lct}(X)$\\
\hline \phantom0   & \\
 $\mathbb{D}_5$ & $\frac{1}{6}$\\
\phantom0  & \\
 \hline \phantom0  & \\
  $\mathbb{A}_3+2\mathbb{A}_1$ & $\frac{1}{3}$\\
 \phantom0  & \\
\hline
\end{tabular}}
\end{table}

\begin{table}[htbp]
\label{table:degree5PicZ} \caption[Table 1]{Del Pezzo surfaces of
degree 5  and $\text{Pic}(X) \cong \mathbb{Z}$}
 \centerline{
\begin{tabular}{|l | c |}
\hline {\bf \phantom0Singularity Type} & $\mathrm{lct}(X)$\\
\hline \phantom0 & \\
 \phantom0 $\mathbb{A}_4$ &  $\frac{1}{6}$\\
\phantom0 & \\
\hline
\end{tabular}}
\end{table}

\begin{table}[htbp]
\label{table:degree6PicZ} \caption[Table 1]{Del Pezzo surfaces of
degree 6  and $\text{Pic}(X) \cong \mathbb{Z}$}
 \centerline{
\begin{tabular}{|l | c |}
\hline {\bf \phantom0Singularity Type} & $\mathrm{lct}(X)$\\
\hline \phantom0 & \\
 \phantom0 $ \mathbb{A}_2 + \mathbb{A}_1$ &  $\frac{1}{6}$\\
\phantom0 & \\
\hline
\end{tabular}}
\end{table}

\newpage

\begin{table}[htbp]
\label{degree1E8DnA8A7A6A5} \caption[Table 1]{Del Pezzo surfaces
of degree 1}
 \centerline{
\begin{tabular}{|l | c |}
\hline {\bf \phantom0Singularity Type} & $\mathrm{lct}(X)$\\
\hline \phantom0 & \\
 \phantom0 $\mathbb{E}_8$ & $\frac{1}{6}$\\
\phantom0 & \\
\hline \phantom0 & \\
\phantom0 $\mathbb{E}_7 \text{, } \mathbb{E}_7+ \mathbb{A}_1$ & $\frac{1}{4}$\\
\phantom0 & \\
\hline \phantom0 & \\
\phantom0 $\mathbb{E}_6 \text{, } \mathbb{E}_6+ \mathbb{A}_2
\text{, } \mathbb{E}_6+ \mathbb{A}_1
 $ & $\frac{1}{3}$\\
\phantom0 & \\
\hline \phantom0 & \\
\phantom0 $\mathbb{D}_8$ & $\frac{1}{3}$\\
\phantom0 & \\
\hline \phantom0 & \\
\phantom0 $\mathbb{D}_7$ & $\frac{2}{5}$\\
\phantom0 & \\
\hline \phantom0 & \\
\phantom0 $\mathbb{D}_6 \text{, }\mathbb{D}_6 + 2 \mathbb{A}_1
\text{, }
\mathbb{D}_6 + \mathbb{A}_1 $ & $\frac{1}{2}$\\
\phantom0 & \\
\hline \phantom0 & \\
\phantom0 $\mathbb{D}_5 \text{, }\mathbb{D}_5 + \mathbb{A}_3
\text{, } \mathbb{D}_5 + \mathbb{A}_2 \text{, } \mathbb{D}_5 + 2
\mathbb{A}_1 \text{, }
\mathbb{D}_5 +  \mathbb{A}_1 $ & $\frac{1}{2}$\\
\phantom0 & \\
\hline \phantom0 & \\
\phantom0 $\mathbb{D}_4 \text{, }\mathbb{D}_4 + \mathbb{D}_4
\text{, } \mathbb{D}_4 + \mathbb{A}_3 \text{, }\mathbb{D}_4 +
\mathbb{A}_2 \text{, }
\mathbb{D}_4 + 4 \mathbb{A}_1 \text{, }$ & \\
 \phantom0
$\mathbb{D}_4 + 3 \mathbb{A}_1 \text{, } \mathbb{D}_4 +
2 \mathbb{A}_1 \text{, } \mathbb{D}_4 + \mathbb{A}_1$ & $\frac{1}{2}$\\
\phantom0 & \\
\hline \phantom0 & \\
 \phantom0 $\mathbb{A}_8$ & $\frac{1}{2}$\\
\phantom0 & \\
\hline \phantom0 & \\
\phantom0 $\mathbb{A}_7 \text{, } \mathbb{A}_7+ \mathbb{A}_1$ & $\frac{1}{2}$\\
\phantom0 & \\
\hline \phantom0 & \\
\phantom0 $\mathbb{A}_7'$ & $\frac{8}{15}$\\
\phantom0 & \\
\hline \phantom0 & \\
\phantom0 $\mathbb{A}_6 \text{, } \mathbb{A}_6 + \mathbb{A}_1$ & $\frac{2}{3}$\\
\phantom0 & \\
\hline \phantom0 & \\
\phantom0 $\mathbb{A}_5 \text{, } \mathbb{A}_5 + \mathbb{A}_1
\text{, } \mathbb{A}_5 + 2 \mathbb{A}_1 \text{, } \mathbb{A}_5 +
\mathbb{A}_2 \text{, } \mathbb{A}_5 +  \mathbb{A}_2 + \mathbb{A}_1
$ & $\frac{2}{3}$\\
\phantom0 & \\
\hline
%\multicolumn{2}{l}
%{\footnotesize Parameter values: $V_0=0.0625$,
%$\mu=0.1$, $\xi=16$, $S_0=100$ and $r=0.0953$.}
\end{tabular}}
\end{table}

\newpage
\begin{table}[htbp]
\label{degree1A4A3}
 \caption[Table 6]{Del Pezzo surfaces of degree 1}
 \centerline{
\begin{tabular}{|l | c |}
\hline {\bf \phantom0Singularity Type} & $\mathrm{lct}(X)$\\
\hline \phantom0 & \\
\phantom0 & \\
\phantom{00} $\mathbb{A}_4 \text{, } \mathbb{A}_4 + \mathbb{A}_4
\text{, }\mathbb{A}_4 + \mathbb{A}_3 $ & $\frac{4}{5}$ \\
\phantom0 & \\
\phantom0 & \\
\hline \phantom0 & \\
\phantom{00} $\mathbb{A}_4 + \mathbb{A}_2 +  \mathbb{A}_1 \text{, }
\mathbb{A}_4 + 2 \mathbb{A}_1  \text{, } \mathbb{A}_4 +
\mathbb{A}_2  \text{, } \mathbb{A}_4 + \mathbb{A}_1$ & \\
\phantom{00} If $|-K_X|$ has no cuspidal curve $C$ such that & $\frac{4}{5}$ \\
\phantom{00} $\mathbb{A}_1= \text{Sing}(C)$ or $\mathbb{A}_2= \text{Sing}(C)$ & \\
\phantom0 & \\
\hline \phantom0 & \\
\phantom{00} $\mathbb{A}_4 + \mathbb{A}_2 +  \mathbb{A}_1 \text{, }
\mathbb{A}_4 + 2 \mathbb{A}_1  \text{, }  \mathbb{A}_4 + \mathbb{A}_1$ & \\
\phantom{00} If $|-K_X|$ has a cuspidal curve $C$ such that $ \mathbb{A}_1 = \text{Sing}(C)$, & $\frac{3}{4}$\\
\phantom{00} but no cuspidal curve $C$ such that $ \mathbb{A}_2 = \text{Sing}(C)$ & \\
\phantom0 & \\
\hline \phantom0 & \\
\phantom{00} $\mathbb{A}_4 + \mathbb{A}_2 +  \mathbb{A}_1 \text{, }
\mathbb{A}_4 +\mathbb{A}_2$ &\\
\phantom{00} If $|-K_X|$ has a cuspidal curve $C$ such that $\mathbb{A}_2 = \text{Sing}(C)$ & $\frac{2}{3}$\\
\phantom0 & \\
\phantom0 & \\
\hline \phantom0 & \\
\phantom0 & \\
\phantom{00} $\mathbb{A}_3 \text{, } 2 \mathbb{A}_3$ & $1$\\
\phantom0 & \\
\phantom0 & \\
\hline \phantom0 & \\
\phantom{00} $\mathbb{A}_3+ 4 \mathbb{A}_1  \text{, } \mathbb{A}_3+ 3
\mathbb{A}_1  \text{, } 2\mathbb{A}_3+ 2 \mathbb{A}_1  \text{, }
\mathbb{A}_3+ 2 \mathbb{A}_1  \text{, }
\mathbb{A}_3+  \mathbb{A}_1  $ & \\
\phantom{00} $ \mathbb{A}_3+ \mathbb{A}_2 \text{, } \mathbb{A}_3+
\mathbb{A}_2 +\mathbb{A}_1 \text{, } \mathbb{A}_3+
\mathbb{A}_2 + 2\mathbb{A}_1$ &  $1$\\
\phantom{00} If $|-K_X|$ has no cuspidal curves & \\
\phantom0 & \\
\hline \phantom0 & \\
\phantom{00} $\mathbb{A}_3+ 4 \mathbb{A}_1  \text{, } \mathbb{A}_3+ 3
\mathbb{A}_1  \text{, } \mathbb{A}_3+ 2 \mathbb{A}_1  \text{, }
\mathbb{A}_3+  \mathbb{A}_1  \text{, } $ & \\
\phantom{00} If $|-K_X|$ has a cuspidal curve such that $\text{Sing}(C) = \mathbb{A}_1 $, & $\frac{3}{4}$\\
\phantom{00} but no  cuspidal curve $C$ such that $ \mathbb{A}_2 = \text{Sing}(C)$ & \\
\phantom0 & \\
\hline \phantom0 & \\
\phantom{00} $\mathbb{A}_3+ 4 \mathbb{A}_1  \text{, } \mathbb{A}_3+ 3
\mathbb{A}_1  \text{, } \mathbb{A}_3+ 2 \mathbb{A}_1  \text{, }
\mathbb{A}_3+  \mathbb{A}_1  \text{, } $ & \\
\phantom{00} $\mathbb{A}_3+ \mathbb{A}_2 \text{, } \mathbb{A}_3+
\mathbb{A}_2 +\mathbb{A}_1 $ &  $\frac{5}{6}$\\
\phantom{00} If $|-K_X|$ has cuspidal curves $C$, but $\text{Sing}(C) \neq \mathbb{A}_1 $ and $\text{Sing}(C) \neq  \mathbb{A}_2 $& \\
\phantom0 & \\
\hline \phantom0 & \\
\phantom{00} $\mathbb{A}_3+  \mathbb{A}_2  \text{, } \mathbb{A}_3+
\mathbb{A}_2 +\mathbb{A}_1 $ &  \\
\phantom{00} If $|-K_X|$ has a cuspidal curve $C$ such that $\text{Sing}(C) = \mathbb{A}_2 $& $\frac{2}{3}$\\
\phantom0 & \\
\phantom0 & \\
\hline
\end{tabular}}
\end{table}

\newpage

\tableofcontents

 \vskip 1cm

        \noindent
           School of Mathematics\\%
           The University of Edinburgh\\%
           Kings Buildings, Mayfield Road\\%
           Edinburgh EH9 3JZ, UK\\%
           \vskip 0.2cm

           \noindent\texttt{D.Kosta@sms.ed.ac.uk}

\end{document}